\documentclass[a4paper,10pt]{article}
\usepackage{amssymb}
\usepackage{latexsym}
\usepackage{amsmath}
\usepackage{amsthm}
\usepackage{amsfonts}
\usepackage{amssymb}
\usepackage{amsmath,amscd}
\usepackage{graphicx}

\newtheorem{Th}{Theorem}

\newtheorem{Lm}{Lemma}

\newtheorem{Dfi}{Definition}
\newtheorem{Rm}{Remark}

\newcommand{\be}{\begin{equation}}
\newcommand{\ee}{\end{equation}}
\newcommand{\bes}{\begin{equation*}}
\newcommand{\ees}{\end{equation*}}

\newcommand{\R}{\mathbb{R}}
\newcommand{\N}{\mathbb{N}}
\newcommand{\C}{\mathbb{C}}
\newcommand{\Z}{\mathbb{Z}}

\newcommand\res{\mathop{\hbox{\vrule height 7pt width .5pt depth 0pt
\vrule height .5pt width 6pt depth 0pt}}\nolimits}

\newcommand{\reset}{\setcounter{equation}{0}\setcounter{Th}{0}\setcounter{Prop}{0}\setcounter{Co}{0}
\setcounter{Lm}{0}\setcounter{Rm}{0}}
\DeclareMathOperator{\arcsinh}{arcsinh}
\def\Xint#1{\mathchoice
{\XXint\displaystyle\textstyle{#1}}%
{\XXint\textstyle\scriptstyle{#1}}%
{\XXint\scriptstyle\scriptscriptstyle{#1}}%
{\XXint\scriptscriptstyle\scriptscriptstyle{#1}}%
\!\int}
\def\XXint#1#2#3{{\setbox0=\hbox{$#1{#2#3}{\int}$}
\vcenter{\hbox{$#2#3$}}\kern-.5\wd0}}
\def\dashint{\Xint-}
\def\La{\Lambda}
\def\La{\Lambda}
\def\ti{\tilde}
\def\lf{\left}
\def\rg{\right}

\def\al{\alpha}
\def\la{\lambda}

\def\ds{\displaystyle}
\def\ov{\overline}
\def\Om{\Omega}
\def\om{\omega}
\def\p{\partial}

\def\res{\mathop{\hbox{\vrule height 7pt width .5pt 
depth 0pt\vrule height .5pt width 6pt depth 0pt}}\nolimits}

\makeatletter
\newcommand{\tpitchfork}{%
  \vbox{
    \baselineskip\z@skip
    \lineskip-.52ex
    \lineskiplimit\maxdimen
    \m@th
    \ialign{##\crcr\hidewidth\smash{$-$}\hidewidth\crcr$\pitchfork$\crcr}
  }%
}
\makeatother

\begin{document}

\title{The Parametric Willmore Flow}

\author{Francesco Palmurella\footnote{Scuola Normale Superiore Piazza dei Cavalieri 7, 56126 Pisa, Italy. }  \, and 
Tristan Rivi\`ere\footnote{Department of Mathematics, ETH Zentrum,
CH-8093 Z\"urich, Switzerland.}}

%\date{ }
\maketitle

{\bf Abstract :}{\it We establish a minimal  positive existence time of the parametric Willmore flow for any smooth initial data  (smooth immersion of a closed oriented surface). The minimal existence time 
is a function exclusively of geometric data which in particular are all well defined for general weak lipschitz $W^{2,2}$ immersions. This fact combined with the conservation law formulation of the equation given by the first author in \cite{Riv-w} opens the possibility for defining the Willmore flow for weak lipschitz $W^{2,2}$ initial data.}

\medskip

\noindent{\bf Math. Class. 53E40, 53E10, 49Q10,  35K46, 53A05, 58E15, 58E30, 35J35, 35J48}
\section{Introduction}
The analysis of Willmore energy takes it's roots in the pioneered work of Leon Simon \cite{Sim} giving the existence of a Willmore Torus minimizing the Willmore energy in any euclidian space. It has been followed by important contributions by Ernst Kuwert and Reiner Sch\"atzle (\cite{KS1},\cite{KS2},\cite{KS3}...) in which in particular the two authors have established a short time existence of the Willmore flow.  The approach adopted in these papers consists in working with the immersed surfaces in its ambient (mostly euclidian) space. The so called ``ambient approach'' - combined with the GMT minmax theory of Fred Almgren and Jon Pitts - has also been successfully adopted by Fernando Cod\'a Marques and Andre Neves in 
their milestone paper \cite{MN} establishing the absolute minimality of the Willmore Torus in ${\R}^3$ among non zero genus surfaces (i.e. proving the Willmore conjecture in 3 dimension).

In parallel to these important contributions to the geometric analysis of the Willmore functional the second author of the present work introduced in \cite{Riv-w} and \cite{Riv-14} an alternative approach to the variational study of this Lagrangian. The main objects in this approach, called ``parametric approach'', are immersions of surfaces into euclidian or more generally Riemannian spaces. One of the motivation was to rely on almost classical functional  analysis of maps at the so called energy level. The advantage of the ``parametric approach''  over the ``ambient approach'' is that one can consider a space of weak  immersions ($W^{2,2}_loc$ Lipschitz branched immersions) which is weakly sequentially closed under Willmore energy control and the framework is getting closer to the one of classical conformally invariant geometric analysis of maps.
The parametric approach has showed  to be successful for instance for performing refined bubble tree analysis and for ``tracking'' the possible loss of energies in neck regions (\cite{BeRi}, \cite{LaRi-w}, \cite{MiRi}...)
or for implementing minmax operation (the ambient approach, by the absence of a weak formulation of the energy, is mostly used for performing minimisation operations or for ``running'' the flow for smooth initial data).

The goal of the present work is to extend the parametric approach to the Willmore flow. Our motivation, as it was originally in the static case, is to consider in future works a precise ``bubble tree analysis'' of the Willmore flow and in particular the energy quantization phenomenon. The other motivation is to have a suitable framework for continuing the flow after the first blow-up time while the asymptotic surface might create branched points
at the blow-up.

\medskip

We follow the notations in \cite{Pa-Ci}. $\vec{\Phi}$ is an immersion from a closed oriented surface $\Sigma$ into ${\R}^m$ for $m>2$. We denote by $g_{\vec{\Phi}}$ (sometimes it will simply be denoted by $g$)
the  {\it first fundamental form} that is the metric induced by the immersion 
\[
\forall\, X,Y\in T_x\Sigma\quad\quad g_{\vec{\Phi}}(X,Y)=\vec{\Phi}^\ast g_{{\R}^m}(X,Y):=\lf<d\vec{\Phi}_x(X),d\vec{\Phi}_x(Y)\rg>\ ,
\]
where $<\cdot,\cdot>$ denotes the scalar product in ${\R}^m$.The \emph{volume form} associated to $g_{\vec{\Phi}}$ on $\Sigma$ is locally given 
by
\[
dvol_{g_{\vec{\Phi}}}:= \sqrt{ det (g(\p_{x_i},\p_{x_j}))}\ dx_1\wedge dx_2,
\]
where $(x_1,x_2)$ are arbitrary local positive coordinates. We denote by $\vec{n}(x)\in \wedge^{m-2}{\R}^m$ thee Gauss unit associated to this immersion giving the normal $m-2$-plane normal to the tangent
plane at $\vec{\Phi}(x)$ and by $\pi_{\vec{n}}$ the orthogonal projection onto this plane. We shall be denoting $\pi_T$ the tangential projection onto $\vec{\Phi}_\ast T_x\Sigma$.

\medskip

The \emph{second fundamental form} at $p\in S$ is the bilinear map which assigns to a pair of vectors ${X},{Y}$ in $T_x\Sigma$ an orthogonal vector to $\vec{\Phi}_\ast T_x\Sigma$ that we shall
denote by $\vec{\mathbb I}({X},{Y})$. This normal vector 
expresses how much the Gauss map varies along these directions ${X}$ and ${Y}$. Precisely, it is given by
 \[
 \begin{array}{rcl}
\ds\vec{\mathbb I}_{\vec{\Phi}}(x)\colon T_x\Sigma\times T_x\Sigma &\ds \longrightarrow & \ds(\vec{\Phi}_\ast T_x\Sigma)^\perp\\[5mm]
 \ds({X},{Y}) &\ds\mapsto &\ds \pi_{\vec{n}}\left(d^2\vec{\Phi}(X,Y)\right),
 \end{array}
 \]
where $X$  and $Y$ are extended smoothly around $x$. The {\it mean curvature vector} is given by
\[
\vec{H}_{\vec{\Phi}}:=\frac{1}{2} \mbox{tr}(g_{\vec{\Phi}}^{-1}\ \vec{\mathbb I}_{\vec{\Phi}})=\frac{1}{2} \sum_{i,j=1}^2 g^{ij}\  \vec{\mathbb I}(\p_{x_i},\p_{x_j}) \ ,
\]
where $(x_1,x_2)$ are arbitrary local coordinates in $\Sigma$ and $(g^{ij})_{ij}$ is the inverse matrix to $(g_{\vec{\Phi}}(\p_{x_i},\p_{x_j}))$. We are interested with the heat flow of the {\it Willmore Energy} of the immersion
given by
\[
W(\vec{\Phi}):=\int_{\Sigma}|\vec{H}|^2\ dvol_g.
\]
Using the fundamental Gauss identity
\[
|\vec{\mathbb I}_{\vec{\Phi}}|_{g_{\vec{\Phi}}}^2=4|\vec{H}_{\vec{\Phi}}|^2-2\,K_{g_{\vec{\Phi}}}.
\]
where $K_g$ is the Gauss curvature of the metric $g$,  and using Gauss-Bonnet theorem, we obtain the following expression of the Willmore energy of an immersion into ${\R}^m$
of an arbitrary closed surface: 
\[
W(\vec{\Phi})=\frac{1}{4}\int_{\Sigma}|\vec{\mathbb I}|_g^2\ dvol_g+\pi\ \chi(\Sigma).
\]
where $\chi(\Sigma)$ is the Euler characteristic of $\Sigma$. One has (see \cite{Pa-Ci} section 1.2.2)
\[
|d\vec{n}|^2_g=|\vec{\mathbb I}|_g^2\ .
\]
Hence the Willmore energy is the homogeneous $W^{1,2}$ pseudo-norm of the Gauss map modulo the subtraction of a multiple of the Euler characteristic. For such an immersion there exists a constant Gauss curvature metric $h$ of value $+1,0,-1$ depending on the topology of $\Sigma$ and a function denoted $\al$ called ``conformal factor'' such that
\[
g=e^{2\al}\,h\ .
\]
When $g(\Sigma)>1$ (i.e. the genus of $\Sigma$ is larger than one) the metric $h$ is unique and its volume is given (thanks to Gauss Bonnet) by
\[
V_h:=\int_\Sigma dvol_h=4\pi\,(g(\Sigma)-1)\ .
\]
For the torus case we shall fix $V_h=1$ and in the sphere case, the volume is equal to $4\pi$ but the uniqueness of $h$ holds true modulo the action of the M\"obius transformations of $S^2$. We shall now consider the complex structure induced by $h$ with respect to which $g$ is conformal. For such complex structure the first variation of the Willmore energy with respect to a variation $\vec{w}\in C^\infty(\Sigma,{\R}^m)$ is given by
\[
\delta W_{\vec{\Phi}}(\vec{w}):=-\,4\,\int_\Sigma\  \Re\lf( \ov{\p}^h\lf[\p^h\vec{H}_{\vec{\Phi}}+|\vec{H}_{\vec{\Phi}}|^2\,\p^h\vec{\Phi}+2\,\vec{H}_{\vec{\Phi}}\cdot\vec{\frak{h}}^{\,0}\res \ov{\p}\vec{\Phi}  \rg]\rg)\cdot\vec{w}
\]
where $\vec{\frak{h}}_{\vec{\Phi}}^{\,0}$ is the Weingarten operator given in local complex coordinates by
\[
\vec{\frak{h}}_{\vec{\Phi}}^{\,0}=\p_z\lf(e^{-2\la}\,\p_z\vec{\Phi}\rg)\ dz\otimes \p_{\ov{z}}\quad \mbox{and }\quad g_{\vec{\Phi}}=2^{-1}\ e^{2\la}\ [dz\otimes d\ov{z}+d\ov{z}\otimes dz]\ ,
\]
and $\res$ is the tautological contraction between tensors. 
We denote by $\delta\vec{W}_{\vec{\Phi}}$ the vector in ${\R}^m$ given by
\[
\delta\vec{W}_{\vec{\Phi}}=-\,4\,  \lf<\Re\lf( \ov{\p}^h\lf[\p^h\vec{H}_{\vec{\Phi}}+|\vec{H}_{\vec{\Phi}}|^2\,\p^h\vec{\Phi}+2\,\vec{H}_{\vec{\Phi}}\cdot\vec{\frak{h}}^{\,0}_{\vec{\Phi}}\res \ov{\p}\vec{\Phi}  \rg]\rg),\ dvol_{g_{\vec{\Phi}}}\rg>
\]
Observe that\footnote{We have the convention
 \[
 \p_{z}:=\frac{1}{2} (\p_{x_1}-i\,\p_{x_2}) \quad\mbox{ and } \quad dz=dx_1+i\,dx_2
 \] in such a way that $dz\res\p_z=1$.}
\[
\vec{\frak{h}}_{\vec{\Phi}}^{\,0}=2^{-1}\p^h\lf(\p^h\vec{\Phi}\res g_{\vec{\Phi}}^{-1}\rg) \quad \mbox{ where }\quad g^{-1}_{\vec{\Phi}}=2\ e^{-2\la}\ [\p_z\otimes \p_{\ov{z}}+\p_{\ov{z}}\otimes \p_z]\ ,
\]
We denote
\[
\vec{\mathcal H}_{\vec{\Phi}}^{\,0}:=\vec{\frak{h}}_{\vec{\Phi}}^{\,0}\res h\ .
\]
We introduce
\begin{Dfi}
\label{df-proj}
 Denote by $P_{h}$ the $L^2$ orthogonal projection (for the $L^2$ scalar product induced by $h$) of quadratic form (i.e. sections of $(T^\ast\Sigma\otimes{\C})^{\otimes^2}$) onto the finite dimensional space of holomorphic quadratic forms\footnote{Holomorphic quadratic forms are sections of $(T^\ast\Sigma\otimes{\C})^{\otimes^2}$ which in local conformal coordinates  are of the form
 \[
 q= f(z)\ dz\otimes dz 
 \]
 where $\p_{\ov{z}}f=0$. The vector space of the holomorphic quadratic  forms identify in a canonical way  to the tangent space at $h$ to the space of constant Gauss curvature metric modulo the action of diffeomorphisms (see \cite{Tro})}.
\end{Dfi}
Finally for any section $U^{0,1}:=u^{0,1}\, \p_{\ov{z}}$ of $\wedge^{0,1}T\Sigma$ we assign the tangent vector-field $\vec{U}$ along the immersion $\vec{\Phi}(\Sigma)$ as follows
\[
\vec{U}:=\vec{\Phi}_\ast \lf(U^{0,1}+\ov{U^{0,1}}\rg)\ .
\]
We introduce the following definition
\begin{Dfi}
\label{df-param-conf-will}
A $C^1$ family $\vec{\Phi}_t$  into the space of immersions solution to the following system
\be
\label{I.flow}
\lf\{
\begin{array}{l}
\ds\frac{\p \vec{\Phi}}{\p t} =\delta \vec{\mathcal W}_{\vec{\Phi}} +\vec{U}\\[5mm]
\p^h\lf( U^{0,1}  \rg)\res h =2\, \lf(I-{P}_{h}\rg)\lf(  \delta\vec{W}   \cdot \vec{\mathcal H}^0\rg)\\[5mm]
\ds\frac{d h}{dt}=-\,4\,\Re \lf[P_{h}\lf(\delta\vec{W}\cdot \vec{\cal H}^0\rg)\rg]\\[5mm]
\ds \delta\vec{W}_{\vec{\Phi}}=-\,4\,  \lf<\Re\lf( \ov{\p}^h\lf[\p^h\vec{H}_{\vec{\Phi}}+|\vec{H}_{\vec{\Phi}}|^2\,\p^h\vec{\Phi}+2\,\vec{H}_{\vec{\Phi}}\cdot\vec{\frak{h}}^{\,0}_{\vec{\Phi}}\res \ov{\p}\vec{\Phi}  \rg]\rg),\ dvol_{g_{\vec{\Phi}}}\rg>\\[5mm]
\ds\vec{U}:=\vec{\Phi}_\ast \lf(U^{0,1}+\ov{U^{0,1}}\rg)\ .
\end{array}
\rg.
\ee
is called  solution to the Parametric Willmore Flow.\hfill $\Box$
\end{Dfi}
Our main result in this work is the following
\begin{Th}
\label{th-I.1}
Let $\vec{\Phi}_0$ be a smooth immersion of a closed oriented surface $\Sigma$ into ${\R}^m$.  There exists $\beta>0$, $\sigma>0$ and $\La>0$ depending exclusively on the topology of $\Sigma$, the dimension $m$, the initial energy $W(\vec{\Phi}(0))$, the length $l_\ast(0)$ of the shortest closed geodesic of $(\Sigma,h_0)$
and the $L^\infty$ norm of the conformal factor $\al(0)$ of $\vec{\Phi}_0$ with respect to $h_0$ such that for any $\sigma>R>0$ statisfying
\be
\label{small-init}
\sup_{x_0\in\Sigma}\int_{B_{R}^{h(0)}(x_0)}|\vec{\mathbb I}_{\vec{\Phi}(0)}|^2_{g_{\vec{\Phi}}}\ \chi^{h(0)}_{R,x_0}(x)\ \ dvol_{g_{\vec{\Phi}(0)}}\le \beta
\ee
and for any $T>0$ satisfying 
\be
\label{time}
\frac{T}{R^4}\le \La\ ,
\ee
there exists a smooth solution to the parametric Willmore Flow on $\Sigma\times [0,T]$ equal to $\vec{\Phi}_0$ at $t=0$. Moreover, if the genus of $\Sigma$ is non zero, the length $l_\ast(t)$ of the shortest closed geodesic of $(\Sigma,h(t))$ satisfies
\[
\inf_{t\in[0,T]}l_\ast(t)\ge l>0\ ,
\]
where
\[
l\ \exp\lf[C(\Sigma,m)\ l^{30}\rg]=l_\ast(0)\ ,
\]
where $C(\Sigma,m)>0$ only depends on the topology of $\Sigma$ and the dimension $m$.
\hfill $\Box$ 
\end{Th}
\begin{Rm}
\label{rm-I.0} 
If $\Sigma$ is a sphere  the statement of our result (as well as the proof of theorem~\ref{th-I.1})  simplifies drastically since there is no issue regarding the control of the length of the shortest geodesic. We have $h(t)\equiv h(0)$ and we are back to the much simpler Willmore parametric flow on $S^2$ considered previously by the authors in \cite{PalRiv} for initial data with small umbilic energy.\hfill $\Box$
\end{Rm}
\begin{Rm}
\label{rm-I.1} The parametric Willmore flow coincides  with the ``normal Willmore flow'' of Kuwert and Sch\"atzle if one ``forgets'' the parametrisation and in our notation it coincides with
\[
\pi_{\vec{n}}\frac{\p\vec{\Phi}}{\p t}=\delta \vec{\mathcal W}_{\vec{\Phi}}
\]
where $\pi_{\vec{n}}$ is the projection onto the normal bundle. It is then natural to compare the small existence time condition in \cite{KS2} with (\ref{time}). The dependence of the existence time with respect to the radius $\rho$ ensuring small energy is also quartic\footnote{Which is normal since the linear part of the flow is the bi-harmonic heat operator.}. However it has to be stressed at this point that the radius $\rho$ in \cite{KS2} is taken with respect to the distance in the ambient space ${\R}^m$ while the radius $R$ is taken with respect to the metric $h(0)$ of the underlying constant Gauss curvature metric. \hfill $\Box$
\end{Rm}
\begin{Rm}
\label{rm-I.2}
It is important to observe at this stage that the quantities taken from the data $\vec{\Phi}_0$ and giving a lower bound of the existence time, that is $l_\ast(0)$, $\|\al(0)\|_\infty$, $W(\vec{\Phi}_0)$ and the concentration radius $R$ are all \underbar{finite}
for any given weak lipschitz $W^{2,2}$ immersion (see \cite{Pa-Ci}) and  continuous with respect to the induced  topology. Hence our existence result with the lower bound on the existence time should be extended
to any initial data belonging to the weak lipschitz $W^{2,2}$ immersion which are approximable by smooth ones. This  question is the subject of a forthcoming work.\hfill $\Box$
\end{Rm}
\begin{Rm}
\label{rm-I.3}
It would be very interesting to extend the short time existence result~\ref{th-I.1} to the situation where the initial data has finitely many isolated branched points. In such a case we still have $l_\ast(0)>0$ but the initial conformal factor $\al(0)$ tends to minus infinity at these points. Hence a delicate analysis has to be carried out around these points. The motivation behind this question is  to extend the flow after the blow-ups until $T=+\infty$ for any initial data in the space of smooth branched immersion or even weak lipschitz $W^{2,2}$ possibly branched immersion which are approximable by smooth ones. .\hfill $\Box$
\end{Rm}
\section{Preliminaries}
\subsection{Conformality  Preserving  Variations : Flowing in the Teichm\"uller Space}
Recall the definition of the ${\p}^{h}$ operator acting on $0-1$ vector-field for $(\Sigma,h)$ . It is given in local coordinates by
\[
\p^{h}\lf(a\ \p_{\ov{z}}\rg):=\p_za\ dz\otimes\p_{\ov{z}}
\]
This definition is clearly independent of the choice of positive conformal coordinates because a change of such coordinates is holomorphic : $w(z)$ with $\ov{\p} w=0$ and we have in one hand
\[
\p_za\ dz\otimes\p_{\ov{z}}=\p_w a\ \ov{(\p_zw)}\ dw\otimes \p_{\ov{w}}=\p_w a\ \ov{(\p_wz)}^{-1}\ dw\otimes \p_{\ov{w}}=\p_w \lf(a\ \ov{(\p_wz)}^{-1}\rg)\ dw\otimes \p_{\ov{w}}\ ,
\]
and in the other hand
\[
a\ \p_{\ov{z}}=a\ \ov{(\p_wz)}^{-1}\ \p_{\ov{w}}\ .
\]
In a similar way ${\p}^{h}$ is also intrinsically defined on sections of $(T^\ast\Sigma)^{(0,1)}\otimes  (T^\ast\Sigma)^{(0,1)}$ locally in conformal coordinates for $h$ by
\[
{\p}^{h}\lf(b\, d\ov{z}\otimes d\ov{z}\rg)={\p}b\otimes d\ov{z}\otimes d\ov{z}\ .
\]
We keep denoting  $\vec{\Phi}$ a conformal immersion from $(\Sigma,h)$ (where $h$ is a constant Gauss curvature metric). For such an immersion, we write $g_{\vec{\Phi}}=e^{2\al}\,h$ and in conformal coordinates
\[
g_{\vec{\Phi}}= e^{2\la}\, [dx_1^2+dx_2^2]=\frac{e^{2\la}}{2}\ \lf[dz\otimes d\ov{z}+d\ov{z}\otimes d{z}\rg]
\]
and
\[
h= e^{2\nu}\, [dx_1^2+dx_2^2]=\frac{e^{2\nu}}{2}\ \lf[dz\otimes d\ov{z}+d\ov{z}\otimes d{z}\rg]\quad\mbox{ where }\quad \nu=\la-\al\ .
\]
We also denote
\[
g^{-1}_{\vec{\Phi}}:= e^{-2\la}\, [\p_{x_1}^2+\p_{x_2}^2]=2\,e^{-2\la}\ \lf[\p_{z}\otimes\p_{\ov{z}}+ \p_{\ov{z}} \otimes\p_{{z}} \rg]
\]
and
\[
h^{-1}:= e^{-2\nu}\, [\p_{x_1}^2+\p_{x_2}^2]=  2\,e^{- 2\,\nu}\ \lf[\p_{z}\otimes\p_{\ov{z}}+ \p_{\ov{z}} \otimes\p_{{z}} \rg]\ .
\]
We now define contractions operators. Let  $A:= E\otimes dz$ and  $B:= \p_{{z}}\otimes F$ where $E$ and $F$ are arbitrary elements from the tensor algebra of $(T\Sigma)^{(1,0)}$ and $(T\Sigma)^{(0,1)}$
\[
A\res B=  E\otimes F\ ,
\]
and in a dual way for $C:=E\otimes\p_z $ and $D:= dz\otimes F$
\[
C\res D=E\otimes F\ .
\]
The same holds replacing $dz$ and $\p_z$ respectively by $d\ov{z}$ and $\p_{\ov{z}}$. For instance, for ${F}$ being a section of $(T^\ast\Sigma)^{(1,0)}\otimes (T\Sigma)^{(0,1)}\otimes{\R}^m$ (i.e. in local conformal coordinates ${F}:=\vec{F}_{\ov{z}}^{z}\ d{z}\otimes \p_{\ov{z}}$) one has
\[
{F}\,\res \,h=2^{-1}\ e^{2\nu}\ {F}_{\ov{z}}^{z}\ d{z}\otimes d{z}
\]
%\[
%\al\res\ov{\beta} =\lf(a\, \p_{{z}}\rg)\res \lf(\ov{b}\ d\ov{z}\otimes d\ov{z}\rg)= a\, \ov{b}\ d\ov{z}\ ,
%\]
%and we have in particular
%\[
%\p(\al\res\ov{\beta} )=d(\al\res\ov{\beta} )=\p_z (a\,\ov{b})\,dz\wedge d\ov{z}\ ,
%\]
%and similarly
%\[
%\p(\ov{\al}\res{\beta} )=d(\ov{\al}\res{\beta} )=\p_z (\ov{a}\,{b})\,dz\wedge d\ov{z}\ .
%\]
We have the following lemma from \cite{Ri1} (lemma 2.2).
\begin{Lm}
\label{lm-delbar}
Let $(\Sigma,h)$ be a Riemann Surface equipped with a compatible constant Gauss curvature $h$. Let ${F}$ be a smooth section of $(T^\ast\Sigma)^{(0,1)}\otimes (T\Sigma)^{(1,0)}$ such that
\be
\label{I.15}
P_{h}\lf({F}\,\res \,h\rg)=0\ ,
\ee
Then there exists a unique ${U}^{0,1}_h$ section of $(T^{(0,1)}\Sigma)\otimes {\R}^n$ orthogonal (for the $L^2$ scalar product induced by $h$) to the finite dimensional vector space of anti-holomorphic vector-fields solution of
\be
\label{I.16}
\p^{h}{U}^{0,1}_h={F}\quad\mbox{ on }\quad(\Sigma,h)
\ee
Moreover for any $1<p<2$ one has
\be
\label{I.17}
\|{U^{0,1}_h}\|_{L^{p^\ast}_{h}(\Sigma)}\le C_{h,p}\ \|{F}\|_{L^{p}_{h}(\Sigma)}\
\ee
where $p^\ast$ is the Sobolev critical exponent given by $(p^\ast)^{-1}=p^{-1}-2^{-1}$ and we have also
\be
\label{I.18}
\|{U}^{0,1}_h\|_{L^{(2,\infty)}_{h}(\Sigma)}\le C_{h}\ \|{F}\|_{L^{1}_{h}(\Sigma)}\
\ee
\end{Lm}
The existence and uniqueness of the solution to (\ref{I.16}) is given in \cite{Ri1} (lemma 2.2). The estimates (\ref{I.17}) and (\ref{I.18}) are consequences of the following Lemma.
\begin{Lm}
\label{lm-potentiel}
Let $U^{0,1}_h$ be a section of $(T^{(0,1)}\Sigma)\otimes {\R}^n$ such that $U^{0,1}_h\res h$ is orthogonal (for the $L^2$ scalar product induced by $h$) to the finite dimensional vector space of holomorphic 1 forms.
Then there exists a unique complex function $\phi$ on $\Sigma$ with average 0 on each component of $\Sigma$ such that
\[
U^{0,1}_h=\p^h\phi\res h^{-1}
\]
The function $\phi$ is called ``potential'' from $U^{0,1}_h$. Moreover, let $F:=\p^h U^{0,1}_h$, we have
\be
\label{phi-F}
4^{-1}\, \Delta_h\lf[\Delta_h\phi+2\, K_h\,\phi\rg]=\ov{\p}^h\lf(h^{-1}\res_2\ov{\p}^h\lf[F\res h\rg]\rg)\res_2 h^{-1}\ .
\ee
where $K_h$ is the Gauss curvature of $h$. In particular the map which to $F$ assigns 2 derivatives of $\phi$ is a Calderon Zygmund operator.
\end{Lm}
\noindent{\bf Proof of lemma~\ref{lm-potentiel}.} The existence of $\phi$ such that
\[
U^{0,1}_h\res h=\p^h\phi\quad\Longleftrightarrow\quad U^{0,1}_h=\p^h\phi\res h^{-1}
\]
is a direct consequence of the fact that the image of $W^{1,2}(\Sigma,{\C})$ by $\p^h$ is equal to the $L^2-$orthogonal to the holomorphic 1-forms thanks to Dolbeault theorem : $H^{1,0}_{\p^h}(\Sigma,{\C})=\mbox{Ker}\,\p^h/\mbox{Im}\,\p^h=L^2(\wedge^{1,0}\Sigma)/\mbox{Im}\,\p^h$ is isomorphic to the space of holomorphic one forms.

In local conformal coordinates for $h$ such that $h=e^{2\nu}\,[dx_1^2+dx_2^2]$ one has
\[
2\,\p_{z}\lf[e^{-2\nu}\p_z\phi\rg] \ dz\otimes \p_{\ov{z}}=F=F_z^{\ov{z}} \ dz\otimes \p_{\ov{z}}\ .
\]
Hence
\[
\p^h\lf[\p^h\phi\res h^{-1}\rg]\res h=e^{2\nu}\,\p_{z}\lf[e^{-2\nu}\p_z\phi\rg] \ dz\otimes dz=F\res h=2^{-1}\,e^{2\nu}\, F_z^{\ov{z}} \ dz\otimes dz\ .
\]
This implies
\[
\ov{\p}^h\lf[\p^h\lf[\p^h\phi\res h^{-1}\rg]\res h\rg]=\p_{\ov{z}}\lf[e^{2\nu}\,\p_{z}\lf[e^{-2\nu}\p_z\phi\rg]\rg] \ d\ov{z}\otimes dz\otimes dz=\ov{\p}^h\lf[F\res h\rg]=2^{-1}\,\p_{\ov{z}}\lf[e^{2\nu}\, F_z^{\ov{z}}\rg] \ d\ov{z}\otimes dz\otimes dz\ .
\]
We deduce
\[
h^{-1}\res_2\ov{\p}^h\lf[\p^h\lf[\p^h\phi\res h^{-1}\rg]\res h\rg]=2\,e^{-2\nu}\,\p_{\ov{z}}\lf[e^{2\nu}\,\p_{z}\lf[e^{-2\nu}\p_z\phi\rg]\rg] \ dz=h^{-1}\res_2\ov{\p}^h\lf[F\res h\rg]=e^{-2\nu}\,\p_{\ov{z}}\lf[e^{2\nu}\, F_z^{\ov{z}}\rg] \ dz\ ,
\]
and then
\[
\begin{array}{l}
\ds\ov{\p}^h\lf(h^{-1}\res_2\ov{\p}^h\lf[\p^h\lf[\p^h\phi\res h^{-1}\rg]\res h\rg]\rg)=2\,\p_{\ov{z}}\lf[e^{-2\nu}\,\p_{\ov{z}}\lf[e^{2\nu}\,\p_{z}\lf[e^{-2\nu}\p_z\phi\rg]\rg]\rg]\ d\ov{z}\otimes dz\\[5mm]
\ds=\ov{\p}^h\lf(h^{-1}\res_2\ov{\p}^h\lf[F\res h\rg]\rg)=\p_{\ov{z}}\lf(e^{-2\nu}\,\p_{\ov{z}}\lf[e^{2\nu}\, F_z^{\ov{z}}\rg]\rg) \ d\ov{z}\otimes dz\ .
\end{array}
\]
Finally 
\[
\begin{array}{l}
\ds\ov{\p}^h\lf(h^{-1}\res_2\ov{\p}^h\lf[\p^h\lf[\p^h\phi\res h^{-1}\rg]\res h\rg]\rg)\res_2 h^{-1}=4\,e^{-2\nu}\,\p_{\ov{z}}\lf[e^{-2\nu}\,\p_{\ov{z}}\lf[e^{2\nu}\,\p_{z}\lf[e^{-2\nu}\p_z\phi\rg]\rg]\rg]\\[5mm]
\ds\ds=\ov{\p}^h\lf(h^{-1}\res_2\ov{\p}^h\lf[F\res h\rg]\rg)\res_2 h^{-1}=2 \,e^{-2\nu}\,  \p_{\ov{z}}\lf(e^{-2\nu}\,\p_{\ov{z}}\lf[e^{2\nu}\, F_z^{\ov{z}}\rg]\rg)\ .
\end{array}
\]
We compute
\[
\begin{array}{l}
\ds e^{-2\nu}\,\p_{\ov{z}}\lf[e^{-2\nu}\,\p_{\ov{z}}\lf[e^{2\nu}\,\p_{z}\lf[e^{-2\nu}\p_z\phi\rg]\rg]\rg]=e^{-2\nu}\,\p_{\ov{z}}\lf[e^{-2\nu}\,\p_{\ov{z}}\lf[\p^2_{z^2}\phi-2\,\p_z\nu\,\p_z\phi\rg]\rg]\\[5mm]
\ds = e^{-2\nu}\,\p_{\ov{z}}\lf[e^{-2\nu}\,\lf[ 4^{-1}\p_z\Delta\phi -2^{-1}\Delta\nu\,\p_z\phi-2^{-1}\p_z\nu\,\Delta\phi    \rg]\rg]\\[5mm]
\ds= e^{-2\nu}\,\p_{\ov{z}}\lf[4^{-1} \p_z\lf[e^{-2\nu}\,\Delta\phi\rg]+2^{-1}\ K_h\,\p_z\phi\rg]=16^{-1}\, \Delta_h\lf[\Delta_h\phi+2\, K_h\,\phi\rg]
\end{array}
\]
where we have used the Liouville equation $e^{-2\nu}\Delta \nu=-\,K_h$ and hence we have

\bigskip

In conformal coordinates the Weingarten operator  is given by
\[
\vec{\mathfrak h}^{\,0}:=\p_z\lf(e^{-2\la}\ \p_z\vec{\Phi}\rg)\ dz\otimes \p_{\ov{z}}\ .
\]
This is a section of $\wedge^{1,0}T^\ast\Sigma\otimes\wedge^{0,1}T\Sigma\otimes N_{\vec{\Phi}}\Sigma$ where $N_{\vec{\Phi}}\Sigma$ is the normal bundle to the immersion by $\vec{\Phi}$ of $\Sigma$. We have in particular
\[
\vec{\mathfrak h}^{\,0}\res g_{\vec{\Phi}}=\frac{1}{2}\,\pi_{\vec{n}}\lf( \p^2_{z^2}\vec{\Phi}\rg)\ dz\otimes dz\ .
\]
We will also denote
\[
\vec{\cal H}^0:=\vec{\mathfrak h}^{\,0}\res h_0=\frac{1}{2}\,e^{-2\al}\,\pi_{\vec{n}}\lf( \p^2_{z^2}\vec{\Phi}\rg)\ dz\otimes dz\ .
\]
Observe that, 
\[
\p^{h}\vec{\Phi}\res g^{-1}_{\vec{\Phi}}=2\,e^{-2\la}\ \p_z\vec{\Phi}\ \p_{\ov{z}}\ .
\]
Hence in particular
\[
\p^h\lf[\p^{h}\vec{\Phi}\res g^{-1}_{\vec{\Phi}}\rg]=2\, \vec{\mathfrak h}^{\,0}\ .
\]
Let $\vec{X}\in \Gamma((\vec{\Phi})_\ast T\Sigma)$. Then there exists a unique section $X\in \Gamma(T\Sigma)$ such that
\[
\vec{X}=\vec{\Phi}_\ast X\ .
\]
For any vector $X\in T\Sigma$ we denote respectively by $X^{1,0}_h$ and by $X^{0,1}_h$ the projection of $X$ on $T^{(1,0)}\Sigma$ respectively $T^{(0,1)}\Sigma$ for the structure induced by $h$. We write
\[
\vec{X}=X_{z}\ \p_z\vec{\Phi}+X_{\ov{z}}\ \p_{\ov{z}}\vec{\Phi}=\vec{\Phi}_\ast\lf(X_z\ \p_z+X_{\ov{z}}\ \p_{\ov{z}}\rg)\ .
\]
This gives
\be
\label{0-1part}
\vec{X}\cdot\p^{h}\vec{\Phi}\res g^{-1}_{\vec{\Phi}}= X_{\ov{z}}\ \p_{\ov{z}}=X^{0,1}_h\ ,
\ee
where we have used
\[
\p_z\vec{\Phi}\cdot \p_z\vec{\Phi}=0\quad\quad\mbox{ and }\quad\quad \p_z\vec{\Phi}\cdot \p_{\ov{z}}\vec{\Phi}=\frac{e^{2\la}}{2}\ .
\]
We are now proving the following result.
\begin{Lm}
\label{lm-conclu}
Let $h_t$ be a smooth path of constant Gauss curvature metrics such that there exists an holomorphic quadratic  form ${q}_0$ form of $(\Sigma,h_0)$ such that\footnote{The constant Gauss curvature metric is evolving ``parallelly'' that is orthogonally to the perturbations of the form ${\mathcal L}_Xh$ (see \cite{Tro}).}
\[
\frac{d h}{dt}(0)=\Re\lf[{q}_0\rg]\ .
\]
Let $\vec{\Phi}_t$ be a smooth path of immersions of $\Sigma$  into ${\R}^m$  such that $u_0$ is conformal from $(\Sigma,h_0)$ into ${\R}^m$. There exists a unique  ${X}\in \Gamma(T\Sigma)$ solution of
\be
\label{I.18-d}
\p^{h_0}\lf(X^{0,1}_{h_0}\rg)= \lf(I-\ti{P}_{h_0}\rg)\lf(\frac{d\vec{\Phi}_t}{dt}\cdot \vec{\frak h}_0\rg)\ .
\ee
where 
\[
\tilde{P}_{h_0}\lf(F\rg):= P_{h_0}\lf(  F\res h_0\rg)\res h_0^{-1}\ .
\]
such that
${X}^{0,1}_h$ is orthogonal (for the $L^2$ scalar product induced by $h$) to the finite dimensional vector space of holomorphic vector-fields 
Then, for any smooth path of diffeomorphisms $\Psi_t$ of $\Sigma$ such that $${X}:= \lf.\frac{d\Psi_t}{dt}\rg|_{t=0}\ .$$ the following holds
\be
\label{I.18-e}
\frac{d}{dt}\lf.\lf[\p^{h_t}(\vec{\Phi}_t\circ\Psi_t)\cdot\p^{h_t}(\vec{\Phi}_t\circ\Psi_t)\rg]\rg|_{t=0}=0\ .
\ee
Moreover
\be
\label{I.18-f}
{q}_0=P_{h_0}\lf[\frac{d\vec{\Phi}_t}{dt}\cdot\vec{\cal H}^0_0\rg]
\ee
\hfill $\Box$
\end{Lm}
\noindent{\bf Proof of Lemma~\ref{lm-conclu}} We establish the proof in the hyperbolic case and we consider
 $h_t$ to be a familly of constant Gauss  curvature metrics equal to $-1$  in such a way that
\be
\label{I.14-dd}
\frac{d h}{dt}=\Re {q}_t
\ee
where ${q}_t$ is an holomorphic quadratic form for $h_t$. We consider a conformal immersion $\vec{\Phi}$ from $(\Sigma, h_0)$ into ${\R}^m$ and consider a perturbation $\vec{\Phi}+t\,\vec{w}$ in such a way that
 \be
 \label{I.1}
 \vec{\Phi}+t\,\vec{w}\ :\ (\Sigma,h_t)\ \longrightarrow\ {\R}^m
 \ee
 is conformal. This gives
 \be
 \label{I.1a}
 \p^{h_t} \lf[\vec{\Phi}+t\,\vec{w}\rg]\,\dot{\otimes}\ \p^{h_t} \lf[\vec{\Phi}+t\,\vec{w}\rg]=0\ ,
 \ee
 where $\p^{h_t}$ the delbar operator associated to the conformal structure defined by $h^t$. We have
 \be
 \label{I.1b}
 \p^{h_t}\varphi:= \frac{1}{2}\lf[d\varphi + i\,\ast_t\, d\varphi\rg]
 \ee
 For any pair of functions $a$ and $b$ 
 \[
 d a\wedge \ast_{h_t} db= <da,db>_{h_t} \ dvol_{h_t}=h^{ij}_t\,\p_{x_i}a\,\p_{x_j}b\ \sqrt{h_{t,11} \,h_{t,22} -h^2_{t,12}}\ dx_1\wedge dx_2
 \]
 We have in local conformal coordinates for $h_0$ where $h_0=e^{2\nu_0}\ [dx^2_1+dx_2^2]$
 \be
 \label{I.1c}
  d a\wedge \frac{d\ast_{h_t}}{dt}(0) db=e^{2\nu_0}\,\frac{d h^{ij}_t}{dt}(0)\ \p_{x_i}a\,\p_{x_j}b\ dx_1\wedge dx_2+ \frac{1}{2}\,\delta_{ij}\ \p_{x_i}a\,\p_{x_j}b\ e^{-2\nu_0}\,\frac{d[h_{11}+h_{22}  ]}{dt}(0)\ dx_1\wedge dx_2
 \ee
 This gives in particular for $a=x_1$ and $b=x_1$
 \be
 \label{I.1dd}
 \begin{array}{l}
 \ds dx_1\wedge\frac{d\ast_{h_t}}{dt}(0) dx_1=e^{-2\nu_0}\,\lf[e^{4\nu_0}\,\frac{d h^{11}_t}{dt}(0)+\frac{1}{2}\,\frac{d[h_{11}+h_{22}  ]}{dt}(0)\rg]\ dx_1\wedge dx_2\ \\[5mm]
 \ds\quad=-\, e^{-2\nu_0}\, \frac{1}{2}\,\frac{d[h_{11}-h_{22}  ]}{dt}(0)\ dx_1\wedge dx_2\
 \end{array}
 \ee
and for $a=x_2$ and $b=x_1$
 \be
 \label{I.1d}
 \begin{array}{l}
 \ds dx_2\wedge\frac{d\ast_{h_t}}{dt}(0) dx_1=e^{2\nu_0}\,\frac{d h^{12}_t}{dt}(0)\  dx_1\wedge dx_2=-e^{-2\nu_0}\,\frac{d h_{12}}{dt}(0)\  dx_1\wedge dx_2
 \end{array}
 \ee
 Thus we deduce in one hand
 \be
 \label{I.1e}
 \frac{d\ast_{h_t}}{dt}(0) dx_1= e^{-2\nu_0}\ \frac{d h_{12}}{dt}(0)\  dx_1-\frac{e^{-2\nu_0}}{2}\ \frac{d(h_{11}-h_{22})}{dt}\ dx_2
 \ee
 Similarly in the other hand
  \be
 \label{I.1f}
 \frac{d\ast_{h_t}}{dt}(0) dx_2= -\frac{e^{-2\nu_0}}{2}\ \frac{d(h_{11}-h_{22})}{dt}\ dx_1-e^{-2\nu_0}\ \frac{d h_{12}}{dt}(0)\  dx_2
 \ee
 Taking the derivative at $t=0$ of (\ref{I.1a}) gives
  \be
 \label{I.1g}
 \p^{h_0}\vec{w}\,\dot{\otimes}\ \p^{h_0} \vec{\Phi}+\p^{h_0} \vec{\Phi}\,\dot{\otimes} \,     \p^{h_0}\vec{w}
 +\frac{i}{2}\frac{d\ast_{h_t}}{dt}(0)\, d\vec{\Phi}\,\dot{\otimes}\ \p^{h_0} \vec{\Phi}
 +\p^{h_0} \vec{\Phi}\,\dot{\otimes}\,\frac{i}{2}\frac{d\ast_{h_t}}{dt}(0)\, d\vec{\Phi}\,=0\ ,
 \ee
 We have in local conformal coordinates for $h_0$
 \be
 \label{I.1h}
 \begin{array}{l}
\ds \frac{i}{2}\frac{d\ast_{h_t}}{dt}(0)\, d\vec{\Phi}\,\dot{\otimes}\ \p^{h_0} \vec{\Phi}= \frac{i}{4}\frac{d\ast_{h_t}}{dt}(0)\, d\vec{\Phi}\,\dot{\otimes}\ (\p_{x_1}\vec{\Phi}-i\,\p_{x_2}\vec{\Phi})\ (dx_1+i\,dx_2)\\[5mm]
\ds = e^{2\la}\,\frac{i}{4}\frac{d\ast_{h_t}}{dt}(0) dx_1\otimes(dx_1+i\,dx_2)- e^{2\la}\,\frac{i}{4}\frac{d\ast_{h_t}}{dt}(0) i\,dx_2\otimes(dx_1+i\,dx_2)\\[5mm]
\ds = e^{-2\nu_0}\,e^{2\la}\,\frac{i}{4}\lf[\frac{d h_{12}}{dt}(0)\  dx_1\otimes dz-\frac{1}{2}\,\frac{d(h_{11}-h_{22})}{dt}\ dx_2\otimes dz + \,\frac{i}{2}\,  \frac{d(h_{11}-h_{22})}{dt}\ dx_1\otimes dz +i\, \frac{d h_{12}}{dt}(0)\  dx_2\otimes dz\rg]\\[5mm]
 \ds =-\frac{e^{2\la}}{8}\ e^{-2\nu_0}\, \frac{d}{dt}\lf[ h_{11}-h_{22}-2\, i\ h_{12}\rg]\ dz\otimes dz\ ,
 \end{array}
 \ee
 and we have similarly
 \be
 \label{I.1ah}
 \p^{h_0} \vec{\Phi}\ \dot{\otimes}\ \frac{i}{2}\frac{d\ast_{h_t}}{dt}(0)\, d\vec{\Phi}=-\frac{e^{2\la}}{8}\ e^{-2\nu_0}\, \frac{d}{dt}\lf[ h_{11}-h_{22}-2\, i\ h_{12}\rg]\ dz\otimes dz\ .
 \ee
 Combining (\ref{I.1g}) and (\ref{I.1h}) gives then
 \be
 \label{I.1i}
   \ds\p^{h_0}\vec{w}\dot{\otimes}\ \p^{h_0} \vec{\Phi}+\p^{h_0} \vec{\Phi}\,\dot{\otimes} \,     \p^{h_0}\vec{w} =\frac{e^{2\al}}{4}\ \frac{d}{dt}\lf[ h_{11}-h_{22}-2\, i\ h_{12}\rg]\ dz\otimes dz\ .
 \ee
Since the familly of constant Gauss curvature metric is satisfying
\be
\label{I.14-d}
\frac{d h}{dt}=\Re {q}_t
\ee
where ${q}_t$ is an holomorphic quadratic form for $h_t$. This gives in particular at $t=0$ in conformal coordinates for $h_0$
\be
\label{I.14-e}
\begin{array}{l}
\ds\frac{d h}{dt}= \Re \lf[(a+ib) dz^{\otimes^2}\rg]= a\, \lf[dx_1\otimes dx_1-dx_2\otimes dx_2\rg]-b\, \lf[dx_1\otimes dx_2+dx_2\otimes dx_1\rg]\\[5mm]
\ds \quad\mbox{ and }\quad \p_{\ov{z}}(a+i\,b)=0
\end{array}
\ee
Hence
\be
\label{I.14-f}
\frac{dh_{11}}{dt}=-\frac{dh_{22}}{dt}\quad\mbox{ and }\quad\p_{\ov{z}}\lf( \frac{d}{dt}\lf[ h_{11}-h_{22}-2\, i\ h_{12}\rg]  \rg)=0\ .
\ee
and 
 \be
 \label{I.3e}
q= \frac{1}{4} \frac{d[h_{11}-h_{22}-\,2\, i\, h_{12}]}{dt}  \ dz\otimes dz
 \ee
 defines an holomorphic quadratic form of $(\Sigma,h_0)$.  This gives 
\be
\label{I.7}
\begin{array}{l}
\p\vec{\Phi}\,\dot{\otimes}\,\p\vec{w}\,\res g^{-1}_{\vec{\Phi}}+\p\vec{w}\,\dot{\otimes}\,\p\vec{\Phi}\,\res g^{-1}_{\vec{\Phi}}\\[5mm]
\ds=\frac{1}{2}\ e^{-2\nu_0}\, \frac{d}{dt}\lf[ h_{11}-h_{22}-2\, i\ h_{12}\rg]\ \ dz\otimes \p_z
=q\res h_0^{-1}\ ,
\end{array}
\ee
where $q$ is the holomorphic quadratic form given by (\ref{I.3e}), where 
\be
\label{I.8}
g^{-1}_{\vec{\Phi}}:=2\,e^{-2\la}\ \lf[\p_{z}\otimes\p_{\ov{z}}+ \p_{\ov{z}} \otimes\p_{{z}} \rg]\quad\mbox{ and }\quad h^{-1}_0:=  2\,e^{- 2\,\nu_0}\ \lf[\p_{z}\otimes\p_{\ov{z}}+ \p_{\ov{z}} \otimes\p_{{z}} \rg]\ .
\ee
and $\res$ denotes the contraction operators between covariant and contravariant tensors. 

\medskip

Observe that
\be
\label{I.9}
\p\vec{\Phi}\,\dot{\otimes}\,\p\vec{w}\,\res g^{-1}_{\vec{\Phi}}= 2\,e^{-2\la}\ \p_z\vec{\Phi}\cdot\p_z\vec{w}\ dz\otimes \p_{\ov{z}}=2\,\p\lf(e^{-2\la}\ \p_z\vec{\Phi}\cdot\vec{w}\rg)\otimes \p_{\ov{z}}-2\,\p_z\lf(e^{-2\la}\ \p_z\vec{\Phi}\rg)\cdot\vec{w}\ \ dz\otimes \p_{\ov{z}}
\ee
Recall the definition of the ${\p}^{h_0}$ operator acting on $0-1$ vector-field for $(\Sigma,h_0)$ given in local coordinates by
\be
\label{I.9a}
\p^{h_0}\lf(a\ \p_{\ov{z}}\rg):=\p_za\ dz\otimes\p_{\ov{z}}
\ee
This definition is clearly independent of the choice of positive conformal coordinates because a change of such coordinates is holomorphic : $w(z)$ with $\ov{\p} w=0$ and we have in one hand
\[
\p_za\ dz\otimes\p_{\ov{z}}=\p_w a\ \ov{(\p_zw)}\ dw\otimes \p_{\ov{w}}=\p_w a\ \ov{(\p_wz)}^{-1}\ dw\otimes \p_{\ov{w}}=\p_w \lf(a\ \ov{(\p_wz)}^{-1}\rg)\ dw\otimes \p_{\ov{w}}\ ,
\]
and in the other hand
\[
a\ \p_{\ov{z}}=a\ \ov{(\p_wz)}^{-1}\ \p_{\ov{w}}\ .
\]
Observe that with this notation we have in one hand
\be
\label{I.10}
2\,\p\lf(e^{-2\la}\ \p_z\vec{\Phi}\cdot\vec{w}\rg)\otimes \p_{\ov{z}}=\p^{h_0}\lf( \p\vec{\Phi}\cdot \vec{w}\res  g^{-1}_{\vec{\Phi}}\rg)
\ee
where $\p\vec{\Phi}\cdot \vec{w}\res  g^{-1}_{\vec{\Phi}}$ is a section of $T^{0,1}(\Sigma,h_0)$. Observe that we have also
\be
\label{I.11}
\p\vec{\Phi}\res g^{-1}_{\vec{\Phi}}=2\ e^{-2\la}\ \p_z\vec{\Phi}\ \p_{\ov{z}}\ ,
\ee
which gives
\be
\label{I.12}
\p\lf(\p\vec{\Phi}\res g^{-1}_{\vec{\Phi}}\rg)=2\ \p_z\lf(e^{-2\la}\ \p_z\vec{\Phi}\rg)\ \ dz\otimes \p_{\ov{z}}=2\,\vec{\frak h}^{\,0}\ ,\ .
\ee
%The Weingarten operator  is given by
%\be
%\label{I.13}
%\vec{\mathfrak h}^{\,0}:=\p_z\lf(e^{-2\la}\ \p_z\vec{\Phi}\rg)\ dz\otimes \p_{\ov{z}}
%\ee
%which is a section of $\wedge^{1,0}T^\ast\Sigma\otimes\wedge^{0,1}T\Sigma\otimes N_{\vec{\Phi}}\Sigma$ where $N_{\vec{\Phi}}\Sigma$ is the normal bundle to the immersion by $\vec{\Phi}$ of $\Sigma$. 
Combining (\ref{I.7}) with (\ref{I.9}) and the following notations we obtain finally
\be
\label{I.14}
\p^{h_0}\lf( \p\vec{\Phi}\cdot \vec{w}\res  g^{-1}_{\vec{\Phi}}\rg)=2\,\vec{\mathfrak h}^{\,0}\cdot\vec{w}\,+\,\frac{q}{2}\res h_0^{-1}\ ,
\ee
We shall now use the following Lemma
\begin{Lm}
\label{lm-delbar-image}
The image of $\Gamma\lf(T^{(1,0)}\Sigma\rg)$ by $\p^{h}$ is exactly the $L^2_{h}$ orthogonal to the contraction with $h^{-1}$ of the space of holomorphic quadratic form.
\hfill $\Box$
\end{Lm}
\noindent{\bf Proof of lemma~\ref{lm-delbar-image}.} Let $A\in \Gamma\lf(T^{(1,0)}\Sigma\rg)$, we have (see the computations in the first lines of the proof of lemma 2.2 of \cite{Ri1}) writing in local coordinates $A=a\ \p_{\ov{z}}$, let $B=b\ dz\otimes dz$ and $h=2^{-1}\,e^{2\nu}\ (dz\otimes d\ov{z}+ d\ov{z}\otimes dz)$
\be
\label{18-aa}
\begin{array}{l}
\ds\lf<\p^{h}A \,,\, B\res h^{-1}\rg>_{h}\ dvol_{h}= 2^{-1}\,\lf<\p_z a\ dz\otimes \p_{\ov{z}}   , e^{-2\nu_0}\,b \ dz\otimes\p_{\ov{z}}  \rg>_{h}\ dvol_{h}\\[5mm]
\ds\ =2^{-1}\,\frac{i}{4}\,\lf[ \p_z a\,\ov{b}+  \p_{\ov{z}} \ov{a}\,{b}  \rg]\ dz\wedge d\ov{z}
\end{array}
\ee
Thus
\be
\label{18-ab}
\begin{array}{l}
\ds\lf<\p^{h}A,  q\res h^{-1}\rg>_{h}\ dvol_{h}=\frac{i}{4} \, d\lf( A\res\ov{B}+\ov{A}\res{B}\rg)-\frac{i}{4}\ \lf[ a\,\p_z\ov{b}+\ov{a}\,\p_{\ov{z}}b\rg]\,dz\wedge d\ov{z}\\[5mm]
\ds=\frac{i}{4} \, d\lf( A\res\ov{B}+\ov{A}\res{B}\rg)-\lf<A,  \lf(((\ov{\p}^{h}B)\res_2 h^{-1}\rg)\res h^{-1}\rg>_{h}\ dvol_{h}\ .
\end{array}
\ee
This finally gives
\be
\label{I.18a}
\forall A\in \Gamma((T\Sigma)^{1,0})\quad\quad\int_\Sigma\lf<\p^{h_0}A,  B\res h_0^{-1}\rg>_{h_0}\ dvol_{h_0}=-\int_\Sigma\lf<A,  \lf(((\ov{\p}^{h_0}B)\res_2 h_0^{-1}\rg)\res h_0^{-1}\rg>_{h_0}\ dvol_{h_0}=0
\ee
Hence $B\res h_0^{-1}$ is $L^2_h$ orthogonal to the image by $\p^h$ if and only if $B$ is holomorphic. This concludes the proof of lemma~\ref{lm-delbar-image}.\hfill $\Box$

\bigskip

Applying lemma~\ref{lm-delbar-image} to (\ref{I.14}) gives then in particular
\be
\label{I.18b}
 0=P_{h_0}\lf(2\,\vec{w}\cdot\vec{\frak{h}}_0\res h_0+2^{-1}\,q\rg)=2\,P_{h_0}\lf(\vec{w}\cdot\vec{\frak{h}}_0\res h_0\rg)+2^{-1}\,q\ .
\ee
This concludes the proof of lemma~\ref{lm-conclu}.\hfill $\Box$

\subsection{The Parametric Willmore Flow System for general Surfaces.}

The first derivative of the Willmore Functional in ${\R}^m$ has been written in conservative form in \cite{Riv-w}
\[
\lf.\frac{d}{dt}\int_\Sigma|\vec{H}_{\vec{\Phi}+t\vec{w}}|^2\ dvol_{g_{\vec{\Phi}+t\vec{w}}}\rg|_{t=0}=\int_\Sigma\vec{w}\cdot d^\ast\lf[d\vec{H}-3\,\pi_{\vec{n}}(d\vec{H})-\star(\vec{H}\wedge \ast d\vec{n})  \rg]\ dvol_{g_{\vec{\Phi}}}
\]
where $\star$ is the standard Hodge Operator on multi-vectors in ${\R}^m$ which to a $m-p$ $\vec{A}$ vector assigns a $p$ vector $\star\vec{A}$ such that
\[
\forall\ \vec{B}\in \wedge^{m-p}{\R}^m\quad\quad \vec{B}\wedge\star\vec{A}=\lf<\vec{B},\vec{A}\rg>\ \star 1
\]
where $<\cdot,\cdot>$ is the standard scalar product in $\wedge^{m-p}{\R}^m$. It is proved in \cite{Riv-notes} that for any immersion $\vec{\Phi}$ there holds
\[
\pi_{\vec{n}}\lf[  d^\ast\lf[d\vec{H}-3\,\pi_{\vec{n}}(d\vec{H})-\star(\vec{H}\wedge \ast d\vec{n})  \rg]  \rg]=d^\ast\lf[d\vec{H}-3\,\pi_{\vec{n}}(d\vec{H})-\star(\vec{H}\wedge \ast d\vec{n})  \rg]
\]
In conformal coordinates there holds (see \cite{Riv-notes} lemma X.3)
\[
\delta\vec{W}:=d^\ast\lf[d\vec{H}-3\,\pi_{\vec{n}}(d\vec{H})-\star(\vec{H}\wedge \ast d\vec{n})  \rg]=-4\, e^{-2\la}\, \Re\lf(  \p_{\ov{z}}\lf[\pi_{\vec{n}}(\p_z\vec{H})+\vec{H}\cdot\vec{H}^0\ \p_{\ov{z}}\vec{\Phi}  \rg]\rg)
\]
where
\[
\vec{H}^0:= 2\, e^{-2\la}\ \pi_{\vec{n}}(\p_z\p_z\vec{\Phi})
\]
We have
\[
\begin{array}{rl}
\ds\pi_{\vec{n}}(\p_z\vec{H})&=\p_{z}\vec{H}-\pi_T(\p_z\vec{H})=\p_{z}\vec{H}-2\, \p_z\vec{H}\cdot\p_z\vec{\Phi}\ e^{-2\la}\,\p_{\ov{z}}\vec{\Phi}-2\, \p_z\vec{H}\cdot\p_{\ov{z}}\vec{\Phi}\ e^{-2\la}\,\p_{{z}}\vec{\Phi}\\[5mm]
 &\ds=\p_{z}\vec{H}+ \vec{H}\cdot\vec{H}^0\,\p_{\ov{z}}\vec{\Phi}+|\vec{H}|^2\,\p_{{z}}\vec{\Phi}
\end{array}
\]
Hence we finally obtain
\be
\label{delW}
\delta\vec{W}:=-4\, e^{-2\la}\, \Re\lf(  \p_{\ov{z}}\lf[\p_z\vec{H}+|\vec{H}|^2\,\p_{{z}}\vec{\Phi}+2\,\vec{H}\cdot\vec{H}^0\ \p_{\ov{z}}\vec{\Phi}  \rg]\rg)
\ee
and in intrinsic notations this gives
\be
\label{delW-int}
\delta\vec{W}\ dvol_{g_{\vec{\Phi}}}:=-4\  \Re\lf( \ov{\p}\lf[\p\vec{H}+|\vec{H}|^2\,\p\vec{\Phi}+2\,\vec{H}\cdot\vec{\frak{h}}^{\,0}\res \ov{\p}\vec{\Phi}  \rg]\rg)
\ee
Let $(\Sigma,h)$ be a Riemann surface (equipped with a CGC metric), let $A^{0,1}_h$ be a section of $\Gamma\lf(T^{(0,1)}\Sigma\rg)$ and let $\vec{\Phi}$ be a conformal immersion of $(\Sigma,h)$ into ${\R}^m$. To $A^{0,1}_h=a\, \p_{\ov{z}}$ we associate the section $A\in \Gamma\lf(T\Sigma\rg)$ given by
\[
A=A^{0,1}_h+\ov{A^{0,1}_h}:=a\, \p_{\ov{z}}+\ov{a}\, \p_{{z}}
\]
and $\vec{A}$ such that
\[
\vec{A}=\vec{\Phi}_\ast A\ .
\]
We recall the definition~\ref{df-param-conf-will} from the introduction.
\begin{Dfi}
\label{df-param-will}
A $C^1$ family $\vec{\Phi}_t$  into the space of immersions solution to the following system
\be
\label{I.25}
\lf\{
\begin{array}{l}
\ds\frac{\p \vec{\Phi}}{\p t} =\delta \vec{\mathcal W} +\vec{U}\\[5mm]
\p^h\lf( U^{0,1}_h  \rg)=2\, \lf(I-\ti{P}_{h}\rg)\lf(  \delta\vec{W}   \cdot \vec{\frak h}_0\rg)\\[5mm]
\ds\frac{d h}{dt}=-\,4\,\Re \lf[P_{h}\lf(\delta\vec{W}\cdot \vec{\cal H}^0\rg)\rg]\\[5mm]
\ds \delta\vec{W}=d^{\ast_h}\lf[d\vec{H}-3\,\pi_{\vec{n}}(d\vec{H})+\star(\vec{H}\wedge d\vec{n})  \rg]
\end{array}
\rg.
\ee
is called  solution to the Parametric Willmore Flow.\hfill $\Box$
\end{Dfi}
\begin{Lm}
\label{lm-conf-cond} Let $\vec{\Phi}_t$ be a smooth solution to (\ref{I.25}) then it satisfies
\be
\label{I.26}
\p^{h_t}\vec{\Phi}_t\cdot\p^{h_t}\vec{\Phi}_t\equiv 0\quad .
\ee
\hfill $\Box$
\end{Lm}
Lemma~\ref{lm-conf-cond} is a direct consequence of lemma~\ref{lm-conclu}.
\subsection{Estimating the $W^{1,1}$ norm of the Green Kernel on degenerating Riemann Surfaces.}
We now establish the following lemma
\begin{Lm}
\label{lm-L1-estim}
Let $l_\ast>0$ be the length of the shortest closed geodesic of an hyperbolic surface $(\Sigma,h)$. There exists a constant $C_{\Sigma}>0$ depending only $\Sigma$ such that
\be
\label{est-green}
\sup_{y\in \Sigma}\int_\Sigma |d_xG_h(x,y)|_h\ dvol_h\le\ \frac{C_{\Sigma}}{l_\ast}\ \quad\mbox{ and }\quad\ \sup_{y\in \Sigma}\int_\Sigma |G_h(x,y)|_h\ dvol_h\le\ \frac{C_{\Sigma}}{l^2_\ast}\ .
\ee
where $G_h(x,y)$ is the Green function associated to the positive Laplace Beltrami operator on $\Sigma$ solving\footnote{We recall that we have adopted the convention that $\Delta_h$ is the negative Laplace Beltrami. The distribution $\delta^h_{x=y}$ is ``acting'' on 2-forms and is defined by 
\[
\forall f\in C^\infty(\Sigma)\quad<\delta^h_{x=y},f\ dvol_h>=f(x)
\] In particular in local conformal coordinates for which $dvol_h=e^{2\nu}\ dx_1\wedge dx_2$ one has
\[
\delta^h_{x=y}= e^{-2\nu(x_1,x_2)} \delta_{(x_1,x_2)=(y_1,y_2)}
\]}
\be
\label{green-def}
-\Delta^h_xG_h(x,y)=\delta^h_{x=y}-\frac{1}{\int_{\Sigma}dvol_h}\  \quad\mbox{ and }\quad \int_\Sigma G_h(x,y)\ dvol_h=0
\ee
\hfill $\Box$
\end{Lm}
\begin{Rm}
\label{rm-Green}
Observe that the estimate (\ref{est-green}) is optimal in the following sense. We consider $\Sigma=S^1\times [-L,L]$ where we identify $S^1\times \{L\}$ and $S^1\times\{-L\}$ equipped with the flat metric
\[
h:=\frac{1}{2\pi\,L}\ (d\theta^2+dt^2)
\]
in such a way that $V_h=1$. Following \cite{LaRi} proof of proposition 1.1 we consider the average of the Green function $g(\theta,t):=(2\pi)^{-1}\,\int_0^{2\pi} G((\theta,t),(\phi,0))\ d\phi$. The computation gives
$|dg|_h(\theta,t)=\sqrt{2\pi L}/4\pi$. Observe that the length of the shortest geodesic $l_\ast=\sqrt{2\pi /L}$. Hence we have
\[
\int_{\Sigma}\lf|\dashint_0^{2\pi}d_xG(x,(\phi,0) \, d\phi\rg|_h\ dvol_h=\frac{1}{2\,l_\ast}\ .
\]
\hfill $\Box$
\end{Rm}
Before establishing this lemma we recall a pointwise lower bound for the Green Function of the positive Laplace Beltrami which can be deduced from the main arguments in \cite{ChLi}
\begin{Lm}
\label{lm-ptw-green}
There exists a universal constant $c_4>0$ such that
\be
\label{ptw-green}
\inf_{(x,y)\in \Sigma}G_h(x,y)\ge - c_4\, \frac{V_h^2}{l_\ast^4}\ .
\ee
\hfill $\Box$
\end{Lm}
\noindent{\bf Proof of Lemma~\ref{lm-ptw-green}.} Let $(\phi_i)_{i\in {\N}}$ be an $L^2_h(\Sigma)$ orthonormal Hilbert Basis of eigenfunctions for the positive Laplace Beltrami operator $-\Delta^h$ and denote by
$\la_i\ge 0$ the corresponding eigenvalues. The first eigenvalues $\la_0=0$ and the corresponding eigenfunction is constant on $\Sigma$. We have
\[
G_h(x,y):=\sum_{i>0}\frac{\phi_i(x)\ \phi_i(y)}{\la_i}
\]
where $\int_\Sigma\phi_i\ dvol_h=0$ for any $i>0$. The Heat Kernel associated to the positive Laplace Beltrami operator $-\Delta^h$ is given by
\[
H_h(x,y,t):=\sum_{i\ge 0} e^{-\,\la_i t}\ \phi_i(x)\ \phi_i(y)
\]
It is known that $H_h>0$ for $t>0$ moreover obviously
\[
G_h(x,y)=\int_{0}^{+\infty}G_h(x,y,t)\ dt\quad\mbox{ where }\quad G_h(x,y,t):=\sum_{i>0} e^{-\la_i t}\ \phi_i(x)\ \phi_i(y)=H_h(x,y,t)-\frac{1}{V_h}
\]
where $V_h=\int_\Sigma dvol_h$. From \cite{ChLi} (corollary 1) there exists a universal constant $C>0$ such that for any $i>0$
\[
\|\phi_i\|^2_{L^\infty(\Sigma)}\le C\, \frac{\la_i}{C_1(\Sigma,h)}\ ,
\]
where $C_1(\Sigma,h)>0$ is the interpolation constant such that
\[
\forall f\in W^{1,2}(\Sigma)\quad\lf(\int_\Sigma f^2\ dvol_h\rg)^2\le \frac{1}{C_1(\Sigma,h)}\ \lf(\int_\Sigma |f|\ dvol_h\rg)^2\  \int_\Sigma |df|^2_h\ dvol_h\ .
\]
We write
\[
\begin{array}{l}
\ds G_h(x,y)=\int_{0}^{T}G_h(x,y,t)\ dt+\int_{T}^{+\infty}G_h(x,y,t)\ dt\ge -\frac{T}{V_h}-\int_T^{+\infty} \sum_{i>0}\, e^{-\la_i t}\ \|\phi_i\|^2_{L^\infty(\Sigma)}\ dt\\[5mm]
\ds\quad\ge -\frac{T}{V_h}-\frac{C}{C_1(\Sigma,h)}\,\int_T^{+\infty} \sum_{i>0}\, \la_i\, e^{-\la_i t}\ dt\ge -\frac{T}{V_h}- \frac{C}{C_1(\Sigma,h)}\,\sum_{i>0}\, e^{-\la_i T}
\end{array}
\]
It is proved in \cite{ChLi} (2.9) that
\[
\sum_{i>0}\, e^{-\la_i T}\le c_2\,\frac{V_h}{T\,C_1(\Sigma,h)}\ ,
\]
where $c_2>0$ is universal. We choose $T>0$ such that
\[
\frac{V_h}{T\,C_1(\Sigma,h)}=1
\]
and we obtain the pointwise lower bound
\be
\label{lower-green}
G_h(x,y)\ge - \frac{c_3}{C_1(\Sigma,h)}
\ee
where $c_3>0$ is universal.

 It is proved in \cite{ChLi} that
\[
C_1(\Sigma,h)\ge c_0\, I^2(\Sigma,h)\ ,
\]
where $c_0>0$ is a universal constant and $I(\Sigma,h)$ is the isoperimetric constant of $(\Sigma,h)$ given by
\[
I(\Sigma,h):=\inf_{\Om\subset \Sigma}\,\frac{({\mathcal H}^1_h(\p\Om))^2}{\min\{|\Om|_h,|\Sigma\setminus\Om|_h\}}\ ,
\]
where $\Om$ is taken among all finite perimeter subsets of $\Sigma$, $|\Om|_h$ and $|\Sigma\setminus\Om|_h$ are the area  respectively of $\Om$ and $\Sigma\setminus\Om$ while ${\mathcal H}^1_h(\p\Om)$ is the perimeter of $\Om$ with respect to the distance induced by $h$. It is proved in \cite{Cro} (proposition 12) that
\[
I(\Sigma,h)\ge \frac{8}{V_h}\, \mbox{inj}_{rad}^2(\Sigma,h)\ ,
\]
where $\mbox{inj}_{rad}(\Sigma,h)$ is the injectivity radius of $(\Sigma,h)$. For an hyperbolic surface the injectivity radius is half the systole (see for instance \cite{Bav})
\[
\mbox{inj}_{rad}(\Sigma,h)=\frac{l_\ast}{2}\ .
\]
Combining the previous we obtain
\be
\label{lower-c1}
 {C_1(\Sigma,h)}\ge 4\, c_0\,\frac{l_\ast^4}{V_h^2}\ .
\ee
Combining (\ref{lower-green}) and (\ref{lower-c1}) is implying finally
\be
\label{lower-green-systo}
G_h(x,y)\ge - c_4\, \frac{V_h^2}{l_\ast^4}
\ee
where $c_4>0$ is a universal constant which implies lemma~\ref{lm-ptw-green}. \hfill $\Box$

\bigskip

\noindent{\bf Proof of lemma~\ref{lm-L1-estim}.} First we assume 
\be
\label{upper-l1-green}
\int_{\Sigma}|G_h|\ dvol_h\le\frac{c_\star}{l_\ast^2} 
\ee
where $c_\star>0$ is a constant depending only on $\Sigma$ which is going to be fixed later.
%\[
%\int_{x\in\Sigma} G_h(x,y)\ dvol_h=0
%\]
%we obtain in particular the upper bound
%\be
%\label{upper-l1-green}
%\int_{x\in\Sigma} |G_h(x,y)|\ dvol_h\le 2\,c_4\, \frac{V_h^3}{l_\ast^4}\ .
%\ee
%which also gives 
%\[
%\|\phi_i\|^2_{L^\infty(\Sigma)}\le c_1\, {\la_i}\,\frac{V_h^2}{l_\ast^4}\ .
%\]
Classical pointwise estimates on the gradient of $G_h$ (see \cite{Aub}) imply that, for a a fixed $h$, the sup in $y\in\Sigma$ of the $L^{2,\infty}$ norm of $d_xG(x,y)$ is in $L^{2,\infty}(\Sigma)$. For almost every $t\in {\R}$, integrating (\ref{green-def}) on the upper-level set $G_h(\cdot,y)>t$ gives
\be
\label{univ-borne-green}
\int_{\{G_h(\cdot,y)=t\}}|d_xG_h(\cdot,y)|\ dl_h=-\int_{\{G_h(\cdot,y)=t\}}\frac{\p G_h(\cdot,y)}{\p\nu}\ dl_h=1- \frac{\int_{\{G_h(\cdot,y)\ge t\}}\ dvol_h}{\int_{\Sigma}dvol_h}\le 1\ . 
\ee
Let $s>0$ to be fixed later. We have
\be
\label{perim}
\begin{array}{l}
\ds\int_{-s}^{s}{\mathcal H}^1_h(\{G_h(\cdot,y)=t\})\ dt= \int_{|G_h|<s}|d_xG_h(\cdot,y)|\ dvol_h\le V_h^{1/2}\, \lf[ \int_{|G_h|<s}|d_xG_h(\cdot,y)|^2\ dvol_h\rg]^{1/2}\\[5mm]
\ds\quad\le V_h^{1/2}\, \lf[\int_0^s\ dt\,\int_{\{G_h(\cdot,y)=t\}}|d_xG_h(\cdot,y)|\ dl_h\rg]^{1/2}\le  V_h^{1/2}\,\sqrt{s}
\end{array}
\ee
From (\ref{upper-l1-green}) we have
\[
\int_0^{s}\lf|\{G_h(\cdot,y)\ge t\}\rg|\ dt\le \frac{c_\star}{l_\ast^2} 
\]
We now choose $s=4\,V_h\, c_\star/l_\ast^2 $.  which gives the existence of $t_\ast\in [2\,V_h\, c_\star/l_\ast^2\, 4\, V_h\, c_\star/l_\ast^2]$ such that
\be
\label{gap}
\forall\ t\ge t_\ast\quad\quad\lf|\{G_h(\cdot,y)\ge t\}\rg|< \frac{V_h}{2}\ .
\ee
Because of (\ref{perim}) we have
\[
\frac{1}{t_\ast}\int_{t_\ast}^{2\,t_\ast}{\mathcal H}^1_h(\{G_h(\cdot,y)=t\})\ dt\le 4\,\frac{\sqrt{c_\star}}{t_\ast}\, \frac{V_h}{l_\ast}\le 2\,\frac{l_\ast}{\sqrt{c_\star}}\ .
\]
We will consider $c_\star$ such that $2\,{l_\ast}/{\sqrt{c_\star}}\le l_\ast/2=\mbox{inj}_{rad}(\Sigma,h)$. That is $c_\star>16$ ($c_\star$ will be fixed definitively later in the argument). Then, using the mean value theorem we have the existence of a regular value $t_{\ast\ast}$ for $G_h(\cdot,y)$ such that $t_{\ast\ast}\in [t_\ast, 2\,t_\ast]$ and
\be
\label{good-slice}
{\mathcal H}^1_h(\{G_h(\cdot,y)=t_{\ast\ast}\})\le \frac{l_\ast}{2}\ .
\ee
Denote by $(\gamma_i)_{i\in I}$ the connected components of $\{G_h(\cdot,y)=t_{\ast\ast}\}$. Because of (\ref{good-slice}) the length of each $\gamma_i$ is below the minimal length for realizing a non-trivial homotopy class and there exists a disc $\om_i\subset \Sigma$ such that
\[
\p \om_i=\gamma_i
\]
Since $\om_i$ is a disc, the isoperimetric inequality by Huber (see \cite{Hub}) gives
\[
4\pi\,|\om_i|_h\le\, {\mathcal H}_h^1(\gamma_i)^2
\]
The closure of these discs are distinct to each-other hence $\Sigma\setminus \cup_{i\in I}\ov{\om}_i$ is connected (this can be proved by induction on Card$(I)$). We have
\[
\forall x\in \,\Sigma\setminus \cup_{i\in I}\ov{\om}_i\quad G_h(x,y)>t_{\ast\ast}\quad\mbox{ or }\quad\forall x\in \,\Sigma\setminus \cup_{i\in I}\ov{\om}_i\quad G_h(x,y)<t_{\ast\ast}\ .
\]
We have
\[
\lf|\cup_{i\in I}\ov{\om}_i\rg|\le (4\pi)^{-1}\,\sum_{i\in I}{\mathcal H}_h^1(\gamma_i)^2\le  (4\pi)^{-1}\, \lf({\mathcal H}^1_h(\{G_h(\cdot,y)=t_{\ast\ast}\})\rg)^2\le \frac{l^2_\ast}{\pi}<\frac{V_h}{2}
\]
hence the first alternative is excluded by (\ref{gap}). Thus we have proved
\[
 \{ G_h(\cdot,y)>t_{\ast\ast}\}\subset \cup_{i\in I}\ov{\om}_i\ ,
\]
and
\[
\lf| \{ G_h(\cdot,y)>t_{\ast\ast}\}\rg|_h\le (4\pi)^{-1}\,\lf({\mathcal H}^1_h(\{G_h(\cdot,y)=t_{\ast\ast}\})\rg)^2
\]
As a consequence, for any regular value $t\ge t_{\ast\ast}$ for $G_h(\cdot,y)$, the union of disjoint circles realizing $\{G_h(\cdot,y)=t\}$ are bounding discs in $\Sigma$ and for the same reasons as above
for any $t>t_{\ast\ast}$ there holds
\be
\label{iso-clas}
\lf| \{ G_h(\cdot,y)>t\}\rg|_h\le (4\pi)^{-1}\,\lf({\mathcal H}^1_h(\{G_h(\cdot,y)=t\})\rg)^2
\ee
We do exactly the same for the domain $ \{ G_h(\cdot,y)<- t_{\ast}\}$  and we also obtain that\footnote{By applying jointly the mean value argument for the regions $[-2t_\ast,t_\ast]$ and $[t_\ast,2t_{\ast}]$ we can always choose $t_{\ast\ast}$ to be the same value for the two intervals.} for $t<-t_{\ast\ast}$ 
\be
\label{iso-clas-low}
\lf| \{ G_h(\cdot,y)<t\}\rg|_h\le (4\pi)^{-1}\,\lf({\mathcal H}^1_h(\{G_h(\cdot,y)=t\})\rg)^2
\ee
Using H\"older inequality we have
\[
\begin{array}{l}
\ds\int_{|G_h|>t_{\ast\ast}}{|d_xG_h(\cdot,y)|}\ dvol_h\le \lf[\int_{|G_h|>t_{\ast\ast}}{|G_h(\cdot,y)|^{3/2}} dvol_h\rg]^{1/2}\ \lf[\int_{|G_h|>t_{\ast\ast}}\frac{|d_xG_h(\cdot,y)|^2}{|G_h(\cdot,y)|^{3/2}} dvol_h\rg]^{1/2}\\[5mm]
\ds\quad\le C\   \lf[\int_{t_{\ast\ast}}^{+\infty} t^{1/2}\ | \{ |G_h|(\cdot,y)>t\} \ dt\rg]^{1/2}\  \lf[\int_{t_{\ast\ast}}^{+\infty} t^{-3/2}\ \int_{\{|G_h|(\cdot,y)=t\}}|d_xG_h(\cdot,y)|\ dl_h\ dt\rg]^{1/2}\\[5mm]
\ds\quad\le C\ \sup_{t>t_{\ast\ast}} t^{1/4}\ | \{ |G_h|(\cdot,y)>t\} |^{1/4}\ \lf[\int_{t_{\ast\ast}}^{+\infty}  | \{ G_h(\cdot,y)>t\} |^{1/2} \ dt\rg]^{1/2}\ t_{\ast\ast}^{-1/4}\\[5mm]
\ds\quad\le C\ \frac{\sqrt{l_\ast}}{V^{1/4}_h\ c_\star^{1/4}}\ \lf[\int_{|G_h|>t_{\ast\ast}}{|G_h(\cdot,y)|}\ dvol_h\rg]^{1/4} \ \lf[ \int_{t_{\ast\ast}}^{+\infty}{\mathcal H}^1_h(\{|G_h|(\cdot,y)=t\})\ dt\rg]^{1/2 }
\end{array}
\]
Hence we deduce
\[
\int_{|G_h|>t_{\ast\ast}}{|d_xG_h(\cdot,y)|}\ dvol_h\le C
\]
Combining this fact with the bound
\[
\int_{|G_h|<t_{\ast\ast}}{|d_xG_h(\cdot,y)|}\ dvol_h\le V_h^{1/2}\ \sqrt{t_{\ast\ast}}\le 2\, C\ \sqrt{c_\star}\,V_h\ l_\ast^{-1}
\]
keeping in mind that $V_h$ is a function depending on the topology of $\Sigma$ (thanks to Gauss Bonnet in the case $g(\Sigma)\ge 2$ and in the torus case we can fix $V_h=1$) we obtain the result 
under the assumption (\ref{upper-l1-green}).  

\medskip  

We assume 
\be
\label{lower-l1-green}
\int_{\Sigma}|G_h|\ dvol_h\ge\frac{c_\star}{l_\ast^2} \ .
\ee
Since
\be
\label{null}
\int_{x\in\Sigma} G_h(x,y)\ dvol_h=0
\ee
we have
\be
\label{low-neg}
\int_{G_h<0}|G_h|\ dvol_h\ge\frac{c_\star}{2\,l_\ast^2} \ .
\ee
Because of (\ref{lower-green-systo}) we have
\be
\label{meas-neg}
\lf|\lf\{x\in\Sigma\ ;\ G_h(x,y)<0\rg\}\rg|_h\ge \frac{c_\star}{2\, c_4\, V_h} \ l_\ast^2 \ .
\ee
Let $t>0$ such that
\[
\lf|\lf\{x\in\Sigma\ ;\ G_h(x,y)>t\rg\}\rg|_h\ge 2\,V_h\, l^{-1}_\ast\ {\mathcal H}^1_h\lf(\lf\{x\in\Sigma\ ;\ G_h(x,y)=t\rg\}\rg)
\]
For such a $t$ we have
\[
{\mathcal H}^1_h\lf(\lf\{x\in\Sigma\ ;\ G_h(x,y)=t\rg\}\rg)<\frac{l_\ast}{2}
\]
Denote by $(\gamma_i)_{i\in I}$ the connected components of $\{G_h(\cdot,y)=t\}$. Because of the previous inequality the length of each $\gamma_i$ is below the minimal length for realizing a non-trivial homotopy class and there exists a disc $\om_i\subset \Sigma$ such that
\[
\p \om_i=\gamma_i
\]
Since $\om_i$ is a disc, the isoperimetric inequality by Huber (see \cite{Hub}) gives
\[
4\pi\,|\om_i|_h\le\, {\mathcal H}_h^1(\gamma_i)^2
\]
The closure of these discs are distinct to each-other hence $\Sigma\setminus \cup_{i\in I}\ov{\om}_i$ is connected (this can be proved by induction on Card$(I)$). We have
\[
\forall x\in \,\Sigma\setminus \cup_{i\in I}\ov{\om}_i\quad G_h(x,y)>t\quad\mbox{ or }\quad\forall x\in \,\Sigma\setminus \cup_{i\in I}\ov{\om}_i\quad G_h(x,y)<t\ .
\]
We have
\[
\lf|\cup_{i\in I}\ov{\om}_i\rg|\le (4\pi)^{-1}\,\sum_{i\in I}{\mathcal H}_h^1(\gamma_i)^2\le  (4\pi)^{-1}\, \lf({\mathcal H}^1_h(\{G_h(x,y)=t\})\rg)^2\le \frac{l^2_\ast}{\pi}
\]
We now fix $c_\star$ definitively such that 
\[
\frac{c_\star}{2\, c_4\, V_h} >\frac{1}{\pi}\ .
\]
Combining this choice with (\ref{meas-neg}) we have
\[
\lf|\lf\{x\in\Sigma\ ;\ G_h(x,y)<t\rg\}\rg|_h>\lf|\cup_{i\in I}\ov{\om}_i\rg|
\]
hence the first alternative is excluded and we deduce that under the assumption (\ref{lower-l1-green}) there holds
\[
\forall t>0\quad\quad \lf|\lf\{x\in\Sigma\ ;\ G_h(x,y)>t\rg\}\rg|_h\le 2\,V_h\, l^{-1}_\ast\ {\mathcal H}^1_h\lf(\lf\{x\in\Sigma\ ;\ G_h(x,y)=t\rg\}\rg)
\]
Integrating over $t\in[0,+\infty)$ gives
\[
\int_{G_h>0}|G_h(\cdot,y)|\ dvol_h\le 2\,V_h\, l^{-1}_\ast\ \int_{G_h>0}|d_xG_h(\cdot,y)|\ dvol_h\ .
\]
Using (\ref{null}) we deduce
\be
\label{poinc}
\int_{\Sigma}|G_h(\cdot,y)|\ dvol_h\le 4\,V_h\, l^{-1}_\ast\ \int_{\Sigma}|d_xG_h(\cdot,y)|\ dvol_h\ .
\ee
Let
\[
t_\ast:= 8\, l^{-1}_\ast\ \int_{\Sigma}|d_xG_h(\cdot,y)|\ dvol_h
\]
Because of (\ref{poinc}) we have
\[
\int_{0}^{t_\ast} \lf| \{ G_h(\cdot,y)>t\} \rg|_h\ dt\le \frac{V_h}{2}
\]
Hence 
\[
\forall t>t_\ast\quad\quad\lf|\{G_h(\cdot,y)\ge t\}\rg|< \frac{V_h}{2}
\]
We have
\[
\frac{1}{t_\ast}\int_{t_\ast}^{2\,t_\ast}{\mathcal H}^1_h(\{|G_h|(\cdot,y)=t\})\ dt\le t_\ast^{-1}\int_{\Sigma}|d_xG_h(\cdot,y)|\ dvol_h=\frac{l_\ast}{8}\ .
\]
there exists $t_{\ast\ast}\in [t_\ast,2\,t_\ast]$ such that
\[
{\mathcal H}^1_h(\{|G_h|(\cdot,y)=t_{\ast\ast}\})\le \frac{l_\ast}{8}
\]
We argue exactly as above and we finally obtain
\[
\begin{array}{l}
\ds\int_{|G_h|>t_{\ast\ast}}{|d_xG_h(\cdot,y)|}\ dvol_h\le \lf[\int_{|G_h|>t_{\ast\ast}}{|G_h(\cdot,y)|^{3/2}} dvol_h\rg]^{1/2}\ \lf[\int_{|G_h|>t_{\ast\ast}}\frac{|d_xG_h(\cdot,y)|^2}{|G_h(\cdot,y)|^{3/2}} dvol_h\rg]^{1/2}\\[5mm]
\ds\quad\le C\   \lf[\int_{t_{\ast\ast}}^{+\infty} t^{1/2}\ | \{ |G_h|(\cdot,y)>t\} \ dt\rg]^{1/2}\  \lf[\int_{t_{\ast\ast}}^{+\infty} t^{-3/2}\ \int_{\{|G_h|(\cdot,y)=t\}}|d_xG_h(\cdot,y)|\ dl_h\ dt\rg]^{1/2}\\[5mm]
\ds\quad\le C\ \sup_{t>t_{\ast\ast}} t^{1/4}\ | \{ |G_h|(\cdot,y)>t\} |^{1/4}\ \lf[\int_{t_{\ast\ast}}^{+\infty}  | \{ G_h(\cdot,y)>t\} |^{1/2} \ dt\rg]^{1/2}\ t_{\ast\ast}^{-1/4}\\[5mm]
\ds\quad\le C\  \lf[\int_{|G_h|>t_{\ast\ast}}{|G_h(\cdot,y)|}\ dvol_h\rg]^{1/4} \ \lf[ \int_{t_{\ast\ast}}^{+\infty}{\mathcal H}^1_h(\{|G_h|(\cdot,y)=t\})\ dt\rg]^{1/2 }\ l_\ast^{1/4}\ \lf[\int_{\Sigma}{|d_xG_h(\cdot,y)|}\ dvol_h\rg]^{-1/4}
\end{array}
\]
Using (\ref{perim}) and (\ref{poinc}) we have
\[
\begin{array}{l}
\ds\int_{\Sigma}{|d_xG_h(\cdot,y)|}\ dvol_h\le C\  \lf[ \int_{\Sigma}{|d_xG_h(\cdot,y)|}\ dvol_h\rg]^{1/2 }+V_h^{1/2}\, l_\ast^{-1}
\end{array}
\]
This concludes the proof of lemma~\ref{lm-L1-estim} in all cases.\hfill $\Box$

\subsection{Explicit constants for the inversion of $(\p^h)^{-1}$ on $L^1_h((T^\ast \Sigma)^{1,0}\otimes (T\Sigma)^{0,1})$}

\begin{Lm}
\label{lm-delbar-inv} Assume the genus of $\Sigma$  is strictly larger than one. There exists a constant $C_\Sigma>0$ depending only on the topology of $\Sigma$ such that for any
hyperbolic metric $h$ on $\Sigma$ such that $V_h=1$ for $g(\Sigma)=1$ and $V_h=4\pi\,(g(\Sigma)-1$ for $g(\Sigma)>1$ and for any $W^{1,1}$-section $U^{0,1}$
of $T^{0,1}\Sigma$ there holds respectively 
\be
\label{delbar-inv}
\int_\Sigma|U^{0,1}|_{h}\ dvol_h\le \frac{C_\Sigma}{l_\ast^4}\ \int_\Sigma|\p^hU^{0,1}|_{h}\ dvol_h\ 
\ee
if $g(\Sigma)>1$ and
\be
\label{delbar-inv-g=1}
\int_\Sigma|U^{0,1}|_{h}\ dvol_h\le \frac{C_\Sigma}{l_\ast^3}\ \int_\Sigma|\p^hU^{0,1}|_{h}\ dvol_h\ 
\ee
if $g(\Sigma)=1$.
\hfill $\Box$
\end{Lm}
\begin{Rm}
\label{rm-delbar}
In the case $g(\Sigma)=1$ the result is optimal\footnote{We consider a rectangular torus ${\R}^2/{\Z}\oplus {ib}{\Z}$ equipped with the flat metric $h=b^{-1}\, [dx^2+dy^2]$. The length of the shortest geodesic is $\sqrt{b^{-1}}$.  Let $U^{0,1}=\min\{y,b-y\}\,\p_{\ov{z}}$.
We have $\int_\Sigma |U^{0,1}|_{h}\ dvol_h=2^{-3/2}b^{3/2}$ and $\int_\Sigma |U^{0,1}|_{h}\ dvol_h=1$.      } while we don't know if (\ref{delbar-inv}) is optimal.\hfill $\Box$ 
\end{Rm}
\noindent{\bf Proof of lemma~\ref{lm-delbar-inv}.} Assume the result wouldn't be true, then we would extract a sequence of degenerating hyperbolic metrics $h^k$ such that the thin-thick parts decomposition
of Keen and Randol applies : We denote by $(\gamma_i^k)_{i\in I}$ the finite set of simple closed geodesics with length $l^i(k)$ going to zero. Set
\[
l_\ast(k):=\inf_{i\in I}l_i(k)\ .
\]
For any $\arcsinh 1>\delta>0$ there exists a set of disjoints neighbourhoods ${\mathcal C}_{i}^k(\delta)$ in $\Sigma$ of $\gamma_i^k$ isometric (for the metric $h_k$) to $(-X_\delta(l_i^k),X_\delta(l_i^k))\times S^1$ equipped with the conformal metric
\[
(\rho^k(s))^2\ \lf[ds^2+d\theta^2\rg]
\]
where
\[
\rho^k_i(s):=\frac{l_i^k}{2\pi\,\cos\lf(\frac{l_i^k\, s}{2\pi}\rg)}\quad\mbox{ and }\quad X_\delta(l_i^k):=\frac{\pi^2}{l_i^k}-\frac{2\,\pi}{l_i^k}\,\arcsin\lf(  \frac{\sinh(l_k^i/2)}{ \sinh\delta} \rg)
\]
The annuli are called ``$\delta-$collars''. Moreover there exists a hyperbolic punctured Riemann surface $(S,\sigma)$ with punctures (i.e. closed Riemann Surface $S_0$, called ``nodal surface, minus finitely many points $S=S_0\setminus\{(a_i,b_i)\}$ such that $\sigma$ is a complete hyperbolic metric) and there exists a sequence of diffeomorphisms
\[
\Phi^k \ :\ \Sigma\setminus \cup_{i\in I} \gamma_i^k\ \longrightarrow\ S
\]
such that
\[
(\Phi^k)^\ast h^k\ \longrightarrow\ \sigma\quad\mbox{ in }C^p_{loc}(S)\quad\forall p\in {\N}\ .
\]
Finally, the union of the pre-image $(\Phi^k)^{-1}({\mathcal C}_{i}^k(\delta))$ converge to a basis of neighbourhood's of the punctures. The part of $\Sigma_{\delta,thick}:=\Sigma\setminus \cup_{i\in I} {\mathcal C}_{i}^k(\delta)$ is called the $\delta-$thick part of $(\Sigma,h_k)$. On this set the metric $h_k$ is converging strongly in any norm and the injectivity radius  of $(\Sigma_{\delta,thick}, h_k)$ is uniformly bounded from below  by a positive constant (depending on $\delta$). The union of the collars $\Sigma_{\delta,thin}:= \cup_{i\in I} {\mathcal C}_{i}^k(\delta)$ is called the $\delta-$thin part of $(\Sigma, h^k)$.

\medskip

Let $\eta^k_i$ be an $L^2-$orthonormal basis of the space of holomorphic (1-0)-forms on $(\Sigma,h^k)$ which is a space of dimension $g$. We introduce the {\it canonical } (1-1)-form given by
\[
\om^k_{can}:= \frac{i}{g}\sum_{i=1}^g  \eta^k_i\wedge \ov{\eta^k_i}\ .
\]
One proves that it is independent of the chosen orthonormal basis. It is a positive (1-1)-form\footnote{Thanks to Riemann Roch theorem all holomorphic $1-1$-forms cannot have a common zero.}. The associated
canonical (or Bergman metric) is given by
\[
h^k_{\it b}:=\frac{1}{g} \sum_{i=1}^g|\eta^k_i|^2_{h^k}\ h^k\ .
\]
In local complex coordinates, if $\eta^k_i=f_i^k(z)\ dz$, there holds
\[
h^k_{\it b}= \frac{1}{2\,g}\sum_{i=1}^g|f_i^k(z)|^2\ [dz^2+d\ov{z}^2]\ .
\]
From \cite{Mas} and \cite{HaJo} (proposition 3.2) we have that
\[
(u^k)^\ast h^k_{\it b} \quad\longrightarrow\quad\sigma_{\it b}\quad\quad\mbox{ in } C^p_{loc}(S)\quad\forall p\in {\N}\ ,
\]
where $\sigma_{\it b}$ is the Bergman metric of $S_0$. As a consequence, for any $\delta>0$ there exists $C_\delta>1$ such that for any $k$
\be
\label{berg-bound-thick}
\forall x\in \Sigma_{\delta,thick}\quad\quad C_\delta^{-1}\le \sum_{i=1}^g|\eta^k_i|^2_{h^k}(x)\le C_\delta
\ee
Let ${F}$ be a smooth section of $(T^\ast\Sigma)^{(0,1)}\otimes (T\Sigma)^{(1,0)}$ such that
\be
\label{I.15-bis}
P_{h^k}\lf({F}\,\res \,h^k\rg)=0\ ,
\ee
and let ${U}^{0,1}_{h^k}$ section of $(T^{(0,1)}\Sigma)\otimes {\R}^n$ solution\footnote{The uniqueness and existence of ${U}^{0,1}_{h^k}$  is given by lemma~\ref{lm-delbar} and is taking into account that the space
of anti-holomorphic vector-fields for $g>1$ is empty.}    of
\be
\label{I.16-bis}
\p^{h^k}{U}^{0,1}_{h^k}={F}\quad\mbox{ on }\quad(\Sigma,h)
\ee
We assume first that $ {U}^{0,1}_{h^k}\res h_k$ is $L^2-$orthogonal to the space of holomorphic one forms and then we have the existence of $\phi$ such that $ {U}^{0,1}_{h^k}=\p\phi\res h^{-1}$.    To any such a pair we consider
\[
F_i^k:=F\res\ov{\eta_i^k}= F_{\ov{z}}^z\, \ov{f_i^k(z)}\, dz\quad\mbox{ and }\quad  {U}^i_k:={U}^{0,1}_{h^k}\res \ov{\eta_i^k}= u_k\,\ov{f_i^k(z)}\ ,
\]
where ${U}^{0,1}_{h^k}=u_k\ \p_{\ov{z}}$. We omit to write explicitly the index $h^k$ for the various operators and we observe that the following equation holds
\[
\p{U}^i_k=F_i^k\quad\mbox{ on }\Sigma\ .
\]
This gives in particular that
\[
\Delta{U}^i_k=\om_{h^k}\cdot \ov{\p} F_i^k\quad\mbox{ on }\Sigma\ .
\]
Observe that in local complex coordinates
\[
U^i_k\, \om_{h^k}=(\p\phi\res h^{-1})\res \ov{\eta}_i\, \om_{h^k}=\p_z\varphi \, \ov{f_i(z)}\ \frac{1}{i}\ dz\wedge \ov{dz}=\frac{1}{i}\, \p\varphi\wedge \ov{\eta_i}
\]
Hence in particular there holds
\be
\label{u-null}
\int_\Sigma U^i_k\, \om_{h^k}=0\ .
\ee
%Using the notation introduced in the previous subsection for the Green function we have
%\[
%{U}^k_i(z)=\int_\Sigma G_{h^k}(z,w)\ \om_{h^k}\cdot \ov{\p} F_i^k\ \om_{h^k}=\int_\Sigma G_{h^k}(z,w)\ \ov{\p} F_i^k=-\int_\Sigma \ov{\p_w} G_{h^k}(z,w)\wedge F_i^k(w)
%\]
Observe that
\[
\int_\Sigma\ov{\eta_j^k}\wedge F_i^k=\int_\Sigma\ov{\eta_j^k}\wedge \p{U}^i_k=0\ ,
\]
Hence the function defined by
\be
\label{def-varphi}
\varphi_{ji}^k:=\int_\Sigma G_{h^k}(z,w)\ \ov{\eta_j^k}\wedge F_i^k\ .
\ee
is satisfying
\[
\p\ov{\p}\varphi_{ji}^k:=\ov{\eta_j^k}\wedge F_i^k=\p\lf( {U}^k_i(z) \,\ov{\eta_j^k}\rg)
\]
Thus
\[
\p\lf[ \ov{\p} \phi_{ji}^k- {U}^k_i(z) \,\ov{\eta_j^k}\rg]=0
\]
This implies the existence of $c^l_{ji}\in {\C}$ such that
\be
\label{u-dev}
\ov{\p} \varphi_{ji}^k- {U}^k_i(z) \,\ov{\eta_j^k}=\sum_{l=1}^g\, c^l_{ji}\,\ov{\eta_l^k}
\ee
and
\[
c^l_{ji}:=-\int_\Sigma {U}^k_i(z) \,\ov{\eta_j^k}\wedge \eta_l^k
\]
For $j\ne l$ we have 
\[
\int_\Sigma\ov{\eta_j^k}\wedge \eta_l^k=0\ .
\]
Hence there exists $\sigma^k_{j,l}$ such that
\[
\p\ov{\p}\sigma^k_{j,l}=\ov{\eta_j^k}\wedge \eta_l^k
\]
Then the coefficients $c^l_{ji}$ for $l\ne j$ are given by
\be
\label{exp-coeff}
c^l_{ji}:=\int_\Sigma F_i^k\wedge \ov{\p}\sigma_{jl}^k\ .
\ee
We have
\[
\ov{\p}\sigma^k_{jl}=\int_\Sigma \ov{\p}_z G_{h^k}(z,w)\ \ov{\eta_j^k}\wedge \eta_l^k(w)\ .
\]
Hence
\be
\label{est-sigma}
\begin{array}{rl}
\ds\||\p\sigma^k_{jl}|_{h^k}\|_{L^\infty(\Sigma)}&\ds\le \int_\Sigma |d_zG_{h^k}(z,w)|\ dvol_{h^k}\ \|\ov{\eta_j^k}\wedge \eta_l^k(w)\cdot\om_{h^k}\|_{L^\infty(\Sigma)}\\[5mm]
\ds\quad&\ds\le \frac{C_\Sigma}{l_\ast(k)}\ \max_{i=1\cdots g}\||\eta_i^k|_{h^k}\|^2_{L^\infty(\Sigma)}\ .
\end{array}
\ee
Combining (\ref{exp-coeff}) and (\ref{est-sigma}) is giving
\be
\label{coef-contr}
\begin{array}{rl}
\ds\forall l\ne j\quad\quad |c^l_{ji}|&\ds\le  \frac{C_\Sigma}{l_\ast(k)}\ \max_{i=1\cdots g}\||\eta_i^k|_{h^k}\|^2_{L^\infty(\Sigma)}\ \int_\Sigma|F^k_i|_{h^k}\ dvol_{h^k}\\[5mm]
 &\ds \le \frac{C_\Sigma}{l_\ast(k)}\ \max_{i=1\cdots g}\||\eta_i^k|_{h^k}\|^3_{L^\infty(\Sigma)}\ \int_\Sigma|F|_{h^k}\ dvol_{h^k}\ .
\end{array}
\ee
For $j=l$ we have thanks to (\ref{u-null})
\[
c^l_{li}:=-\int_\Sigma {U}^k_i(z) \,\ov{\eta_l^k}\wedge \eta_l^k=-\int_\Sigma {U}^k_i(z) \,\lf[\ov{\eta_l^k}\wedge \eta_l^k- V^{-1}\,{\om_{h^k}}\rg]
\]
where $V$ is the volume of $(\Sigma,h^k)$ and thanks to Gauss-Bonnet $V=4\pi\,g-4\pi$. Since $\int_\Sigma\ov{\eta_l^k}\wedge \eta_l^k-V^{-1}\,{\om_{h^k}}=0$ there exists $\hat{\sigma}_l^k\in C^\infty(\Sigma)$ such that
\[
\p\ov{\p}\hat{\sigma}_l^k=\ov{\eta_l^k}\wedge \eta_l^k-V^{-1}\,{\om_{h^k}}\ .
\]
We have
\[
\ov{\p}\hat{\sigma}^k_{l}=\int_\Sigma \ov{\p}_z G_{h^k}(z,w)\ \lf[\ov{\eta_l^k}\wedge \eta_l^k(w)-V^{-1}\,{\om_{h^k}}\rg]\ .
\]
Hence
\be
\label{est-sigma-j=l}
\begin{array}{rl}
\ds\||\p\hat{\sigma}^k_{l}|_{h^k}\|_{L^\infty(\Sigma)}&\ds\le \int_\Sigma |d_zG_{h^k}(z,w)|\ dvol_{h^k}\ \lf[V^{-1}+\|\ov{\eta_l^k}\wedge \eta_l^k(w)\cdot\om_{h^k}\|_{L^\infty(\Sigma)}\rg]\\[5mm]
\ds\quad&\ds\le \frac{C_\Sigma}{l_\ast(k)}\ \lf[   V^{-1}+\max_{i=1\cdots g}\||\eta_i^k|_{h^k}\|^2_{L^\infty(\Sigma)}\rg]\ .
\end{array}
\ee
We deduce similarly as before
\be
\label{coef-contr-bis}
\begin{array}{l}
\ds\forall l\ne j\quad\quad |c^l_{li}|\ds\le  \frac{C_\Sigma}{l_\ast(k)}\ \lf[   V^{-1}+\max_{i=1\cdots g}\||\eta_i^k|_{h^k}\|^2_{L^\infty(\Sigma)}\rg]\ \int_\Sigma|F^k_i|_{h^k}\ dvol_{h^k}\\[5mm]
 \ds\quad \le \frac{C_\Sigma}{l_\ast(k)}\ \lf[V^{-1}\,\max_{i=1\cdots g}\||\eta_i^k|_{h^k}\|_{L^\infty(\Sigma)}+\max_{i=1\cdots g}\||\eta_i^k|_{h^k}\|^3_{L^\infty(\Sigma)}\rg]\ \int_\Sigma|F|_{h^k}\ dvol_{h^k}\ .
\end{array}
\ee
Combining (\ref{u-dev}) with (\ref{def-varphi}), (\ref{coef-contr}) and (\ref{coef-contr-bis}) gives for any $i,j=1\cdots g$
\[
\begin{array}{l}
\ds\int_\Sigma|U^k_i|\,|\eta_j^k|_{h^k}\ dvol_{h^k}\le \int_\Sigma\sup_{w} |d_zG(z,w)|_{h^k}\ dvol_{h^k}\ \int_\Sigma|F^k_i|_{h^k}\, |\eta_j^k|_{h^k}\ dvol_{h^k}\\[5mm]
\ds\quad\quad+ \,g\,\max_{l=1\cdots g}|c^l_{ji}|\,\int_\Sigma|\eta_l^k|\,dvol_{h^k}
\end{array}
\]
Combining the above we obtain for any $i$ and $j$ between 1 and $g$
\[
\begin{array}{l}
\ds\int_\Sigma|U^{0,1}|_{h^k}\,|\eta_i^k|_{h^k}\,|\eta_j^k|_{h^k}\ dvol_{h^k}\\[5mm]
\ds\le\frac{C_\Sigma}{l_\ast(k)}\ \lf[\max_{i=1\cdots g}\||\eta_i^k|_{h^k}\|^2_{L^\infty(\Sigma)}+\max_{i=1\cdots g}\||\eta_i^k|_{h^k}\|^3_{L^\infty(\Sigma)}\rg]\ \int_\Sigma|F|_{h^k}\ dvol_{h^k}\ .
\end{array}
\]
This gives in particular
\be
\label{ui}
\begin{array}{l}
\ds\int_\Sigma|U^{0,1}|_{h^k}\, dvol_{h^k_{\it b}}
%\\[5mm]
\ds\le\frac{C_\Sigma}{l_\ast(k)}\ \lf[\max_{i=1\cdots g}\||\eta_i^k|_{h^k}\|^2_{L^\infty(\Sigma)}+\max_{i=1\cdots g}\||\eta_i^k|_{h^k}\|^3_{L^\infty(\Sigma)}\rg]\ \int_\Sigma|F|_{h^k}\ dvol_{h^k}\ .
\end{array}
\ee   
Assume now $U^{0,1}_{h^k}=\eta_i^k\res h^{-1}_k$ we have
\[
\int_\Sigma|U^{0,1}|^2_{h^k}\, dvol_{h^k}=1
\]
and
\[
\ov{\eta_i^k}\res\p\lf(\eta_i^k\res h^{-1}_k\rg)=\p\lf(  \ov{\eta_i^k}\res\eta_i^k\res h^{-1}_k \rg)=2\,\p|\eta_i^k|_{h^k}^2
\]
wedging with $\ov{\eta}_i^k$ gives finally
\[
2\, \p\lf[|\eta_i^k|_{h^k}^2\,\ov{\eta_i^k}\rg]=\lf[\ov{\eta_i^k}\res\p\lf(\eta_i^k\res h^{-1}_k\rg)\rg]\wedge\ov{\eta_i^k}
\]
Let $\xi_i^k$ of average zero on $\Sigma$ such that
\be
\label{equ-xi}
\p\ov{\p}\xi_i^k=2\, \p\lf[|\eta_i^k|_{h^k}^2\,\ov{\eta_i^k}\rg]
\ee
Hence
\[
\ov{\p}\xi_i^k=\int_\Sigma \ov{\p}_zG(z,w)\ \lf[\ov{\eta_i^k}\res\p\lf(\eta_i^k\res h^{-1}_k\rg)\rg]\wedge\ov{\eta_i^k}\ dz
\]
We have in particular
\be
\label{est-xi}
\begin{array}{l}
\ds\int_\Sigma|\ov{\p}\xi_i^k|_{h^k}\ dvol_{h^k}\le \int_\Sigma\max_{w}|\ov{\p}_zG(z,w)|\ \int_\Sigma|\eta_i^k|_{h^k}^2\, |\p\lf(\eta_i^k\res h^{-1}_k\rg)|_{h^k}\ dvol_{h^k}\\[5mm]
\ds\quad\le\ \frac{C_\Sigma}{l_\ast(k)}\ \||\eta_i^k|_{h^k}\|^2_{L^\infty(\Sigma)}\  \int_\Sigma |\p\lf(\eta_i^k\res h^{-1}_k\rg)|_{h^k}\ dvol_{h^k}
\end{array}
\ee
Because of (\ref{equ-xi}) there exists $d_{ij}^k\in{\C}$ such that
\[
|\eta_i^k|_{h^k}^2\,\ov{\eta_i^k}=\ov{\p}\xi_i^k+\sum_{j=1}^{g} d_{ij}^k\ \ov{\eta_j^k}\ .
\]
This implies in particular
\[
d_{ij}^k=\int_{\Sigma}|\eta_i^k|_{h^k}^2\,\ov{\eta_i^k}\wedge \eta_j^k
\]
For $i\ne j$ we recall the notation $\p\ov{\p}\sigma_{ij}=\ov{\eta_i}\wedge\eta_j$. Hence
\[
\forall\, i\ne j\quad\quad d_{ij}^k=\int_{\Sigma}|\eta_i^k|_{h^k}^2\,\p\ov{\p}\sigma_{ij}=-\, \int_{\Sigma}\p|\eta_i^k|_{h^k}^2\wedge\ov{\p}\sigma_{ij}\ .
\]
Thus, using (\ref{est-sigma})
\be
\label{est-dij}
\begin{array}{rl}
\ds\forall\, i\ne j\quad\quad |d_{ij}^k|&\ds\le \||\p{\sigma}^k_{ij}|_{h^k}\|_{L^\infty(\Sigma)}\ \int_\Sigma\lf|\p|\eta_i^k|_{h^k}^2\rg|\ dvol_{h^k} \\[5mm]
\ds\quad&\ds\le \frac{C_\Sigma}{l_\ast(k)}\ \max_{j=1\cdots g}\||\eta_j^k|_{h^k}\|^2_{L^\infty(\Sigma)}\ \int_\Sigma |\p\lf(\eta_i^k\res h^{-1}_k\rg)|_{h^k}\ dvol_{h^k}
\end{array}
\ee
We have moreover
\[
\begin{array}{l}
\ds d_{ii}^k=V^{-1}+\int_{\Sigma}|\eta_i^k|_{h^k}^2\,\ov{\eta_i^k}\wedge \eta_i^k-V^{-1}\om_{h^k}=V^{-1}+\int_{\Sigma}|\eta_i^k|_{h^k}^2\,\p\ov{\p}\hat{\sigma}_{i}\\[5mm]
\ds\quad =V^{-1}-\int_{\Sigma}\p|\eta_i^k|_{h^k}^2\wedge\ov{\p}\hat{\sigma}_{i}
\end{array}
\]
Using (\ref{est-sigma-j=l}) we have then
\be
\label{est-dii}
|d^k_{ii}-V^{-1}|\le  \frac{C_\Sigma}{l_\ast(k)}\ \lf[   V^{-1}+\max_{i=1\cdots g}\||\eta_i^k|_{h^k}\|^2_{L^\infty(\Sigma)}\rg]\ \int_\Sigma |\p\lf(\eta_i^k\res h^{-1}_k\rg)|_{h^k}\ dvol_{h^k}
\ee
In other words
\be
\label{lower-1}
\begin{array}{l}
\ds V^{-1}\lf[\int_\Sigma dvol_{h^k}\ \int_\Sigma|\eta_i^k|^4_{h^k}\ dvol_{h^k}-\lf(\int_\Sigma|\eta_i^k|^2_{h^k}\ dvol_{h^k}\rg)^2\rg] \\[5mm]
\ds\quad \le\frac{C_\Sigma}{l_\ast(k)}\ \lf[   V^{-1}+\max_{i=1\cdots g}\||\eta_i^k|_{h^k}\|^2_{L^\infty(\Sigma)}\rg]\ \int_\Sigma |\p\lf(\eta_i^k\res h^{-1}_k\rg)|_{h^k}\ dvol_{h^k}
\end{array}
\ee
We shall now be using the following lemma
\begin{Lm}
\label{lm-lower-1}
Under the previous notations we have
\be
\label{lower-2}
\liminf_{k\rightarrow +\infty}\min_{i=1\cdots g}\int_\Sigma dvol_{h^k}\ \int_\Sigma|\eta_i^k|^4_{h^k}\ dvol_{h^k}-\lf(\int_\Sigma|\eta_i^k|^2_{h^k}\ dvol_{h^k}\rg)^2>0
\ee
\end{Lm}
\noindent{\bf proof of lemma~\ref{lm-lower-1}.}
We argue by contradiction. Assume 
\[
\liminf_{k\rightarrow +\infty}\min_{i=1\cdots g}\,\int_\Sigma dvol_{h^k}\ \int_\Sigma|\eta_i^k|^4_{h^k}\ dvol_{h^k}-\lf(\int_\Sigma|\eta_i^k|^2_{h^k}\ dvol_{h^k}\rg)^2=0
\]
then there exists $i$ and a subsequence that we keep denoting $k$ such that
\[
\lim_{k\rightarrow +\infty}\int_\Sigma dvol_{h^k}\ \, \int_\Sigma|\eta_i^k|^4_{h^k}\ dvol_{h^k}-\lf(\int_\Sigma|\eta_i^k|^2_{h^k}\ dvol_{h^k}\rg)^2=0\ ,
\]
or in other words
\[
\lim_{k\rightarrow +\infty} \int_\Sigma\lf||\eta_i^k|^2_{h^k}-\dashint_\Sigma|\eta_i^k|^2_{h^k}\ dvol_{h^k}\rg|^2\  dvol_{h^k}=\lim_{k\rightarrow +\infty} \int_\Sigma\lf||\eta_i^k|^2_{h^k}-V^{-1}\rg|^2\  dvol_{h^k}\ .
\]
For any $\delta>0$ we have
\be
\label{squiz-1}
\lim_{k\rightarrow +\infty}\int_{\Sigma_{\delta,thick}} \lf||\eta_i^k|^2_{h^k}\circ\Phi^k-V^{-1}\rg|^2\  dvol_{(u^k)^\ast h^k}=0
\ee
Since $\sigma^k:=(\Phi^k)^\ast h^k$ is pre-compact in any $C^l$ norm, $(\Phi^k)^\ast \eta_i^k$ is converging strongly locally (in any $C^l_{loc}(\Sigma_{\delta,thick})$ topology) towards
a limiting holomorphic one form on $S$ with norm (w.r.t. $\sigma$) identically equal to $V^{-1}$. Hence, in particular we have
\be
\label{squiz-2}
\forall \delta>0\quad\exists\, k_\delta\in{\N}\quad\forall\,k\ge k_\delta\quad \inf_{x\in\Sigma_{\delta,thick}} |\eta_i^k|^2_{h^k}(x)>2^{-1}\, V^{-1}
\ee
Inside a given collar $(-X_\delta(l_i^k),X_\delta(l_i^k))\times S^1$ equipped with the conformal metric
\[
(\rho_i^k)^2(s)\ \lf[ds^2+d\theta^2\rg]\quad \mbox{ where }
\]
we write $\eta_i^k=f_i^k(z)\ dz$ where $z=s+i\theta$. hence $|\eta_i^k|^2_{h^k}=(\rho^k)^{-2}\,|f_i^k(z)|^2$ and we have because of (\ref{squiz-2}) for any $R>1$ 
\be
\label{conv-cons}
\lim_{k\rightarrow +\infty}\||f_i^k(z)|^2-V^{-1}\,(\rho^k(s))^2\|_{L^\infty(-X(l_i^k),-X(l_i^k)+R)\cup(X(l_i^k)-R,X(l_i^k))}=0\ .
\ee
The Fourier development of $f_i^k$ on $(-X(l_i^k),X(l_i^k))\times S^1$ gives the existence of $(a^k_n)_{n\in{\Z}}\in{\C}^{\Z}$ and $b^k$ such that
\[
f_i^k(z)=a^k_0+b^k_0\ s+\sum_{n\in {\Z}^\ast} a_n^k\, e^{i\,n\,\theta}\, e^{n\, s}\ .
\]
We have
\[
\dashint_0^{2\pi}|f_i^k(z)|^2\ d\theta=|a^k_0+b^k_0\ s|^2+\sum_{n\in {\Z}^\ast} |a_n^k|^2\, e^{2\,n\, s}
\]
Because of (\ref{conv-cons}), we have for $t\in[0,R]$
\be
\label{coll-est-1}
\lf\{
\begin{array}{l}
\ds\lim_{k\rightarrow +\infty}\lf\| |a^k_0+b^k_0\ (t-X(l_i^k))|^2+\sum_{n\in {\Z}^\ast} |a_n^k|^2\, e^{-2\,n\, X(l_i^k)}\, e^{2\,n\,t}-V^{-1}\,(\rho^k_i(s))^2\rg\|_{L^\infty([0,R])}= 0\\[8mm]
\ds\lim_{k\rightarrow +\infty}\lf\| |a^k_0+b^k_0\ (-t+X(l_i^k))|^2+\sum_{n\in {\Z}^\ast} |a_n^k|^2\, e^{2\,n\, X(l_i^k)}\, e^{-\,2\,n\,t}-V^{-1}\,(\rho^k_i(s))^2\rg\|_{L^\infty([0,R])}= 0
\end{array}
\rg.
\ee
For $n\ne0$ we denote $\hat{a}_{n,\pm}^k:= a_n^k e^{\pm\,n\, X(l_i^k)}$ and $\hat{a}_{0,\pm}^k:= a_0^k\pm X(l_i^k)$. With these notations we have
\be
\label{coll-est-2}
\lf\{
\begin{array}{l}
\ds\lim_{k\rightarrow +\infty}\lf\| |\hat{a}^k_{0,-}+b^k_0\ t|^2+\sum_{n\in {\Z}^\ast} |\hat{a}_{n,-}^k|^2\, e^{2\,n\,t}-V^{-1}\,(\rho^k_i(s))^2\rg\|_{L^\infty([0,R])}= 0\\[8mm]
\ds\lim_{k\rightarrow +\infty}\lf\| |\hat{a}^k_{0,+}-b^k_0\ t|^2+\sum_{n\in {\Z}^\ast} |\hat{a}_{n,+}^k|^2\, e^{-\,2\,n\,t}-V^{-1}\,(\rho^k_i(s))^2\rg\|_{L^\infty([0,R])}= 0
\end{array}
\rg.
\ee
Recall
\[
\rho^k_i(s):=\frac{l_i^k}{2\pi\,\cos\lf(\frac{l_i^k\, s}{2\pi}\rg)}\quad\mbox{ and }\quad X(l_i^k):=\frac{\pi^2}{l_i^k}-\frac{2\,\pi}{l_i^k}\,\arcsin\lf(  \frac{\sinh(l_k^i/2)}{ \sinh 1} \rg)
\]
Hence for $s=X_\delta(l_i^k)-t$ and $t\in [0,R]$ one has
\[
 X_\delta(l_i^k)=\frac{\pi^2}{l_i^k}-\frac{\pi}{\sinh 1}+O(l_k^i)\quad
 \]
and
\[ 
\rho^k_i(s)=\frac{l_i^k}{2\pi\,\cos\lf(\frac{\pi}{2} -\frac{l_i^k\, }{2\sinh 1}-l_i^k\, \frac{t}{2\pi}+O((l_i^k)^2)\rg)}=\frac{1}{\frac{\pi}{\sinh 1}+t+O(l_i^k)}
\]
Let $k$ large enough (i.e. $k\ge k_0$) such that
\[
 |\hat{a}^k_{0,-}+b^k_0\ t|^2+\sum_{n\in {\Z}^\ast} |\hat{a}_{n,-}^k|^2\, e^{2\,n\,t}\le 2 V^{-1}
\]
Consider $n>0$, we have
\[
 |\hat{a}_{n,+}^k|^2\,\le 2 V^{-1}
\]
from which we deduce
\[
 |a_n^k| \le 2 \,V^{-1}\,e^{-n\,X(l_i^k)}
\]
similarly, consider $n<0$, we have
\[
|\hat{a}_{n,-}^k|^2\,\le 2\, V^{-1}
\]
from which we deduce
\[
|a_n^k|\le 2\, V^{-1}\, \,e^{-|n|\,X(l_i^k)}
\]
Since $X(l_i^k)={\pi^2}/l_i^k+\pi/\sinh\,1+O(l_i^k)$ and $X(l_i^k)\ge X_\delta(l_i^k)+ \pi/\sinh\delta$, we deduce for $k\ge k_0$
\be
\label{control-reste}
\begin{array}{rl}
\ds \forall s\in [-X_\delta(l_i^k),X_\delta(l_i^k)]\quad\quad\sum_{n\in {\Z}^\ast} |a_n^k|^2\, e^{2\,n\, s}&\ds\le  2\, V^{-1}\ \sum_{n\in {\Z}^\ast} e^{-2\,|n|\,X(l_i^k)} \ e^{2\,n\, s}\\[5mm]
 &\ds\le 4\, V^{-1}\, e^{-2\, (X(l_i^k)-s)}
 \end{array}
\ee
For  $s=\pm (X(l_i^k)-t)$, $t\in [0,R]$ and $k$ large enough
\[
\lf||a^k_0+b^k_0\ s|^2-V^{-1}\,(\rho^k_i(s))^2\rg|\le\,4\, V^{-1}\, e^{-2\, (X(l_i^k)-s)}
\]
This gives
\[
\lf|\rho^k_i(s)^{-2}\,|a^k_0+b^k_0\ s|^2-V^{-1}\rg|\le\,4\, V^{-1}\, e^{-2\, (X(l_i^k)-s)}\,\rho^k_i(s)^{-2}
\]
Since $\rho^k_i(s)=\rho^k_i(-s)$ we have for $s=\pm (X(l_i^k)-t)$
\[
\lf|\rho^k_i(s)^{-2} |a^k_0\pm b^k_0\ s|^2-V^{-1}\rg|\le\, 4\, V^{-1}\, e^{-2\, t}\ \lf[\frac{\pi}{\sinh 1}+t+O(l_i^k)\rg]^2
\]
This gives
\[
\lf|\rho^k_i(s)^{-1} |a^k_0\pm b^k_0\ s|-V^{-1/2}\rg|\le\, 4\, V^{-1/2}\, e^{-2\, t}\ \lf[\frac{\pi}{\sinh 1}+2\,t\rg]^2
\]
Hence for any $R.0$ and $k$ large enough
\[
|a^k_0\pm b^k_0\ (X(l_i^k)-R)|\le 4\, V^{-1/2}\, e^{-2\, R}+\frac{V^{-1/2}}{{\frac{\pi}{\sinh 1}+2\,R}}=o_R(1)
\]
This gives
\[
|a^k_0|+|b^k_0\ (X(l_i^k)-R)|=o_R(1)
\]
Hence for any $s\in [X(l_j^k)-R,X(l_j^k)]$ there holds
\[
 |a^k_0\pm b^k_0\ s|^2=o_R(1)
\]
This gives
\[
\lf|\lf(\frac{\pi}{\sinh 1}+2\,t\rg)^2\, o_R(1)-V^{-1}\rg|\le\, 4\, V^{-1}\, e^{-2\, t}\ \lf[\frac{\pi}{\sinh 1}+2\,t\rg]^2
\]
Let $t_\ast>0$ such that
\[
4\, e^{-2\, t_\ast}\ \lf[\frac{\pi}{\sinh 1}+2\,t_\ast \rg]^2=\frac{1}{2}
\]
We then have
\[
V^{-1}\le\, 2^{-1} V^{-1}+o_R(1)\ \lf(\frac{\pi}{\sinh 1}+2\,t_\ast\rg)
\]
For $R$ large enough and $k$ large enough we obtain a contradiction and lemma~\ref{lm-lower-1} is proved. \hfill $\Box$

\medskip
 
 Before continuing the proof of lemma~\ref{lm-delbar-inv} we first estimate $\||\eta_i^k|_{h^k}\|_{L^\infty(\Sigma)}$. We prove the following lemma.
 
 \begin{Lm}
\label{lm-linfty-holom}
Let $\Sigma$ be a closed Riemann surface of non zero genus equipped with a constant gauss curvature metric $h$ with volume $V_h$ equal to $1$ if $g(\Sigma)=1$ and  $V_h=4\pi\,(g(\Sigma)-1)$ otherwize (i.e. $K_h=-1$). Let  $\eta$ be an holomorphic 1-form on $\Sigma$ with $L^2_h$ norm  equal to  1. Then there exists a constant $C_\Sigma>0$ depending only on the topology of $\Sigma$ such that
\be
\label{linfty-holom}
\||\eta^k|_{h^k}\|_{L^\infty(\Sigma)}\le C_\Sigma\ {l_\ast(k)^{-1}}\ .
\ee
\hfill $\Box$
 \end{Lm}
\noindent{\bf Proof of lemma~\ref{lm-linfty-holom}.} First we consider the case $g(\Sigma)=1$. Hence, in this case we have $(\Sigma,h)=({\R}^2/{\Z}\oplus (a+ib){\Z}, b^{-1} \,dz^2)$ and $|a+ib|>1$. The length of the shortest geodesic is equal to $\sqrt{b^{-1}}$ and the $L^2$ unit holomorphic 1-forms are given by $e^{i\theta}\ dz$. Observe that $|dz|_h=\sqrt{b}$. Hence the result is proved in that case.

For the more general case where $g(\Sigma)>1$ we argue by contradiction. Assume 
\[
\limsup_{k\rightarrow+\infty}\max_{i=1\cdots g}\ {l_\ast(k)}\ \||\eta^k|_{h^k}\|_{L^\infty(\Sigma)}=+\infty\ .
\]
We extract a subsequence that we keep denoting $k$ such that $\eta^k$ weakly $L^2$ converges on the thick part and for which
\[
\lim_{k\rightarrow+\infty}\ {l_\ast(k)}\ \||\eta^k|_{h^k}\|_{L^\infty(\Sigma)}=+\infty
\]
Since $\sigma^k:=(\Phi^k)^\ast h^k$ is pre-compact in any $C^l$ norm, $(\Phi^k)^\ast \eta^k$ is converging strongly locally (in any $C^l_{loc}(\Sigma_{\delta,thick})$ topology) towards
a limiting holomorphic one form on $S$. We write $\eta^k=f_i^k(z)\ dz$ where $z=s+i\theta$. hence $|\eta_j^k|^2_{h^k}=(\rho^k)^{-2}\,|f_j^k(z)|^2$. We extend $f^k_j$ periodically on $(-X_\delta(l_i^k),X_\delta(l_i^k))\times {\R}$ by taking $f_j^k(s+i (\theta+2\pi n))=f_j^k(s+i \theta)$ for any $n\in{\Z}$ and $\theta\in [0,2\pi]$. There holds, for any $z_0\in (-X(l_i^k),X(l_i^k))\times [0,2\pi]$
\[
|f_j^k(z_0)|^2\le \frac{i}{2}\dashint_{B_{2\pi}(z_0)}|f_j^k(z)|^2\, dz\wedge d\ov{z}\le \frac{1}{4\pi^2}\int_\Sigma|\eta_j^k|^2_{h^k}\ dvol_{h^k}=\frac{1}{4\pi^2}\ .
\]
Hence we obtain
\[
\||\eta_j^k|_{h^k}\|_{L^{\infty}(\Sigma_{\delta,thin})}\le\frac{1}{4\pi^2}\lf \|{(\rho_i^k(s))^{-1}}\rg\|_{L^\infty(-X(l_i^k),X(l_i^k))}\le \frac{1}{2\pi}\,{\frac{1}{l_i^k}}
\]
which leads to a contradiction. This concludes the proof of lemma~\ref{lm-linfty-holom} for the genus strictly larger than one. We now prove (\ref{delbar-inv-g=1}).

\medskip

\hfill $\Box$

\medskip

\noindent{\bf Proof of Lemma~\ref{lm-delbar-inv} continued.} Combining (\ref{lower-1}) and the previous lemma, we deduce the existence of a positive universal constant $C_\Sigma>0$ such that for any $i=1\cdots g$  such that
\[
\int_{\Sigma}|\eta_i^k\res h^{-1}_k|_{h^k}\ dvol_{h^k}\le\sqrt{V}\le \frac{C_\Sigma}{l_\ast(k)}\ \lf[ 1+V\,\max_{i=1\cdots g}\||\eta_i^k|_{h^k}\|^2_{L^\infty(\Sigma)}\rg]\ \int_\Sigma |\p\lf(\eta_i^k\res h^{-1}_k\rg)|_{h^k}\ dvol_{h^k}
\]
Combining this inequality with (\ref{ui}) and (\ref{berg-bound-thick}) we deduce that for any $\delta>0$ and any $U^{0,1}_{h^k}$ there holds
\be
\label{est-thick}
\int_{\Sigma_{\delta,thick}}|U^{0,1}_{k} |_{h^k}\ dvol_{h^k}\le \frac{C_\Sigma}{l_\ast(k)}\ \lf[ 1+V\,\max_{i=1\cdots g}\||\eta_i^k|_{h^k}\|^3_{L^\infty(\Sigma)}\rg]\ \int_\Sigma |\p U^{0,1}_{k}|_{h^k}\ dvol_{h^k}
\ee
Assume again the inequality would not hold, we then assume that there exists $U^{0,1}_{k}$ with $$\int_{\Sigma_{\delta}}|U^{0,1}_{k} |_{h^k}\ dvol_{h^k}=1\ ,$$ such that
\[
\limsup_{k\rightarrow +\infty}\ (l_\ast^k)^{-4}\, \int_\Sigma |\p U^{0,1}_{h^k}|_{h^k}\ dvol_{h^k}=0\ .
\]
We extract a subsequence such that the above thin-thick decomposition applies and such that
\[
\lim_{k\rightarrow +\infty}\ (l_\ast^k)^{-4}\, \int_\Sigma |\p U^{0,1}_{k}|_{h^k}\ dvol_{h^k}=0\ .
\]
On a given collar $(-X(l_i^k),X(l_i^k))\times S^1$ we write $U^{0,1}_{k}=u^k(z,\ov{z})\ \ov{\p_{z}}$ where $z=s+i\theta$.
Observe that 
\[
|U^{0,1}|_{h^k}\, dvol_{h^k}= |u^k(z,\ov{z})|\,{(\rho_i^k(s))^3}\ ds\wedge d\theta\quad\mbox{and }\quad |\p U^{0,1}|_{h^k}\, dvol_{h^k}= |\p_z\,u^k(z,\ov{z})|\,({\rho_i^k(s)})^2\ ds\wedge d\theta
\] 
%Let
%\[
%u^k(z,\ov{z}):= {\rho_i^k(s)}\, u^k(z,\ov{z}) \quad\mbox{ and }\quad u^k_0(s):=\dashint_0^{2\pi}u^k(s,\theta)\ d\theta={\rho_i^k(s)}\, u^k_0(s)\ .
%\]
%Observe that
%\[
%\begin{array}{l}
%\ds|\p_zu^k(z,\ov{z})|\ds=\lf|{\rho_i^k(s)}\, \p_z u^k(z,\ov{z})+\frac{\dot{\rho}_i^k(s)}{{\rho_i^k(s)}}\,{\rho_i^k(s)}\, u^k(z,\ov{z})\rg|\\[5mm]
%\ds\ds\le{\rho_i^k(s)}\,  \lf|\p_z u^k(z,\ov{z})\rg|+\lf|\frac{\dot{\rho}_i^k(s)}{{\rho_i^k(s)}}\rg|\,|u^k(z,\ov{z})|\ .
%\end{array}
%\]
%and, explicit computations give
%\[
%\lf|\frac{\dot{\rho}_i^k(s)}{{\rho_i^k(s)}}\rg|=\frac{\cos(s\,l_i^k/2\pi)}{2\,\pi\, l^k_i} \, \frac{(l_i^k)^2}{\cos^2(s\,l_i^k/2\pi)}\ |\sin(s\,l_i^k/2\pi)|\le \frac{\rho_i^k(s)}{2\pi}\ .
%\]
%Hence we have finally for any $R>0$
%\be
%\label{phi-u}
%\begin{array}{rl}
%\ds \int_{(-X_\delta(l_i^k)+R,X(l_i^k)-R)\times S^1} |\p_zu^k(z,\ov{z})|\  ds\wedge d\theta &\ds\le\int_{(-X(l_i^k)+R,X(l_i^k)-R)\times S^1}(\rho^k_i)^{-1}\,|\p U^{0,1}|_{h^k}\, dvol_{h^k}\  \\[5mm]
% &\ds+\ \int_{(-X(l_i^k)+R,X(l_i^k)-R)\times S^1} \frac{\rho_i^k(s)}{2\pi}\ |u^k(z,\ov{z})|\  ds\wedge d\theta
% \end{array}
%\ee
We decompose $u^k(s,\theta)$ in its Fourier series and write
\[
u^k(s,\theta)=u^k_-(s,\theta)+u^k_0(s,\theta)+u^k_+(s,\theta)\quad\quad\mbox{where}\quad\quad
\lf\{
\begin{array}{l}
u^k_-(s,\theta)=\sum_{n<0}u^k_n(s)\,e^{i\,n\,\theta}\\[3mm]
u^k_0(s,\theta)=u^k_0(s)\\[3mm]
u^k_+(s,\theta)=\sum_{n>0}u^k_n(s)\,e^{i\,n\,\theta}
\end{array}
\rg.
\]
We have for $n\in{\Z}$
\[
\p_{z}u^k_n(s,\theta)=\p_z\lf[ u^k_n(s)\,e^{i\,n\,\theta}  \rg]=2^{-1}\,e^{i\,n\theta}\, \lf[\dot{u}^k_n(s)+n\,u^k_n(s)  \rg]
\]
Hence
\[
\dashint_0^{2\pi} e^{-i\,n\,\theta}\,\p_{z}u^k(s,\theta)\ d\theta=\dashint_0^{2\pi} e^{-i\,n\,\theta}\,\p_{z}u^k_n(s,\theta)\ d\theta=2^{-1}\, \lf[\dot{u}^k_n(s)+n\,u^k_n(s)  \rg]
\]
We first consider $n>0$ and we take $R>0$ independent of $k$. We have
\[
\int_{-X(l_i^k)+R}^s e^{n\,t}\dashint_0^{2\pi} e^{-i\,n\,\theta}\,\p_{z}u^k(t,\theta)\ d\theta \ dt=   e^{n\,s}\,  u^k_n(s) -e^{n\,(-X(l_i^k)+R)}\, u^k_n(-X(l_i^k)+R)
\]
This gives
\[
u^k_n(s)=\int_{-X(l_i^k)+R}^s e^{n\,(t-s)}\dashint_0^{2\pi} e^{-i\,n\,\theta}\,\p_{z}u^k(t,\theta)\ d\theta \ dt+e^{n\,(-X(l_i^k)+R-s)}\,u^k_n(-X(l_i^k)+R)
\]
from which we deduce for any $(s,\sigma)\in (-X(l_i^k)+R,X(l_i^k)-R)\times S^1$
\[
\begin{array}{l}
\ds u_+^k(s,\sigma)=\int_{-X(l_i^k)+R}^s\dashint_0^{2\pi}\sum_{n>0}e^{n\,((t-i\theta)-(s-i\sigma))}\ \p_{z}u^k(t,\theta)\ d\theta \ dt\\[5mm]
\ds\quad+\dashint_0^{2\pi} \sum_{n>0} e^{n\,(-X(l_i^k)+R-s)}\,e^{-i\,n\,(\theta-\sigma)}\ u^k(-X(l_i^k)+R)\ d\theta   \\[5mm]
\ds=\int_{-X(l_i^k)+R}^s\dashint_0^{2\pi}\frac{e^{\,((t-i\theta)-(s-i\sigma))}}{1-e^{\,((t-i\theta)-(s-i\sigma))}}\,\p_{z}u^k(t,\theta)\ d\theta \ dt\\[5mm]
\ds +\dashint_0^{2\pi} \frac{e^{-X(l_i^k)+R-s-i\,(\theta-\sigma)}}{1-e^{-X(l_i^k)+R-s-i\,(\theta-\sigma)}}\ u^k(-X(l_i^k)+R,\theta)\ d\theta 
\end{array}
\]
Hence
\be
\label{dashint-phi+}
\begin{array}{l}
\ds\dashint_0^{2\pi}|u_+^k(s,\sigma)|\ d\sigma\le \int_{-X(l_i^k)+R}^s\dashint_0^{2\pi}\dashint_0^{2\pi}\frac{e^{t-s}}{|1-e^{\,t-s-i(\theta-\sigma)}|}\ d\sigma\,\lf|\p_{z}u^k(t,\theta)\rg|\ d\theta \ dt\\[5mm]
\ds +\dashint_0^{2\pi} \dashint_0^{2\pi}\frac{e^{-X(l_i^k)+R-s}}{|1-e^{-X(l_i^k)+R-s-i\,(\theta-\sigma)}|}\ d\sigma\, u^k(-X(l_i^k)+R)\ d\theta 
\end{array}
\ee
We have
\be
\label{elliptic-int}
\begin{array}{l}
\ds\dashint_0^{2\pi}\frac{e^{t-s}}{|1-e^{\,t-s-i(\theta-\sigma)}|}\ d\sigma=\dashint_0^{2\pi}\frac{e^{t-s}}{\sqrt{1+e^{2\,(t-s)}-2\, e^{\,t-s}\, \cos(\theta-\sigma)}|}\ d\sigma\\[5mm]
\ds= \frac{e^{t-s}}{\sqrt{1+e^{2\,(t-s)}}}\ \ \dashint_0^{2\pi}\frac{1}{\ds\sqrt{1-\frac{2\, e^{\,t-s}}{1+e^{2\,(t-s)}}\, \cos(\sigma)}|}\ d\sigma
\end{array}
\ee
Using estimates on complete  elliptic integrals of the first kind we obtain for $t-s<0$
\be
\label{elliptic-int+}
\begin{array}{l}
\ds\dashint_0^{2\pi}\frac{e^{t-s}}{|1-e^{\,t-s-i(\theta-\sigma)}|}\ d\sigma=\dashint_0^{2\pi}\frac{e^{t-s}}{\sqrt{1+e^{2\,(t-s)}-2\, e^{\,t-s}\, \cos(\theta-\sigma)}|}\ d\sigma\\[5mm]
\ds=\frac{e^{t-s}}{1+e^{t-s}}\dashint_{0}^{\pi}\frac{d\psi}{\sqrt{1-\frac{4\,e^{t-s}}{(1+e^{t-s})^2}\, \sin^2\psi}}\\[5mm]
\ds=\frac{2}{\pi}\,\frac{e^{t-s}}{1+e^{t-s}}\, \lf[\log\lf[ 1- \frac{4\,e^{t-s}}{(1+e^{t-s})^2}\rg]+ O(1)    \rg]=\frac{1}{\pi}\,\frac{e^{t-s}}{1+e^{t-s}} \lf[2\,\log\lf[ 1- e^{\,t-s} \rg]+ O(1)\rg]
\end{array}
\ee
where $O(1)$ is uniformly bounded independently of $t-s<0$. Combining (\ref{dashint-phi+}), (\ref{elliptic-int}) and (\ref{elliptic-int+}) gives
\be
\label{pos-freq}
\begin{array}{l}
\ds\int_{-X(l_i^k)+R}^{X(l_i^k)-R}\dashint_0^{2\pi}|u_+^k(s,\sigma)|\ d\sigma\, ds\\[5mm]
\ds\le \int_{-X(l_i^k)+R}^{X(l_i^k)-R}\dashint_0^{2\pi} \int_{t}^{X(l_i^k)-R}\frac{e^{t-s}}{{1+e^{(t-s)}}}\  \lf[\frac{2}{\pi}\,\log\lf[ 1- e^{\,t-s} \rg]+ O(1)\rg]\ ds\ \lf|\p_{z}u^k(t,\theta)\rg|\ d\theta \ dt\\[5mm]
\ds +\dashint_0^{2\pi}\int_{-X(l_i^k)+R}^{X(l_i^k)-R} \frac{e^{-s-X(l_i^k)+R}}{{1+e^{\,(-s-X(l_i^k)+R)}}}\  \lf[\frac{2}{\pi}\,\log\lf[ 1- e^{-s-X(l_i^k)+R} \rg]+ O(1)\rg]\ ds\
|u^k(-X(l_i^k)+R,\theta)|\ d\theta \\[5mm]
\ds\le  \int_{-X(l_i^k)+R}^{X(l_i^k)-R}\dashint_0^{2\pi} \int_{e^{t-X(l_i^k)+R}}^1\lf[\frac{2}{\pi}\, \log(1-\xi)+O(1)\rg]\ d\xi\ \lf|\p_{z}u^k(t,\theta)\rg|\ d\theta \ dt\\[5mm]
\ds +\dashint_0^{2\pi} \ \int_{e^{-2\,X(l_i^k)+2\,R}}^1\lf[\frac{2}{\pi}\, \log(1-\xi)+O(1)\rg]\ d\xi\ |u^k(-X(l_i^k)+R,\theta)|\ d\theta\\[5mm]
\ds\le C\  \int_{-X(l_i^k)+R}^{X(l_i^k)-R}\dashint_0^{2\pi}  \lf|\p_{z}u^k(t,\theta)\rg|\ d\theta \ dt+C\,\dashint_0^{2\pi} \ |u^k(-X(l_i^k)+R,\theta)|\ d\theta
\end{array}
\ee
where $C>0$ is universal. Similarly, integrating this time from $X(l_i^k)-R$ instead from $-X(l_i^k)+R$ we obtain the corresponding inequality for $u_-^k(s,\sigma)$
\be
\label{neg-freq}
\begin{array}{l}
\ds\int_{-X(l_i^k)+R}^{X(l_i^k)-R}\dashint_0^{2\pi}|u_-^k(s,\sigma)|\ d\sigma\, ds\\[5mm]
\ds\le C\  \int_{-X(l_i^k)+R}^{X(l_i^k)-R}\dashint_0^{2\pi}  \lf|\p_{z}u^k(t,\theta)\rg|\ d\theta \ dt+C\,\dashint_0^{2\pi} \ |u^k(X(l_i^k)-R,\theta)|\ d\theta
\end{array}
\ee
Regarding now the zero frequency of $u^k$ we simply write
\[
u_0^k(s)=2\,\int_{-X(l_i^k)+R}^s\dashint_0^{2\pi} \p_zu^k(t,\theta)\ d\theta\ dt+u_0^k(-X(l_i^k)+R)
\]
Recall 
\[
{\dot{\rho}_i^k(s)}=\frac{1}{2\,\pi} \, \frac{(l_i^k)^2}{\cos^2(s\,l_i^k/2\pi)}\ \sin(s\,l_i^k/2\pi)
\]
This gives that for $s\in [0,X(l_i^k)]$ $\rho_i^k$ is increasing and its extremal value on $[-X(l_i^k),X(l_i^k)]$ is equal to
\[
\rho_i^k(X(l_i^k))=\frac{l_i^k}{\cos(\pi/2-\arcsin(\sinh(l_i^k/2)/\sinh 1))}=\sinh 1\frac{l_i^k}{\sinh(l_i^k/2)}\le 2\,\sinh 1\le e
\]
Hence\footnote{This inequality is not optimal, but there is no need of improving it in view of (\ref{est-thick}) which is worst.}
\be
\label{zero-freq}
\begin{array}{l}
\ds\int_{-X(l_i^k)+R}^{X(l_i^k)-R}|u_0^k(s)|\,(\rho_i^k(s))^3\ ds\le {C_\Sigma}\int_{-X(l_i^k)+R}^{X(l_i^k)-R}\int_{-X(l_i^k)+R}^{X(l_i^k)-R} (\rho_i^k(s))^2\ ds \ \dashint_0^{2\pi} \lf|\p_zu^k(t,\theta)\rg|\ d\theta\ dt\\[5mm]
\ds +e\,|u_0^k(-X(l_i^k)+R)|\ 
\int_{-X(l_i^k)+R}^{X(l_i^k)-R}(\rho_i^k(s))^2\ ds\\[5mm]
\ds\le \frac{C_\Sigma}{(l_i^k)^2} \int_{-X(l_i^k)+R}^{X(l_i^k)-R} \dashint_0^{2\pi} \lf|\p_zu^k(t,\theta)\rg|\ (\rho_i^k(s))^2\, d\theta\ dt+C_\Sigma\ |u_0^k(-X(l_i^k)+R)|
\end{array}
\ee
Combining (\ref{pos-freq}), (\ref{neg-freq}) and (\ref{zero-freq}) gives finally
\be
\label{l1-bound-phi}
\begin{array}{l}
\ds\int_{-X(l_i^k)+R}^{X(l_i^k)-R}\dashint_0^{2\pi}|u^k(s,\sigma)|\,(\rho_i^k(s))^3\ d\sigma\, ds\le \frac{C_\Sigma}{(l_i^k)^2}\,\int_{-X(l_i^k)+R}^{X(l_i^k)-R}\dashint_0^{2\pi} \lf|\p_zu^k(t,\theta)\rg|\,(\rho_i^k(s))^2\ d\theta\ dt\\[5mm]
\ds +C_\Sigma\,\dashint_0^{2\pi} \ |u^k(-X(l_i^k)+R,\theta)|\ d\theta+C_\Sigma\,\dashint_0^{2\pi} \ |u^k(-X(l_i^k)+R,\theta)|\ d\theta+C_\Sigma\, |u_0^k(-X(l_i^k)+R)|
\end{array}
\ee
We average over $R$ in $[0,1]$ and we can find a ``good slice'' $\rho\in [0,1]$ such that
\[
\begin{array}{l}
\ds\int_{-X(l_i^k)+\rho}^{X(l_i^k)-\rho}\dashint_0^{2\pi}|u^k(s,\sigma)|\,(\rho_i^k(s))^3\ d\sigma\, ds\le \frac{C_\Sigma}{(l_i^k)^2}\,\int_{-X(l_i^k)}^{X(l_i^k)}\dashint_0^{2\pi} \lf|\p_zu^k(t,\theta)\rg|\,(\rho_i^k(s))^2\ d\theta\ dt\\[5mm]
\ds +C_\Sigma\,\dashint_0^{2\pi} \ \int_{0}^1\,|u^k(-X(l_i^k)+R,\theta)|\ d\theta\, dR+C_\Sigma\,\dashint_0^{2\pi} \ \int_{0}^1\,|u^k(-X(l_i^k)+R,\theta)|\ d\theta\, dR\\[5mm]
\ds+{C_\Sigma}\, \int_{0}^1\,|u_0^k(-X(l_i^k)+R)|\ dR\ .
\end{array}
\]
Combining this inequality with  (\ref{est-thick}) gives finally
\[
\liminf_{k\rightarrow +\infty} (l_\ast^k)^{-4}\ \int_\Sigma |\p U^{0,1}_{k}|_{h^k}\ dvol_{h^k}>0
\]
which is a contradiction and lemma~\ref{lm-delbar-inv} is proved.\hfill $\Box$
%%%%%%%%%%%%%%%%%%%%%%%%%%%%%%%%%%%%%%%%%%%%%%%%%%%%%%%%%%%%%%%%%%%%%%%%%%%%%%%%%%%%%%%%%%%%%%%%%%%%%%%%
\section{A-priori global controls}
\reset
\subsection{A-priori global control of the $L^1$ norm of the time derivative of the underlying constant Gauss curvature metric}
In this part we are still considering a smooth solution of the parametrized Willmore flow and we aim at giving an a-priori control of the time derivative of the flow of immersions. We consider the flow over an interval $[0,T]$. We denote by $l_\ast(t)$ the length of the shortest geodesic of $(\Sigma,h(t))$ and we introduce 
\[
l_{\ast\ast}(T):=\min_{t\in[0,T]}l_\ast(t)\ .
\]
First we will prove an $L^1$ control of the time derivative of the underlying constant Gauss curvature metric $h$. We have the following.
\begin{Lm}
\label{lm-L1-h} Leet $\vec{\Phi}(t)$ be a solution of the parametrized Willmore flow, under the previous notations there holds
\be
\label{L1-h}
\begin{array}{l}
\ds\int_0^T\int_\Sigma\lf|\frac{dh}{dt} \rg|_{h(t)}\ dvol_{h(t)}\ dt\\[5mm]
%=\int_0^T\int_\Sigma\, \lf|P_{h(t)}\lf(\frac{d\vec{\Phi}}{dt}\cdot\vec{\mathcal H}^{0}\rg)\rg|_{h(t)}\ dvol_{h(t)}\ dt\\[5mm]
\ds\le C_\Sigma\ \frac{e^{2\,\|\al\|_{\infty,T}}}{l_{\ast,\ast}^2(T)}\ \sqrt{T}\ \sqrt{W(\vec{\Phi}(0))+2\pi\,(g(\Sigma)-1)}\,\sqrt{W(\vec{\Phi}(0))-W(\vec{\Phi}(T))}\ .
\end{array}
\ee
\hfill $\Box$
\end{Lm}
Before proving lemma~\ref{lm-L1-h}, we shall be using the following $L^\infty$ bound on the $L^2_h$ unit holomorphic quadratic forms.
 \begin{Lm}
\label{lm-linfty-holom-quad}
Let $q$ be an holomorphic quadratic form of $(\Sigma,h)$ and $g(\Sigma)\ge1$. In the case $g(\Sigma)=1$ we assume $V_h=\int_\Sigma dvol_{h}=1$. There exists  a constant $C_\Sigma>0$ depending only on the topology of $\Sigma$ such that
\be
\label{linfty-holom-quad}
||q|_{h}\|_{L^\infty(\Sigma)}\le C_\Sigma\ {l_\ast^{-2}}\, \|q\|_{L^2_h(\Sigma)}\ .
\ee
where $l_\ast$ is the length of the shortest closed geodesic on $(\Sigma,h)$.
\hfill $\Box$
 \end{Lm}
\noindent{\bf Proof of lemma~\ref{lm-linfty-holom-quad}.} First we consider the case $g(\Sigma)=1$. Hence, in this case we have $(\Sigma,h)=({\R}^2/{\Z}\oplus (a+ib){\Z}, b^{-1} \,dz^2)$ and $|a+ib|>1$. The length of the shortest geodesic is equal to $\sqrt{b^{-1}}$ and the $L^2$ unit holomorphic quadratic forms are given by $b^{-1}\,e^{i\theta}\, dz^{\otimes^2}$. The result follows in that case.

For the case $g(\Sigma)>1$ we argue by contradiction. Assume $l_\ast^k\rightarrow 0$, $\|q^k\|_{L^2_h(\Sigma)}=1$ and
\[
\limsup_{k\rightarrow+\infty} {(l_\ast^k)^2}\ \||q^k|_{h^k}\|_{L^\infty(\Sigma)}=+\infty\ .
\]
We extract a subsequence that we keep denoting $k$ such that $q^k$ weakly $L^2$ converges on the thick part.  Since $\sigma^k:=(\Phi^k)^\ast h^k$ is pre-compact in any $C^l$ norm, $(\Phi^k)^\ast \eta_i^k$ is converging strongly locally (in any $C^l_{loc}(\Sigma_{\delta,thick})$ topology) towards
a limiting holomorphic quadratic form on $S$. we write $q^k=f^k_i(z)\ dz\otimes dz$ where $z=s+i\theta$ in a degenerating collar $(-X(l_i^k),X(l_i^k))\times [0,2\pi]$ in which we take the coordinates as in the previous section . Hence $|q^k|^2_{h^k}=(\rho_i^k)^{-4}\,|f^k(z)|^2$. We extend $f^k_i$ periodically on $(-X(l_i^k),X(l_i^k))\times {\R}$ by taking $f_j^k(s+i (\theta+2\pi n))=f_j^k(s+i \theta)$ for any $n\in{\Z}$ and $\theta\in [0,2\pi]$. There holds, for any $z_0\in (-X(l_i^k),X(l_i^k))\times [0,2\pi]$
\[
|f_j^k(z_0)|^2\le \frac{i}{2}\,\dashint_{B_{2\pi}(z_0)}|f_j^k(z)|^2\, dz\wedge d\ov{z}\le \frac{1}{4\pi^2}\int_\Sigma  (\rho_i^k)^2 \,| q^k|^2_{h^k}\ dvol_{h^k}
\]
We have 
\[
\frac{\dot{\rho}_i^k}{\rho_i^k}=\frac{\rho_i^k}{2\pi}\,\sin\lf(\frac{s\,l_i^k}{2\pi}\rg)\ .
\]
Hence in particular 
\[
\lf\|\frac{d}{ds}\log(\rho_i^k)\rg\|_{L^\infty((-X(l_i^k),X(l_i^k)\times [0,2\pi])}\le \frac{e}{2\pi}\ .
\]
we deduce
\[
\limsup_{k\rightarrow+\infty} \sup_{z_0\in (-X(l_i^k),X(l_i^k))\times {\R}}\ \frac{\max_{z\in B_{2\pi}(z_0)}{\rho_i^k(z)}}{\max_{z\in B_{2\pi}(z_0)}{\rho_i^k(z)}}<+\infty
\]
Hence we obtain
\[
|q^k|^2_{h^k}(z_0)= (\rho_i^k)^{-4}(z_0)\,|f_j^k(z_0)|^2\le C\, (\rho_i^k)^{-2}(z_0)\,\int_\Sigma  \,| q^k|^2_{h^k}\ dvol_{h^k}\le \frac{C}{(l_i^k)^2}\ 
\]
which leads to a contradiction. This concludes the proof of lemma~\ref{lm-linfty-holom-quad}.\hfill $\Box$

\medskip

\noindent{\bf Proof of lemma~\ref{lm-L1-h}.} Let $(q_i(t))_{i=1\cdots Q}$ be an $L^2_{h(t)}$-orthonormal basis of holomorphic quadratic forms. We have
\[
P_{h(t)}\lf(\frac{d\vec{\Phi}}{dt}\cdot\vec{\mathcal H}^{0}\rg)=\sum_{i=1}^Q\int_\Sigma \lf<\frac{d\vec{\Phi}}{dt}\cdot\vec{\mathcal H}^{0},q_i(t)\rg>_{h(t)}\ dvol_{h(t)}\ q_i(t)\ .
\]
 Using lemma~\ref{linfty-holom-quad} we have
\be
\label{der-conf-213}
\begin{array}{l}
\ds\int_0^T\int_\Sigma\, \lf|P_{h(t)}\lf(\frac{d\vec{\Phi}}{dt}\cdot\vec{\mathcal H}^{0}\rg)\rg|_{h(t)}\ dvol_{h(t)}\ dt\\[5mm]
\ds\le  \sum_{i=1}^Q\sup_{[0,T]} \|q_i(t)\|_{L^\infty(\Sigma)}\ \int_0^T\ \int_\Sigma\, \lf|\pi_{\vec{n}}\frac{d\vec{\Phi}}{dt}\rg|\ \lf|\vec{\mathcal H}^{0}\rg|_{h(t)}\ dvol_{h(t)}\ \sup_{[0,T]}\|q_i(t)\|_{L^1(\Sigma)}\\[5mm]
\ds\le\frac{C_\Sigma}{l^2_{\ast\ast}}\int_0^T\int_\Sigma\lf|\pi_{\vec{n}}\frac{d\vec{\Phi}}{dt}\rg|\ \lf|\vec{\frak h}^{0}\rg|_{h(t)}\ dvol_{h(t)}\ dt\\[5mm]
\ds\le C_\Sigma\ \frac{e^{2\,\|\al\|_\infty}}{l^2_{\ast\ast}}\ \lf[\int_0^T\int_{\Sigma}|\vec{\mathbb I}_{\vec{\Phi}(t)}|^2\ dvol_{g_{\vec{\Phi}(t)}}\ dt\rg]^{1/2} \  \lf[\int_0^T\int_{\Sigma}\lf|\pi_{\vec{n}}\frac{d\vec{\Phi}}{dt}\rg|^2\ dvol_{g_{\vec{\Phi}(t)}}\ dt\rg]^{1/2} \\[5mm]
\ds\le  C_\Sigma\ \frac{e^{2\,\|\al\|_\infty}}{l^2_{\ast\ast}}\ \sqrt{T}\ \lf[\int_{\Sigma}|\vec{\mathbb I}_{\vec{\Phi}(0)}|^2\ dvol_{g_{\vec{\Phi}(0)}}\rg]^{1/2}\,\lf(W(\vec{\Phi}(0))-W(\vec{\Phi}(t))\rg)^{1/2}\ ,
\end{array}
\ee
where we have used that $\lf|\vec{\mathcal H}^{0}\rg|_{h(t)}=\lf|\vec{\frak h}^{0}\res h(t)\rg|_{h(t)}=\lf|\vec{\frak h}^{0}\rg|_{h(t)}$. This concludes the proof of the lemma~\ref{lm-L1-h}.\hfill $\Box$

\subsection{A-priori global control of the $L^1$ norm of the time derivative of the flow of immersions}
Precisely we are proving the following lemma.
\begin{Lm}
\label{lm-L1-control-tangent}
Let $\vec{\Phi}(t)$ be a solution of the parametrized Willmore flow on a time interval $[0,T]$. We have
\be
\label{L1-control-tangent}
\begin{array}{l}
\ds\int_0^T\int_\Sigma\lf|\pi_T\frac{d\vec{\Phi}}{dt}  \rg|\ dvol_{h(t)}\\[5mm]
%\ds \le C_\Sigma\ \frac{e^{3\,\|\al\|_{\infty,T}}}{l_{\ast,\ast}^6(T)}\ \sqrt{T}\ \lf[\int_{\Sigma}|\vec{\mathbb I}_{\vec{\Phi}(0)}|^2\ dvol_{g_{\vec{\Phi}(0)}}\rg]^{1/2}\,\sqrt{W(\vec{\Phi}(0))-W(\vec{\Phi}(T))}\\[5mm]
\ds\le C_\Sigma\ \frac{e^{3\,\|\al\|_{\infty,T}}}{l_{\ast,\ast}^6(T)}\ \sqrt{T}\ \sqrt{W(\vec{\Phi}(0))+2\pi\,(g(\Sigma)-1)}\,\sqrt{W(\vec{\Phi}(0))-W(\vec{\Phi}(T))}\ .
\end{array}
\ee
where $\pi_T$ is the tangential projection, $\vec{\Phi}(t)^\ast g_{{\R}^m}=e^{2\al}\ h(t)$ and
\[
\|\al\|_{\infty,T}:=\|\al\|_{L^\infty([0,T]\times \Sigma)}
\]
\hfill $\Box$
\end{Lm}
\noindent{\bf Proof of lemma~\ref{lm-L1-control-tangent}.} Observe that we have $\pi_T\frac{d\vec{\Phi}}{dt} =\vec{U}=\vec{\Phi}(t)_\ast U$. Hence
\[
|\vec{U}|=|U|_{g_{\vec{\Phi}(t)}}=e^{\al}\,|U|_{h(t)}=\sqrt{2}\,e^{\al}\,|U^{0,1}|_{h(t)}\ .
\]
Hence using lemma~\ref{lm-delbar-inv}
\be
\label{der-conf-212-p1}
\begin{array}{l}
\ds\int_0^T\int_\Sigma\lf|\pi_T\frac{d\vec{\Phi}}{dt}  \rg|\ dvol_{h(t)}=\ds\int_0^T\int_\Sigma\lf|\vec{U}(t)  \rg|\ dvol_{h(t)}=  \sqrt{2}\,\int_0^T\int_\Sigma\lf|{U}^{0,1}_{h(t)}\rg|_{h(t)}\ e^\al\, dvol_{h(t)}\\[5mm]
\ds\le \frac{C_\Sigma}{l^4_{\ast\ast}(T)}\ e^{\|\al\|_{\infty,T}}\ \int_0^T\int_\Sigma |\p^{h(t)}{U}^{0,1}|_{h(t)}\, dvol_{h(t)}\\[5mm]
\ds\le \frac{C_\Sigma}{l^4_{\ast\ast}(T)}\ e^{\|\al\|_{\infty,T}}\ \int_0^T\int_\Sigma \lf|(I-\ti{P}_{h(t)})\lf(\frac{d\vec{\Phi}}{dt}\cdot\vec{\frak h}^{0}\rg)\rg|_{h(t)}\ dvol_{h(t)}\ dt\\[5mm]
\ds\le  \frac{C_\Sigma}{l^4_{\ast\ast}(T)}\ e^{\|\al\|_{\infty,T}}\ \int_0^T\int_\Sigma\lf|\pi_{\vec{n}}\frac{d\vec{\Phi}}{dt}\rg|\ \lf|\vec{\frak h}^{0}\rg|_{h(t)}\ dvol_{h(t)}\ dt\\[5mm]
\ds+ \frac{C_\Sigma}{l^4_{\ast\ast}(T)}\ e^{\|\al\|_{\infty,T}}\ \int_0^T\int_\Sigma\, \lf|P_{h(t)}\lf(\frac{d\vec{\Phi}}{dt}\cdot\vec{\mathcal H}^{0}\rg)\rg|_{h(t)}\ dvol_{h(t)}\ dt
\end{array}
\ee
Observe first that $\lf|\vec{\frak h}^{0}\rg|_{h(t)}=\lf|\vec{\frak h}^{0}\rg|_{g_{\vec{\Phi}(t)}}$. Thus
\be
\label{der-conf-213-p2}
\begin{array}{l}
\ds\frac{C_\Sigma}{l^4_{\ast\ast}(T)}\ e^{\|\al\|_{\infty,T}}\ \int_0^T\int_\Sigma\lf|\pi_{\vec{n}}\frac{d\vec{\Phi}}{dt}\rg|\ \lf|\vec{\frak h}^{0}\rg|_{h(t)}\ dvol_{h(t)}\ dt\\[5mm]
\ds\le C_\Sigma\ \frac{e^{3\,\|\al\|_\infty}}{l^4_{\ast\ast}(T)}\ \lf[\int_0^T\int_{\Sigma}|\vec{\mathbb I}_{\vec{\Phi}(t)}|^2\ dvol_{g_{\vec{\Phi}(t)}}\ dt\rg]^{1/2} \  \lf[\int_0^T\int_{B_R^{h(t)}(x_0)}\lf|\pi_{\vec{n}}\frac{d\vec{\Phi}}{dt}\rg|^2\ dvol_{g_{\vec{\Phi}(t)}}\ dt\rg]^{1/2} \\[5mm]
\ds\le  C_\Sigma\ \frac{e^{3\,\|\al\|_\infty}}{l^4_{\ast\ast}(T)}\ \sqrt{T}\  \sqrt{W(\vec{\Phi}(0))+2\pi\,(g(\Sigma)-1)}\,\lf(W(\vec{\Phi}(0))-W(\vec{\Phi}(t))\rg)^{1/2}\ ,
\end{array}
\ee
Using lemma~\ref{lm-L1-h} we have
\be
\label{der-conf-213-bis}
\begin{array}{l}
\ds\frac{C_\Sigma}{l^4_{\ast\ast}(T)}\ e^{\|\al\|_{\infty,T}}\ \int_0^T\int_\Sigma\, \lf|P_{h(t)}\lf(\frac{d\vec{\Phi}}{dt}\cdot\vec{\mathcal H}^{0}\rg)\rg|_{h(t)}\ dvol_{h(t)}\ dt\\[5mm]
\ds\le  C_\Sigma\ \frac{e^{3\,\|\al\|_\infty}}{l_{\ast\ast}^6(T)}\ \sqrt{T}\     \sqrt{W(\vec{\Phi}(0))+2\pi\,(g(\Sigma)-1)}    \,\lf(W(\vec{\Phi}(0))-W(\vec{\Phi}(t))\rg)^{1/2}\ ,
\end{array}
\ee
where we have used that $\lf|\vec{\mathcal H}^{0}\rg|_{h(t)}=\lf|\vec{\frak h}^{0}\res h(t)\rg|_{h(t)}=\lf|\vec{\frak h}^{0}\rg|_{h(t)}$.

%%%%%%%%%%%%%%%%%%%%%%%%%%%%%%%%%%%%%%%%%%%%%%%%%%%%%%%%%%%%%%%%%%%%%%%%%%%%%%%%%%%%%%%%%%%%%%
\subsection{Controlled conformal charts for Constant Gauss Curvature metrics}
The following lemma is a consequence of the classical geometry of constant Gauss curvature surfaces and the fact that half of the length of the shortest closed geodesic is equal to the systole.
\begin{Lm}
\label{lm-controled-charts} 
%Let $\vec{\Phi}_t$ be a $L^2_tW^{4,2}_x([0,T]\times\Sigma)$ solution to the conformal Parametric Willmore Flow (\ref{I.25}) and 

Let $(\Sigma,h)$ be closed oriented constant Gauss curvature Riemannian surface with volume equal to $1$ if $g(\Sigma)=1$ or $4\pi\,|g-1|$ if $g(\Sigma)\ne 1$. Denote by $l_\ast$ the length of the shortest closed geodesic. Then for any $p\in \Sigma$ and any $R\le l_\ast/2$ there exist a conformal diffeomorphism
\[
x_{p}\ :\ B_{R}^{h}(p)\ \longrightarrow\ B_{R}(0)\subset {\C}
\]
where $B_{R}^{h}(p)$ is denoting the geodesic ball for the metric $h$ equipped with this metric, such that
\[
(x^{-1})^\ast h=e^{2\nu}\, [dx_1^2+dx_2^2]\ ,
\]
moreover
\[
\forall\ k\in {\N}\quad\quad \|\nabla^k \nu\|_{L^\infty(B_R(0))}\le C_k\ ,
\]
where $C_k$ only depends on $k$, not on $h$ and $l_\ast$ in particular.
\hfill $\Box$
\end{Lm}
\subsection{Control of the tangential part of the flow in the weak $L^2-$space}
We shall now establish some straightforward consequences of the previous lemma.
\begin{Lm}
\label{lm-L2-infty}
Let $U^{0,1}$ be an $L^1_h(\Sigma)$ section of $T^{0,1}\Sigma$ such that $\p^hU^{0,1}\in L^1_h(\Sigma)$ then $U^{0,1}\in L^{2,\infty}(\Sigma)$ and there holds
\be
\label{L2-infty}
\forall \,p\in \Sigma\quad\quad \|U^{0,1}\|_{L^{2,\infty}(B^h_{l^\ast/4}(p))}\le C \,\|\p^hU^{0,1}\|_{L^{1}(B^h_{l^\ast/2}(p))}+C\,l_\ast^{-1}\,\|U^{0,1}\|_{L^{1}(B^h_{l^\ast/2}(p))}\ .
\ee
\hfill $\Box$
\end{Lm}
\noindent{\bf Proof of lemma~\ref{lm-L2-infty}.} Thanks to the chart given by lemma~\ref{lm-controled-charts} we reduce the problem to the flat case on euclidian discs of radii comparable to $l_\ast$ and (\ref{L2-infty}) is a direct consequence of classical elliptic estimates.\hfill $\Box$

\subsection{Local Controls of the conformal factor in $L^\infty$, $W^{1,2}$ and $W^{1,(2,1)}$}
Another consequence of the existence of controlled  conformal charts given by lemma~\ref{lm-controled-charts} is the following lemma.

\begin{Lm}
\label{lm-L2-1}
Let $\vec{\Phi}$ be a conformal immersion into ${\R}^m$ of $(\Sigma,h)$ a closed oriented constant Gauss curvature Riemannian surface with volume equal to $1$ if $g(\Sigma)=1$ or $4\pi\,|g-1|$ if $g(\Sigma)\ne 1$. Denote by $l_\ast$ the length of the shortest closed geodesic on $(\sigma,h)$ and by $\al$ the conformal factor satisfying
\[
\vec{\Phi}^{\,\ast}g_{{\R}^m}=e^{2\al}\, h\ .
\]
Let $0<R<l_\ast/2$ and $p\in \Sigma$ such that
\[
\int_{B^h_{R}(p)}|\vec{\mathbb I}|^2_{g_{\vec{\Phi}}}\ dvol_{g_{\vec{\Phi}}}<\frac{8\pi}{3}\ .
\]
Then there exists $C>0$ depending only on $\Sigma$ and $m$ such that for any $\sigma<1/4$
\be
\label{L2-1}
\begin{array}{rl}
\ds\lf\|\al-\dashint_{B^h_{\sigma\,R}(p))}\al\ dvol_{h}\rg\|_{L^{\infty}(B^h_{\sigma\,R}(p))}&\ds\le C\, \|d\al\|_{L^{2,1}(B^h_{\sigma\,R}(p))}\\[5mm]
\ds\quad&\ds\le C\ \int_{B^h_{R}(p)}|\vec{\mathbb I}|^2_{g_{\vec{\Phi}}}\ dvol_{g_{\vec{\Phi}}}+ C\ \sigma\, [1+\|d\al\|_{L^{2,\infty}(B^h_{R}(p))}]
\end{array}
\ee
We also have\footnote{The $L^{2,1}$ space is the Lorentz space pre-dual of the weak Marcinkiewicz space $L^{2\infty}$ (see \cite{Gra1})}
\be
\label{L2-2}
\|d\al\|_{L^{2,1}(B^h_{R/4}(p))}\le \frac{C_{\Sigma,m}}{ l_\ast^{2}}\ \int_{\Sigma}|\vec{\mathbb I}|^2_{g_{\vec{\Phi}}}\ dvol_{g_{\vec{\Phi}}}
\ee
\hfill $\Box$
\end{Lm}
\noindent{\bf Proof of lemma~\ref{lm-L2-1}.} Thanks again lemma~\ref{lm-controled-charts} we have a controlled conformal chart in the ball $B^h_{R}(p)$ in which
\[
(x^{-1})^\ast h= e^{2\nu}\ [dx_1^2+dx_2^2]
\]
In this chart the assumption of the lemma imply
\[
\int_{B_{R}(0)}|\nabla\vec{n}|^2 dx^2<\frac{8\pi}{3}\ .
\]
Using the existence of controlled moving frames by H\'elein (see \cite{Hel} lemma 5.1.4), we have the existence of an orthonormal frame\footnote{The space $V_2({\R}^m)$ classically denotes the Stiefel manifold of orthonormal two frames in ${\R}^m$.} $(\vec{e}_1,\vec{e}_2)\in V_2({\R}^m)$ such that $\vec{n}:=\vec{e}_1\wedge\vec{e}_2$ and
\[
\|\nabla\vec{e}_i\|_{L^2(B_R(0))}\le C_m\,\|\nabla\vec{n}\|_{L^2(B_R(0))}=C_m\, \sqrt{\int_{B^h_{R}(p)}|\vec{\mathbb I}|^2_{g_{\vec{\Phi}}}\ dvol_{g_{\vec{\Phi}}}}
\]
In $B_R(0)$ we have
\[
(x^{-1})^\ast g_{\vec{\Phi}}=e^{2\al+2\nu}\ [dx_1^2+dx_2^2]\ .
\]
and we have
\[
\Delta (\al+\nu)= \nabla^\perp\vec{e}_1\cdot\nabla\vec{e}_2\ \quad\quad\mbox{ in }B_R(0)\ .
\]
Let $\mu$ be the solution of
\[
\lf\{
\begin{array}{l}
\ds\Delta \mu= \nabla^\perp\vec{e}_1\cdot\nabla\vec{e}_2\ \quad\quad\mbox{ in }B_R(0)\\[5mm]
\ds\mu=0\quad\quad\mbox{ in }\quad \p B_R(0)
\end{array}
\rg.
\]
Classical integrability by compensation result (see again \cite{Hel}) gives then 
\[
\|\nabla\mu\|_{L^{2,1}(B_{R}(0))}\le C\,\int_{B^h_{R}(p)}|\vec{\mathbb I}|^2_{g_{\vec{\Phi}}}\ dvol_{g_{\vec{\Phi}}}\ .
\]
Since $\al+\nu-\mu$ is harmonic, the monotonicity formula gives  for any $\sigma<1/4$
\[
\begin{array}{l}
\ds\,C^{-1}\,\|\nabla(\al+\nu-\mu)\|^2_{L^{2,1}(B_{\sigma\,R}(0))}\le \int_{B_{2\,\sigma\,R}}|\nabla(\al+\nu-\mu)|^2 dx^2\\[5mm]
\ds\le 4\, \sigma^2\, \int_{B_{2^{-1}\,R}}|\nabla(\al+\nu-\mu)|^2 dx^2\le C\,\sigma^2\ \|\nabla(\al+\nu-\mu)\|^2_{L^{2,\infty}(B_R(0))}\ .
\end{array}
\]
Combining the previous gives then
\[
\|\nabla(\al+\nu)\|_{L^{2,1}(B_{\sigma\,R}(0))}
\le C\,\int_{B^h_{R}(p)}|\vec{\mathbb I}|^2_{g_{\vec{\Phi}}}\ dvol_{g_{\vec{\Phi}}}+C\ \sigma\ \|\nabla(\al+\nu)\|_{L^{2,\infty}(B_R(0))}\ .
\]
Since $\|\nabla \nu\|_\infty$ is universally bounded on this local conformal chart we obtain (\ref{L2-1}).

\medskip

The function $\al$ satisfies the Liouville equation
\[
-\Delta_h\al=e^{2\al}\,K_{g_{\vec{\Phi}}}-K_h
\]
We have in particular
\[
d\al(x)=\int_\Sigma d_xG(x,y)\ e^{2\al}\,K_{g_{\vec{\Phi}}}\ dvol_h\ .
\]
Thanks to lemma~\ref{lm-L1-estim} we deduce
\be
\label{daal-L1}
\|d\al\|_{L^1(\Sigma)}\le \sup_{y\in\Sigma}\int_\Sigma |d_xG(x,y)|\ dvol_{h}(x)\ \int_\Sigma e^{2\al}\,|K_{g_{\vec{\Phi}}}|\ dvol_h\le \frac{C_\Sigma}{l_\ast} \int_{\Sigma}|\vec{\mathbb I}|^2_{g_{\vec{\Phi}}}\ dvol_{g_{\vec{\Phi}}}\ .
\ee
Writing the Liouville equation in the chart centred at $p$ for the radius $l_\ast/2$ gives
\[
-\Delta \al= e^{2(\al+\nu)}\,K_{g_{\vec{\Phi}}}- K_h\,e^{2\nu}\ .
\]
Hence, thanks to classical elliptic estimates we have for any $t<1$
\[
\|\nabla\al\|_{L^{2,\infty}(B_{t\,l_\ast/2}(0))}\le \int_{B_{l_\ast/2}(0)}\lf|e^{2(\al+\nu)}\,K_{g_{\vec{\Phi}}}- K_h\,e^{2\nu}\rg| \ dx^2+C\, l_\ast^{-1}\,\|\nabla\al\|_{L^1(B_{l_\ast/2}(0)}
\]
This estimate brought back on $\Sigma$ thanks to the control on $\nu$ is giving for $R<l_\ast/2$ and thanks to the fact that $ \int_{\Sigma}|\vec{\mathbb I}|^2_{g_{\vec{\Phi}}}$ is bounded from below by a universal constant we obtain
\[
\|d\al\|_{L^{2,\infty}(B^h_{R})}\le C \int_{\Sigma}|\vec{\mathbb I}|^2_{g_{\vec{\Phi}}}\ dvol_{g_{\vec{\Phi}}}+C\, l_\ast^{-1}\,\|d\al\|_{L^1(\Sigma)}
\]
Combining this inequality with (\ref{daal-L1}) and (\ref{L2-1}) we obtain (\ref{L2-2}) and lemma~\ref{lm-L2-1} is proved.\hfill $\Box$
%%%%%%%%%%%%%%%%%%%%%%%%%%%%%%%%%%%%%%%%%%%%%%%%%%%%%%%%%%%%%%%%%%%%%%%%%%%%%%%%%%%%%%%%%%%%%%%%%%%%%%
\subsection{A-priori global control of the conformal factor}
In this part we are still considering a smooth solution of the parametrized Willmore flow and we aim at giving an a-priori control of the time dependent conformal factor $\al(t)$ in $L^\infty$ norm uniformly on $\Sigma\times[0,T]$ for some $T>0$ chosen smaller than some data. We recall that $\al_t$ is given by
\[
\vec{\Phi}_t^{\, \ast}g_{{\R}^m}=e^{2\,\al(t)}\,h(t)\ .
\]
We first start with the following lemma
We have the following lemma
\begin{Lm}
\label{lm-param-conf-var} Under the previous notations we have
\be
\label{param-conf-var}
\frac{d\al}{dt}=\,\frac{d\vec{\Phi}}{dt}\cdot \vec{H}+e^{-\,2\al(t)}\frac{1}{2}\,d^{\ast_{h(t)}}\lf[\frac{d\vec{\Phi}}{dt}\cdot d\vec{\Phi}(t)\rg]\ .
\ee
\hfill$\Box$
\end{Lm}
\noindent{\bf  Proof of lemma~\ref{lm-param-conf-var}. } We have
\[
\vec{H}:=2^{-1}\,\Delta_{g_{\vec{\Phi}(t)}}\vec{\Phi}(t)
\]
Taking the scalar product with $d\vec{\Phi}/dt$ and multiplying by $e^{2\al(t)}$ gives
\[
2\,e^{2\al(t)}\,\frac{d\vec{\Phi}}{dt}\cdot \vec{H}=\frac{d\vec{\Phi}}{dt}\cdot\Delta_{h(t)}\vec{\Phi}(t)=-d^{\ast_{h(t)}}\lf[\frac{d\vec{\Phi}}{dt}\cdot d\vec{\Phi}(t)\rg]-\, d\frac{d\vec{\Phi}}{dt}\cdot_{h(t)} d\vec{\Phi}(t)
\]
We have
\[
\begin{array}{l}
\ds d\frac{d\vec{\Phi}}{dt}\cdot_{h(t)} d\vec{\Phi}(t)=\lf(d\frac{d\vec{\Phi}}{dt}\,\dot{\otimes }\,d\vec{\Phi}\rg)\res_2 h^{-1}(t)=\frac{1}{2}\lf(d\frac{d\vec{\Phi}}{dt}\,\dot\,{\otimes }d\vec{\Phi}+d\vec{\Phi}\,\dot{\otimes}\,d\frac{d\vec{\Phi}}{dt}\rg)\res_2 h^{-1}(t)\\[5mm]
\ds=\frac{1}{2}\frac{d}{dt}\lf[ d\vec{\Phi}\,\dot{\otimes}\, d\vec{\Phi}\rg]\res_2 h^{-1}(t)-\frac{1}{2}\,d\vec{\Phi}\,\dot{\otimes}\, d\vec{\Phi}\res_2 \frac{d h^{-1}}{dt}\\[5mm]
\ds=\frac{1}{2}\frac{d}{dt}\lf[ g_{\vec{\Phi}(t)}\rg]\res_2 h^{-1}(t)-\frac{1}{2}\,g_{\vec{\Phi}(t)}\res_2 \frac{d h^{-1}}{dt}=\frac{1}{2}\frac{d}{dt}\lf[ e^{2\al(t)}\,h(t)\rg]\res_2 h^{-1}(t)-\frac{1}{2}\,e^{2\al(t)}\, h(t)\res_2 \frac{d h^{-1}}{dt}\\[5mm]
\ds=2\,\frac{d\al}{dt}\,e^{2\al}+e^{2\al}\,\frac{dh}{dt}\res_2 h^{-1}
\end{array}
\]
where we have used the fact that $h\res_2 h^{-1}\equiv 2$.  Observe that $dh/dt$ being the real part of a quadratic form it is traceles. Hence
\[
\frac{dh}{dt}\res_2 h^{-1}=0
\] 
Combining the three last identities  gives the result and lemma~\ref{lm-param-conf-var} is proved.\hfill $\Box$

\medskip

Now we come to the main result of the present subsection.
\begin{Lm}
\label{lm-conf-fact-infty}
Under the previous notations, let 
$$l_{\ast\ast}(T):=\inf_{[0,T]}l_\ast(t)\ , $$ and $R<l_{\ast\ast}(T)/4$ such that
\be
\label{smallness}
\sup_{x\in\Sigma}\ \sup_{[0,T]}\ \int_{B_{2\,R}^{h(t)}(x)}|\vec{\mathbb I}_{\vec{\Phi}(t)}|^2\ dvol_{g_{\vec{\Phi}(t)}}<\frac{8\pi}{3}\ .
\ee
Then the following inequality holds
\be
\label{conf-fact-infty}
\begin{array}{l}
\ds\int_0^T\lf|\frac{d}{dt}\int_\Sigma\al\ \chi^{h(t)}_{R,x_0}(x)\ dvol_{h(t)}\rg|\ dt\\[5mm]
\ds\le C_\Sigma\ \frac{e^{2\,\|\al\|_\infty}}{l_{\ast\ast}^9(T)}\ \sqrt{T}\     E(0)^{3/2}    \,\sqrt{W(\vec{\Phi}(0))-W(\vec{\Phi}(t))}\\[5mm]
\ds +C_\Sigma\ \frac{e^{2\,\|\al\|_\infty}}{l^4_{\ast\ast}(T)}\ \sqrt{T}\  \sup_{[0,T]}\,\int_{B_R^{h(t)}(x_0)}\ |\al|\ dvol_{h(t)}\   \sqrt{E(0)}\,\sqrt{W(\vec{\Phi}(0))-W(\vec{\Phi}(t))} 
\end{array}
\ee
where $E(0):=W(\vec{\Phi}(0))+2\pi\,(g(\Sigma)-1)$.
\hfill $\Box$
\end{Lm}
\noindent{\bf Proof of lemma~\ref{lm-conf-fact-infty}.} Let $x_0\in\Sigma$ and denote by $\chi_{R,x_0}(x)$ a cut-off function given by $$\chi_{R,x_0}(x):=\chi\lf(\frac{|x-x_0|_{h(t)}}{R}\rg)$$ where $\chi\in C^\infty_0([-1,1])$ and $\chi\equiv 1$ on $[-1/2,1/2]$. Multiplying (\ref{param-conf-var}) by $\chi_{R,x_0}(x)$ and integrating over $\Sigma$ gives
\[
\begin{array}{l}
\ds\int_{\Sigma}\frac{d\al}{dt}\ \chi^{h(t)}_{R,x_0}(x)\ dvol_{h(t)}\\[5mm]
\ds=\frac{d}{dt}\int_\Sigma\al\ \chi^{h(t)}_{R,x_0}(x)\ dvol_{h(t)}-\int_\Sigma\al\ \chi^{h(t)}_{R,x_0}(x)\ \frac{d}{dt}dvol_{h(t)}-\frac{1}{R}\,\int_\Sigma\al\ \frac{d|x-x_0|_{h(t)}}{dt}\,\chi'\lf( \frac{|x-x_0|_{h(t)}}{R}  \rg) dvol_{h(t)}\\[5mm]
\ds=\int_\Sigma \chi^{h(t)}_{R,x_0}(x)\ \frac{d\vec{\Phi}}{dt}\cdot \vec{H}\ dvol_{h(t)}+\frac{1}{2}\int_\Sigma \chi^{h(t)}_{R,x_0}(x)\ e^{-\,2\al(t)}\,d^{\ast_{h(t)}}\lf[\frac{d\vec{\Phi}}{dt}\cdot d\vec{\Phi}(t)\rg]\ dvol_{h(t)}\end{array}
\]
Observe that
\[
\frac{d}{dt}dvol_{h(t)}=\frac{1}{2}\frac{dh}{dt}\res_2 h^{-1}(t)\ dvol_{h(t)}=0\ .
\]
Hence we have
\be
\label{deri-conf}
\begin{array}{l}
\ds\frac{d}{dt}\int_\Sigma\al\ \chi^{h(t)}_{R,x_0}(x)\ dvol_{h(t)}\\[5mm]
\ds=\int_\Sigma \chi^{h(t)}_{R,x_0}(x)\ \frac{d\vec{\Phi}}{dt}\cdot \vec{H}\ dvol_{h(t)}+\frac{1}{2}\int_\Sigma d\lf[\chi^{h(t)}_{R,x_0}(x)\ e^{-\,2\al(t)}\rg]\wedge \lf[\frac{d\vec{\Phi}}{dt}\cdot d\vec{\Phi}(t)\rg]\\[5mm]
\ds+\frac{1}{R}\int_\Sigma\al\ \frac{d|x-x_0|_{h(t)}}{dt}\, \,\chi'\lf( \frac{|x-x_0|_{h(t)}}{R}  \rg) \ dvol_{h(t)}\
\end{array}
\ee
We bound the first term in the r.h.s of (\ref{deri-conf}) as follows
\be
\label{der-conf-1}
\begin{array}{l}
\ds\int_0^T\lf|\int_\Sigma \chi^{h(t)}_{R,x_0}(x)\ \frac{d\vec{\Phi}}{dt}\cdot \vec{H}\ dvol_{h(t)}\rg|\, dt\\[5mm]
\ds\le e^{2\,\|\al\|_\infty}\ \lf[\int_0^T\int_{B_R^{h(t)}(x_0)}|\vec{H}|^2\ dvol_{g_{\vec{\Phi}(t)}}\ dt\rg]^{1/2} \  \lf[\int_0^T\int_{B_R^{h(t)}(x_0)}\lf|\pi_{\vec{n}}\frac{d\vec{\Phi}}{dt}\rg|^2\ dvol_{g_{\vec{\Phi}(t)}}\ dt\rg]^{1/2} \\[5mm]
\ds\le e^{2\,\|\al\|_\infty}\ \sqrt{T}\ \lf[\sup_{x\in\Sigma}\ \sup_{[0,T]}\ \int_{B_R^{h(t)}(x)}|\vec{\mathbb I}_{\vec{\Phi}(t)}|^2\ dvol_{g_{\vec{\Phi}(t)}}\rg]^{1/2}\,\lf(W(\vec{\Phi}(0))-W(\vec{\Phi}(t))\rg)^{1/2}\ ,
\end{array}
\ee
where we use the notation
\[
\|\al\|_\infty:=\|\al\|_{L^\infty([0,T]\times\Sigma)}\ .
\]
We bound the second term in the r.h.s of (\ref{deri-conf}) as follows
\be
\label{der-conf-2}
\begin{array}{l}
\ds\int_0^T\ds\lf|\frac{1}{2}\int_\Sigma d\lf[\chi^{h(t)}_{R,x_0}(x)\ e^{-\,2\al(t)}\rg]\wedge \lf[\frac{d\vec{\Phi}}{dt}\cdot d\vec{\Phi}(t)\rg]\rg|\ dt\\[5mm]
\ds\le R^{-1}\, \int_0^T\ds\lf|\frac{1}{2}\int_\Sigma \chi'\lf( \frac{|x-x_0|_{h(t)}}{R}  \rg)  e^{-2\al}\, d|x-x_0|_{h(t)}\wedge\lf[\frac{d\vec{\Phi}}{dt}\cdot d\vec{\Phi}(t)\rg]\rg|\ dt\\[5mm]
\ds+\int_0^T\ds\lf|\int_\Sigma \chi^{h(t)}_{R,x_0}(x)\ e^{-\,2\al(t)}\ d\al\wedge \lf[\frac{d\vec{\Phi}}{dt}\cdot d\vec{\Phi}(t)\rg]\rg|\ dt
\end{array}
\ee
We now consider each of the two terms in the r-h-s of (\ref{der-conf-2}). We first have
\be
\label{der-conf-21}
\begin{array}{l}
\ds R^{-1}\,\int_0^T\ds\lf|\frac{1}{2}\int_\Sigma \chi'\lf( \frac{|x-x_0|_{h(t)}}{R}  \rg)  e^{-2\al}\, d|x-x_0|_{h(t)}\wedge\lf[\frac{d\vec{\Phi}}{dt}\cdot d\vec{\Phi}(t)\rg]\rg|\ dt\\[5mm]
\ds\le C\,R^{-1}\,\int_0^T\ds\int_{B_R^{h(t)}(x)} e^{-2\al}\, |d\vec{\Phi}|_{h(t)}\ |\vec{U}|\ dvol_{h(t)}
\end{array}
\ee
We have $\vec{U}=\vec{\Phi}(t)_\ast U$. Hence 
\be
\label{u}
|\vec{U}|=|U|_{g_{\vec{\Phi}(t)}}=e^{\al}\,|U|_{h(t)}=\sqrt{2}\,e^{\al}\,|U^{0,1}|_{h(t)}\ , 
\ee
and $|d\vec{\Phi}|_{h(t)}=e^\al\, |d\vec{\Phi}|_{g_{\vec{\Phi}(t)}}=\sqrt{2}\, e^\al$. Hence we have
\be
\label{der-conf-211}
\begin{array}{l}
\ds R^{-1}\,\int_0^T\ds\lf|\frac{1}{2}\int_\Sigma \chi'\lf( \frac{|x-x_0|_{h(t)}}{R}  \rg)  e^{-2\al}\, d|x-x_0|_{h(t)}\wedge\lf[\frac{d\vec{\Phi}}{dt}\cdot d\vec{\Phi}(t)\rg]\rg|\ dt\\[5mm]
\ds\le C\, R^{-1}\,\int_0^T\int_{B_R^{h(t)}(x)}|U^{0,1}|_{h(t)}\ dvol_{h(t)}\ dt\le C\,\int_0^T\|U^{0,1}\|_{L^{2,\infty}(B_R^{h(t)}(x))}\ dt
\end{array}
\ee
Using lemma~\ref{lm-L2-infty} and lemma~\ref{lm-delbar-inv} we obtain
\be
\label{der-conf-212-a}
\begin{array}{l}
\ds R^{-1}\,\int_0^T\ds\lf|\frac{1}{2}\int_\Sigma \chi'\lf( \frac{|x-x_0|_{h(t)}}{R}  \rg)  e^{-2\al}\, d|x-x_0|_{h(t)}\wedge\lf[\frac{d\vec{\Phi}}{dt}\cdot d\vec{\Phi}(t)\rg]\rg|\ dt\\[5mm]
\ds\le \frac{C}{l^5_{\ast\ast}(T)}\,\int_0^T\int_\Sigma |\p^{h(t)}U^{0,1}|_{h(t)}\ dvol_{h(t)}\ dt\ .
\end{array}
\ee
Using the explicit expression of $\p^{h(t)}U^{0,1}$ given by the flow we obtain the following.
\be
\label{der-conf-212-b}
\begin{array}{l}
\ds R^{-1}\,\int_0^T\ds\lf|\frac{1}{2}\int_\Sigma \chi'\lf( \frac{|x-x_0|_{h(t)}}{R}  \rg)  e^{-2\al}\, d|x-x_0|_{h(t)}\wedge\lf[\frac{d\vec{\Phi}}{dt}\cdot d\vec{\Phi}(t)\rg]\rg|\ dt\\[5mm]
\ds\le \frac{C_\Sigma}{l^5_{\ast\ast}(T)}\int_0^T\int_\Sigma \lf|(I-\ti{P}_{h(t)})\lf(\frac{d\vec{\Phi}}{dt}\cdot\vec{\frak h}^{0}\rg)\rg|_{h(t)}\ dvol_{h(t)}\ dt\\[5mm]
\ds\le  \frac{C_\Sigma}{l^5_{\ast\ast}(T)}\int_0^T\int_\Sigma\lf|\pi_{\vec{n}}\frac{d\vec{\Phi}}{dt}\rg|\ \lf|\vec{\frak h}^{0}\rg|_{h(t)}\ dvol_{h(t)}\ dt\\[5mm]
\ds+ \frac{C_\Sigma}{l^5_{\ast\ast}(T)}\int_0^T\int_\Sigma\, \lf|P_{h(t)}\lf(\frac{d\vec{\Phi}}{dt}\cdot\vec{\mathcal H}^{0}\rg)\rg|_{h(t)}\ dvol_{h(t)}\ dt
\end{array}
\ee
Using (\ref{der-conf-213-p2}) and ((\ref{der-conf-213-bis}) we obtain
\be
\label{der-conf-212-c}
\begin{array}{l}
\ds R^{-1}\,\int_0^T\ds\lf|\frac{1}{2}\int_\Sigma \chi'\lf( \frac{|x-x_0|_{h(t)}}{R}  \rg)  e^{-2\al}\, d|x-x_0|_{h(t)}\wedge\lf[\frac{d\vec{\Phi}}{dt}\cdot d\vec{\Phi}(t)\rg]\rg|\ dt\\[5mm]
 \ds \le C_\Sigma\ \frac{e^{2\,\|\al\|_\infty}}{l_{\ast\ast}^7(T)}\ \sqrt{T}\     \sqrt{W(\vec{\Phi}(0))+2\pi\,(g(\Sigma)-1)}    \,\lf(W(\vec{\Phi}(0))-W(\vec{\Phi}(t))\rg)^{1/2}\ ,
\end{array}
\ee
We are now controlling the second term in the r.h.s. of (\ref{der-conf-2}). We have using (\ref{u}),  (\ref{L2-infty}), (\ref{L2-2}), (\ref{der-conf-213-p2}) and ((\ref{der-conf-213-bis})
\be
\label{der-conf-23-c}
\begin{array}{l}
\ds\int_0^T\ds\lf|\int_\Sigma \chi^{h(t)}_{R,x_0}(x)\ e^{-\,2\al(t)}\ d\al\wedge \lf[\frac{d\vec{\Phi}}{dt}\cdot d\vec{\Phi}(t)\rg]\rg|\ dt\\[5mm]
\ds\le \int_0^T\ds\int_{B_R^{h(t)}(x)} e^{-2\al}\,\ |d\al|_{h(t)}\  |d\vec{\Phi}|_{h(t)}\ |\vec{U}|\ dvol_{h(t)}\\[5mm]
\ds\le\,\int_0^T \|d\al\|_{L^{2,1}(B^{h(t)}_{R}(p))}\ \|U^{0,1}\|_{L^{2,\infty}(B_R^{h(t)}(x))}\ dt\\[5mm]
\ds\le \frac{C_{\Sigma,m}}{l^7_{\ast\ast}(T)}\ \lf[W(\vec{\Phi}(0))+2\pi\,(g(\Sigma)-1)\rg]\ \int_0^T\int_\Sigma |\p^{h(t)}U^{0,1}|_{h(t)}\ dvol_{h(t)}\ dt\ \\[5mm]
\ds \le C_\Sigma\ \frac{e^{2\,\|\al\|_\infty}}{l_{\ast\ast}^9(T)}\ \sqrt{T}\     \lf[W(\vec{\Phi}(0))+2\pi\,(g(\Sigma)-1)\rg]^{3/2}    \,\lf(W(\vec{\Phi}(0))-W(\vec{\Phi}(t))\rg)^{1/2}\ ,
\end{array}
\ee
Finally we are estimating the third term in the r.h.s. of (\ref{deri-conf}). We observe the following. Let $s\in[0,1]\rightarrow\gamma^{h(t)}_{x_0}(s)$ be the unique geodesic bounding $x$ and $x_0$ inside the convex ball $B_R^{h(t)}(x_0)$.
We have
\[
|x-x_0|_{h(t)}=\int_0^1\lf|\frac{d\gamma^{h(t)}_{x,x_0}}{ds}\rg|_{h(t)}\ ds\ .
\]
for $t+\delta(t)$ we obviously have
\[
|x-x_0|_{h(t+\delta)}\le \int_0^1\lf|\frac{d\gamma^{h(t)}_{x,x_0}}{ds}\rg|_{h(t+\delta)}\ ds
\]
Hence
\[
\frac{|x-x_0|_{h(t+\delta)}-|x-x_0|_{h(t)}}{\delta}\le  \int_0^1\frac{1}{\delta}\,\lf[\lf|\frac{d\gamma^{h(t)}_{x,x_0}}{ds}\rg|_{h(t+\delta)}-\lf|\frac{d\gamma^{h(t)}_{x,x_0}}{ds}\rg|_{h(t)}\rg]\ ds
\]
We have respectively
\[
\lf|\frac{d\gamma^{h(t)}_{x,x_0}}{ds}\rg|_{h(t+\delta)}=\sqrt{ h(t+\delta)\res_2\frac{d\gamma^{h(t)}_{x,x_0}}{ds}\otimes \frac{d\gamma^{h(t)}_{x,x_0}}{ds}}
\]
and
\[
\lf|\frac{d\gamma^{h(t)}_{x,x_0}}{ds}\rg|_{h(t+\delta)}=\sqrt{ h(t)\res_2\frac{d\gamma^{h(t)}_{x,x_0}}{ds}\otimes \frac{d\gamma^{h(t)}_{x,x_0}}{ds}}
\]
Hence
\[
\begin{array}{l}
\ds\lim_{\delta\rightarrow 0}\frac{1}{\delta}\,\lf[\lf|\frac{d\gamma^{h(t)}_{x,x_0}}{ds}\rg|_{h(t+\delta)}-\lf|\frac{d\gamma^{h(t)}_{x,x_0}}{ds}\rg|_{h(t)}\rg]\\[8mm]
\ds\quad=\frac{1}{2}\frac{dh}{dt}(t)\res_2\frac{d\gamma^{h(t)}_{x,x_0}}{ds}\otimes \frac{d\gamma^{h(t)}_{x,x_0}}{ds}\lf|\frac{d\gamma^{h(t)}_{x,x_0}}{ds}\rg|_{h(t)}^{-1}
\end{array}
\]
Observe that
\[
\begin{array}{l}
\ds\lf|\frac{dh}{dt}(t)\res_2\frac{d\gamma^{h(t)}_{x,x_0}}{ds}\otimes \frac{d\gamma^{h(t)}_{x,x_0}}{ds}\rg|\ \lf|\frac{d\gamma^{h(t)}_{x,x_0}}{ds}\rg|_{h(t)}^{-1}\\[5mm]
\ds\le \lf|\frac{dh}{dt}(t)\rg|_{h(t)}\,\lf|\frac{d\gamma^{h(t)}_{x,x_0}}{ds}\rg|_{h(t)}
\end{array}
\]
Hence we deduce
\be
\label{dist-cont}
\limsup_{\delta\rightarrow 0}\frac{|x-x_0|_{h(t+\delta)}-|x-x_0|_{h(t)}}{\delta}\le  \lf\| \lf|\frac{dh}{dt}(t)\rg|_{h(t)}\rg\|_{L^\infty(B_R^{h(t)}(x_0))}\ |x-x_0|_{h(t)}
\ee
Exchanging the role of $t+\delta$ and $t$ gives finally
\be
\label{der-c-1}
\limsup_{\delta\rightarrow 0}\lf|\frac{|x-x_0|_{h(t+\delta)}-|x-x_0|_{h(t)}}{\delta}\rg|\le R\,  \lf\| \lf|\frac{dh}{dt}(t)\rg|_{h(t)}\rg\|_{L^\infty(B_R^{h(t)}(x_0))}
\ee
Inserting this fact in  the third term in the r.h.s. of (\ref{deri-conf}) gives
\be
\label{der-conf-23-cc}
\begin{array}{l}
\ds\frac{1}{R}\int_\Sigma\al\ \frac{d|x-x_0|_{h(t)}}{dt}\, \,\chi'\lf( \frac{|x-x_0|_{h(t)}}{R}  \rg) \ dvol_{h(t)}\\[5mm]
\ds\le C\ \sup_{[0,T]}\,\int_{B_R^{h(t)}(x_0))}\ |\al|\ dvol_{h(t)}\  \int_0^T\lf\| \lf|\frac{dh}{dt}(t)\rg|_{h(t)}\rg\|_{L^\infty(B_R^{h(t)}(x_0))}\  dt
\end{array}
\ee
We have
\be
\label{der-conf-23-d}
\begin{array}{l}
\ds \int_0^T\lf\| \lf|\frac{dh}{dt}(t)\rg|_{h(t)}\rg\|_{L^\infty(B_R^{h(t)}(x_0))}\  dt\\[5mm]
\ds\le \int_0^T \lf\|\lf|P_{h(t)}\lf(\frac{d\vec{\Phi}}{dt}\cdot\vec{\mathcal H}^{0}\rg)\rg|_{h(t)}\rg\|_{L^\infty(B_R^{h(t)}(x_0))}\ dt\\[5mm]
\ds\le  \sum_{i=1}^Q\sup_{[0,T]} \|q_i(t)\|^2_{L^\infty(\Sigma)}\ \int_0^T\ \int_\Sigma\, \lf|\pi_{\vec{n}}\frac{d\vec{\Phi}}{dt}\rg|\ \lf|\vec{\mathcal H}^{0}\rg|_{h(t)}\ dvol_{h(t)}\\[5mm]
\ds\le\frac{C_\Sigma}{l^4_{\ast\ast}(T)}\int_0^T\int_\Sigma\lf|\pi_{\vec{n}}\frac{d\vec{\Phi}}{dt}\rg|\ \lf|\vec{\frak h}^{0}\rg|_{h(t)}\ dvol_{h(t)}\ dt\\[5mm]
\ds\le C_\Sigma\ \frac{e^{2\,\|\al\|_\infty}}{l^4_{\ast\ast}(T)}\ \lf[\int_0^T\int_{\Sigma}|\vec{\mathbb I}_{\vec{\Phi}(t)}|^2\ dvol_{g_{\vec{\Phi}(t)}}\ dt\rg]^{1/2} \  \lf[\int_0^T\int_{\Sigma}\lf|\pi_{\vec{n}}\frac{d\vec{\Phi}}{dt}\rg|^2\ dvol_{g_{\vec{\Phi}(t)}}\ dt\rg]^{1/2} \\[5mm]
\ds\le  C_\Sigma\ \frac{e^{2\,\|\al\|_\infty}}{l^4_{\ast\ast}(T)}\ \sqrt{T}\ \lf[\int_{\Sigma}|\vec{\mathbb I}_{\vec{\Phi}(0)}|^2\ dvol_{g_{\vec{\Phi}(0)}}\rg]^{1/2}\,\lf(W(\vec{\Phi}(0))-W(\vec{\Phi}(t))\rg)^{1/2}\ ,
\end{array}
\ee
where we have used again that $\lf|\vec{\mathcal H}^{0}\rg|_{h(t)}=\lf|\vec{\frak h}^{0}\res h(t)\rg|_{h(t)}=\lf|\vec{\frak h}^{0}\rg|_{h(t)}$. Combining (\ref{deri-conf}), (\ref{der-conf-1}), (\ref{der-conf-212-c}), (\ref{der-conf-23-c}) and (\ref{der-conf-23-d})
\be
\label{deri-conf-fin}
\begin{array}{l}
\ds\int_0^T\lf|\frac{d}{dt}\int_\Sigma\al\ \chi^{h(t)}_{R,x_0}(x)\ dvol_{h(t)}\rg|\ dt\\[5mm]
\ds\le C_\Sigma\ \frac{e^{2\,\|\al\|_\infty}}{l_{\ast\ast}^9(T)}\ \sqrt{T}\     E(0)^{3/2}    \,\sqrt{W(\vec{\Phi}(0))-W(\vec{\Phi}(t))}\\[5mm]
\ds +C_\Sigma\ \frac{e^{2\,\|\al\|_\infty}}{l^4_{\ast\ast}(T)}\ \sqrt{T}\  \sup_{[0,T]}\,\int_{B_R^{h(t)}(x_0)}\ |\al|\ dvol_{h(t)}\   \sqrt{E(0)}\,\sqrt{W(\vec{\Phi}(0))-W(\vec{\Phi}(t))} 
\end{array}
\ee
where $E(0):=W(\vec{\Phi}(0))+2\pi\,(g(\Sigma)-1)$. This concludes the proof of lemma~\ref{lm-conf-fact-infty}. \hfill $\Box$
\begin{Lm}
\label{lm-conf-fact-shift} Under the above notations 
Let $\sigma<1/8$, $x_0\in \Sigma$  and $R< 2^{-1}\,l_{\ast\ast}(T)$ such that
\[
\sup_{t\in[0,T]}\sup_{x_0\in\Sigma}\,\int_{B_{R}^{h(t)}(x)}|\vec{\mathbb I}_{\vec{\Phi}(t)}|^2\ dvol_{g_{\vec{\Phi}(t)}}<\frac{8\pi}{3}\ .
\]
then there exists a constant $C_{\Sigma,m}>0$ depending only on the topology of $\Sigma$ and the dimension $m$ such that
\be
\label{conf-fact-shift}
\begin{array}{l}
\ds\sup_{t\in[0,T]}\lf\|\al(t)-\al(0)\rg\|_{L^\infty(\Sigma)}\le\,   C_{\Sigma,m}\,\Sigma\, \sup_{t\in[0,T]}\sup_{x_0\in\Sigma}\int_{B^{h(t)}_{R}(p)}|\vec{\mathbb I}_{\vec{\Phi}(t)}|^2_{g_{\vec{\Phi}(t)}}\ dvol_{g_{\vec{\Phi}(t)}}+\,\frac{C_{\Sigma,m}}{  l_{\ast\ast}^{2}(T)}\,\sigma\ E(0)\ \\[5mm]
\ds\quad\quad+\,C_{\Sigma,m}\ \frac{e^{2\,\|\al\|_\infty}}{l_{\ast\ast}^9(T)}\ \frac{\sqrt{T}}{\sigma^2\,R^2}\     E(0)^{3/2}    \,\sqrt{W(\vec{\Phi}(0))-W(\vec{\Phi}(t))}\\[5mm]
\ds\quad\quad+C_{\Sigma,m}\ \frac{e^{2\,\|\al\|_\infty}}{l^4_{\ast\ast}(T)}\ {\sqrt{T}}\ \log\sqrt{\frac{\mbox{Vol}(\vec{\Phi}(0))}{V_{h(0)}}}\ \   \sqrt{E(0)}\,\sqrt{W(\vec{\Phi}(0))-W(\vec{\Phi}(t))}\\[5mm]
\ds\quad\quad+\,C_{\Sigma,m}\,  \sup_{x_0\in\Sigma}\, \lf\|\al(0)(x)-\al(0)(y)\rg\|_{ L^\infty\lf(B_{\La(T)\,\sigma\,R}^{h(0)}(x_0)\times B_{\La(T)\,\sigma\,R}^{h(0)}(x_0)\rg)}\ .
\end{array}
\ee
where 
\[
E(0):=W(\vec{\Phi}(0))+2\pi\,(g(\Sigma)-1)
\]
moreover
\[
\|\al\|_\infty:=\|\al\|_{L^\infty([0,T]\times\Sigma)}
\]
and
\[
\log\,\La(T):= C_\Sigma\  \frac{e^{2\,\|\al\|_\infty}}{l^4_{\ast\ast}(T)}\ \sqrt{T}\ \sqrt{E(0)}\,\lf(W(\vec{\Phi}(0))-W(\vec{\Phi}(T))\rg)^{1/2}\ .
\]
\hfill $\Box$
\end{Lm}
\noindent{\bf Proof of lemma~\ref{lm-conf-fact-shift}.}
Let $\sigma<1/8$, $x_0\in \Sigma$  and $R< 2^{-1}\,l_{\ast\ast}(T)$ such that
\[
\sup_{t\in[0,T]}\sup_{x_0\in\Sigma}\,\int_{B_{R}^{h(t)}(x)}|\vec{\mathbb I}_{\vec{\Phi}(t)}|^2\ dvol_{g_{\vec{\Phi}(t)}}<\frac{8\pi}{3}\ .
\]
We recall from (\ref{L2-1}) and (\ref{L2-2}) we have for any $t\in[0,T]$
\be
\label{shift-1}
\begin{array}{rl}
\ds\lf\|\al(t)-\dashint_{B^h_{\sigma\,R}(p)}\al(t)\ dvol_{h}\rg\|_{L^{\infty}(B^{h(t)}_{\sigma\,R}(x_0))}&\ds\le   C\, \int_{B^{h(t)}_{R}(x_0))}|\vec{\mathbb I}(t)|^2_{g_{\vec{\Phi}(t)}}\ dvol_{g_{\vec{\Phi}(t)}}\\[5mm]
\ds &\ds+\frac{C_{\Sigma,m}}{ l_\ast^{2}}\,\sigma\, \lf[W(\vec{\Phi}(0))+2\pi\,(g(\Sigma)-1)\rg]
\end{array}
\ee
This gives in particular
\be
\label{shift-2}
\begin{array}{l}
\ds\lf\|\al(t)\dashint\chi^{h(t)}_{\sigma\, R,x_0}-\dashint_{B^{h(t)}_{\sigma\,R}(p)}\al(t)\ \chi^{h(t)}_{\sigma\, R,x_0}\ dvol_{h(t)}\rg\|_{L^{\infty}(B^{h(t)}_{\sigma\,R}(x_0))}\\[5mm]
\ds\le \lf\|\lf[\al(t)-  \dashint_{B^{h(t)}_{\sigma\,R}(p)}\al(t)\ dvol_{h(t)} \rg]\,\dashint\chi^{h(t)}_{\sigma\, R,x_0}\ dvol_{h(t)}\rg\|_{L^{\infty}(B^{h(t)}_{\sigma\,R}(x_0))}\\[5mm]
\ds+ \lf|\dashint_{B^{h(t)}_{\sigma\,R}(p)}\lf[\al(t)-\dashint_{B^{h(t)}_{\sigma\,R}(p)}\al(t)\ dvol_{h(t)}\rg] \chi^{h(t)}_{\sigma\, R,x_0}\ dvol_{h}\rg|\\[5mm]
\ds\le   C\, \int_{B^{h(t)}_{R}(p)}|\vec{\mathbb I}(t)|^2_{g_{\vec{\Phi}(t)}}\ dvol_{g_{\vec{\Phi}(t)}}+\frac{C_{\Sigma,m}}{ l_\ast^{2}}\,\sigma\, \lf[W(\vec{\Phi}(0))+2\pi\,(g(\Sigma)-1)\rg]
\end{array}
\ee
This implies
\be
\label{shift-3}
\begin{array}{l}
\ds\lf\|\al(t)\int\chi^{h(t)}_{\sigma\, R,x_0}-\int_{B^{h(t)}_{\sigma\,R}(p)}\al(t)\ \chi^{h(t)}_{\sigma\, R,x_0}\ dvol_{h(t)}\rg\|_{L^{\infty}(B^{h(t)}_{\sigma\,R}(x_0))}\\[5mm]
\ds\le   C\,  \int_{B^{h(t)}_{R}(p)}|\vec{\mathbb I}(t)|^2_{g_{\vec{\Phi}(t)}}\ dvol_{g_{\vec{\Phi}(t)}}\ \int_{B^{h(t)}_{\sigma\,R}(p)}\ dvol_{h(t)}\\[5mm]
\ds+\frac{C_{\Sigma,m}}{ l_\ast^{2}}\,\sigma\, \lf[W(\vec{\Phi}(0))+2\pi\,(g(\Sigma)-1)\rg]\ \int_{B^{h(t)}_{\sigma\,R}(p)}\ dvol_{h(t)}\ 
\end{array}
\ee
The explicit expression of $h(t)$ is giving that below the injectivity radius, for $R<l_\ast(t)/2$ there holds $\int_{B^{h(t)}_{\sigma\,R}(p)}\ dvol_{h(t)}\simeq \sigma^2\,R^2$. Hence we have finally
\be
\label{shift-4}
\begin{array}{l}
\ds\lf\|\al(t)\int_\Sigma\chi^{h(t)}_{\sigma\, R,x_0}-\int_{\Sigma}\al(t)\ \chi^{h(t)}_{\sigma\, R,x_0}\ dvol_{h(t)}\rg\|_{L^{\infty}(B^{h(t)}_{\sigma\,R}(x_0))}\\[5mm]
\ds\le   C\, \sigma^2\,R^2\ \int_{B^{h(t)}_{R}(p)}|\vec{\mathbb I}(t)|^2_{g_{\vec{\Phi}(t)}}\ dvol_{g_{\vec{\Phi}(t)}}\ \\[5mm]
\ds+\frac{C_{\Sigma,m}}{ l_{\ast\ast}^{2}(T)}\,\sigma^3\, R^2\, \lf[W(\vec{\Phi}(0))+2\pi\,(g(\Sigma)-1)\rg]\ 
\end{array}
\ee
This implies in particular, for any $t\in[0,T]$
\be
\label{shift-4-a}
\begin{array}{l}
\ds\lf\|\chi^{h(t)}_{\sigma\, R,x_0}\,\al(t)\int_\Sigma\chi^{h(t)}_{\sigma\, R,x_0}-\chi^{h(t)}_{\sigma\, R,x_0}\,\int_{\Sigma}\al(t)\ \chi^{h(t)}_{\sigma\, R,x_0}\ dvol_{h(t)}\rg\|_{L^{\infty}(\Sigma)}\\[5mm]
\ds\le   C\, \sigma^2\,R^2\ \int_{B^{h(t)}_{R}(p)}|\vec{\mathbb I}(t)|^2_{g_{\vec{\Phi}(t)}}\ dvol_{g_{\vec{\Phi}(t)}}\ \\[5mm]
\ds+\frac{C_{\Sigma,m}}{ l_{\ast\ast}^{2}(T)}\,\sigma^3\, R^2\, \lf[W(\vec{\Phi}(0))+2\pi\,(g(\Sigma)-1)\rg]\ \int_{B^{h(t)}_{\sigma\,R}(p)}\ dvol_{h(t)}\ 
\end{array}
\ee
From lemma~\ref{lm-conf-fact-infty} we deduce
\be
\label{shift-5}
\begin{array}{l}
\ds\sup_{t\in[0,T]}\lf|\int_{\Sigma}\al(t)\ \chi^{h(t)}_{\sigma\, R,x_0}\ dvol_{h(t)}-\int_{\Sigma}\al(0)\ \chi^{h(0)}_{\sigma\, R,x_0}\ dvol_{h(0)}\rg|\\[5mm]
\ds\le C_\Sigma\ \frac{e^{2\,\|\al\|_\infty}}{l_{\ast\ast}^9(T)}\ \sqrt{T}\     E(0)^{3/2}    \,\sqrt{W(\vec{\Phi}(0))-W(\vec{\Phi}(t))}\\[5mm]
\ds +C_\Sigma\ \frac{e^{2\,\|\al\|_\infty}}{l^4_{\ast\ast}(T)}\ \sqrt{T}\  \sup_{[0,T]}\,\int_{B_{\sigma\,R}^{h(t)}(x_0)}\ |\al|\ dvol_{h(t)}\   \sqrt{E(0)}\,\sqrt{W(\vec{\Phi}(0))-W(\vec{\Phi}(t))} 
\end{array}
\ee
Now we have
\be
\label{shift-5-a}
\begin{array}{l}
\ds\lf\|\chi^{h(t)}_{\sigma\, R,x_0}\,\al(0)\int_\Sigma\chi^{h(0)}_{\sigma\, R,x_0}-\chi^{h(t)}_{\sigma\, R,x_0}\,\int_{\Sigma}\al(0)\ \chi^{h(0)}_{\sigma\, R,x_0}\ dvol_{h(0)}\rg\|_{L^{\infty}(\Sigma)}\\[5mm]
\ds \le C\, R^2\ \lf\|\al(0)(x)-\al(0)(y)\rg\|_{L^\infty(B_{\sigma\,R}^{h(t)}(x_0)\times B_{\sigma\,R}^{h(t)}(x_0))}
\end{array}
\ee
We have also using (\ref{dist-cont})
\[
\int_0^T\lf| \frac{d}{ds}\log |x-x_0|_{h(s)} \rg|\ ds\le \int_0^T \lf|\lf\|\frac{dh}{dt}(t)\rg|_{h(t)}\rg\|_{L^\infty(\Sigma)}\ ds
\]
Using (\ref{der-conf-23-d}) we obtain in particular
\be
\label{def-control}
 \lf|\log\frac{|x-x_0|_{h(0)}}{|x-x_0|_{h(s)}}\rg|\le C_\Sigma\ \frac{e^{2\,\|\al\|_\infty}}{l^4_{\ast\ast}(T)}\ \sqrt{T}\ \lf[\int_{\Sigma}|\vec{\mathbb I}_{\vec{\Phi}(0)}|^2\ dvol_{g_{\vec{\Phi}(0)}}\rg]^{1/2}\,\lf(W(\vec{\Phi}(0))-W(\vec{\Phi}(T))\rg)^{1/2}\ .
\ee
Introduce
\[
\log\,\La(T):= C_\Sigma\  \frac{e^{2\,\|\al\|_\infty}}{l^4_{\ast\ast}(T)}\ \sqrt{T}\ \lf[\int_{\Sigma}|\vec{\mathbb I}_{\vec{\Phi}(0)}|^2\ dvol_{g_{\vec{\Phi}(0)}}\rg]^{1/2}\,\lf(W(\vec{\Phi}(0))-W(\vec{\Phi}(T))\rg)^{1/2}\ ,
\]
then (\ref{shift-5-a}) and (\ref{def-control}) give
\be
\label{shift-5-b}
\begin{array}{l}
\ds\lf\|\chi^{h(t)}_{\sigma\, R,x_0}\,\al(0)\int_\Sigma\chi^{h(0)}_{\sigma\, R,x_0}-\chi^{h(t)}_{\sigma\, R,x_0}\,\int_{\Sigma}\al(0)\ \chi^{h(0)}_{\sigma\, R,x_0}\ dvol_{h(0)}\rg\|_{L^{\infty}(\Sigma)}\\[5mm]
\ds \le C\, \sigma^2\,R^2\ \lf\|\al(0)(x)-\al(0)(y)\rg\|_{ L^\infty\lf(B_{\La(T)\,\sigma\,R}^{h(0)}(x_0)\times B_{\La(T)\, \sigma\,R}^{h(0)}(x_0)\rg)}\ .
\end{array}
\ee
We have using (\ref{der-conf-23-c}) and (\ref{deri-conf-fin})
\be
\label{shift-6}
\begin{array}{l}
\ds\sup_{t\in[0,T]}\lf| \int_\Sigma \chi^{h(t)}_{\sigma\, R,x_0}\ dvol_{h(t)}-\int_\Sigma \chi^{h(0)}_{\sigma\, R,x_0}\ dvol_{h(0)}\rg|\\[5mm]
\ds\le\int_0^T\lf|\frac{d}{dt}\int_\Sigma \chi^{h(t)}_{\sigma\, R,x_0}\ dvol_{h(t)}\rg|\ dt\ds=\sigma^{-1}\,R^{-1}\,\int_0^T\lf|\int_\Sigma \frac{d|x-x_0|_{h(t)}}{dt}\chi'\lf(\frac{|x-x_0|_{h(t)}}{\sigma\,R} \rg)\ dvol_{h(t)}\rg|\ dt\\[5mm]
\ds\le C\, \sigma^2\,R^2\,\int_0^T\lf\| \lf|\frac{dh}{dt}(t)\rg|_{h(t)}\rg\|_{L^\infty(B_{\sigma\,R}^{h(t)}(x_0))}\  dt\\[5mm]
 \ds\le  C_\Sigma\ \frac{e^{2\,\|\al\|_\infty}}{l^4_{\ast\ast}(T)}\ \sigma^2\,R^2\, \sqrt{T}\ \lf[\int_{\Sigma}|\vec{\mathbb I}_{\vec{\Phi}(0)}|^2\ dvol_{g_{\vec{\Phi}(0)}}\rg]^{1/2}\,\lf(W(\vec{\Phi}(0))-W(\vec{\Phi}(t))\rg)^{1/2}\ ,
\end{array}
\ee
Combining (\ref{shift-4-a}), (\ref{shift-5}), (\ref{shift-5-b}) and (\ref{shift-6}) gives
\be
\label{shift-7}
\begin{array}{l}
\ds\lf\|\chi^{h(t)}_{\sigma\, R,x_0}\,[\al(t)-\al(0)]\ \int_\Sigma\chi^{h(t)}_{\sigma\, R,x_0}\rg\|_{L^\infty(\Sigma)}\le\,   C\, \sigma^2\,R^2\ \int_{B^{h(t)}_{R}(x_0)}|\vec{\mathbb I}(t)|^2_{g_{\vec{\Phi}(t)}}\ dvol_{g_{\vec{\Phi}(t)}}\ \\[5mm]
\ds+\,\frac{C_{\Sigma,m}}{ l_\ast^{2}}\,\sigma^3\, R^2\, \lf[W(\vec{\Phi}(0))+2\pi\,(g(\Sigma)-1)\rg]\ \\[5mm] 
\ds+\,C_\Sigma\ \frac{e^{2\,\|\al\|_\infty}}{l_{\ast\ast}^9(T)}\ \sqrt{T}\     E(0)^{3/2}    \,\sqrt{W(\vec{\Phi}(0))-W(\vec{\Phi}(t))}\\[5mm]
\ds +\,C_\Sigma\ \frac{e^{2\,\|\al\|_\infty}}{l^4_{\ast\ast}(T)}\ \sqrt{T}\  \sup_{[0,T]}\,\int_{B_{\sigma\,R}^{h(t)}(x_0)}\ |\al|\ dvol_{h(t)}\   \sqrt{E(0)}\,\sqrt{W(\vec{\Phi}(0))-W(\vec{\Phi}(t))} \\[5mm]
\ds +\, C_\Sigma\ \frac{e^{2\,\|\al\|_\infty}}{l^4_{\ast\ast}(T)}\ \sqrt{T}\ \|\al(0)\|_\infty\ \sigma^2\, R^2\ \lf[\int_{\Sigma}|\vec{\mathbb I}_{\vec{\Phi}(0)}|^2\ dvol_{g_{\vec{\Phi}(0)}}\rg]^{1/2}\,\lf(W(\vec{\Phi}(0))-W(\vec{\Phi}(t))\rg)^{1/2}\ ,\\[5mm]
\ds +\,C_\Sigma\, \sigma^2\,R^2\ \lf\|\al(0)(x)-\al(0)(y)\rg\|_{ L^\infty\lf(B_{\La(T)\,\sigma\,R}^{h(0)}(x_0)\times B_{\La(T)\,\sigma\,R}^{h(0)}(x_0)\rg)}\ .
\end{array}
\ee
Assume now that for $R<l_{\ast\ast}/2$ there holds
\[
\sup_{t\in[0,T]}\sup_{x_0\in\Sigma}\,\int_{B_{R}^{h(t)}(x)}|\vec{\mathbb I}_{\vec{\Phi}(t)}|^2\ dvol_{g_{\vec{\Phi}(t)}}<\frac{8\pi}{3}\ .
\]
Then, the inequality (\ref{shift-7}) holds true for any $x_0\in\Sigma$, any $t\in[0,T]$ and for any $\sigma<1/8$. This implies 
\be
\label{shift-8}
\begin{array}{l}
\ds\sup_{t\in[0,T]}\lf\|\al(t)-\al(0)\rg\|_{L^\infty(\Sigma)}\le\,   C\, \sup_{t\in[0,T]}\sup_{x_0\in\Sigma}\int_{B^{h(t)}_{R}(p)}|\vec{\mathbb I}_{\vec{\Phi}(t)}|^2_{g_{\vec{\Phi}(t)}}\ dvol_{g_{\vec{\Phi}(t)}}\ \\[5mm]
\ds+\,\frac{C_{\Sigma,m}}{ l_\ast^{2}}\,\sigma\ \lf[W(\vec{\Phi}(0))+2\pi\,(g(\Sigma)-1)\rg]\ \\[5mm] 
\ds+\,C_\Sigma\ \frac{e^{2\,\|\al\|_\infty}}{l_{\ast\ast}^9(T)}\ \frac{\sqrt{T}}{\sigma^2\,R^2}\     E(0)^{3/2}    \,\sqrt{W(\vec{\Phi}(0))-W(\vec{\Phi}(t))}\\[5mm]
\ds +\,C_\Sigma\ \frac{e^{2\,\|\al\|_\infty}}{l^4_{\ast\ast}(T)}\ \frac{\sqrt{T}}{\sigma^2\,R^2}\  \sup_{[0,T]}\,\sup_{x_0\in\Sigma}\,\int_{B_{\sigma R}^{h(t)}(x_0)}\ |\al|\ dvol_{h(t)}\   \sqrt{E(0)}\,\sqrt{W(\vec{\Phi}(0))-W(\vec{\Phi}(t))} \\[5mm]
\ds +\, C_\Sigma\ \frac{e^{2\,\|\al\|_\infty}}{l^4_{\ast\ast}(T)}\ \frac{\sqrt{T}}{\sigma^2\,R^2}\   \|\al(0)\|_\infty\ \sigma^2\,R^2\ \lf[\int_{\Sigma}|\vec{\mathbb I}_{\vec{\Phi}(0)}|^2\ dvol_{g_{\vec{\Phi}(0)}}\rg]^{1/2}\,\lf(W(\vec{\Phi}(0))-W(\vec{\Phi}(t))\rg)^{1/2}\ ,\\[5mm]
\ds +\,C_\Sigma\,  \sup_{x_0\in\Sigma}\, \lf\|\al(0)(x)-\al(0)(y)\rg\|_{ L^\infty\lf(B_{\La(T)\,\sigma\,R}^{h(0)}(x_0)\times B_{\La(T)\,\sigma\,R}^{h(0)}(x_0)\rg)}\ .
\end{array}
\ee
We now estimate $\|\al(0)\|_\infty$. The mean value theorem gives the existence of $\ov{x}\in\Sigma$ such that
\be
\label{conf-init}
\al(0)(\ov{x})=\frac{1}{2}\,\log\frac{\ds\int_{\Sigma}\,e^{2\al(0)}\ dvol_{h(0)}}{\ds\int_{\Sigma} dvol_{h(0)}}
=\log\sqrt{\frac{\mbox{Vol}(\vec{\Phi}(0))}{V_{h(0)}}}\ .
\ee
Moreover for any $x_0\in \Sigma$, (\ref{shift-1}) is implying
\be
\label{shift-9}
\begin{array}{rl}
\ds\lf\|\al(0)(x)-\al(0)(y)\rg\|_{L^{\infty}\lf(B^{h(0)}_{8^{-1}\,R}(x_0)\times B^{h(0)}_{8^{-1}\,R}(x_0)\rg)}&\ds\le   C\, \int_{B^{h(0)}_{R}(x_0)}|\vec{\mathbb I}_{\vec{\Phi}(0)}|^2_{g_{\vec{\Phi}(0)}}\ dvol_{g_{\vec{\Phi}(0)}}\\[5mm]
\ds &\ds+\frac{C_{\Sigma,m}}{ l_\ast(0)^{2}}\, \lf[W(\vec{\Phi}(0))+2\pi\,(g(\Sigma)-1)\rg]
\end{array}
\ee
We consider a good covering of $\Sigma$ by balls $(B^{h(0)}_{8^{-1}\,R}(x_i))_{i\in I}$ (which are all convex since $R<l_\ast(0)/2$) in such a way that each point of $\Sigma$ is covered by a universally bounded number of balls. Obviously we have
\be
\label{shift-10}
\begin{array}{rl}
\ds\lf\|\al(0)(x)-\al(0)(y)\rg\|_{L^{\infty}(\Sigma\times \Sigma)}&\ds\le C\,   \sum_{i\in I}\int_{B^{h(0)}_{R}(x_i)}|\vec{\mathbb I}_{\vec{\Phi}(0)}|^2_{g_{\vec{\Phi}(0)}}\ dvol_{g_{\vec{\Phi}(0)}}\\[5mm]
\ds&\ds\ +\frac{C_{\Sigma,m}}{ l_\ast(0)^{2}}\, \mbox{Card}(I)\ \lf[W(\vec{\Phi}(0))+2\pi\,(g(\Sigma)-1)\rg]
\end{array}
\ee
We have also 
\be
\label{shift-11}
\mbox{Card}(I) \, R^2\le C_\Sigma\ V_{h(0)}\ .
\ee
Thus combining (\ref{shift-9}), (\ref{shift-10}) and (\ref{shift-11}) gives 
\be
\label{shift-12}
\|\al(0)\|_{L^\infty(\Sigma)}\le \frac{C_{\Sigma,m}}{ l_\ast(0)^{2}}\, R^{-2}\, E(0)+\log\sqrt{\frac{\mbox{Vol}(\vec{\Phi}(0))}{V_{h(0)}}}
\ee
Inserting (\ref{shift-12}) in (\ref{shift-8}) gives for any $\sigma<1/8$
\be
\label{shift-13}
\begin{array}{l}
\ds\sup_{t\in[0,T]}\lf\|\al(t)-\al(0)\rg\|_{L^\infty(\Sigma)}\le\,   C_{\Sigma,m}\, \sup_{t\in[0,T]}\sup_{x_0\in\Sigma}\int_{B^{h(t)}_{R}(p)}|\vec{\mathbb I}_{\vec{\Phi}(t)}|^2_{g_{\vec{\Phi}(t)}}\ dvol_{g_{\vec{\Phi}(t)}}+\,\frac{C_{\Sigma,m}}{ l_{\ast\ast}^{2}(T)}\,\sigma\ E(0)\ \\[5mm]
%\ds+\,\frac{C_{\Sigma,m}}{ l_\ast^{2}}\,\sigma\ \lf[W(\vec{\Phi}(0))+2\pi\,(g(\Sigma)-1)\rg]\ \\[5mm] 
\ds+\,C_{\Sigma,m}\ \frac{e^{2\,\|\al\|_\infty}}{l_{\ast\ast}^9(T)}\ \frac{\sqrt{T}}{\sigma^2\,R^2}\     E(0)^{3/2}    \,\sqrt{W(\vec{\Phi}(0))-W(\vec{\Phi}(t))}\\[5mm]
\ds +C_{\Sigma,m}\ \frac{e^{2\,\|\al\|_\infty}}{l^4_{\ast\ast}(T)}\ \sqrt{T}\ \log\sqrt{\frac{\mbox{Vol}(\vec{\Phi}(0))}{V_{h(0)}}}\ \   \sqrt{E(0)}\,\sqrt{W(\vec{\Phi}(0))-W(\vec{\Phi}(t))}\\[5mm]
\ds +\,C_{\Sigma,m}\,  \sup_{x_0\in\Sigma}\, \lf\|\al(0)(x)-\al(0)(y)\rg\|_{ L^\infty\lf(B_{\La(T)\,\sigma\,R}^{h(0)}(x_0)\times B_{\La(T)\,\sigma\,R}^{h(0)}(x_0)\rg)}\ .
\end{array}
\ee
This concludes the proof of lemma~\ref{lm-conf-fact-shift} .\hfill $\Box$

\medskip

\begin{Lm}
\label{lm-conf-fact-cnt} Let $R>0$ such that
\[
\sup_{t\in[0,T]}\sup_{x_0\in\Sigma}\int_{B^{h(t)}_{R}(p)}|\vec{\mathbb I}_{\vec{\Phi}(t)}|^2_{g_{\vec{\Phi}(t)}}\ dvol_{g_{\vec{\Phi}(t)}}<\frac{8\pi}{3}
\]
There exists a constant $\hat{C}_{\Sigma,m}>0$ depending only on $m$ and the topology of $\Sigma$ such that if
\[
\frac{e^{2\,\|\al(0)\|_\infty}}{R^2}\,\sqrt{T}<\hat{C}_{\Sigma,m}\, l^4_{\ast\ast}(T)\, \inf\lf\{1,\frac{R^{-2}}{\lf|\log\sqrt{\frac{\mbox{Vol}(\vec{\Phi}(0))}{V_{h(0)}}}\rg|}, \frac{l^{9}_{\ast\ast}(T)}{ E(0)^{4}} \rg\}
\]
then
\[
\sup_{t\in[0,T]}\lf\|\al(t)-\al(0)\rg\|_{L^\infty(\Sigma)}\le \hat{C}_{\Sigma,m}\ .
\]
\hfill $\Box$
\end{Lm}
\noindent{\bf Proof of lemma~\ref{lm-conf-fact-cnt}.} Let 
\[
A= 16\ C_{\Sigma,m}\, \frac{8\,\pi}{3}
\]
and
\[
\delta:=C_{\Sigma,m}\,e^{2A}\,\frac{e^{2\,\|\al(0)\|_\infty}}{l_{\ast\ast}^9(T)}\ \frac{\sqrt{T}}{\sigma^2\,R^2}\   
\]
where $ C_{\Sigma,m}$ is the constant in (\ref{conf-fact-shift}). Assume 
\be
\label{hyp}
\sup_{t\in[0,T]}\lf\|\al(t)-\al(0)\rg\|_{L^\infty(\Sigma)}\le A\ .
\ee
Choose  $\sigma$ such that
\[
 \sigma:=\frac{ 8\,\pi\,l_{\ast\ast}^{2}(T)}{3\,C_{\Sigma,m}\,E(0)}\ \min\{C_{\Sigma,m},1\}
\]
We have
\[
0<\log\,\La(T)\le C_\Sigma\  \delta\ \sigma^2\,R^2\, {l^5_{\ast\ast}(T)}
\]
Choose $T$ such that
\be
\label{cont-1}
\delta\ \sigma^2<C^{-1}_\Sigma\ {l^5_{\ast\ast}(T)}^{-1}\,\log^{-1}2\ \Longleftrightarrow\ \frac{e^{2\,\|\al(0)\|_\infty}}{R^2}\,\sqrt{T}<\ti{C}_{\Sigma,m}\ l^4_{\ast\ast}(T)
\ee
for some well chosen constant $\ti{C}_{\Sigma,m}$ depending only on the topology of $\Sigma$ and $m$. Hence
\[
0<\log\,\La(T)\le C_\Sigma\ R^2\, \log^{-1}2\
\]Using a covering of $B_{\La(T)\,\sigma\,R}^{h(0)}(x_0)$ by $\simeq\La(T)^{-2}$ balls of radius $\sigma\,R$, using (\ref{L2-1}) and (\ref{L2-2}) gives
\be
\label{cnt-1}
\begin{array}{l}
\ds\,C_{\Sigma,m}\lf\|\al(0)(x)-\al(0)(y)\rg\|_{ L^\infty\lf(B_{\La(T)\,\sigma\,R}^{h(0)}(x_0)\times B_{\La(T)\,\sigma\,R}^{h(0)}(x_0)\rg)}\\[5mm]
\ds\le \frac{C_{\Sigma,m}}{\La(T)^{2}}\,\lf[\int_{B^{h(0)}_{R}(x_0)}|\vec{\mathbb I}_{\vec{\Phi}(0)}|^2_{g_{\vec{\Phi}(0)}}\ dvol_{g_{\vec{\Phi}(t)}}  +\frac{C_{\Sigma,m}}{l_\ast^2(0)}\ \sigma\, E(0) \rg]\\[5mm]
\ds\le 4\, {C_{\Sigma,m}}\,\lf[ \frac{8\pi}{3}+\frac{C_{\Sigma,m}}{l_\ast^2(0)}\ \sigma\, E(0)\rg]\le \frac{A}{2}
\end{array}
\ee
Inserting these estimates and notations in (\ref{conf-fact-shift}) gives
\[
\sup_{t\in[0,T]}\lf\|\al(t)-\al(0)\rg\|_{L^\infty(\Sigma)}\le \frac{3\,A}{4}+\delta\, E(0)^2+\delta\, \sigma^2\,R^2\, {l^5_{\ast\ast}(T)}\,\log\sqrt{\frac{\mbox{Vol}(\vec{\Phi}(0))}{V_{h(0)}}}\
\]
Choosing $T$ small enough such that
\be
\label{cont-2}
\delta\, E(0)^2<\frac{A}{16}\quad\mbox{ and }\quad\delta\, \sigma^2\,R^2\, {l^5_{\ast\ast}(T)}\,\log\sqrt{\frac{\mbox{Vol}(\vec{\Phi}(0))}{V_{h(0)}}}<\frac{A}{16}\ ,
\ee
gives that
\[
\sup_{t\in[0,T]}\lf\|\al(t)-\al(0)\rg\|_{L^\infty(\Sigma)}\le \frac{7\,A}{8}\ .
\]
Since $\lf\|\al(t)-\al(0)\rg\|_{L^\infty(\Sigma)}$ is a continuous function of $t$ and equal to zero at zero we deduce that as long $T$ satisfies (\ref{cont-1}) and (\ref{cont-2}) we have that (\ref{hyp}) holds true and the lemma is proved.
%%%%%%%%%%%%%%%%%%%%%%%%%%%%%%%%%%%%%%%%%%%%%%%%%%%%%%%%%%%%%%%%%%%%%%%%%%%%%%%%%%%%%%%%%%%%%%%%%%%%%%%%%%
\section{A priori control of conformal class degenerations}
In this section we give a control of the conformal class of the evolving surface under the flow. The torus case is considered separately from the hyperbolic one and is treaded in the first subsection.
\subsection{The genus one case}
\begin{Lm}
\label{lm-length-control-torus}
Let $\Sigma=T^2$ be a closed oriented genus one surface. Under the above notations, there exists a constant  $ {C}(m)>0$ depending only on $m$ such that for any $l>0$ satisfying
\be
\label{length-control-torus}
 l^9<E(0)^4\quad,\quad l^2<\lf|\log\sqrt{\frac{\mbox{Vol}(\vec{\Phi}(0))}{V_{h(0)}}}\rg|^{-1}\quad\mbox{ and }\quad l \exp\lf[ {C}(m)\ l^{30} \rg]\le  {l_\ast(0)} ,
\ee
For any $T>0$ and $0<R<l$ satisfying
\[
\sup_{t\in[0,T]}\sup_{x_0\in\Sigma}\int_{B^{h(t)}_{R}(p)}|\vec{\mathbb I}_{\vec{\Phi}(t)}|^2_{g_{\vec{\Phi}(t)}}\ dvol_{g_{\vec{\Phi}(t)}}<\frac{8\pi}{3}
\]
and 
\[
\frac{e^{4\|\al(0)\|_\infty}}{R^4} \ T< \hat{C}_{\Sigma,m}^2\, \frac{l^{26}}{E(0)^8}\ .
\]
where $\hat{C}_{m}$ is given by lemma~\ref{lm-conf-fact-cnt}, then there holds
\[
l_{\ast\ast}(T):=\inf_{t\in[0,T]}l_\ast(t)>l\ \quad\mbox{ and }\quad\sup_{t\in[0,T]}\lf\|\al(t)-\al(0)\rg\|_{L^\infty(\Sigma)}\le \hat{C}_{\Sigma,m}\ .
\]
\hfill $\Box$
\end{Lm}
\noindent{\bf Proof of lemma~\ref{lm-length-control-torus}.} The torus $(T^2,h(t))$ is isometric to the flat torus $T^2_{a,b}:={\C}/{\Z}\oplus(a(t)+ib(t)){\Z}$ with $|a+ib|\ge 1$ equipped with $h(t)=b(t)^{-1}\ [dx_1^2+dx_2^2]$. The holomorphic quadratic form $dh/dt$ has the form
\[
\frac{dh}{dt}=\Re\lf[(c(t)+id(t))\ dz\otimes dz\rg]\ .
\]
and
\[
\lf|\frac{d\log b}{dt}\rg|=\lf|\frac{dh}{dt}\rg|_{h(t)}=b\ |c+id|=\lf|P_{h(t)}\lf(\pi_{\vec{n}}\frac{d\vec{\Phi}}{dt}\cdot\vec{\mathcal H}^0\rg)\rg|_{h(t)}
\]
Let $q(t)= {b(t)^{-1}} dz\otimes dz$. We have $|q(t)|^2_{h(t)}\equiv 1$. Hence $\|q(t)\|_{L^2(T^2_{a,b})}=1$ and
\[
\begin{array}{l}
\ds \lf|P_{h(t)}\lf(\pi_{\vec{n}}\frac{d\vec{\Phi}}{dt}\cdot\vec{\mathcal H}^0\rg)\rg|_{h(t)}=\lf|q(t)\ \int_{T^2_{a,b}}\lf<\pi_{\vec{n}}\frac{d\vec{\Phi}}{dt}\cdot\vec{\mathcal H}^0, q(t)\rg>_{h(t)}\ dvol_{h(t)}\rg|_{h(t)}\\[5mm]
\ds\le \int_{T^2_{a,b}}\ \lf|\pi_{\vec{n}}\frac{d\vec{\Phi}}{dt}\cdot\vec{\frak h}^0\rg|_{h(t)}\ dvol_{h(t)}\le e^{2\,\|\al(t)\|_\infty}\ \sqrt{E(0)}\ \sqrt{W(\vec{\Phi}(0))-W(\vec{\Phi}(T))}
\end{array}
\]
We have
\[
l_\ast(t)=\sqrt{b(t)^{-1}}
\]
Hence we finally obtain
\be
\label{genus-length}
\lf|\frac{d\log l_\ast(t)}{dt}\rg|\le e^{2\,\|\al(t)\|_\infty}\ \sqrt{E(0)}\ \sqrt{W(\vec{\Phi}(0))-W(\vec{\Phi}(T))}\ .
\ee
Let $T>0$, $l>0$ and $l>R>0$ such that
\[
\sup_{t\in[0,T]}\sup_{x_0\in\Sigma}\int_{B^{h(t)}_{R}(p)}|\vec{\mathbb I}_{\vec{\Phi}(t)}|^2_{g_{\vec{\Phi}(t)}}\ dvol_{g_{\vec{\Phi}(t)}}<\frac{8\pi}{3}
\]
and 
\[
\frac{e^{2\,\|\al(0)\|_\infty}}{R^2}\,\sqrt{T}<\hat{C}_{\Sigma,m}\, l^4\, \inf\lf\{1,\frac{R^{-2}}{\lf|\log\sqrt{\frac{\mbox{Vol}(\vec{\Phi}(0))}{V_{h(0)}}}\rg|}, \frac{l^{9}}{ E(0)^{4}} \rg\}
\]
Assume
\[
l<l_{\ast\ast}(T)\quad ,\quad l^9<E(0)^4\quad\mbox{ and }\quad l^2<\lf|\log\sqrt{\frac{\mbox{Vol}(\vec{\Phi}(0))}{V_{h(0)}}}\rg|^{-1}\ .
\]
From lemma~\ref{lm-conf-fact-cnt} we deduce
\[
\sup_{t\in[0,T]}\|\al(t)\|_{L^\infty(\Sigma)}\le \|\al(0)\|_\infty+ \hat{C}_{\Sigma,m}\ .
\]
Integrating (\ref{genus-length})  gives
\[
\begin{array}{rl}
\ds l_\ast(t)&\ds\ge l_\ast(0)\ \exp\lf[ -\,C(\Sigma)\ T\  l^{-3}\ e^{2\|\al(0)\|_\infty}\ e^{2\hat{C}_{\Sigma,m}}\ E(0)   \rg]\\[5mm]
\ds &\ds\ge  l_\ast(0)\ \exp\lf[ -\,C(\Sigma)\ l^{-3}\ R^4\ e^{-2\|\al(0)\|_\infty}\ e^{2\hat{C}_{\Sigma,m}}\ E(0)^{-7} \ l^{26} \rg]\\[5mm]
\ds &\ds\ge l_\ast(0)\  \exp\lf[ -\,\ti{C}(\Sigma,m)\  \ E(0)^{-7} \ l^{30} \rg]\\[5mm]
\ds &\ds\ge l_\ast(0)\  \exp\lf[ -\,\ti{C}(\Sigma,m)\  \ (4\pi)^{-7} \ l^{30} \rg]
\end{array}
\]
where we used the fact that $E(0)^{-7}<(4\pi)^{-7}$ thanks to Willmore inequality. Hence choosing $l$ satisfying respectively
\[
 l^9<E(0)^4\quad,\quad l^2<\lf|\log\sqrt{\frac{\mbox{Vol}(\vec{\Phi}(0))}{V_{h(0)}}}\rg|^{-1}\quad\mbox{ and }\quad l \exp\lf[ \ti{C}(\Sigma,m)\  (4\pi)^{-7}  \ l^{30} \rg]\le  {l_\ast(0)}
\]
and $R<l$ and $T$ such that respectively
\[
\frac{e^{4\|\al(0)\|_\infty}}{R^4} \ T< \hat{C}_{\Sigma,m}^2\, \frac{l^{26}}{E(0)^8}\,
\]
and
\[
\sup_{x_0\in\Sigma}\sup_{t\in[0,T]}\int_{B^{h(t)}_{2\,R}(x_0)}|\vec{\mathbb I}_{\vec{\Phi}(t)}|^2_{g_{\vec{\Phi}(t)}}\ dvol_{g_{\vec{\Phi}(t)}}<\frac{8\pi}{3}\ ,
\]
we have $$l_{\ast,\ast}(T)>l$$. The proof of lemma~\ref{lm-length-control-torus} is complete.\hfill $\Box$
\hfill $\Box$

\subsection{The hyperbolic case}
Thee purpose of the present subsection is to establish the following lemma
\begin{Lm}
\label{lm-length-control}
Let $\Sigma$ be a closed oriented surface of genus larger than one. Under the above notations, there exists a constant  $ {C}(\Sigma,m)>0$ depending only on $m$ and the topology of $\Sigma$ such that for any $l>0$ satisfying
\be
\label{length-control}
 l^9<E(0)^4\quad,\quad l^2<\lf|\log\sqrt{\frac{\mbox{Vol}(\vec{\Phi}(0))}{V_{h(0)}}}\rg|^{-1}\quad\mbox{ and }\quad l \exp\lf[ {C}(\Sigma,m)\ l^{27} \rg]\le  {l_\ast(0)}\ .
\ee
For any $T>0$ and $0<R<l$ satisfying
\[
\sup_{t\in[0,T]}\sup_{x_0\in\Sigma}\int_{B^{h(t)}_{R}(p)}|\vec{\mathbb I}_{\vec{\Phi}(t)}|^2_{g_{\vec{\Phi}(t)}}\ dvol_{g_{\vec{\Phi}(t)}}<\frac{8\pi}{3}
\]
and 
\[
\frac{e^{4\|\al(0)\|_\infty}}{R^4} \ T< \hat{C}_{\Sigma,m}^2\, \frac{l^{26}}{E(0)^8}\ .
\]
where $\hat{C}_{\Sigma,m}$ is given by lemma~\ref{lm-conf-fact-cnt}, then there holds
\[
l_{\ast\ast}(T):=\inf_{t\in[0,T]}l_\ast(t)>l\quad\mbox{ and }\quad\sup_{t\in[0,T]}\lf\|\al(t)-\al(0)\rg\|_{L^\infty(\Sigma)}\le \hat{C}_{\Sigma,m}
\]
\hfill $\Box$
\end{Lm}
Before proving lemma~\ref{lm-length-control} we recall the lemma 2.2 from \cite{RuTo}.
\begin{Lm}
\label{lm-ruto}
Let $\Sigma$ be a closed oriented surface of genus larger or equal than 2. Let $h(t)$ be a $C^1$ path of metric in a neighbourhood of $t=0$ such that at $t=0$
\be
\label{I.27}
\frac{dh}{dt}=\Re\lf[ (I-P_{h_t})(\Psi)  \rg]
\ee
where $\Psi$ is a smooth section of $\wedge^{(1,0)}T\Sigma\otimes\wedge^{(1,0)}T\Sigma$. Assume we have a Collar ${\mathcal C}$ in $(\Sigma, h(0))$ around a simple closed geodesic
of length $l<2\,\arcsinh(1)$. Then, at $t=0$
\be
\label{I.28}
\lf| \frac{d l}{dt}+ \frac{l^2}{ 16\pi^3}\,\Re\lf<\Psi, dz^2\rg>_{L^2({\mathcal C},h)}\rg|\le C\  l^2\ \|\Psi\|_{L^1(\Sigma,h)}\ .
\ee
where $C$ depends only on $\Sigma$ and $$<\Psi, dz^2 \ >_{L^2({\mathcal C},h)}=\int_{\mathcal C} <\Psi, dz\otimes dz>_{h}\ dvol_h$$ and where the collar ${\mathcal C}$ denotes the cylinder
$(-X(l),X(l))\times S^1$ with coordinates $(s,e^{i\,\theta})$ where
\[
X(l)=\frac{2\pi}{l}\,\lf(\frac{\pi}{2}-\arctan\lf(\sinh\lf(\frac{l}{2}\rg) \rg)\rg)
\]
and the metric $h$ in this collar is given by
\[
h=\rho^2\ [ds^2+d\theta^2]\ \quad\mbox{ with }\quad\rho(s)=\frac{l}{2\,\pi\, \cos\lf(\frac{ls}{2\pi}\rg)}
\]
and $z=s+\,i\,\theta$. \hfill $\Box$
\end{Lm}
\noindent{\bf Proof of lemma-\ref{lm-length-control}.}
Let $\vec{\Phi}_t$ be a smooth solution to the conformal parametrized Willmore Flow. We apply lemma~\ref{lm-ruto}
\be
\label{I.29}
\Psi:=\vec{\frak{h}}_t^{\,0}\cdot\delta\vec{W}\res h_t
\ee
We have
\[
\vec{\frak{h}}_t^{\,0}=\p_z\lf(e^{-2\la}\p_z\vec{\Phi}\rg)\ dz\otimes\p_{\ov{z}}
\]
Hence
\be
\label{I.30}
\Psi=\rho^2\,\p_z\lf(e^{-2\la}\p_z\vec{\Phi}\rg)\cdot\delta\vec{W}\ dz\otimes dz\ ,
\ee
We have
\[
<dz\otimes dz,dz\otimes dz>_h=\rho^{-4}
\]
from which we deduce in one hand
\be
\label{I.31}
\begin{array}{l}
\ds<\Psi, dz^2 \ >_{L^2({\mathcal C},h)}=\int_{\mathcal C}\p_z\lf(e^{-2\la}\p_z\vec{\Phi}\rg)\cdot\delta\vec{W}\ \,ds\,d\theta=-\frac{1}{2i}\,\int_{\mathcal C}\,\p_z\lf(e^{-2\la}\p_z\vec{\Phi}\rg)\cdot\frac{d\vec{\Phi}}{dt}\ 
dz\wedge d\ov{z}\\[5mm]
\ds= -\frac{1}{2i}\,\int_{\mathcal C}\frac{d\vec{\Phi}}{dt}\cdot d\lf(e^{-2\la}\p_z\vec{\Phi}\rg)\wedge d\ov{z}=-\frac{1}{2i}\,\int_{\p {\mathcal C}}\frac{d\vec{\Phi}}{dt}\cdot e^{-2\la}\p_z\vec{\Phi} \ d\ov{z}
+\frac{1}{2i}\,\int_{\mathcal C}e^{-2\la}\p_z\vec{\Phi}\cdot\frac{d\p_z\vec{\Phi}}{dt}\ 
dz\wedge d\ov{z}\\[5mm]
\ds=-\frac{1}{2i}\,\int_{\p {\mathcal C}}\frac{d\vec{\Phi}}{dt}\cdot e^{-2\la}\p_z\vec{\Phi} \ d\ov{z}=-\frac{1}{2i}\,\int_{\p {\mathcal C}}\vec{U}\cdot e^{-2\la}\p_z\vec{\Phi} \ d\ov{z}\\[5mm]
\end{array}
\ee
Hence
\[
\lf|<\Psi, dz^2 \ >_{L^2({\mathcal C},h)}\rg|\le 
\]
where we have used $\p_z\vec{\Phi}\cdot\p_z\vec{\Phi}=0$  and hence $\p_z\vec{\Phi}\cdot\frac{d\p_z\vec{\Phi}}{dt}=0$.  In the other hand
\be
\label{I.32}
\|\Psi\|_{L^1(\Sigma,h)}=\int_{\mathcal C} \, \lf|\vec{\cal H}_t^{\,0}\cdot\delta\vec{W}\rg|_h \ dvol_{h}
\ee
We have in the above local conformal coordinates in the collar
\[
\vec{\cal H}_t^{\,0}\cdot\delta\vec{W}=\frac{e^{-2\al}}{2}\, \p_{z^2}^2\vec{\Phi}\cdot\pi_{\vec{n}}\lf[\frac{d\vec{\Phi}}{dt}\rg]\ dz\otimes dz
\]
where $e^{2\al}=\rho^{-2}\, e^{2\la}$. Hence
\[
\lf|\vec{\cal H}_t^{\,0}\cdot\delta\vec{W}\rg|_h=\frac{e^{-2\la}}{2}\, \lf|\p_{z^2}^2\vec{\Phi}\cdot\pi_{\vec{n}}\lf[\frac{d\vec{\Phi}}{dt}\rg]\rg|=\lf|\vec{\frak h}_t^0\cdot\pi_{\vec{n}}\lf[\frac{d\vec{\Phi}}{dt}\rg]\rg|_{g_{\vec{\Phi}_t}}
\]
\be
\label{I.33}
\lf|\frac{d\log l}{dt}\rg|\le C\, l\ \int_{\p{\mathcal C}}\lf|\vec{U}\cdot e^{-2\la}\p_z\vec{\Phi} \rg|\ dl_h+C\ l\ \int_{\mathcal C} \, \lf|\vec{\frak h}_t^{\,0}\rg|_{g_{\vec{\Phi}_t}}\ 
\lf|\pi_{n_t}\frac{d\vec{\Phi}}{dt}\rg|\  dvol_{h}
\ee
Recall 
\[
\vec{U}\cdot e^{-2\la}\p_z\vec{\Phi}\ \p_{\ov{z}}=U^{0,1}_h
\]
We can choose a ``good cut''  between $(-X(l), +X(l))$ and $(-X(l)+1, +X(l)-1)$ such that
\be
\label{I.34}
 \int_{\p{\mathcal C}}\lf|\vec{U}\cdot e^{-2\la}\p_z\vec{\Phi} \rg|\ dl_h\le \int_{S^1}\int_{-X(l)}^{-X(l)+1} |U^{0,1}_h|\ dvol_h+ \int_{S^1}\int_{X(l)-1}^{X(l)} |U^{0,1}_h|\ dvol_h
\ee
Combining (\ref{I.33}), (\ref{I.34}) and lemma~\ref{lm-delbar-inv} gives
\be
\label{I.35}
\begin{array}{l}
\ds\lf|\frac{d\log l}{dt}\rg|\le C(\Sigma)\, l(t)\ \|U^{0,1}\|_{L^1(\Sigma)}+C(\Sigma)\ l\ e^{2\,\|\al(t)\|_\infty}\ \sqrt{E(0)}\ \sqrt{W(\vec{\Phi}(0))-W(\vec{\Phi}(t))}\\[5mm]
\ds\le  C(\Sigma)\,  l(t)^{-3}\ l\ e^{2\,\|\al(t)\|_\infty}\ \sqrt{E(0)}\ \sqrt{W(\vec{\Phi}(0))-W(\vec{\Phi}(t))}
\end{array}
\ee
Let $T>0$, $l>0$ and $l>R>0$ such that
\[
\sup_{t\in[0,T]}\sup_{x_0\in\Sigma}\int_{B^{h(t)}_{R}(p)}|\vec{\mathbb I}_{\vec{\Phi}(t)}|^2_{g_{\vec{\Phi}(t)}}\ dvol_{g_{\vec{\Phi}(t)}}<\frac{8\pi}{3}
\]
and 
\[
\frac{e^{2\,\|\al(0)\|_\infty}}{R^2}\,\sqrt{T}<\hat{C}_{\Sigma,m}\, l^4\, \inf\lf\{1,\frac{R^{-2}}{\lf|\log\sqrt{\frac{\mbox{Vol}(\vec{\Phi}(0))}{V_{h(0)}}}\rg|}, \frac{l^{9}}{ E(0)^{4}} \rg\}
\]
Assume
\[
l<l_{\ast\ast}(T)\quad ,\quad l^9<E(0)^4\quad\mbox{ and }\quad l^2<\lf|\log\sqrt{\frac{\mbox{Vol}(\vec{\Phi}(0))}{V_{h(0)}}}\rg|^{-1}\ .
\]
From lemma~\ref{lm-conf-fact-cnt} we deduce
\[
\sup_{t\in[0,T]}\|\al(t)\|_{L^\infty(\Sigma)}\le \|\al(0)\|+ \hat{C}_{\Sigma,m}
\]
Integrating (\ref{I.35}) gives
\[
\begin{array}{rl}
\ds l_\ast(t)&\ds\ge l_\ast(0)\ \exp\lf[ -\,C(\Sigma)\ T\  l^{-3}\ e^{2\|\al(0)\|_\infty}\ e^{2\hat{C}_{\Sigma,m}}\ E(0)   \rg]\\[5mm]
\ds &\ds\ge  l_\ast(0)\ \exp\lf[ -\,C(\Sigma)\ l^{-3}\ R^4\ e^{-2\|\al(0)\|_\infty}\ e^{2\hat{C}_{\Sigma,m}}\ E(0)^{-7} \ l^{26} \rg]\\[5mm]
\ds &\ds\ge l_\ast(0)\  \exp\lf[ -\,\ti{C}(\Sigma,m)\  \ E(0)^{-7} \ l^{27} \rg]\\[5mm]
\ds &\ds\ge l_\ast(0)\  \exp\lf[ -\,\ti{C}(\Sigma,m)\  \ (4\pi)^{-7} \ l^{27} \rg]
\end{array}
\]
where we used the fact that $E(0)^{-7}<(4\pi)^{-7}$ thanks to Willmore inequality. Hence choosing $l$ satisfying respectively
\[
 l^9<E(0)^4\quad,\quad l^2<\lf|\log\sqrt{\frac{\mbox{Vol}(\vec{\Phi}(0))}{V_{h(0)}}}\rg|^{-1}\quad\mbox{ and }\quad l \exp\lf[ \ti{C}(\Sigma,m)\  (4\pi)^{-7}  \ l^{27} \rg]\le  {l_\ast(0)}
\]
and $R<l$ and $T$ such that respectively
\[
\frac{e^{4\|\al(0)\|_\infty}}{R^4} \ T< \hat{C}_{\Sigma,m}^2\, \frac{l^{26}}{E(0)^8}\,
\]
and
\[
\sup_{x_0\in\Sigma}\sup_{t\in[0,T]}\int_{B^{h(t)}_{2\,R}(x_0)}|\vec{\mathbb I}_{\vec{\Phi}(t)}|^2_{g_{\vec{\Phi}(t)}}\ dvol_{g_{\vec{\Phi}(t)}}<\frac{8\pi}{3}\ ,
\]
we have $$l_{\ast,\ast}(T)>l$$. The proof of lemma~\ref{lm-length-control} is complete.\hfill $\Box$
%%%%%%%%%%%%%%%%%%%%%%%%%%%%%%%%%%%%%%%%%%%%%%%%%%%%%%%%%%%%%%%%%%%%%%%%%%%%%%%%%%%%%%%%%%%%%%%%%%%%%%%%%%%%%
\section{A priori energy concentration evolution and higher order derivatives control}
\reset
\begin{Lm}
\label{lm-conc-evol}
Let $T>0$ and $0<R<l_{\ast\ast}(T)/2$ such that
\[
E(2R,x_0,t):=\sup_{t\in[0,T]}\int_{B^{h(t)}_{2\,R}(x_0)}|\vec{\mathbb I}_{\vec{\Phi}(t)}|^2_{g_{\vec{\Phi}(t)}}\ dvol_{g_{\vec{\Phi}(t)}}<\frac{8\pi}{3}
\]
Then for any $\delta\in (0,1)$ we have
\be
\label{conc-evol}
\begin{array}{l}
\ds\int_{B_{R}^{h(T)}(x_0)}|\vec{H}(T)|^2 \ \chi^{h(T)}_{R,x_0}(x)\ \ dvol_{g_{\vec{\Phi}(T)}} +\int_0^T\int_{B_{R/2}^{h(t)}(x_0)}\ |\Delta_{g_{\vec{\Phi}(t)}}\vec{H}(t)|^2\  \ dvol_{g_{\vec{\Phi}(T)}}\ dt\\[5mm]
\ds\le\int_{B_{R}^{h(0)}(x_0)}|\vec{H}(0)|^2 \  \chi^{h(0)}_{R,x_0}(x)\ dvol_{h(0)}+ [\delta+E(2R,x_0,T)^{1/2}]\  \int_0^T\int_{B_{2\,R}^{h(t)}(x_0)}\ |\Delta_{h(t)}\vec{H}|^2\  \ dvol_{h(t)}\ dt\\[5mm]
\ds\quad+C(\Sigma,l_{\ast\ast}(T),E(0))\ \lf[\frac{e^{4\|\al\|_\infty}}{\delta^3\,R^4}\ T+\sqrt{\frac{e^{4\|\al\|_\infty}}{2\,\delta^2\,R^4}\ T}+\frac{e^{4\|\al\|_\infty}}{\delta^3}\ \lf[ \frac{T}{R^4} \rg]^{1/3}\rg]\ E(2R,x_0,T)\\[5mm]
 \ds\quad+C(\Sigma,l_{\ast\ast}(T),E(0))\ \lf[\sqrt{E(2R,x_0,T)}\ \sqrt{W(\vec{\Phi}(0)-W(\vec{\Phi}(T)}\rg]\ E(2R,x_0,T)\\[5mm]
\ds\quad+C(\Sigma,l_{\ast\ast}(T),E(0))\ \lf[ e^{2\|\al\|_\infty}\ \sqrt{T}\ \sqrt{W(\vec{\Phi}(0)-W(\vec{\Phi}(T)} \rg]\ E(2R,x_0,T)\\[5mm]
\ds\quad+ C(\Sigma,l_{\ast\ast}(T),E(0))\ \lf[\frac{e^{8\|\al\|_\infty}}{\delta}\ \lf[W(\vec{\Phi}(0)-W(\vec{\Phi}(T)\rg] \rg]\ E(2R,x_0,T)\ .
\end{array}
\ee

\hfill $\Box$
\end{Lm}
\noindent{\bf Proof of lemma~\ref{lm-conc-evol}} For any $x_0\in \Sigma$, we take the scalar product in ${\R}^m$ of the first equation of the Willmore flow system (\ref{I.25}) with
 $$
 \Delta_{g_{\vec{\Phi}}}\lf[\chi^{h(t)}_{R,x_0}(x)\ e^{2\al}\,\vec{H}\rg]=2^{-1}\,\Delta_{g_{\vec{\Phi}}}\lf[\chi^{h(t)}_{R,x_0}(x)\ \Delta_{h(t)}\vec{\Phi}\rg]
 $$ 
 at an arbitrary time $t\in[0,T]$ where $\chi^{h(t)}_{R,x_0}(x)$ is the non negative  cut-off function supported on $B^{h(t)}_{R}(p)$ and equal to $1$ on $B^{h(t)}_{R/2}(p)$. 
 
%%%%%%%%%%%%%%%%%%%%%%%%%%%%%%%%%%%%%%%%%%%%%%%%%%%%%%%%%%%%%%%%%%%%%%%%%%%%%%%%%%%%%%%%%%%%%%%%%%%%%%%%%%%%% 
 We take the scalar product in $\R^m$ of the first equation of the Willmore parametrised flow with  $\Delta_{g_{\vec{\Phi}(t)}}\lf[\chi^{h(t)}_{R,x_0}(x)\  \vec{H}\rg]$ and we integrate with respect to
 $dvol_{g_{\vec{\Phi}(t)}}$. Since for any function $f$
\[
 \Delta_{g_{\vec{\Phi}(t)}}f\ dvol_{g_{\vec{\Phi}(t)}}= \Delta_{h(t)}f\ dvol_{h(t)}
 \]
this gives
 \be
\label{n-ene-con-1}
\begin{array}{l}
\ds\int_\Sigma \Delta_{h(t)}\lf[\chi^{h(t)}_{R,x_0}(x)\  \vec{H}\rg]\cdot\frac{\p\vec{\Phi}}{\p t}\ dvol_{h(t)}\\[5mm]
\ds=-\, 4\, \Re \int_\Sigma\Delta_{g_{\vec{\Phi}}}\lf[\chi^{h(t)}_{R,x_0}(x)\ \vec{H}\rg]\cdot \ov{\p}\lf[\p\vec{H}+|\vec{H}|^2\,\p\vec{\Phi}+2\,\vec{H}\cdot\vec{\frak{h}}^{\,0}\res \ov{\p}\vec{\Phi}  \rg]\\[5mm]
\ds+\int_\Sigma\Delta_{h(t)}\lf[\chi^{h(t)}_{R,x_0}(x)\ \vec{H}\rg]\cdot\vec{U}\, \ dvol_{h(t)}
\end{array}
\ee
 We have using local coordinates in $B^{h(t)}_{R}(x_0)$ (which are conformal at time $t$ but not necessary conformal in a neighbourhood of $t$ of course) and the fact that $dh/dt$ is trace free
and in particular $d/dt(\mbox{det}(h^{kl}))=0$ as well as $d/dt(\mbox{det}(h_{kl}))=0$
\be
\label{n-ene-con-2}
\begin{array}{l}
\ds\int_\Sigma \Delta_{h(t)}\lf[\chi^{h(t)}_{R,x_0}(x)\ \vec{H}\rg]\cdot\frac{\p\vec{\Phi}}{\p t}\ dvol_{h(t)}\ds=\frac{1}{2}\int_\Sigma \chi^{h(t)}_{R,x_0}(x)\ e^{-2\al}\ \Delta_{h(t)}\vec{\Phi}\cdot\Delta_{h(t)}\lf[\frac{\p\vec{\Phi}}{\p t}\rg]\ dvol_{h(t)}\\[5mm]
\ds=\frac{1}{4}\int_\Sigma\ e^{-2\al}\  \chi^{h(t)}_{R,x_0}(x)\frac{d}{dt} \lf|\Delta_{h(t)}\vec{\Phi}\rg|^2\ dvol_{h(t)}\\[5mm]
\ds\ -\frac{1}{2}\int_\Sigma \ e^{-2\al}\ \chi^{h(t)}_{R,x_0}(x)\ \Delta_{h(t)}\vec{\Phi}\cdot\, \p_{x_i}\lf(\frac{dh^{ij}}{dt} \sqrt{det(h_{kl})}\ 
\p_{x_j}\vec{\Phi}\rg)\ dx^2\\[5mm]
\ds =\frac{1}{4}\frac{d}{dt}\int_\Sigma \chi^{h(t)}_{R,x_0}(x)\,\ e^{-2\al}\  \lf|\Delta_{h(t)}\vec{\Phi}\rg|^2\ dvol_{h(t)}-\frac{1}{4}\int_\Sigma \frac{d\chi^{h(t)}_{R,x_0}(x)}{dt}\, \ e^{-2\al}\ \lf|\Delta_{h(t)}\vec{\Phi}\rg|^2\ dvol_{h(t)}\\[5mm]
\ds+\frac{1}{2}\,\int_\Sigma \chi^{h(t)}_{R,x_0}(x)\,\ e^{-2\al}\  \frac{d\al}{dt}\ \lf|\Delta_{h(t)}\vec{\Phi}\rg|^2\ dvol_{h(t)}\\[5mm]
\ds+\int_\Sigma\ e^{-2\al}\ \frac{dh}{dt}\res_2 d\lf[ \chi^{h(t)}_{R,x_0}(x)\ e^{2\al}\,\vec{H}\rg]\,\dot{\otimes}\,d\vec{\Phi}\ dvol_{h(t)}
\end{array}
\ee
We write
\be
\label{n-ene-con-3}
\begin{array}{l}
\ds -\, 4\, \Re \int_\Sigma\Delta_{g_{\vec{\Phi}}}\lf[\chi^{h(t)}_{R,x_0}(x)\,\vec{H}\rg]\cdot \ov{\p}\lf[\p\vec{H}+|\vec{H}|^2\,\p\vec{\Phi}+2\,\vec{H}\cdot\vec{\frak{h}}^{\,0}\res \ov{\p}\vec{\Phi}  \rg]\\[5mm]
\ds =-\, 4\, \Re \int_\Sigma\chi^{h(t)}_{R,x_0}(x)\ \Delta_{g_{\vec{\Phi}}}\vec{H}\cdot \ov{\p}\lf[\p\vec{H}+|\vec{H}|^2\,\p\vec{\Phi}+2\,\vec{H}\cdot\vec{\frak{h}}^{\,0}\res \ov{\p}\vec{\Phi}  \rg]\\[5mm]
\ds \ -\, 4\, \Re \int_\Sigma\Delta_{g_{\vec{\Phi}}}\lf[\chi^{h(t)}_{R,x_0}(x)\ \rg]\,\vec{H}\cdot \ov{\p}\lf[\p\vec{H}+|\vec{H}|^2\,\p\vec{\Phi}+2\,\vec{H}\cdot\vec{\frak{h}}^{\,0}\res \ov{\p}\vec{\Phi}  \rg]\\[5mm]
\ds-\, 8\, \Re \int_\Sigma \lf<d\lf[\chi^{h(t)}_{R,x_0}(x)\ \rg],d\vec{H}\rg>_{g_{\vec{\Phi}}}\cdot \ov{\p}\lf[\p\vec{H}+|\vec{H}|^2\,\p\vec{\Phi}+2\,\vec{H}\cdot\vec{\frak{h}}^{\,0}\res \ov{\p}\vec{\Phi}  \rg]
\end{array}
\ee
We finally use the fact that $\vec{U}\cdot\p\vec{\Phi}=U^{0,1}\res g_{\vec{\Phi}}$ and $\vec{U}\cdot\ov{\p}\vec{\Phi}=\ov{U^{0,1}}\res g_{\vec{\Phi}}$ and write
\be
\label{uuu}
\begin{array}{rl}
\ds\vec{U}&\ds=\lf({\ov{\p}}\vec{\Phi}\otimes \p\vec{\Phi}\cdot\vec{U}\rg)\res g^{-1}_{\vec{\Phi}}+\lf({{\p}}\vec{\Phi}\otimes {\ov{\p}}\vec{\Phi}\cdot\vec{U}\rg)\res g^{-1}_{\vec{\Phi}}\\[5mm]
&\ds=\lf({\ov{\p}}\vec{\Phi}\otimes U^{0,1}\res g_{\vec{\Phi}}\rg)\res g^{-1}_{\vec{\Phi}}+\lf({{\p}}\vec{\Phi}\otimes \ov{U^{0,1}}\res g_{\vec{\Phi}}\rg)\res g^{-1}_{\vec{\Phi}}\\[5mm]
&\ds={\ov{\p}}\vec{\Phi}\res U^{0,1}+{{\p}}\vec{\Phi}\res \ov{U^{0,1}}
\end{array}
\ee
Hence
\be
\label{n-ene-con-4}
\begin{array}{rl}
\ds\int_\Sigma\Delta_{h(t)}\lf[\chi^{h(t)}_{R,x_0}(x)\ \vec{H}\rg]\cdot\vec{U}\, \ dvol_{h(t)}&\ds=4\,\int_\Sigma \p\ov{\p} \lf[\chi^{h(t)}_{R,x_0}(x)\ \vec{H}\rg]\cdot\vec{U}\\[5mm]
&\ds=8\,\Re\int_\Sigma \p\ov{\p} \lf[\chi^{h(t)}_{R,x_0}(x)\ \vec{H}\rg]\cdot{\ov{\p}}\vec{\Phi}\res U^{0,1}\\[5mm]
&\ds=-\,8\,\Re\int_\Sigma \ov{\p} \lf[\chi^{h(t)}_{R,x_0}(x)\ \vec{H}\rg]\dot{\wedge}\p\lf[{\ov{\p}}\vec{\Phi}\res U^{0,1}\rg]
\end{array}
\ee
We have
\be
\label{n-ene-con-5}
\begin{array}{l}
\ds\p\lf[{\ov{\p}}\vec{\Phi}\res U^{0,1}\rg]=\p U^{0,1}\res\ov{\p}\vec{\Phi}+\vec{H}\,U^{0,1}\res g_{\vec{\Phi}}
\ds=2\,(I-\ti{P}_h)\lf(\pi_{\vec{n}}\frac{\p\vec{\Phi}}{\p t}\cdot\vec{\frak{h}}_0\rg)\res\ov{\p}\vec{\Phi}+\vec{H}\,U^{0,1}\res g_{\vec{\Phi}}
\end{array}
\ee
This gives
\be
\label{n-ene-con-6}
\begin{array}{l}
\ds\int_\Sigma\Delta_{h(t)}\lf[\chi^{h(t)}_{R,x_0}(x)\ \vec{H}\rg]\cdot\vec{U}\, \ dvol_{h(t)}\\[5mm]
\ds=-\,16\,\Re\int_\Sigma\chi^{h(t)}_{R,x_0}(x)\ \pi_T\lf(\ov{\p}\vec{H}\rg)\dot{\wedge}\lf[\ov{\p}\vec{\Phi}\res(I-\ti{P}_h)\lf(\pi_{\vec{n}}\frac{\p\vec{\Phi}}{\p t}\cdot\vec{\frak{h}}_0\rg)\rg]\\[5mm]
\ds\ \ \ -\,8\,\Re\int_\Sigma\vec{H}\cdot\ov{\p} \lf[\chi^{h(t)}_{R,x_0}(x)\ \vec{H}\rg]\wedge\lf(U^{0,1}\res g_{\vec{\Phi}}\rg)
\end{array}
\ee
Combining (\ref{n-ene-con-2}), (\ref{n-ene-con-3}) and (\ref{n-ene-con-6}) gives
\be
\label{re-n-ene-con-7}
\begin{array}{l}
\ds\frac{1}{4}\frac{d}{dt}\int_\Sigma \chi^{h(t)}_{R,x_0}(x)\, e^{-2\al}\,\lf|\Delta_{h(t)}\vec{\Phi}\rg|^2\ dvol_{h(t)}+4\, \int_\Sigma\chi^{h(t)}_{R,x_0}(x)\ \, e^{-2\al}\,|\Delta_{h(t)}\vec{H}|^2\ dvol_{h(t)}\\[5mm]
\ds=\frac{1}{4}\int_\Sigma \frac{d\chi^{h(t)}_{R,x_0}(x)}{dt}\,\, e^{-2\al}\, \lf|\Delta_{h(t)}\vec{\Phi}\rg|^2\ dvol_{h(t)}-\int_\Sigma\frac{dh}{dt}\res_2 d\lf[ \chi^{h(t)}_{R,x_0}(x)\ \,\vec{H}\rg]\,\dot{\otimes}\,d\vec{\Phi}\ dvol_{h(t)}
\\[5mm]
\ds-\frac{1}{2}\,\int_\Sigma \chi^{h(t)}_{R,x_0}(x)\,\ e^{-2\al}\  \frac{d\al}{dt}\ \lf|\Delta_{h(t)}\vec{\Phi}\rg|^2\ dvol_{h(t)}\\[5mm]
\ds-\, 4\, \Re \int_\Sigma\chi^{h(t)}_{R,x_0}(x)\ \, e^{-2\al}\,\Delta_{h(t)}\vec{H}\cdot \ov{\p}\lf[|\vec{H}|^2\,\p\vec{\Phi}+2\,\vec{H}\cdot\vec{\frak{h}}^{\,0}\res \ov{\p}\vec{\Phi}  \rg]\\[5mm]
\ds \ -\, 4\, \Re \int_\Sigma\Delta_{g_{\vec{\Phi}}}\lf[\chi^{h(t)}_{R,x_0}(x)\rg]\,\vec{H}\cdot \ov{\p}\lf[\p\vec{H}+|\vec{H}|^2\,\p\vec{\Phi}+2\,\vec{H}\cdot\vec{\frak{h}}^{\,0}\res \ov{\p}\vec{\Phi}  \rg]\\[5mm]
\ds-\, 8\, \Re \int_\Sigma \lf<d\chi^{h(t)}_{R,x_0}(x)\ ,d\vec{H}\rg>_{g_{\vec{\Phi}}}\cdot \ov{\p}\lf[\p\vec{H}+|\vec{H}|^2\,\p\vec{\Phi}+2\,\vec{H}\cdot\vec{\frak{h}}^{\,0}\res \ov{\p}\vec{\Phi}  \rg]\\[5mm]
\ds-\,16\,\Re\int_\Sigma\chi^{h(t)}_{R,x_0}(x)\ \pi_T\lf(\ov{\p}\vec{H}\rg)\dot{\wedge}\lf[\ov{\p}\vec{\Phi}\res(I-\ti{P}_h)\lf(\pi_{\vec{n}}\frac{\p\vec{\Phi}}{\p t}\cdot\vec{\frak{h}}_0\rg)\rg]\\[5mm]
\ds -\,8\,\Re\int_\Sigma\vec{H}\cdot\ov{\p} \lf[\chi^{h(t)}_{R,x_0}(x)\ \vec{H}\rg]\wedge\lf(U^{0,1}\res g_{\vec{\Phi}}\rg)
\end{array}
\ee
Now, from lemma~\ref{lm-param-conf-var}
\be
\label{param-conf-var-re}
\frac{d\al}{dt}=\,\frac{d\vec{\Phi}}{dt}\cdot \vec{H}+e^{-\,2\al(t)}\frac{1}{2}\,d^{\ast_{h(t)}}\lf[\frac{d\vec{\Phi}}{dt}\cdot d\vec{\Phi}(t)\rg]\ .
\ee
Hence
\be
\label{re-n-1}
\begin{array}{l}
\ds-\frac{1}{2}\,\int_\Sigma \chi^{h(t)}_{R,x_0}(x)\,\ e^{-2\al}\  \frac{d\al}{dt}\ \lf|\Delta_{h(t)}\vec{\Phi}\rg|^2\ dvol_{h(t)}\\[5mm]
\ds=-2\,\int_\Sigma \chi^{h(t)}_{R,x_0}(x)\,\pi_{\vec{n}}\frac{d\vec{\Phi}}{dt}\cdot \vec{H}\ |\vec{H}|^2\ \ dvol_{g_{\vec{\Phi}(t)}}\\[5mm]
\ds-\int_\Sigma \chi^{h(t)}_{R,x_0}(x)\,e^{-\,2\al(t)}\,d^{\ast_{h(t)}}\lf[\frac{d\vec{\Phi}}{dt}\cdot d\vec{\Phi}(t)\rg]\  |\vec{H}|^2\ \ dvol_{g_{\vec{\Phi}(t)}}\ .
\end{array}
\ee
We have
\be
\label{re-n-2}
\begin{array}{l}
\ds-\int_\Sigma \chi^{h(t)}_{R,x_0}(x)\,e^{-\,2\al(t)}\,d^{\ast_{h(t)}}\lf[\frac{d\vec{\Phi}}{dt}\cdot d\vec{\Phi}(t)\rg]\  |\vec{H}|^2\ \ dvol_{g_{\vec{\Phi}(t)}}\\[5mm]
\ds=-\int_\Sigma \ \lf<d\chi^{h(t)}_{R,x_0}(x),\lf[\frac{d\vec{\Phi}}{dt}\cdot d\vec{\Phi}(t)\rg]\rg>_{h(t)}\  |\vec{H}|^2\ \ dvol_{h(t)}\\[5mm]
\ds-\int_\Sigma \ \chi^{h(t)}_{R,x_0}(x)\  \lf<d|\vec{H}|^2,\lf[\vec{U}\cdot d\vec{\Phi}(t)\rg]\rg>_{h(t)} \ dvol_{h(t)}
\end{array}
\ee
From (\ref{u}) we have $\vec{U}={\ov{\p}}\vec{\Phi}(t)\res U^{0,1}+{{\p}}\vec{\Phi}(t)\res \ov{U^{0,1}}$ hence
\[
\vec{U}\cdot d\vec{\Phi}(t)=\vec{U}\cdot \lf[\p\vec{\Phi}+\ov{\p}\vec{\Phi}\rg]=2\,\Re\lf[U^{0,1}\res g_{\vec{\Phi}}\rg].
\]
this gives
\be
\label{re-n-2-a}
\begin{array}{l}
\ds-\int_\Sigma \ \chi^{h(t)}_{R,x_0}(x)\  \lf<d|\vec{H}|^2,\lf[\vec{U}\cdot d\vec{\Phi}(t)\rg]\rg>_{h(t)} \ dvol_{h(t)}\\[5mm]
\ds=-4\, \Re\lf[ \int_\Sigma \ \chi^{h(t)}_{R,x_0}(x)\  \vec{H}\cdot\lf<d \vec{H},(U^{0,1}\res g_{\vec{\Phi}})\rg>_{h(t)} \ dvol_{h(t)}\rg]\ .
\end{array}
\ee
Hence
\be
\label{re-n-3}
\begin{array}{l}
\ds-\int_\Sigma \chi^{h(t)}_{R,x_0}(x)\,e^{-\,2\al(t)}\,d^{\ast_{h(t)}}\lf[\frac{d\vec{\Phi}}{dt}\cdot d\vec{\Phi}(t)\rg]\  |\vec{H}|^2\ \ dvol_{g_{\vec{\Phi}(t)}}\\[5mm]
\ds =-\int_\Sigma \ \lf<d\chi^{h(t)}_{R,x_0}(x),\lf[\frac{d\vec{\Phi}}{dt}\cdot d\vec{\Phi}(t)\rg]\rg>_{h(t)}\  |\vec{H}|^2\ \ dvol_{h(t)}\\[5mm]
\ds-4\, \Re\lf[ \int_\Sigma \ \chi^{h(t)}_{R,x_0}(x)\  \vec{H}\cdot\lf<d \vec{H},(U^{0,1}\res g_{\vec{\Phi}})\rg>_{h(t)} \ dvol_{h(t)}\rg]\ .
\end{array}
\ee
Combining (\ref{re-n-1}), (\ref{re-n-2}) and (\ref{re-n-3}) with (\ref{re-n-ene-con-7}) and integrating between $0$ and $T$ is giving
\be
\label{n-ene-con-7}
\begin{array}{l}
\ds\frac{1}{4}\int_\Sigma \chi^{h(T)}_{R,x_0}(x)\, e^{-2\al}\,\lf|\Delta_{h(T)}\vec{\Phi}\rg|^2\ dvol_{h(T)}+4\, \int_0^T\int_\Sigma\chi^{h(t)}_{R,x_0}(x)\ \, e^{-2\al}\,|\Delta_{h(t)}\vec{H}|^2\ dvol_{h(t)}\\[5mm]
\ds=\frac{1}{4}\int_\Sigma \chi^{h(0)}_{R,x_0}(x)\, e^{-2\al}\,\lf|\Delta_{h(0)}\vec{\Phi}\rg|^2\ dvol_{h(0)}\ dt
\ds+\frac{1}{4}\int_0^T\int_\Sigma \frac{d\chi^{h(t)}_{R,x_0}(x)}{dt}\,\, e^{-2\al}\, \lf|\Delta_{h(t)}\vec{\Phi}\rg|^2\ dvol_{h(t)}\ dt\\[5mm]
\ds-\int_0^T\int_\Sigma\frac{dh}{dt}\res_2 d\lf[ \chi^{h(t)}_{R,x_0}(x)\ \,\vec{H}\rg]\,\dot{\otimes}\,d\vec{\Phi}\ dvol_{h(t)}\ dt\\[5mm]
\ds-2\,\int_0^T\int_\Sigma \chi^{h(t)}_{R,x_0}(x)\,\pi_{\vec{n}}\frac{d\vec{\Phi}}{dt}\cdot \vec{H}\ |\vec{H}|^2\ \ dvol_{g_{\vec{\Phi}(t)}}\ dt\\[5mm]
\ds-\int_0^T\int_\Sigma \ \lf<d\chi^{h(t)}_{R,x_0}(x),\lf[\frac{d\vec{\Phi}}{dt}\cdot d\vec{\Phi}(t)\rg]\rg>_{h(t)}\  |\vec{H}|^2\ \ dvol_{h(t)}\ dt\\[5mm]
\ds-\int_0^T\int_\Sigma \ \chi^{h(t)}_{R,x_0}(x)\  \lf<d|\vec{H}|^2,\lf[\vec{U}\cdot d\vec{\Phi}(t)\rg]\rg>_{h(t)} \ dvol_{h(t)}\ dt\\[5mm]
%\ds-\int_\Sigma \ \lf<d\chi^{h(t)}_{R,x_0}(x),\lf[\frac{d\vec{\Phi}}{dt}\cdot d\vec{\Phi}(t)\rg]\rg>_{h(t)}\  |\vec{H}|^2\ \ dvol_{h(t)}\\[5mm]
%\ds-\frac{1}{2}\,\int_0^T\int_\Sigma \chi^{h(t)}_{R,x_0}(x)\,\ e^{-2\al}\  \frac{d\al}{dt}\ \lf|\Delta_{h(t)}\vec{\Phi}\rg|^2\ dvol_{h(t)}\ dt\\[5mm]
\ds-\, 4\, \Re \int_0^T\int_\Sigma\chi^{h(t)}_{R,x_0}(x)\ \, e^{-2\al}\,\Delta_{h(t)}\vec{H}\cdot \ov{\p}\lf[|\vec{H}|^2\,\p\vec{\Phi}\,+2\,\vec{H}\cdot\vec{\frak{h}}^{\,0}\res \ov{\p}\vec{\Phi}  \rg]\ dt\\[5mm]
\ds \ -\, 4\, \Re\int_0^T\ \int_\Sigma\Delta_{g_{\vec{\Phi}}}\lf[\chi^{h(t)}_{R,x_0}(x)\rg]\,\vec{H}\cdot \ov{\p}\lf[\p\vec{H}+|\vec{H}|^2\,\p\vec{\Phi}+2\,\vec{H}\cdot\vec{\frak{h}}^{\,0}\res \ov{\p}\vec{\Phi}  \rg]\ dt\\
\end{array}
\ee
We shall now be bounding  each of the terms in the r.h.s. of (\ref{n-ene-con-7}). First, using (\ref{der-conf-23-d}) we obtain
\be
\label{n-ene-con-8}
\begin{array}{l}
\ds\lf|\int_0^T\frac{1}{4}\int_\Sigma \frac{d\chi^{h(t)}_{R,x_0}(x)}{dt}\, \ e^{-2\al}\ \lf|\Delta_{h(t)}\vec{\Phi}\rg|^2\ dvol_{h(t)}\rg|\\[5mm]
\ds\le\frac{1}{4\,R}\lf|\int_0^T\int_\Sigma\frac{d|x-x_0|_{h(t)}}{dt}\chi'\lf(\frac{ |x-x_0|_{h(t)}}{R}  \rg)\, \ e^{-2\al}\  \lf|\Delta_{h(t)}\vec{\Phi}\rg|^2\ dvol_{h(t)}\rg|\\[5mm]
\ds\le \sup_{[0,T]}\int_{B_R^{h(t)}(x_0)} \lf|\Delta_{g_{\vec{\Phi}}}\vec{\Phi}\rg|^2\ dvol_{g_{\vec{\Phi}}}\  \int_0^T\lf\| \lf|\frac{dh}{dt}(t)\rg|_{h(t)}\rg\|_{L^\infty(B_R^{h(t)}(x_0))}\  dt\\[5mm]
\ds\le  C_\Sigma\ \frac{e^{2\,\|\al\|_\infty}}{l^4_{\ast\ast}(T)}\ \sqrt{T}\ \sup_{[0,T]}\int_{B_R^{h(t)}(x_0)}  \lf|\Delta_{g_{\vec{\Phi}}}\vec{\Phi}\rg|^2\ dvol_{g_{\vec{\Phi}}}\  \sqrt{E(0)}\,\sqrt{W(\vec{\Phi}(0))-W(\vec{\Phi}(t))}\ .
\end{array}
\ee
Then, observing that $\vec{H}\,\dot{\otimes}\,d\vec{\Phi}=0$ and, using local conformal coordinates
\be
\label{n-dh-tan}
\begin{array}{l}
|d\vec{H}\dot{\otimes}\,d\vec{\Phi}|^2_{h(t)}\le e^{-4\nu} \sum_{i,j=1}^2|\p_{x_i}\vec{H}\cdot\p_{x_j}\vec{\Phi}|^2= e^{-4\nu} \sum_{i,j=1}^2|\vec{H}\cdot\p_{x_i}\p_{x_j}\vec{\Phi}|^2\\[5mm]
\ds=e^{-4\nu} \sum_{i,j=1}^2|\vec{H}\cdot\vec{\mathbb I}_{ij}|^2=|\vec{H}\cdot\vec{\mathbb I}_{\vec{\Phi}}|^2_{h(t)}
\end{array}
\ee
 we have
\be
\label{n-ene-con-9}
\begin{array}{l}
\ds\lf|\int_0^T\int_\Sigma\frac{dh}{dt}\res_2 d\lf[ \chi^{h(t)}_{R,x_0}(x)\ \vec{H}\rg]\,\dot{\otimes}\,d\vec{\Phi}\ dvol_{h(t)}\ dt\rg|\\[5mm]
\ds\le\int_0^T\int_{\Sigma}\lf|\frac{dh}{dt}\rg|_{h(t)} |\chi^{h(t)}_{R,x_0}(x)|\ |\vec{H}|\ |\vec{\mathbb I}_{\vec{\Phi}(t)}|_{g_{\vec{\Phi}}}\ dvol_{g_{\vec{\Phi}(t)}}\\[5mm]
\ds\le\,\sup_{[0,T]}\int_{B_R^{h(t)}(x_0)}  |\vec{H}|\ |\vec{\mathbb I}_{\vec{\Phi}(t)}|_{g_{\vec{\Phi}(t)}}\ dvol_{g_{\vec{\Phi}(t)}}\ \int_0^T\lf\| \lf|\frac{dh}{dt}(t)\rg|_{h(t)}\rg\|_{L^\infty(B_R^{h(t)}(x_0))}\  dt\\[5mm]
\ds\le  C_\Sigma\ \frac{e^{2\,\|\al\|_\infty}}{l^4_{\ast\ast}(T)}\ \sqrt{T}\ \sup_{[0,T]}\int_{B_R^{h(t)}(x_0)}   |\vec{H}|\ |\vec{\mathbb I}_{\vec{\Phi}(t)}|_{g_{\vec{\Phi}(t)}}\ dvol_{g_{\vec{\Phi}(t)}}\ \sqrt{E(0)}\,\sqrt{W(\vec{\Phi}(0))-W(\vec{\Phi}(t))}\ .
\end{array}
\ee
%Now, from lemma~\ref{lm-param-conf-var}
%\be
%\label{param-conf-var-re}
%\frac{d\al}{dt}=\,\frac{d\vec{\Phi}}{dt}\cdot \vec{H}+e^{-\,2\al(t)}\frac{1}{2}\,d^{\ast_{h(t)}}\lf[\frac{d\vec{\Phi}}{dt}\cdot d\vec{\Phi}(t)\rg]\ .
%\ee
%Hence
%\be
%\label{n-1}
%\begin{array}{l}
%\ds\frac{1}{2}\,\int_\Sigma \chi^{h(t)}_{R,x_0}(x)\,\ e^{-2\al}\  \frac{d\al}{dt}\ \lf|\Delta_{h(t)}\vec{\Phi}\rg|^2\ dvol_{h(t)}\\[5mm]
%\ds=2\,\int_\Sigma \chi^{h(t)}_{R,x_0}(x)\,\pi_{\vec{n}}\frac{d\vec{\Phi}}{dt}\cdot \vec{H}\ |\vec{H}|^2\ \ dvol_{g_{\vec{\Phi}(t)}}\\[5mm]
%\ds+\int_\Sigma \chi^{h(t)}_{R,x_0}(x)\,e^{-\,2\al(t)}\,d^{\ast_{h(t)}}\lf[\frac{d\vec{\Phi}}{dt}\cdot d\vec{\Phi}(t)\rg]\  |\vec{H}|^2\ \ dvol_{g_{\vec{\Phi}(t)}}\ .
%\end{array}
%\ee
We have next
\be
\label{n-2}
\begin{array}{l}
\ds\int_\Sigma \chi^{h(t)}_{R,x_0}(x)\,\pi_{\vec{n}}\frac{d\vec{\Phi}}{dt}\cdot \vec{H}\ |\vec{H}|^2\ \ dvol_{g_{\vec{\Phi}(t)}}\\[5mm]
\ds\le \sqrt{\int_\Sigma \lf|\pi_{\vec{n}}\frac{d\vec{\Phi}}{dt}\rg|^2\  dvol_{g_{\vec{\Phi}(t)}}}\ \lf[\int_{B_R^{h(t)}(x_0)} |\vec{H}|^6\ dvol_{g_{\vec{\Phi}(t)}}\rg]^{1/2}    \\[5mm]
\ds\le \sqrt{\int_\Sigma \lf|\pi_{\vec{n}}\frac{d\vec{\Phi}}{dt}\rg|^2\  dvol_{g_{\vec{\Phi}(t)}}}\  e^{\al_{R,x_0}(t)}\ \lf[\int_{B_R^{h(t)}(x_0)} |\vec{H}|^6\ dvol_{h(t)}\rg]^{1/2}  \ ,
 \end{array}
\ee
where
\[
\al_{R,x_0}(t)=\max_{x\in B_R^{h(t)}(x_0)}\al(t)(x)\ .
\]
Interpolation estimates give
\be
\label{n-ene-con-31-re}
\begin{array}{l}
\ds\lf[\int_{B_R^{h(t)}(x_0)} |\vec{H}|^6\  \ dvol_{h(t)}\rg]^{1/2}\\[5mm]
\ds\le C\ \int_{B_{2\,R}^{h(t)}(x_0)}\ |\vec{H}|^2\  \ dvol_{h(t)}\ \ \lf[\int_{B_{2\,R}^{h(t)}(x_0)}\ |\Delta_{h(t)}\vec{H}|^2\  \ dvol_{h(t)}\rg]^{1/2}\\[5mm]
\ds+ C\ \lf[R^{-4/3}\,\int_{B_{2\,R}^{h(t)}(x_0)}\ |\vec{H}|^2\  \ dvol_{h(t)}\rg]^{3/2}\ 
\end{array}
\ee
 Hence
 \be
\label{n-3}
\begin{array}{l}
\ds C^{-1}\ \int_\Sigma \chi^{h(t)}_{R,x_0}(x)\,\pi_{\vec{n}}\frac{d\vec{\Phi}}{dt}\cdot \vec{H}\ |\vec{H}|^2\ \ dvol_{g_{\vec{\Phi}(t)}}\\[5mm]
\ds\le\sqrt{\int_\Sigma \lf|\pi_{\vec{n}}\frac{d\vec{\Phi}}{dt}\rg|^2\  dvol_{g_{\vec{\Phi}(t)}}}\  e^{\al_{R,x_0}(t)}\ \int_{B_{2\,R}^{h(t)}(x_0)}\ |\vec{H}|^2\  \ dvol_{h(t)}\ \ \sqrt{\int_{B_{2\,R}^{h(t)}(x_0)}\ |\Delta_{h(t)}\vec{H}|^2\  \ dvol_{h(t)}}\\[5mm]
\ds +\sqrt{\int_\Sigma \lf|\pi_{\vec{n}}\frac{d\vec{\Phi}}{dt}\rg|^2\  dvol_{g_{\vec{\Phi}(t)}}}\  e^{\al_{R,x_0}(t)}\ R^{-2} \lf[\int_{B_{2\,R}^{h(t)}(x_0)}\ |\vec{H}|^2\  \ dvol_{h(t)}\rg]^{3/2}\ 
\end{array}
\ee
We denote
\[
E(2R,x_0,t):=\sup_{[0,T]}\int_{B_{2\,R}^{h(t)}(x_0)}\ |\vec{\mathbb I}|^2_{g_{\vec{\Phi}}}\ e^{2\al_{R,x_0}(t)} dvol_{h(t)}
\] 
Observe that, thanks to (\ref{L2-1}) there holds
\be
\label{energy}
\sup_{[0,T]}\int_{B_{2\,R}^{h(t)}(x_0)}\ |\vec{\mathbb I}|^2_{g_{\vec{\Phi}}}\  dvol_{g_{\vec{\Phi}}}\le E(2R,x_0,t)\le C(\Sigma,l_{\ast\ast}(T),E(0))\ \sup_{[0,T]}\int_{B_{2\,R}^{h(t)}(x_0)}\ |\vec{\mathbb I}|^2_{g_{\vec{\Phi}}}\  dvol_{g_{\vec{\Phi}}}\ .
\ee
With this notation (\ref{n-3}) implies
 \be
\label{n-4}
\begin{array}{l}
\ds C^{-1}\ \int_\Sigma \chi^{h(t)}_{R,x_0}(x)\,\pi_{\vec{n}}\frac{d\vec{\Phi}}{dt}\cdot \vec{H}\ |\vec{H}|^2\ \ dvol_{g_{\vec{\Phi}(t)}}\\[5mm]
\ds\le 2^{-1}\,\sqrt{\int_\Sigma \lf|\pi_{\vec{n}}\frac{d\vec{\Phi}}{dt}\rg|^2\  dvol_{g_{\vec{\Phi}(t)}}}\ E(2R,x_0,t)\ \sqrt{\int_{B_{2\,R}^{h(t)}(x_0)}\ e^{-2\al_{R,x_0}(t)}\ |\Delta_{h(t)}\vec{H}|^2\  \ dvol_{h(t)}}\\[5mm]
 \ds+ 2^{-1}\,\sqrt{\int_\Sigma \lf|\pi_{\vec{n}}\frac{d\vec{\Phi}}{dt}\rg|^2\  dvol_{g_{\vec{\Phi}(t)}}}\  {E(2R,x_0,t)}^{3/2}\ e^{2\|\al\|_\infty}\ R^{-2}\\[5mm]
 \end{array}
\ee 
Integrating with respect to $t$ and using Cauchy Schwartz gives
 \be
\label{n-5}
\begin{array}{l}
\ds C^{-1}\ \int_0^T\int_\Sigma \chi^{h(t)}_{R,x_0}(x)\,\pi_{\vec{n}}\frac{d\vec{\Phi}}{dt}\cdot \vec{H}\ |\vec{H}|^2\ \ dvol_{g_{\vec{\Phi}(t)}}\ dt\\[5mm]
\ds\le \sqrt{E(2R,x_0,t)}\ \int_{B_{2\,R}^{h(t)}(x_0)}\ e^{-2\al_{R,x_0}(t)}\ |\Delta_{h(t)}\vec{H}|^2\  \ dvol_{h(t)}\ dt\\[5mm]
\ds+{E(2R,x_0,t)}^{3/2}\ \sqrt{\int_0^T\int_\Sigma \lf|\pi_{\vec{n}}\frac{d\vec{\Phi}}{dt}\rg|^2\  dvol_{g_{\vec{\Phi}(t)}}\ dt}\ [\sqrt{T}\ e^{2\|\al\|_\infty}\ R^{-2}+1]\ .
\end{array}
\ee 
%We estimate now the second term in the r.h.s. of (\ref{n-1})
%\be
%\label{n-6}
%\begin{array}{l}
%\ds\int_\Sigma \chi^{h(t)}_{R,x_0}(x)\,e^{-\,2\al(t)}\frac{1}{2}\,d^{\ast_{h(t)}}\lf[\frac{d\vec{\Phi}}{dt}\cdot d\vec{\Phi}(t)\rg]\  |\vec{H}|^2\ \ dvol_{g_{\vec{\Phi}(t)}} \\[5mm]
%\ds=\int_\Sigma \chi^{h(t)}_{R,x_0}(x)\,\frac{1}{2}\,d^{\ast_{g_{\vec{\Phi}(t)}}}\lf[\frac{d\vec{\Phi}}{dt}\cdot d\vec{\Phi}(t)\rg]\  |\vec{H}|^2\ \ dvol_{g_{\vec{\Phi}(t)}}\\[5mm]
%\ds =\frac{1}{2}\,\int_\Sigma \lf<d\chi^{h(t)}_{R,x_0}(x),\vec{U}\cdot d\vec{\Phi}(t)\rg>_{g_{\vec{\Phi}(t)}}\  |\vec{H}|^2\  \ dvol_{g_{\vec{\Phi}(t)}}\\[5mm]
%\ds+\frac{1}{2}\,\int_\Sigma \ \chi^{h(t)}_{R,x_0}(x)\ \lf<d|\vec{H}|^2,\vec{U}\cdot d\vec{\Phi}(t)\rg>_{g_{\vec{\Phi}(t)}}\  dvol_{g_{\vec{\Phi}(t)}}
%\end{array}%\ee 
%From (\ref{u}) we have $\vec{U}={\ov{\p}}\vec{\Phi}(t)\res U^{0,1}+{{\p}}\vec{\Phi}(t)\res \ov{U^{0,1}}$ hence
%\[
%\vec{U}\cdot d\vec{\Phi}(t)=\vec{U}\cdot \lf[\p\vec{\Phi}+\ov{\p}\vec{\Phi}\rg]=2\,\Re\lf[U^{0,1}\res g_{\vec{\Phi}}\rg]=2\, e^{2\al}\,\Re\lf[U^{0,1}\res h\rg]\ .
%\]
We have next
\be
\label{n-7}
\begin{array}{l}
\ds\lf|\int_\Sigma e^{-2\al}\ \lf<d\chi^{h(t)}_{R,x_0}(x),\frac{d\vec{\Phi}}{dt}\cdot d\vec{\Phi}(t)\rg>_{h(t)}\  |\vec{H}|^2\  \ dvol_{g_{\vec{\Phi}(t)}}\rg|\\[5mm]
\ds\le C\,R^{-1}\ \int_{B_{R}^{h(t)}(x_0)}\  e^{2\al}\ |U^{0,1}|_{h(t)}\ |\vec{H}|^2\  \ dvol_{h(t)}\\[5mm]
\ds\le C\, R^{-1}\ e^{2\al_{R,x_0}(t)} \|U^{0,1}\|_{L^{2,\infty}_{h(t)}(B_{R}^{h(t)}(x_0))}\ \||\vec{H}|^2\ \|_{L^{2,1}_{h(t)}(B_{R}^{h(t)}(x_0))}\
%\ds\le C\, R^{-1}\  e^{2\al_{R,x_0}(t)} \ \lf[\int_{B_{R}^{h(t)}(x_0)}\ |U^{0,1}|^{3/2}_{h(t)}\ dvol_{h(t)}\rg]^{2/3}\  \lf[\int_{B_{2\,R}^{h(t)}(x_0)}\ |\vec{H}|^6\  \ dvol_{h(t)}\rg]^{1/3}\
\end{array}
\ee
Lorentz-Sobolev embedding  inequalities give
\be
\label{n-8}
\begin{array}{l}
\ds\||\vec{H}|^2\ \|_{L^{2,1}_{h(t)}(B_{R}^{h(t)}(x_0))}\le C\ \|\nabla(\vec{H}\cdot\vec{H})\|_{L^1_{h(t)}(B_{3R/2}^{h(t)}(x_0))}+C\,R^{-1}\,\int_{B_{3R/2}^{h(t)}(x_0)}|\vec{H}|^2\ dvol_{h(t)}\\[5mm]
\ds\le C\ \lf[\int_{B_{3R/2}^{h(t)}(x_0)}|\vec{H}|^2\ dvol_{h(t)}\rg]^{1/2}\ \lf[\int_{B_{3R/2}^{h(t)}(x_0)}|d\vec{H}|^2\ dvol_{h(t)}\rg]^{1/2}\\[5mm]
\ds\quad+C\,R^{-1}\,\int_{B_{3R/2}^{h(t)}(x_0)}|\vec{H}|^2\ dvol_{h(t)}\\[5mm]
\end{array}
\ee
and classical interpolation inequality imply
\be
\label{n-9}
\begin{array}{l}
\ds\lf[\int_{B_{3R/2}^{h(t)}(x_0)}|d\vec{H}|^2\ dvol_{h(t)}\rg]^{1/2}\le C\ \lf[\int_{B_{2R}^{h(t)}(x_0)}|\vec{H}|^2\ dvol_{h(t)}\rg]^{1/4}\ \lf[\int_{B_{2R}^{h(t)}(x_0)}|\Delta_{h(t)}\vec{H}|^2\ dvol_{h(t)}\rg]^{1/4}\\[5mm]
\ds \quad+C\ R^{-1}\lf[\int_{B_{2R}^{h(t)}(x_0)}|\vec{H}|^2\ dvol_{h(t)}\rg]^{1/2}
\end{array}
\ee
Combining (\ref{n-8}) and (\ref{n-9}) is giving
\be
\label{n-9-b}
\begin{array}{l}
\ds\||\vec{H}|^2\ \|_{L^{2,1}_{h(t)}(B_{R}^{h(t)}(x_0))}\le  C\ \lf[\int_{B_{2R}^{h(t)}(x_0)}|\vec{H}|^2\ dvol_{h(t)}\rg]^{3/4}\ \lf[\int_{B_{2R}^{h(t)}(x_0)}|\Delta_{h(t)}\vec{H}|^2\ dvol_{h(t)}\rg]^{1/4}\\[5mm]
\ds\quad+C\,R^{-1}\,\int_{B_{2R}^{h(t)}(x_0)}|\vec{H}|^2\ dvol_{h(t)}
\end{array}
\ee
Using a combination of lemma~\ref{lm-delbar-inv} and lemma~\ref{lm-L2-infty} gives
 \be
 \label{n-9-c}
 \begin{array}{l}
 \ds\|U^{0,1}\|_{L^{2,\infty}_{h(t)}(B_{R}^{h(t)}(x_0))}\le C(l_{\ast\ast}(T))\ \|\p^hU^{0,1}\|_{L^1_{h(t)}(\Sigma)}
 \end{array}
 \ee
%Classical H\"older-Lorentz inequality (\cite{Gra1}) combined with lemma~\ref{lm-delbar-inv} and lemma~\ref{lm-L2-infty}
%\be
%\label{n-8}
%\begin{array}{l}
%\ds R^{-1} \lf[ \int_{B_{R}^{h(t)}(x_0)}\ |U^{0,1}|^{3/2}_{h(t)}\ dvol_{h(t)}\rg]^{2/3}\le C\ R^{-2/3}\ \|U^{0,1}\|_{L^{2,\infty}_{h(t)}(B_{R}^{h(t)}(x_0))}\\[5mm]
%\ds\le C\ l_{\ast\ast}(T)^{-5}\ R^{-2/3}\ \|\p^hU^{0,1}\|_{L^1_{h(t)}(\Sigma)}
%\end{array}
%\ee
%where $C$ is universal. Recall
%\be
%\label{n-ene-con-31-bis}
%\begin{array}{l}
%\ds\lf[\int_{B_R^{h(t)}(x_0)} |\vec{H}|^6\  \ dvol_{h(t)}\rg]^{1/3}\\[5mm]
%\ds\le C\ \lf[\int_{B_{2\,R}^{h(t)}(x_0)}\ |\vec{H}|^2\  \ dvol_{h(t)}\rg]^{2/3}\ \ \lf[\int_{B_{2\,R}^{h(t)}(x_0)}\ |\Delta_{h(t)}\vec{H}|^2\  \ dvol_{h(t)}\rg]^{1/3}\\[5mm]
%\ds+ C\ R^{-4/3}\,\int_{B_{2\,R}^{h(t)}(x_0)}\ |\vec{H}|^2\  \ dvol_{h(t)}\ 
%\end{array}
%\ee
%Hence
%\newpage
%\be
%\label{n-9}
%\begin{array}{l}
%\ds\lf|\int_\Sigma e^{-2\al}\ \lf<d\chi^{h(t)}_{R,x_0}(x),\vec{U}\cdot d\vec{\Phi}(t)\rg>_{h(t)}\  |\vec{H}|^2\  \ dvol_{g_{\vec{\Phi}(t)}}\rg|\\[5mm]
%\ds\le \frac{C( l_{\ast\ast}(T))}{R^{2/3}}\ \  e^{2\al_{R,x_0}(t)} \  \|\p^hU^{0,1}\|_{L^1_{h(t)}(\Sigma)}\ \lf[\int_{B_{2\,R}^{h(t)}(x_0)}\ |\vec{H}|^2\  \ dvol_{h(t)}\rg]^{2/3}\ \ \lf[\int_{B_{2\,R}^{h(t)}(x_0)}\ |\Delta_{h(t)}\vec{H}|^2\  \ dvol_{h(t)}\rg]^{1/3}\\[5mm]
%\ds+\frac{C( l_{\ast\ast}(T))}{R^{2}}\  \|\p^hU^{0,1}\|_{L^1_{h(t)}(\Sigma)}\ \int_{B_{2\,R}^{h(t)}(x_0)}\ |\vec{H}|^2\ \  e^{2\al_{R,x_0}(t)} \  \ dvol_{h(t)}\ 
%\end{array}
%\ee
We have
\be
\label{n-10}
\begin{array}{rl}
\ds\|\p^hU^{0,1}\|_{L^1_{h(t)}(\Sigma)}&\ds=\int_\Sigma\lf|(I-\ti{P}_{h(t)})\lf(\frac{d\vec{\Phi}}{dt}\cdot\vec{\frak h}^{0}\rg)\rg|_{h(t)}\ dvol_{h(t)}\\[5mm]
\ds&\ds\le \int_\Sigma\lf|\frac{d\vec{\Phi}}{dt}\cdot\vec{\frak h}^{0}\rg|_{h(t)}\ dvol_{h(t)}+ \int_\Sigma\lf|\ti{P}_{h(t)}\lf(\frac{d\vec{\Phi}}{dt}\cdot\vec{\frak h}^{0}\rg)\rg|_{h(t)}\ dvol_{h(t)}
\end{array}
\ee
As in the previous section we bound
\be
\label{n-11}
\begin{array}{l}
\ds\int_\Sigma\lf|\frac{d\vec{\Phi}}{dt}\cdot\vec{\frak h}^{0}\rg|_{h(t)}\ dvol_{h(t)}
\le e^{2\|\al\|_\infty}\ \sqrt{\int_\Sigma \lf|\pi_{\vec{n}}\frac{d\vec{\Phi}}{dt}\rg|^2\  dvol_{g_{\vec{\Phi}(t)}}}\ \sqrt{ E(0)}\\[5mm]

\end{array}
\ee
and
\be
\label{n-12}
\begin{array}{l}
\ds
\int_\Sigma\lf|\ti{P}_{h(t)}\lf(\frac{d\vec{\Phi}}{dt}\cdot\vec{\frak h}^{0}\rg)\rg|_{h(t)}\ dvol_{h(t)}\le \sum_{i=1}^Q\|q_i(t)\|_{L^\infty(\Sigma)}\sqrt{V_{h(t)}}\ {\int_\Sigma\lf|\frac{d\vec{\Phi}}{dt}\cdot\vec{\frak h}^{0}\rg|_{h(t)}\ dvol_{h(t)}}\\[5mm]
\ds\le C_\Sigma\  l_{\ast\ast}(T)^{-2}\ e^{2\|\al\|_\infty}\ \sqrt{\int_\Sigma \lf|\pi_{\vec{n}}\frac{d\vec{\Phi}}{dt}\rg|^2\  dvol_{g_{\vec{\Phi}(t)}}}\ \sqrt{ E(0)}\ .
\end{array}
\ee
This gives
\be
\label{n-13}
\ds\|\p^hU^{0,1}\|_{L^1_{h(t)}(\Sigma)}\le \frac{C_\Sigma}{l_{\ast\ast}(T)^{2}}\ e^{2\|\al\|_\infty}\ \sqrt{\int_\Sigma \lf|\pi_{\vec{n}}\frac{d\vec{\Phi}}{dt}\rg|^2\  dvol_{g_{\vec{\Phi}(t)}}}\ \sqrt{ E(0)}\ .
\ee
Combining (\ref{n-7}), (\ref{n-9-b}), (\ref{n-9-c}) and (\ref{n-13}) is giving
\be
\label{n-14}
\begin{array}{l}
\ds\int_0^T\lf|\int_\Sigma e^{-2\al}\ \lf<d\chi^{h(t)}_{R,x_0}(x),\vec{U}\cdot d\vec{\Phi}(t)\rg>_{h(t)}\  |\vec{H}|^2\  \ dvol_{g_{\vec{\Phi}(t)}}\rg|\ dt\\[5mm]
\ds\le \frac{C_\Sigma\,\sqrt{ E(0)}}{l_{\ast\ast}(T)^{2}}\ \frac{e^{2\|\al\|_\infty}}{R}\ \int_0^T\sqrt{\int_\Sigma \lf|\pi_{\vec{n}}\frac{d\vec{\Phi}}{dt}\rg|^2\  dvol_{g_{\vec{\Phi}(t)}}}\  e^{\al_{R,x_0}(t)} \  \\  
\ds\times \  \lf[\int_{B_{2\,R}^{h(t)}(x_0)}\ |\vec{H}|^2\  e^{2\,\al_{R,x_0}(t)} \ dvol_{h(t)}\rg]^{3/4}\ \ \lf[\int_{B_{2\,R}^{h(t)}(x_0)}\ |\Delta_{h(t)}\vec{H}|^2\  e^{-2\al_{R,x_0}(t)} \ dvol_{h(t)}\rg]^{1/4}\\[5mm]
\ds+ \frac{C_\Sigma\,\sqrt{ E(0)}}{l_{\ast\ast}(T)^{2}}\ \frac{e^{2\|\al\|_\infty}}{R^{2}}\ \int_0^T\sqrt{\int_\Sigma \lf|\pi_{\vec{n}}\frac{d\vec{\Phi}}{dt}\rg|^2\  dvol_{g_{\vec{\Phi}(t)}}}\  \int_{B_{2\,R}^{h(t)}(x_0)}\ |\vec{H}|^2\  \ 
e^{2\al_{R,x_0}(t)} \ dvol_{h(t)}\ .
\end{array}
\ee
% using that the time derivative of $h$ is trace free and consequently $d(dvol_{h(t)})/dt=0$
%\be
%\label{n-ene-con-60-7}
%\frac{d}{dt}\int_\Sigma\ dvol_{g_{\vec{\Phi}(t)}}=2\, \int_{\Sigma}\frac{d\al}{dt}\ e^{2\al}\ dvol_{h(t)}=\int_{\Sigma}\pi_{\vec{n}}\frac{d\vec{\Phi}}{dt}\cdot \vec{H}\ dvol_{g_{\vec{\Phi}(t)}}
%\ee
%Thus
%\be
%\label{n-vol-evol}
%|V(\vec{\Phi}(T))-V(\vec{\Phi}(0))|\le\int_0^T\lf|\frac{d}{dt}\int_\Sigma\ dvol_{g_{\vec{\Phi}(t)}}\rg|\ dt\le \sqrt{E(0)}\ \sqrt{W(\vec{\Phi}(0))-W(\vec{\Phi}(T)}
%\ee
%Let $x(t)\in\Sigma$ such that
%\[
%\al(t)(x(t))=\max_{x_0\in\Sigma}\al_{R,x_0}(t)\ .
%\]
%On $B_{R}(x_(t))$ we have
%\[
%\|\al(t)(x)-\al(t)(x(t))\|_{L^\infty(B_{R}(x_(t)))}\le C(\Sigma, l_{\ast\ast}(T),E(0))\ .
%\]
%Hence 
%\be
%\label{n-15}
%C(\Sigma, l_{\ast\ast}(T),E(0))^{-1}\ R^2\ \max_{x_0\in\Sigma}\,e^{2\al_{R,x_0}(t)}\le\ V(\vec{\Phi}(T))\le V(\vec{\Phi}(0))+\sqrt{E(0)}\ \sqrt{W(\vec{\Phi}(0))-W(\vec{\Phi}(T)}\ .
%\ee
%Combining (\ref{n-14}) and (\ref{n-15}) gives
This implies
\be
\label{n-16}
\begin{array}{l}
\ds C(\Sigma, l_{\ast\ast}(T),E(0))^{-1}\ \int_0^T\lf|\int_\Sigma e^{-2\al}\ \lf<d\chi^{h(t)}_{R,x_0}(x),\vec{U}\cdot d\vec{\Phi}(t)\rg>_{h(t)}\  |\vec{H}|^2\  \ dvol_{g_{\vec{\Phi}(t)}}\rg|\ dt\\[5mm]
\ds\le  \frac{e^{3\|\al\|_\infty}}{R}\ \sqrt{\int_0^T\int_\Sigma \lf|\pi_{\vec{n}}\frac{d\vec{\Phi}}{dt}\rg|^2\  dvol_{g_{\vec{\Phi}(t)}} dt}\\  
\ds\times \  E(2R,x_0,t)^{3/4}\ \ \sqrt{\int_0^T\lf[\int_{B_{2\,R}^{h(t)}(x_0)}\ |\Delta_{h(t)}\vec{H}|^2\  e^{-2\al_{R,x_0}(t)} \ dvol_{h(t)}\rg]^{1/2}\ dt}\\[5mm]
\ds+C\ \frac{e^{2\|\al\|_\infty}}{R^{2}}\ \int_0^T\sqrt{\int_\Sigma \lf|\pi_{\vec{n}}\frac{d\vec{\Phi}}{dt}\rg|^2\  dvol_{g_{\vec{\Phi}(t)}}}\  \int_{B_{2\,R}^{h(t)}(x_0)}\ |\vec{H}|^2\  \ 
e^{2\al_{R,x_0}(t)} \ dvol_{h(t)}
\end{array}
\ee
Hence
\be
\label{n-17}
\begin{array}{l}
\ds C(\Sigma, l_{\ast\ast}(T),E(0))^{-1}\ \int_0^T\lf|\int_\Sigma e^{-2\al}\ \lf<d\chi^{h(t)}_{R,x_0}(x),\vec{U}\cdot d\vec{\Phi}(t)\rg>_{h(t)}\  |\vec{H}|^2\  \ dvol_{g_{\vec{\Phi}(t)}}\rg|\ dt\\[5mm]
\ds\le\,\delta\ \int_0^T \int_{B_{2\,R}^{h(t)}(x_0)}\ |\Delta_{h(t)}\vec{H}|^2\  e^{-2\al_{R,x_0}(t)} \ dvol_{h(t)}\\[5mm]
\ds+\ E(2R,x_0,t)  \ \frac{{e^{4\|\al\|_\infty}}}{\delta^3}\ \lf(\frac{T}{R^{4}}\rg)^{1/3}\ \lf(\int_0^T\int_\Sigma \lf|\pi_{\vec{n}}\frac{d\vec{\Phi}}{dt}\rg|^2\  dvol_{g_{\vec{\Phi}(t)}}\ dt\rg)^{2/3}\\[5mm]
\ds+\  E(2R,x_0,t)\ {e^{2\|\al\|_\infty}}\sqrt{\frac{{T}}{R^{4}}}\ \ \sqrt{\int_0^T\int_\Sigma \lf|\pi_{\vec{n}}\frac{d\vec{\Phi}}{dt}\rg|^2\  dvol_{g_{\vec{\Phi}(t)}}\ dt}\  
\end{array}
\ee
We bound now using (\ref{re-n-2-a})
\be
\label{n-18}
\begin{array}{l}
\ds\int_0^\infty\lf|\int_\Sigma \ \chi^{h(t)}_{R,x_0}(x)\ \lf<d|\vec{H}|^2,\vec{U}\cdot d\vec{\Phi}(t)\rg>_{g_{\vec{\Phi}(t)}}\  dvol_{g_{\vec{\Phi}(t)}}\rg|\ dt\\[5mm]
\ds\le\int_0^\infty e^{2\al_{R,x_0}(t)}\ \|U^{0,1}\|_{L^{2,\infty}_{h(t)}(B_{R}^{h(t)}(x_0))}\ \|d|\vec{H}|^2\ \|_{L^{2,1}_{h(t)}(B_{R}^{h(t)}(x_0))}\ dt\ .
\end{array}
\ee 
We have
\be
\label{n-19}
\begin{array}{l}
\ds\|d|\vec{H}|^2\ \|_{L^{2,1}_{h(t)}(B_{R}^{h(t)}(x_0))}\le C\ \| \nabla^2_{h(t)}|\vec{H}|^2\|_{L^1_{h(t)}(B_{3R/2}^{h(t)}(x_0))}+R^{-1}\ \|d|\vec{H}|^2\ \|_{L^{1}_{h(t)}(B_{R}^{h(t)}(x_0))}\\[5mm]
\ds\le C\  \| \nabla^2_{h(t)}\vec{H}\|_{L^2_{h(t)}(B_{3R/2}^{h(t)}(x_0))}\ \| \vec{H}\|_{L^2_{h(t)}(B_{3R/2}^{h(t)}(x_0))}+C\ \| d\vec{H}\|^2_{L^2_{h(t)}(B_{3R/2}^{h(t)}(x_0))}\\[5mm]
\ds+C\ R^{-2} \int_{B_{3R/2}^{h(t)}(x_0)}|\vec{H}|^2\ dvol_{h(t)}\\[5mm]
\ds\le C\  \| \Delta_{h(t)}\vec{H}\|_{L^2_{h(t)}(B_{2R}^{h(t)}(x_0))}\ \| \vec{H}\|_{L^2_{h(t)}(B_{2R}^{h(t)}(x_0))}+C\ R^{-2} \int_{B_{2R}^{h(t)}(x_0)}|\vec{H}|^2\ dvol_{h(t)}
\end{array}
\ee
Combining (\ref{n-9-c}), (\ref{n-13}) (\ref{n-18}) and (\ref{n-19}) gives
\be
\label{n-20}
\begin{array}{l}
\ds C(\Sigma, l_{\ast\ast}(T),E(0))^{-1}\ \int_0^\infty\lf|\int_\Sigma \ \chi^{h(t)}_{R,x_0}(x)\ \lf<d|\vec{H}|^2,\vec{U}\cdot d\vec{\Phi}(t)\rg>_{g_{\vec{\Phi}(t)}}\  dvol_{g_{\vec{\Phi}(t)}}\rg|\ dt\\[5mm]
\ds\le \frac{C_\Sigma\,\sqrt{ E(0)}}{l_{\ast\ast}(T)^{2}}\ e^{2\|\al\|_\infty}\ \int_0^T\sqrt{\int_\Sigma \lf|\pi_{\vec{n}}\frac{d\vec{\Phi}}{dt}\rg|^2\  dvol_{g_{\vec{\Phi}(t)}}}\  e^{2\al_{R,x_0}(t)} \  \\  
\ds\times \  \lf[\int_{B_{2\,R}^{h(t)}(x_0)}\ |\vec{H}|^2\   \ dvol_{h(t)}\rg]^{1/2}\ \ \lf[\int_{B_{2\,R}^{h(t)}(x_0)}\ |\Delta_{h(t)}\vec{H}|^2\  \ dvol_{h(t)}\rg]^{1/2}\\[5mm]
\ds+ \frac{C_\Sigma\,\sqrt{ E(0)}}{l_{\ast\ast}(T)^{2}}\ \frac{e^{2\|\al\|_\infty}}{R^{2}}\ \int_0^T\sqrt{\int_\Sigma \lf|\pi_{\vec{n}}\frac{d\vec{\Phi}}{dt}\rg|^2\  dvol_{g_{\vec{\Phi}(t)}}}\  \int_{B_{2\,R}^{h(t)}(x_0)}\ |\vec{H}|^2\  \ 
e^{2\al_{R,x_0}(t)} \ dvol_{h(t)}\ .
\end{array}
\ee
This gives
\be
\label{n-21}
\begin{array}{l}
\ds C(\Sigma, l_{\ast\ast}(T),E(0))^{-1}\ \int_0^\infty\lf|\int_\Sigma \ \chi^{h(t)}_{R,x_0}(x)\ \lf<d|\vec{H}|^2,\vec{U}\cdot d\vec{\Phi}(t)\rg>_{g_{\vec{\Phi}(t)}}\  dvol_{g_{\vec{\Phi}(t)}}\rg|\ dt\\[5mm]
\ds\le  e^{4\|\al\|_\infty}\ \sqrt{E(2R,x_0,t)}\ \sqrt{\int_0^T\int_\Sigma \lf|\pi_{\vec{n}}\frac{d\vec{\Phi}}{dt}\rg|^2\  dvol_{g_{\vec{\Phi}(t)}} dt}\    \\  
\ds\times \ \ \lf[\int_0^T\int_{B_{2\,R}^{h(t)}(x_0)}\ |\Delta_{h(t)}\vec{H}|^2\  \ dvol_{h(t)}\ dt\rg]^{1/2}\\[5mm]
\ds+ \frac{C_\Sigma\,\sqrt{ E(0)}}{l_{\ast\ast}(T)^{2}}\ \frac{e^{2\|\al\|_\infty}}{R^{2}}\ \int_0^T\sqrt{\int_\Sigma \lf|\pi_{\vec{n}}\frac{d\vec{\Phi}}{dt}\rg|^2\  dvol_{g_{\vec{\Phi}(t)}}}\  \int_{B_{2\,R}^{h(t)}(x_0)}\ |\vec{H}|^2\  \ 
e^{2\al_{R,x_0}(t)} \ dvol_{h(t)}\ .
\end{array}
\ee
hence
\be
\label{n-22}
\begin{array}{l}
\ds C(\Sigma, l_{\ast\ast}(T),E(0))^{-1}\ \int_0^\infty\lf|\int_\Sigma \ \chi^{h(t)}_{R,x_0}(x)\ \lf<d|\vec{H}|^2,\vec{U}\cdot d\vec{\Phi}(t)\rg>_{g_{\vec{\Phi}(t)}}\  dvol_{g_{\vec{\Phi}(t)}}\rg|\ dt\\[5mm]
\ds\le\,\delta\ \int_0^T \int_{B_{2\,R}^{h(t)}(x_0)}\ |\Delta_{h(t)}\vec{H}|^2\  e^{-2\al_{R,x_0}(t)} \ dvol_{h(t)}\\[5mm]
\ds+\frac{e^{8\|\al\|_\infty}}{\delta}\ {E(2R,x_0,t)}\ {\int_0^T\int_\Sigma \lf|\pi_{\vec{n}}\frac{d\vec{\Phi}}{dt}\rg|^2\  dvol_{g_{\vec{\Phi}(t)}} dt}\\[5mm] 
\ds+\  E(2R,x_0,t)\ {e^{2\|\al\|_\infty}}\sqrt{\frac{{T}}{R^{4}}}\ \ \sqrt{\int_0^T\int_\Sigma \lf|\pi_{\vec{n}}\frac{d\vec{\Phi}}{dt}\rg|^2\  dvol_{g_{\vec{\Phi}(t)}}\ dt}\  
\end{array}
\ee
%Combining the previous gives finally
%\be
%\label{n-23}
%\begin{array}{l}
%\ds\frac{1}{2}\,\int_0^T\lf|\int_\Sigma \chi^{h(t)}_{R,x_0}(x)\,\ e^{-2\al}\  \frac{d\al}{dt}\ \lf|\Delta_{h(t)}\vec{\Phi}\rg|^2\ dvol_{h(t)}\rg|\ dt\\[5mm]
%\ds\le \lf[\sqrt{E(2R,x_0,t)}+\delta\rg]\ \int_0^T\int_{B_{2\,R}^{h(t)}(x_0)}\ e^{-2\al_{R,x_0}(t)}\ |\Delta_{h(t)}\vec{H}|^2\  \ dvol_{h(t)}\ dt\\[5mm]
%\ds+{E(2R,x_0,t)}^{3/2}\ \int_0^T\int_\Sigma \lf|\pi_{\vec{n}}\frac{d\vec{\Phi}}{dt}\rg|^2\  dvol_{g_{\vec{\Phi}(t)}}\ dt+{E(2R,x_0,t)}^{3/2}\  e^{4\|\al\|_\infty}\ R^{-4}\ T\ \\[5mm]
%\ds+\ E(2R,x_0,t)  \ \frac{{e^{4\|\al\|_\infty}}}{\delta^3}\ \lf(\frac{T}{R^{4}}\rg)^{1/3}\ \lf(W(\vec{\Phi}(0))-W(\vec{\Phi}(T))\rg)^{2/3}\\[5mm]
%\ds+\  E(2R,x_0,t)\ {e^{2\|\al\|_\infty}}\sqrt{\frac{{T}}{R^{4}}}\ \ \sqrt{W(\vec{\Phi}(0))-W(\vec{\Phi}(T))}\  \\[5mm]
%\ds+\frac{e^{8\|\al\|_\infty}}{\delta}\ {E(2R,x_0,t)}\ \lf(W(\vec{\Phi}(0))-W(\vec{\Phi}(T))\rg)\\[5mm] 
%\ds+{E(2R,x_0,t)}^{3/2}\ e^{2\|\al\|_\infty}\ \sqrt{\frac{T}{R^4}}\ \sqrt{W(\vec{\Phi}(0))-W(\vec{\Phi}(T))}\ .
%\ds+\  E(2R,x_0,t)\ {e^{2\|\al\|_\infty}}\sqrt{\frac{{T}}{R^{4}}}\ \ \sqrt{\int_0^T\int_\Sigma \lf|\pi_{\vec{n}}\frac{d\vec{\Phi}}{dt}\rg|^2\  dvol_{g_{\vec{\Phi}(t)}}\ dt}\  
%\end{array}
%\ee
Then we have
\be
\label{n-ene-con-10}
\begin{array}{l}
\ds\lf|\int_0^T\int_\Sigma\chi^{h(t)}_{R,x_0}(x)\ e^{-2\al}\ \Delta_{h(t)}\vec{H}\cdot \ov{\p}\lf[|\vec{H}|^2\,\p\vec{\Phi}\rg]\ dt\rg|\\[5mm]
\ds\le\int_0^T\lf|\int_\Sigma\chi^{h(t)}_{R,x_0}(x)\ e^{-2\al}\ \Delta_{h(t)}\vec{H}\cdot \p\vec{\Phi}\wedge\ov{\p}\lf[|\vec{H}|^2\rg]\ dt\rg|\\[5mm]
\ds +2^{-1}\,\lf|\int_0^T\int_\Sigma\chi^{h(t)}_{R,x_0}(x)\ e^{-2\al}\  \Delta_{h(t)}\vec{H}\cdot \Delta_{h(t)}\vec{\Phi}\, |\vec{H}|^2\ \ dvol_{h(t)}\ dt\rg|
\end{array}
\ee
We control the two terms in the r.h.s. of (\ref{n-ene-con-10}). We  now consider the first term.
\be
\label{n-ene-con-10-1}
\begin{array}{l}
\ds\lf|\int_0^T\int_\Sigma\chi^{h(t)}_{R,x_0}(x)\ \ e^{-2\al}\  \Delta_{h(t)}\vec{H}\cdot \p\vec{\Phi}\wedge\ov{\p}\lf[|\vec{H}|^2\rg]\ dt\rg|\\[5mm]
\ds\le\int_0^T\lf[\int_\Sigma\chi^{h(t)}_{R,x_0}(x)\ \ e^{-2\al}\   |\Delta_{h(t)}\vec{H}|^2\ \ dvol_{h(t)}\rg]^{1/2}\ F_{1}(R,x_0,t) \ dt
\end{array}
\ee
where
\be
\label{n-ene-con-10-2}
\begin{array}{l}
\ds F_{1}(R,x_0,t) :=e^{-\al_{R,x_0}(t)}\ \lf[\int_\Sigma\chi^{h(t)}_{R,x_0}(x)\ |d\vec{H}|^4_{h(t)}\  \ dvol_{h(t)}\rg]^{1/4}\ \lf[\int_\Sigma\chi^{h(t)}_{R,x_0}(x)\ |\vec{H}|^4_{h(t)}\  \ dvol_{h(t)}\rg]^{1/4}
\end{array}
\ee
%where $$\al_{R,x_0}(t)=\max_{x\in B_R^{h(t)}(x_0)}\al(t)(x)$$. Classical interpolation inequalities give respectively
%\be
%\label{ene-con-10-21}
%\begin{array}{l}
%\ds\lf[\int_{B_R^{h(t)}(x_0)}\ |dH|^4_{h(t)}\  \ dvol_{h(t)}\rg]^{1/4}\le C\ \lf[\int_{B_{2\,R}^{h(t)}(x_0)}\ |dH|^2_{h(t)}\  \ dvol_{h(t)}\rg]^{1/4}\ \ \lf[\int_{B_{2\,R}^{h(t)}(x_0)}\ |\Delta H|^2_{h(t)}\  \ dvol_{h(t)}\rg]^{1/4}\\[5mm]
%\ds+ C\ \lf[R^{-2}\,\int_{B_{2\,R}^{h(t)}(x_0)}\ |dH|^2_{h(t)}\  \ dvol_{h(t)}\rg]^{1/4}\ 
%\end{array}
%\ee
Recall
\be
\label{n-ene-con-10-21}
\begin{array}{l}
\ds\lf[\int_{B_R^{h(t)}(x_0)}\ |d\vec{H}|^4_{h(t)}\  \ dvol_{h(t)}\rg]^{1/4}\\[5mm]
\ds\le C\ \lf[\int_{B_{2\,R}^{h(t)}(x_0)}\ |\vec{H}|^2\  \ dvol_{h(t)}\rg]^{1/8}\ \ \lf[\int_{B_{2\,R}^{h(t)}(x_0)}\ |\Delta_{h(t)} H|^2\  \ dvol_{h(t)}\rg]^{3/8}\\[5mm]
\ds+ C\ \lf[R^{-3}\,\int_{B_{2\,R}^{h(t)}(x_0)}\ |\vec{H}|^2  \ dvol_{h(t)}\rg]^{1/2} 
\end{array}
\ee
and 
\be
\label{n-ene-con-10-22}
\begin{array}{l}
\ds\lf[\int_{B_R^{h(t)}(x_0)} |\vec{H}|^4\  \ dvol_{h(t)}\rg]^{1/4}\\[5mm]
\ds\le C\ \lf[\int_{B_{2\,R}^{h(t)}(x_0)}\ |\vec{H}|^2\  \ dvol_{h(t)}\rg]^{3/8}\ \ \lf[\int_{B_{2\,R}^{h(t)}(x_0)}\ |\Delta_{h(t)}\vec{H}|^2\  \ dvol_{h(t)}\rg]^{1/8}\\[5mm]
\ds+ C\ \lf[R^{-1}\,\int_{B_{2\,R}^{h(t)}(x_0)}\ |\vec{H}|^2\  \ dvol_{h(t)}\rg]^{1/2}\ 
\end{array}
\ee
Hence combining (\ref{n-ene-con-10-1})...(\ref{n-ene-con-10-22}) gives
\be
\label{n-ene-con-10-23}
\begin{array}{l}
\ds\lf|\int_0^T\int_\Sigma\chi^{h(t)}_{R,x_0}(x)\ e^{-2\al}\ \Delta_{h(t)}\vec{H}\cdot \p\vec{\Phi}\wedge\ov{\p}\lf[|\vec{H}|^2\rg]\ dt\rg|\\[5mm]
\ds\le\,C\,\int_0^T\int_{B_{2\,R}^{h(t)}(x_0)}\ e^{-2\,\al_{R,x_0}(t)} |\Delta_{h(t)}\vec{H}|^2\  \ dvol_{h(t)}\   \lf[E(2R,x_0,T)\rg]^{1/2}\ dt\\[5mm]
\ds+\, C\,\int_0^T\lf[\int_{B_{2\,R}^{h(t)}(x_0)}\ e^{-2\,\al_{R,x_0}(t)}\  |\Delta_{h(t)} H|^2\  \ dvol_{h(t)}\rg]^{7/8}\ \frac{e^{-\al_{R,x_0}(t)/2}}{\sqrt{R}}\ \lf[E(2R,x_0,T)\rg]^{5/8}\ dt\\[5mm]
\ds+\, C\,\int_0^T\lf[\int_{B_{2\,R}^{h(t)}(x_0)}\ e^{-2\,\al_{R,x_0}(t)}\ |\Delta_{h(t)} H|^2  \ dvol_{h(t)}\rg]^{5/8}\ \frac{e^{-3\al_{R,x_0}(t)/2}}{R^{3/2}}\  \lf[E(2R,x_0,T)\rg]^{7/8}\ dt\\[5mm]
\ds+\, C\,\int_0^T\lf[\int_{B_{2\,R}^{h(t)}(x_0)}\  e^{-2\,\al_{R,x_0}(t)}\ |\Delta_{h(t)} H|^2\  \ dvol_{h(t)}\rg]^{1/2}\  e^{-2\,\al_{R,x_0}(t)}\ \frac{1}{R^2}\  E(2R,x_0,T)\ dt\\[5mm]
\ds\le \sup_{[0,T]} \lf[\int_{B_{2\,R}^{h(t)}(x_0)}\ |\vec{H}|^2\  e^{2\,\al_{R,x_0}(t)}\ dvol_{h(t)}\rg]^{1/2}\ \lf[A+\frac{A^{7/8}\, B^{1/8}}{\sqrt{R}}+\frac{A^{5/8}\, B^{3/8}}{R^{3/2}}+\frac{A^{1/2}\,B^{1/2}}{R^2}\rg]
\end{array}
\ee
Where we introduced the notations respectively
\[
A:=\int_0^T\int_{B_{2\,R}^{h(t)}(x_0)}\ e^{-2\,\al_{R,x_0}(t)}\  |\Delta_{h(t)}\vec{H}|^2\  \ dvol_{h(t)}\ \quad\mbox{ and }\quad\ B:=\int_0^T\int_{B_{2\,R}^{h(t)}(x_0)}\ |\vec{H}|^2  \ e^{2\,\al_{R,x_0}(t)}\  dvol_{h(t)}\ 
\]
We bound respectively 
\[
\frac{A^{7/8}\, B^{1/8}}{\sqrt{\rho}}\le \frac{7\,A+\rho^{-4}\,B}{8}\quad,\quad\frac{A^{5/8}\, B^{3/8}}{\rho^{3/2}}\le \frac{4\,A+3\,\rho^{-4}\,B}{8}\quad\mbox{ and }\quad\frac{A^{1/2}\,B^{1/2}}{\rho^2}\le \frac{A+\rho^{-4}B}{2}\ .
\]
where $\rho:=\ R\ e^{2\,\al_{R,x_0}(t)}$.   This gives finally
\be
\label{n-ene-con-10-24}
\begin{array}{l}
\ds\lf|\int_0^T\int_\Sigma\chi^{h(t)}_{R,x_0}(x)\ \Delta_{h(t)}\vec{H}\cdot \p\vec{\Phi}\wedge\ov{\p}\lf[|\vec{H}|^2\rg]\ dt\rg|\\[5mm]
\ds\le C\, \sup_{[0,T]} \lf[\int_{B_{2\,R}^{h(t)}(x_0)}\ |\vec{H}|^2\  e^{2\,\al_{R,x_0}(t)}\ dvol_{h(t)}\rg]^{1/2}\ \int_0^T\int_{B_{2\,R}^{h(t)}(x_0)}\  e^{-2\,\al_{R,x_0}(t)}\  |\Delta_{h(t)}\vec{H}|^2\  \ dvol_{h(t)}\ dt\\[5mm]
\ds+C \, e^{4\|\al\|_\infty}\ \frac{T}{R^4}\  \sup_{[0,T]} \lf[\int_{B_{2\,R}^{h(t)}(x_0)}\ |\vec{H}|^2\  e^{2\,\al_{R,x_0}(t)}\ dvol_{h(t)}\rg]^{3/2}\ 
\end{array}
\ee
We control the second term in the r.h.s. of (\ref{n-ene-con-10})
\be
\label{n-ene-con-30}
\begin{array}{l}
\ds +2^{-1}\,\lf|\int_0^T\int_\Sigma\chi^{h(t)}_{R,x_0}(x)\, \  e^{-\,2\,\al_{R,x_0}(t)}\ \Delta_{h(t)}\vec{H}\cdot \Delta_{h(t)}\vec{\Phi}\, |\vec{H}|^2\ \ dvol_{h(t)}\ dt\rg|\\[5mm]
\ds\le C\,\int_0^T\lf[\int_{B_{2\,R}^{h(t)}(x_0)} \ e^{-\,2\,\al_{R,x_0}(t)}\  |\Delta_{h(t)}\vec{H}|^2\  \ dvol_{h(t)}\rg]^{1/2}\,e^{\al_{R,x_0}(t)}\,\lf[\int_{B_{2\,R}^{h(t)}(x_0)}\ |\vec{H}|^6\  \ dvol_{h(t)}\rg]^{1/2}\ dt
\end{array}
\ee
Classical interpolation inequalities gives
\be
\label{n-ene-con-31}
\begin{array}{l}
\ds\lf[\int_{B_R^{h(t)}(x_0)} |\vec{H}|^6\  \ dvol_{h(t)}\rg]^{1/2}\\[5mm]
\ds\le C\ \int_{B_{2\,R}^{h(t)}(x_0)}\ |\vec{H}|^2\  \ dvol_{h(t)}\ \ \lf[\int_{B_{2\,R}^{h(t)}(x_0)}\ |\Delta_{h(t)}\vec{H}|^2\  \ dvol_{h(t)}\rg]^{1/2}\\[5mm]
\ds+ C\ \lf[R^{-4/3}\,\int_{B_{2\,R}^{h(t)}(x_0)}\ |\vec{H}|^2\  \ dvol_{h(t)}\rg]^{3/2}\ 
\end{array}
\ee
Hence we obtain
\be
\label{n-ene-con-32}
\begin{array}{l}
\ds +2^{-1}\,\lf|\int_0^T\int_\Sigma\chi^{h(t)}_{R,x_0}(x)\  e^{-\,2\,\al_{R,x_0}(t)}\  \Delta_{h(t)}\vec{H}\cdot \Delta_{h(t)}\vec{\Phi}\, |\vec{H}|^2\ \ dvol_{h(t)}\ dt\rg|\\[5mm]
\ds\le C\,\sup_{[0,T]} \lf[\int_{B_{2\,R}^{h(t)}(x_0)}\ |\vec{H}|^2\  e^{2\,\al_{R,x_0}(t)}\ dvol_{h(t)}\rg]\,\int_0^T\int_{B_{2\,R}^{h(t)}(x_0)} \ e^{-\,2\,\al_{R,x_0}(t)}\  |\Delta_{h(t)}\vec{H}|^2\  \ dvol_{h(t)}\ dt\\[5mm]
\ds+C\,\int_0^T\lf[\int_{B_{2\,R}^{h(t)}(x_0)}\  e^{-\,2\,\al_{R,x_0}(t)}\ |\Delta_{h(t)}\vec{H}|^2\  \ dvol_{h(t)}\rg]^{1/2}\,\frac{e^{\al_{R,x_0}(t)}}{R^2}\ \lf[\,\int_{B_{2\,R}^{h(t)}(x_0)}\ |\vec{H}|^2\  \ dvol_{h(t)}\rg]^{3/2}\ dt
\end{array}
\ee
We write
\be
\label{n-ene-con-33}
\begin{array}{l}
\ds\int_0^T\lf[\int_{B_{2\,R}^{h(t)}(x_0)}\  e^{-\,2\,\al_{R,x_0}(t)}\ |\Delta_{h(t)}\vec{H}|^2\  \ dvol_{h(t)}\rg]^{1/2}\,\frac{e^{\al_{R,x_0}(t)}}{R^2}\ \lf[\,\int_{B_{2\,R}^{h(t)}(x_0)}\ |\vec{H}|^2\  \ dvol_{h(t)}\rg]^{3/2}\ dt\\[5mm]
\ds\le \sup_{[0,T]} \lf[\int_{B_{2\,R}^{h(t)}(x_0)}\ |\vec{H}|^2\  e^{2\,\al_{R,x_0}(t)}\ dvol_{h(t)}\rg] \int_0^T\int_{B_{2\,R}^{h(t)}(x_0)}\  e^{-\,2\,\al_{R,x_0}(t)}\ |\Delta_{h(t)}\vec{H}|^2\  \ dvol_{h(t)}\ dt\\[5mm]
\ds+ R^{-4}\  \int_0^T\lf[\int_{B_{2\,R}^{h(t)}(x_0)}\ |\vec{H}|^2\  \ dvol_{h(t)}\rg]^2\ dt\ .
\end{array}
\ee
This gives finally
\be
\label{n-ene-con-34}
\begin{array}{l}
\ds +2^{-1}\,\lf|\int_0^T\int_\Sigma\chi^{h(t)}_{R,x_0}(x)\  e^{-\,2\,\al}\  \Delta_{h(t)}\vec{H}\cdot \Delta_{h(t)}\vec{\Phi}\, |\vec{H}|^2\ \ dvol_{h(t)}\ dt\rg|\\[5mm]
\ds\le C\,\sup_{[0,T]} \lf[\int_{B_{2\,R}^{h(t)}(x_0)}\ |\vec{H}|^2\  e^{2\,\al_{R,x_0}(t)}\ dvol_{h(t)}\rg] \int_0^T\int_{B_{2\,R}^{h(t)}(x_0)}\  e^{-\,2\,\al_{R,x_0}(t)}\ |\Delta_{h(t)}\vec{H}|^2\  \ dvol_{h(t)}\ dt\\[5mm]
\ds\quad+\frac{T}{R^4}\,e^{4\|\al\|_\infty}\ \sup_{[0,T]} \lf[\int_{B_{2\,R}^{h(t)}(x_0)}\ |\vec{H}|^2\  e^{2\,\al_{R,x_0}(t)}\ dvol_{h(t)}\rg]^2\ 
\end{array}
\ee
We are now controlling
\be
\label{n-ene-con-35}
\ds\lf|\int_0^T \int_\Sigma\chi^{h(t)}_{R,x_0}(x)\  e^{-\,2\,\al}\  \Delta_{h(t)}\vec{H}\cdot \ov{\p}\lf[\vec{H}\cdot\vec{\frak{h}}^{\,0}\res \ov{\p}\vec{\Phi}  \rg]\\[5mm] \ dt\rg|
\ee
%In local conformal coordinates we compute
%\[
%\Delta_{h(t)}\vec{H}\cdot \p\vec{\Phi}=2\, e^{-2\nu}\,\p_{{z}}\p_{\ov{z}}\vec{H}\cdot\p_z\vec{\Phi}\ dz=2\, e^{-2\nu}\,\p_{z}\lf[\p_{\ov{z}}\vec{H}\cdot\p_z\vec{\Phi}\rg]-2\, e^{-2\nu}\,\p_{\ov{z}}\vec{H}\cdot\p^2_{z^2}\vec{\Phi}
%\]
To that aim we will make use of the following lemma  based on Codazzi Mainardi identity.
\begin{Lm}
\label{lm-n-cod-main} We have respectively
\be
\label{n-cod-main-1}
\ov{\p}\lf(\vec{H}\cdot\vec{\frak{h}}^{\,0}\res g_{\vec{\Phi}}\rg)=\ov{\p}\vec{H}\,\dot{\otimes}\,\vec{\frak{h}}^{\,0}\res g_{\vec{\Phi}}+2^{-1}\,g_{\vec{\Phi}}\otimes \vec{H}\cdot\p\vec{H}
\ee
and
\be
\label{n-cod-main-2}
\begin{array}{rl}
\ds\ov{\p}\lf[\vec{H}\cdot\vec{\frak{h}}^{\,0}\res \ov{\p}\vec{\Phi}\rg]&\ds=  \lf(\ov{\p}\vec{H}\,\dot{\otimes}\,\vec{\frak{h}}^{\,0}\res g_{\vec{\Phi}}+2^{-1}\,\hat{g}_{\vec{\Phi}}\otimes \vec{H}\cdot\p\vec{H} \rg)\res\lf(g_{\vec{\Phi}}^{-1} \res\ov{\p}\vec{\Phi}\rg)\\[5mm]
\ds&\ds+2\,\lf(\vec{H}\cdot\vec{\frak{h}}^{\,0}\res g_{\vec{\Phi}}\rg)\res\ov{\vec{\frak h}^{\,0}}
%\ov{\p}\lf(g_{\vec{\Phi}}^{-1} \res\ov{\p}\vec{\Phi}\rg)
\end{array}
\ee
whre $\hat{g}_{\vec{\Phi}}$ is the component of ${g}_{\vec{\Phi}}$ in $\wedge^{0,1}\Sigma\otimes\wedge^{1,0}\Sigma$. This implies in particular
\be
\label{n-cod-main-3}
\lf|\ov{\p}\lf[\vec{H}\cdot\vec{\frak{h}}^{\,0}\res \ov{\p}\vec{\Phi}\rg]\rg|_{h(t)}\le e^{\al}\, \lf[|\p\vec{H}|_{h(t)}\ |\vec{\frak{h}}^{\,0}|_{h(t)}+|\vec{H}|\ |\p\vec{H}|_{h(t)}+2\, e^\al\, |\vec{H}|\ |\vec{\frak{h}}^{\,0}|_{h(t)}^2\rg]\ .
\ee
\hfill $\Box$
\end{Lm}
\noindent{\bf Proof of Lemma~\ref{lm-n-cod-main}} We have
\[
\begin{array}{l}
\ds\ov{\p}\lf[\vec{H}\cdot\vec{\frak{h}}^{\,0}\res \ov{\p}\vec{\Phi}\rg]=\ov{\p}\lf[(\vec{H}\cdot\vec{\frak{h}}^{\,0}\res g_{\vec{\Phi}})\res\lf(g_{\vec{\Phi}}^{-1} \res\ov{\p}\vec{\Phi}\rg)\rg]\\[5mm]
\ds=\ov{\p}\lf(\vec{H}\cdot\vec{\frak{h}}^{\,0}\res g_{\vec{\Phi}}\rg)\res\lf(g_{\vec{\Phi}}^{-1} \res\ov{\p}\vec{\Phi}\rg)+\lf(\vec{H}\cdot\vec{\frak{h}}^{\,0}\res g_{\vec{\Phi}}\rg)\res\ov{\p}\lf(g_{\vec{\Phi}}^{-1} \res\ov{\p}\vec{\Phi}\rg)
\end{array}
\]
In local conformal coordinates 
\[
\vec{H}\cdot\vec{\frak{h}}^{\,0}\res g_{\vec{\Phi}}= e^{-2\la}\ \p_{z}\p_{\ov{z}}\vec{\Phi}\cdot\pi_{\vec{n}}\p_{z}\p_{{z}}\vec{\Phi} \ dz\otimes dz
\]
Hence
\[
\begin{array}{rl}
\ds\ov{\p}\lf(\vec{H}\cdot\vec{\frak{h}}^{\,0}\res g_{\vec{\Phi}}\rg)&\ds= \p_{\ov{z}}\lf(e^{-2\la}\ \p_{z}\p_{\ov{z}}\vec{\Phi}\rg)\cdot\pi_{\vec{n}}\p_{z}\p_{{z}}\vec{\Phi} \ d\ov{z}\otimes dz\otimes dz\\[5mm]
&\ds+ e^{-2\la}\ \p_{z}\p_{\ov{z}}\vec{\Phi}\cdot \pi_{\vec{n}}\p_{\ov{z}}\lf( \pi_{\vec{n}}\p_{z}\p_{{z}}\vec{\Phi}  \rg)\ d\ov{z}\otimes dz\otimes dz
\end{array}
\]
We have (using $\p_{z}\vec{\Phi}\cdot\p_{z}\vec{\Phi}=0$)
\[
\pi_{\vec{n}}\p_{z}\p_{{z}}\vec{\Phi}=\p_z\p_z\vec{\Phi}-2\, \p_{\ov{z}}\vec{\Phi}\cdot \p_{z}\p_{{z}}\vec{\Phi}\ e^{-2\la}\,\p_z\vec{\Phi}
\]
Thus
\[
\begin{array}{l}
\ds\pi_{\vec{n}}\p_{\ov{z}}\lf(\pi_{\vec{n}}\p_{z}\p_{{z}}\vec{\Phi}\rg)=\pi_{\vec{n}}\p_{\ov{z}}\p_z\p_z\vec{\Phi}-2\, \p_{\ov{z}}\vec{\Phi}\cdot \p_{z}\p_{{z}}\vec{\Phi}\ e^{-2\la}\,\p_{\ov{z}}\p_z\vec{\Phi}\\[5mm]
\ds=\pi_{\vec{n}}\lf[\p_{z}\lf[2^{-1}\, e^{2\la}\,\vec{H}\rg] -\p_{z}|\p_z\vec{\Phi}|^2 \ \vec{H}\rg]=2^{-1}\, e^{2\la}\, \pi_{\vec{n}}\p_z\vec{H}
\end{array}
\]
Finally we have obtained
\[
\ov{\p}\lf(\vec{H}\cdot\vec{\frak{h}}^{\,0}\res g_{\vec{\Phi}}\rg)=\ov{\p}\vec{H}\,\dot{\otimes}\,\vec{\frak{h}}^{\,0}\res g_{\vec{\Phi}}+2^{-1}\,\hat{g}_{\vec{\Phi}}\otimes \vec{H}\cdot\p\vec{H}
\]
which is the desired identity and  lemma~\ref{lm-n-cod-main} is proved.\hfill $\Box$.

\medskip

\noindent{\bf Proof of lemma~\ref{lm-conc-evol} continued} Using Lemma~\ref{lm-n-cod-main}
\be
\label{n-ene-con-40}
\begin{array}{l}
\ds\lf|\int_0^T \int_\Sigma\chi^{h(t)}_{R,x_0}(x)\   e^{-\,2\,\al}\ \Delta_{h(t)}\vec{H}\cdot \ov{\p}\lf[\vec{H}\cdot\vec{\frak{h}}^{\,0}\res \ov{\p}\vec{\Phi}  \rg] \ dt\rg|\\[5mm]
\ds\le 2\,\int_0^T \int_\Sigma\chi^{h(t)}_{R,x_0}(x)\ |\Delta_{h(t)}\vec{H}|\ e^{-\al}\, \lf[|\p\vec{H}|_{h(t)}\ |\vec{\frak{h}}^{\,0}|_{h(t)}+|\vec{H}|\ |\p\vec{H}|_{h(t)}+\,e^\al\, |\vec{H}|\ |\vec{\frak{h}}^{\,0}|_{h(t)}^2\rg]\ dvol_{h(t)}\ dt
\end{array}
\ee
Observe first that the term in the middle has already been considered in (\ref{n-ene-con-10-1}). We bound now the first one in thee following
\be
\label{n-ene-con-41}
\begin{array}{l}
\ds\int_0^T \int_\Sigma\chi^{h(t)}_{R,x_0}(x)\ |\Delta_{h(t)}\vec{H}|\ e^{-\al}\, |\p\vec{H}|_{h(t)}\ |\vec{\frak{h}}^{\,0}|_{h(t)}\ dvol_{h(t)}\  dt\\[5mm]
\ds\ds\le\int_0^T\lf[\int_\Sigma\chi^{h(t)}_{R,x_0}(x)\   e^{-\,2\,\al_{R,x_0}(t)}\   |\Delta_{h(t)}\vec{H}|^2\ \ dvol_{h(t)}\rg]^{1/2}\ F_{2}(R,x_0,t) \ dt
\end{array}
\ee
where
\be
\label{n-ene-con-42}
\begin{array}{l}
\ds F_{2}(R,x_0,t) :=e^{\al_{R,x_0}(t)}\ \lf[\int_\Sigma\chi^{h(t)}_{R,x_0}(x)\ |dH|^4_{h(t)}\  \ dvol_{h(t)}\rg]^{1/4}\ \lf[\int_\Sigma\chi^{h(t)}_{R,x_0}(x)\ |\vec{\frak{h}}^{\,0}|^4_{h(t)}\  \ dvol_{h(t)}\rg]^{1/4}
\end{array}
\ee
where $\al_{R,x_0}(t)=\max_{x\in B_R^{h(t)}(x_0)}\al(t)(x)$. We choose the local conformal coordinates given by lemma~\ref{lm-controled-charts}
\[
|\vec{\frak{h}}^{\,0}|_{h(t)}\le e^{-2\al_{R,x_0}(t)}\ |\p_{z}\p_{z}\vec{\Phi}|
\]
Classical Calderon-Zygmund theory gives then
\be
\label{n-CZ}
\int_{B^{h(t)}_R(x_0)}|\p_{z}\p_{z}\vec{\Phi}|^4\ \ dx^2\le C\ \int_{B^{h(t)}_{3\,R/2}(x_0)}|\p_{z}\p_{\ov{z}}\vec{\Phi}|^4\ dx^2+R^{-2}\,\lf[\int_{B^{h(t)}_{3R/2}(x_0)}|\p_{z}\p_{z}\vec{\Phi}|^2\ dx^2\rg]^2
\ee
Recall that from lemma~\ref{lm-L2-1} we have
\be
\label{n-osc-cont-ceff-conf}
\|\al(t)(x)-\al(t)(y)\|_{L^{\infty}(B^{h(t)}_{2R}(x_0))\times B^{h(t)}_{2R}(x_0))}\le C\ \lf[(1+l_{\ast\ast}^{-2})\,E(0)+1\rg]
\ee
Multiplying inequality (\ref{n-CZ}) by $e^{-8\al_{R,x_0}(t)}$  and using (\ref{n-osc-cont-ceff-conf}) we obtain
\be
\label{n-h-0-H-L4}
\begin{array}{l}
\ds\lf[\int_\Sigma\chi^{h(t)}_{R,x_0}(x)\ |\vec{\frak{h}}^{\,0}|^4_{h(t)}\  \ dvol_{h(t)}\rg]^{1/4}\le  e^{C\ (1+l_{\ast\ast}^{-2}(T))\,E(0)}\ \lf[\int_{B_{3\,R/2}^{h(t)}(x_0)}\ |\vec{H}|^4\ dvol_{h(t)}\rg]^{1/4}\\[5mm]
\ds +e^{C\ (1+l_{\ast\ast}^{-2}(T))\,E(0)}\ \sqrt{R^{-1}}\ \lf[\int_{B_{2\,R}^{h(t)}(x_0)}\ |\vec{\frak{h}^{\,0}}|^2\ dvol_{h(t)}\rg]^{1/2}
\end{array}
\ee
Arguing as for the estimation of (\ref{n-ene-con-10-1}) we obtain
\be
\label{n-ene-con-3-bis}
\begin{array}{l}
\ds e^{-C\ (1+l_{\ast\ast}^{-2}(T))\,E(0)} \ \int_0^T \int_\Sigma\chi^{h(t)}_{R,x_0}(x)\ |\Delta_{h(t)}\vec{H}|\ e^{-\al}\, |\p\vec{H}|_{h(t)}\ |\vec{\frak{h}}^{\,0}|_{h(t)}\ dvol_{h(t)}\  dt\\[5mm]
\ds\le  \sup_{[0,T]} \lf[\int_{B_{2\,R}^{h(t)}(x_0)}\  |\vec{\mathbb I}|_{g_{\vec{\Phi}(t)}}^2\    e^{2\,\al_{R,x_0}(t)}\ dvol_{h(t)}\rg]^{1/2}\ \int_0^T\int_{B_{2\,R}^{h(t)}(x_0)}\  e^{-\,2\,\al_{R,x_0}(t)}\   |\Delta_{h(t)}\vec{H}|^2\  \ dvol_{h(t)}\ dt\\[5mm]
\ds\quad\quad\,+ \,e^{4\|\al\|_\infty}\ \frac{T}{R^4}\  \sup_{[0,T]} \lf[\int_{B_{2\,R}^{h(t)}(x_0)}\ |\vec{\mathbb I}|_{g_{\vec{\Phi}(t)}}^2\  e^{2\,\al_{R,x_0}(t)}\ dvol_{h(t)}\rg]^{3/2}\ 
\end{array}
\ee
Now for the third contribution in (\ref{n-ene-con-40}) we proceed as follows
\be
\label{n-ene-con-43}
\begin{array}{l}
\ds\int_0^T \int_\Sigma\chi^{h(t)}_{R,x_0}(x)\ |\Delta_{h(t)}\vec{H}|\  |\vec{H}|\ |\vec{\frak{h}}^{\,0}|_{h(t)}^2\ dvol_{h(t)}\ dt\\[5mm]
\ds\ds\le\int_0^T\lf[\int_\Sigma\chi^{h(t)}_{R,x_0}(x)\  e^{-\,2\,\al_{R,x_0}(t)}\   |\Delta_{h(t)}\vec{H}|^2\ \ dvol_{h(t)}\rg]^{1/2}\ \ F_{3}(R,x_0,t) \ dt
\end{array}
\ee
where
\be
\label{n-ene-con-44}
\begin{array}{l}
\ds F_{3}(R,x_0,t) :=e^{\al_{R,x_0}(t)}\ \lf[\int_\Sigma\chi^{h(t)}_{R,x_0}(x)\ |\vec{H}|^6\  \ dvol_{h(t)}\rg]^{1/6}\ \lf[\int_\Sigma\chi^{h(t)}_{R,x_0}(x)\ |\vec{\frak{h}}^{\,0}|^6_{h(t)}\  \ dvol_{h(t)}\rg]^{1/3}
\end{array}
\ee
The first quantity involving the $L^6$ norm of the mean curvature in (\ref{n-ene-con-44}) is estimated following (\ref{n-ene-con-31}) while the last quantity, the use of the fact that the map
which to $\p^2_{z\ov{z}}\vec{\Phi}$ assigns $\p^2_{z{z}}\vec{\Phi}$ is Calderon Zygmund as we did in (\ref{n-CZ}) for the $L^4$ norm, for the $L^6$ norm is giving
 \be
\label{n-h-0-H-L6}
\begin{array}{l}
\ds\lf[\int_\Sigma\chi^{h(t)}_{R,x_0}(x)\ |\vec{\frak{h}}^{\,0}|^6_{h(t)}\  \ dvol_{h(t)}\rg]^{1/3}\le  e^{C\ (1+l_{\ast\ast}^{-2}(T))\,E(0)}\ \lf[\int_{B_{3\,R/2}^{h(t)}(x_0)}\ |\vec{H}|^6\ dvol_{h(t)}\rg]^{1/3}\\[5mm]
\ds +e^{C\ (1+l_{\ast\ast}^{-2}(T))\,E(0)}\  R^{-4/3}\,\int_{B_{2\,R}^{h(t)}(x_0)}\ |\vec{\frak{h}^{\,0}}|^2\ dvol_{h(t)}
\end{array}
\ee
This gives
\be
\label{n-ene-con-44-aa}
\begin{array}{l}
\ds F_{3}(R,x_0,t) :=e^{\al_{R,x_0}(t)}\ \lf[\int_{B_{3\,R/2}^{h(t)}(x_0)}\ |\vec{H}|^6\ dvol_{h(t)}\rg]^{1/2}\\[5mm]
\ds +e^{\al_{R,x_0}(t)}\ \lf[\int_{B_R^{h(t)}(x_0)}\ |\vec{H}|^6\ dvol_{h(t)}\rg]^{1/6}\ R^{-4/3}\,\int_{B_{2\,R}^{h(t)}(x_0)}\ |\vec{\frak{h}^{\,0}}|^2\ dvol_{h(t)}\\[5mm]
\ds\le C\,e^{\al_{R,x_0}(t)}\ \lf[\int_{B_{3\,R/2}^{h(t)}(x_0)}\ |\vec{H}|^6\ dvol_{h(t)}\rg]^{1/2}\\[5mm]
\ds\quad+ C\,e^{\al_{R,x_0}(t)}\ \lf[R^{-4/3}\,\int_{B_{2\,R}^{h(t)}(x_0)}\ |\vec{\frak{h}^{\,0}}|^2\ dvol_{h(t)}\rg]^{3/2}
\end{array}
\ee
Recall from  (\ref{n-ene-con-31})
\be
\label{n-ene-con-31-bis}
\begin{array}{l}
\ds\lf[\int_{B_R^{h(t)}(x_0)} |\vec{H}|^6\  \ dvol_{h(t)}\rg]^{1/2}\\[5mm]
\ds\le C\ \int_{B_{2\,R}^{h(t)}(x_0)}\ |\vec{H}|^2\  \ dvol_{h(t)}\ \ \lf[\int_{B_{2\,R}^{h(t)}(x_0)}\ |\Delta_{h(t)}\vec{H}|^2\  \ dvol_{h(t)}\rg]^{1/2}\\[5mm]
\ds+ C\ \lf[R^{-4/3}\,\int_{B_{2\,R}^{h(t)}(x_0)}\ |\vec{H}|^2\  \ dvol_{h(t)}\rg]^{3/2}\ 
\end{array}
\ee
Inserting (\ref{n-ene-con-31-bis}) into  (\ref{n-ene-con-44-aa}) gives
\be
\label{n-ene-con-45-a}
\begin{array}{l}
 e^{-\,C\ (1+l_{\ast\ast}^{-2}(T))\,E(0)}\ F_{3}(R,x_0,t) \\[5mm]
 \ds\le C\  e^{\al_{R,x_0}(t)}\  \int_{B_{2\,R}^{h(t)}(x_0)}\ |\vec{H}|^2\  \ dvol_{h(t)} \ \lf[\int_{B_{2\,R}^{h(t)}(x_0)}\ |\Delta_{h(t)}\vec{H}|^2\  \ dvol_{h(t)}\rg]^{1/2}\\[5mm]
\ds+ C\ e^{\al_{R,x_0}(t)}\ \lf[R^{-4/3}\,\int_{B_{2\,R}^{h(t)}(x_0)}\ |\vec{\mathbb I}|_{g_{\vec{\Phi}(t)}}^2\  \ dvol_{h(t)}\rg]^{3/2}\ 
\end{array}
\ee
Inserting (\ref{n-ene-con-45-a}) into (\ref{n-ene-con-43})  gives
\be
\label{n-ene-con-46}
\begin{array}{l}
\ds e^{-\,C\ (1+l_{\ast\ast}^{-2}(T))\,E(0)}\ \int_0^T \int_\Sigma\chi^{h(t)}_{R,x_0}(x)\ |\Delta_{h(t)}\vec{H}|\ |\vec{H}|\ |\vec{\frak{h}}^{\,0}|_{h(t)}^2\ dvol_{h(t)}\ dt\\[5mm]
\ds\le C\,\sup_{[0,T]} \lf[\int_{B_{2\,R}^{h(t)}(x_0)}\  |\vec{\mathbb I}|_{g_{\vec{\Phi}(t)}}^2\  e^{2\,\al_{R,x_0}(t)}\ dvol_{h(t)}\rg] \int_0^T\int_{B_{2\,R}^{h(t)}(x_0)}\ e^{-2\al_{R,x_0}(t)}\ |\Delta_{h(t)}\vec{H}|^2\  \ dvol_{h(t)}\ dt\\[5mm]
\ds+\int_0^T\lf[\int_\Sigma\chi^{h(t)}_{R,x_0}(x)\ e^{-2\al_{R,x_0}(t)}\ |\Delta_{h(t)}\vec{H}|^2\ \ dvol_{h(t)}\rg]^{1/2}\ e^{\al_{R,x_0}(t)}\ \lf[R^{-4/3}\,\int_{B_{2\,R}^{h(t)}(x_0)}\ |\vec{\mathbb I}|_{g_{\vec{\Phi}(t)}}^2\  \ dvol_{h(t)}\rg]^{3/2}\ .
\end{array}
\ee
We bound
\be
\label{n-ene-con-46-bis}
\begin{array}{l}
\ds\int_0^T\lf[\int_\Sigma\chi^{h(t)}_{R,x_0}(x)\ e^{-2\al_{R,x_0}(t)}\ |\Delta_{h(t)}\vec{H}|^2\ \ dvol_{h(t)}\rg]^{1/2}\ e^{\al_{R,x_0}(t)}\ \lf[R^{-4/3}\,\int_{B_{2\,R}^{h(t)}(x_0)}\ |\vec{\mathbb I}|_{g_{\vec{\Phi}(t)}}^2\  \ dvol_{h(t)}\rg]^{3/2}\\[5mm]
\ds \le \sup_{[0,T]}\ \lf[\int_{B_{2\,R}^{h(t)}(x_0)}\  |\vec{\mathbb I}|_{g_{\vec{\Phi}(t)}}^2\  e^{2\,\al_{R,x_0}(t)}\ dvol_{h(t)}\rg] \int_0^T\int_\Sigma\chi^{h(t)}_{R,x_0}(x)\ |\Delta_{h(t)}\vec{H}|^2\ \ dvol_{h(t)}\ dt\\[5mm]
\ds\quad+\  e^{2\|\al\|_\infty}\ \,\frac{T}{R^4}\ \sup_{[0,T]}\ \lf[\int_{B_{2\,R}^{h(t)}(x_0)}\  |\vec{\mathbb I}|_{g_{\vec{\Phi}(t)}}^2\  e^{2\,\al_{R,x_0}(t)}\ dvol_{h(t)}\rg]^2\ .
\end{array}
\ee
Combining (\ref{n-ene-con-46}) and (\ref{n-ene-con-46-bis}) gives finally
\be
\label{n-ene-con-46-ter}
\begin{array}{l}
\ds e^{-\,C\ (1+l_{\ast\ast}^{-2}(T))\,E(0)}\ \int_0^T \int_\Sigma\chi^{h(t)}_{R,x_0}(x)\ |\Delta_{h(t)}\vec{H}|\  |\vec{H}|\ |\vec{\frak{h}}^{\,0}|_{h(t)}^2\ dvol_{h(t)}\ dt\\[5mm]
\ds \le E(2R,x_0,t)\  \int_0^T\int_\Sigma\chi^{h(t)}_{R,x_0}(x)\ e^{-2\al_{R,x_0}(t)}\ |\Delta_{h(t)}\vec{H}|^2\ \ dvol_{h(t)}\ dt\\[5mm]
\ds\quad+\  e^{4\|\al\|_\infty}\ \,\frac{T}{R^4}\ E(2R,x_0,t)^2\  .
\end{array}
\ee
We now estimate the next term in (\ref{n-ene-con-7})
\be
\label{n-ene-con-47}
\begin{array}{l}
\ds  \int_0^T\lf|\int_\Sigma\Delta_{g_{\vec{\Phi}}}\lf[\chi^{h(t)}_{R,x_0}(x)\ \rg]\,\vec{H}\cdot \ov{\p}\p\vec{H}\rg|\ dt\\[5mm]
\ds\le \frac{C}{R^2}\, \int_0^T\int_{B_{R}^{h(t)}(x_0)} \ e^{-2\al}\ |\vec{H}|\, |\Delta_{\vec{h(t)}}\vec{H}|\ dvol_{h(t)}\ dt
%\ds+\frac{C}{R}\, \int_0^T\int_{B_{R}^{h(t)}(x_0)} \ |d\al|_{h(t)}\ |\vec{H}|\, |\Delta_{\vec{h(t)}}\vec{H}|\ dvol_{h(t)}\ dt\\[5mm]
%\ds+C\, \int_0^T\int_{B_{R}^{h(t)}(x_0)} \,|d\al|^2_{h(t)}\ |\vec{H}|\, |\Delta_{\vec{h(t)}}\vec{H}|\ dvol_{h(t)}\ dt\\[5mm]
%\ds+C\, \int_0^T\int_{B_{R}^{h(t)}(x_0)} \,|\Delta_{h(t)}\al|\ |\vec{H}|\, |\Delta_{{h(t)}}\vec{H}|\ dvol_{h(t)}\ dt
\end{array}
\ee
We first have for ay $\delta>0$ to be fixed later
\be
\label{n-ene-con-48}
\begin{array}{l}
\ds\int_0^T\int_\Sigma \,R^{-2}\ e^{-2\al}\ |\vec{H}|\, |\Delta_{\vec{h(t)}}\vec{H}|\ dvol_{h(t)}\ dt\\[5mm]
\ds\le C \frac{e^{2\,\|\al\|_\infty}}{R^2}\  \sqrt{E(2R,x_0,t)}\ \int_0^T\lf[ \int_{B_{R}^{h(t)}(x_0)}\ e^{-2\al_{R,x_0}(t)}\ |\Delta_{h(t)}\vec{H}|^2\ dvol_{h(t)} \rg]^{1/2}\, dt
\\[5mm]
\ds\le C \frac{e^{2\,\|\al\|_\infty}\sqrt{T}}{R^2}\  \sqrt{E(2R,x_0,t)}\  \sqrt{\int_0^T \int_{B_{R}^{h(t)}(x_0)}\ e^{-2\al_{R,x_0}(t)}\  |\Delta_{h(t)}\vec{H}|^2\ dvol_{h(t)}\ dt}\\[5mm]
\ds\le C\ \delta\ \int_0^T \int_{B_{R}^{h(t)}(x_0)}\ e^{-2\al_{R,x_0}(t)}\  |\Delta_{h(t)}\vec{H}|^2\ dvol_{h(t)}\ dt\\[5mm]
\ds\quad+ \frac{C}{\delta}\ \frac{e^{4\,\|\al\|_\infty}\,{T}}{R^4}\ E(2R,x_0,t)\ .
\end{array}
\ee
Next we have
\be
\label{n-ene-con-49}
\begin{array}{l}
\ds\int_0^T\lf|\int_\Sigma\Delta_{g_{\vec{\Phi}}}\lf[\chi^{h(t)}_{R,x_0}(x)\rg]\,\vec{H}\cdot \ov{\p}\lf[|\vec{H}|^2\,\p\vec{\Phi} \rg]\rg|\ dt\\[5mm]
\ds\le \frac{C}{R^2}\,\int_0^T\int_{B_{R}^{h(t)}(x_0)}\ e^{-\al}\ |\vec{H}|^2\ |\p\vec{H}|_{h(t)}\ dvol_{h(t)}\ dt\\[5mm]
\ds+\frac{C}{R^2}\,\int_0^T\int_{B_{R}^{h(t)}(x_0)}\  |\vec{H}|^4\ dvol_{h(t)}\ dt
\end{array}
\ee
We recall from (\ref{n-ene-con-10-22})
\be
\label{n-ene-con-49-1}
\begin{array}{l}
\ds \frac{C}{R^2}\ \int_0^T\int_{B_R^{h(t)}(x_0)} |\vec{H}|^4\  \ dvol_{h(t)}\ dt\\[5mm]
\ds\le  \frac{C}{R^2}\ \int_0^T\lf[\int_{B_{2\,R}^{h(t)}(x_0)}\ |\vec{H}|^2\  \ dvol_{h(t)}\rg]^{3/2}\ \ \lf[\int_{B_{2\,R}^{h(t)}(x_0)}\ |\Delta_{h(t)}\vec{H}|^2\  \ dvol_{h(t)}\rg]^{1/2}\ dt\\[5mm]
\ds + C\, \frac{e^{4\,\|\al\|_\infty}\,{T}}{R^4}\ E(2R,x_0,t)^2\\[5mm]
\ds \le E(2R,x_0,t)\ \int_0^T \int_{B_{2R}^{h(t)}(x_0)}\ e^{-2\al_{R,x_0}(t)}\  |\Delta_{h(t)}\vec{H}|^2\ dvol_{h(t)}\ dt+ C\, \frac{e^{4\,\|\al\|_\infty}\,{T}}{R^4}\ E(2R,x_0,t)^2
\end{array}
\ee
Interpolation inequalities give also
\be
\label{n-ene-con-49-2}
\begin{array}{rl}
\ds\int_{B_R^{h(t)}(x_0)} |\p\vec{H}|^2\  \ dvol_{h(t)}&\ds\le C\,\sqrt{\int_{B_{2\,R}^{h(t)}(x_0)}\ |\vec{H}|^2\  \ dvol_{h(t)}}\ \sqrt{\int_{B_{2\,R}^{h(t)}(x_0)}\ |\Delta_{h(t)}\vec{H}|^2\  \ dvol_{h(t)}}\\[5mm]
\ds &\ds\ +\frac{C}{R^2}\,\int_{B_{2R}^{h(t)}(x_0)} |\vec{H}|^2\  \ dvol_{h(t)}
\end{array}
\ee
Hence for any $\delta>0$
\be
\label{n-ene-con-49-3}
\begin{array}{l}
\ds\frac{C}{R^2}\,\int_0^T\int_{B_{R}^{h(t)}(x_0)}\ e^{-2\al}\ |\p\vec{H}|^2_{h(t)}\ dvol_{h(t)}\ dt\\[5mm]
\ds\le  C\ \delta\ \int_0^T \int_{B_{2R}^{h(t)}(x_0)}\ e^{-2\al_{R,x_0}(t)}\  |\Delta_{h(t)}\vec{H}|^2\ dvol_{h(t)}\ dt\\[5mm]
\ds\quad+ \frac{C}{\delta}\ \frac{e^{4\,\|\al\|_\infty}\,{T}}{R^4}\ E(2R,x_0,t)\ .
\end{array}
\ee
Next we have using lemma~\ref{lm-n-cod-main}
\be
\label{n-ene-con-50}
\begin{array}{l}
\ds\int_0^T\lf|\int_\Sigma\Delta_{g_{\vec{\Phi}}}\lf[\chi^{h(t)}_{R,x_0}(x)\rg]\,\vec{H}\cdot \ov{\p}\lf[\vec{H}\cdot\vec{\frak{h}}^{\,0}\res \ov{\p}\vec{\Phi}  \rg]\rg|\ dt\\[5mm]
\ds  \le \frac{C}{R^2}\, \int_0^T\int_\Sigma e^{-\al}\, \lf[|\p\vec{H}|_{h(t)}\ |\vec{\frak{h}}^{\,0}|_{h(t)}+|\vec{H}|\ |\p\vec{H}|_{h(t)}\rg]+\, |\vec{H}|\ |\vec{\frak{h}}^{\,0}|_{h(t)}^2\ dvol_{h(t)}\ dt
\end{array}
\ee
These integrals are estimated as the previous ones using in particular (\ref{n-h-0-H-L4}) in addition to (\ref{n-ene-con-10-22}). Next we estimate
\be
\label{n-ene-con-51}
\begin{array}{l}
\ds\int_0^T\lf|\int_\Sigma \lf<d\chi^{h(t)}_{R,x_0}(x)\ ,d\vec{H}\rg>_{g_{\vec{\Phi}}}\cdot \ov{\p}\p\vec{H}  \rg|\ dt\\[5mm]
\ds\le \frac{C}{R}\ \int_0^T \int_{B_{R}^{h(t)}(x_0)}\ e^{-2\al}\ |d\vec{H}|_{h(t)}\ |\Delta_{h(t)}\vec{H}|\ dvol_{h(t)}\ dt
\end{array}
\ee
We have
\be
\label{n-ene-con-51-b}
\begin{array}{l}
\ds\int_{B_{R}^{h(t)}(x_0)}\ e^{-2\al}\ |d\vec{H}|_{h(t)}\ |\Delta_{h(t)}\vec{H}|\ dvol_{h(t)}\ dt\\[5mm]
\ds\le \ e^{-\al_{R,x_0}(t)}\ \sqrt{\int_{B_{R}^{h(t)}(x_0)}\ |d\vec{H}|^2_{h(t)}\ dvol_{h(t)}}\  \sqrt{\int_{B_{R}^{h(t)}(x_0)}\ e^{-2\al_{R,x_0}(t)}\ |\Delta_{h(t)}\vec{H}|^2_{h(t)}\ dvol_{h(t)}}\\[5mm]
\end{array}
\ee
We use the fact
\be
\label{n-ene-con-51-c}
\begin{array}{l}
\ds  \sqrt{\int_{B_{R}^{h(t)}(x_0)}\ |d\vec{H}|^2_{h(t)}\ dvol_{h(t)}}\ \\[5mm]
\ds\le C\,\lf[\int_{B_{2\,R}^{h(t)}(x_0)}\ |\vec{H}|^2\  \ dvol_{h(t)}\rg]^{1/4}\ \lf[\int_{B_{2\,R}^{h(t)}(x_0)}\ |\Delta_{h(t)}\vec{H}|^2\  \ dvol_{h(t)}\rg]^{1/4}\\[5mm]
\ds\ +\frac{C}{R}\,\sqrt{\int_{B_{2R}^{h(t)}(x_0)} |\vec{H}|^2\  \ dvol_{h(t)}}\ ,
\end{array}
\ee
to deduce
\be
\label{re-n-ene-con-51}
\begin{array}{l}
\ds\int_0^T\lf|\int_\Sigma \lf<d\chi^{h(t)}_{R,x_0}(x)\ ,d\vec{H}\rg>_{g_{\vec{\Phi}}}\cdot \ov{\p}\p\vec{H}  \rg|\ dt\\[5mm]
\ds\le \frac{C}{R}\, \int_0^T e^{-\al_{R,x_0}(t)}\ \lf[\int_{B_{2\,R}^{h(t)}(x_0)}\ |\vec{H}|^2\ e^{2\al_{R,x_0}(t)}  \ dvol_{h(t)}\rg]^{1/4}\  \lf[\int_{B_{2\,R}^{h(t)}(x_0)}\ e^{-2\al_{R,x_0}(t)}\ |\Delta_{h(t)}\vec{H}|^2_{h(t)}\ dvol_{h(t)}\rg]^{3/4}\ dt\\[5mm]
\ds+\frac{C}{R^2}\, \int_0^T e^{-2\al_{R,x_0}(t)}\ \sqrt{\int_{B_{2\,R}^{h(t)}(x_0)}\ |\vec{H}|^2\ e^{2\al_{R,x_0}(t)}\ dvol_{h(t)}}\  \sqrt{\int_{B_{2\,R}^{h(t)}(x_0)}\ e^{-2\al_{R,x_0}(t)}\ |\Delta_{h(t)}\vec{H}|^2_{h(t)}\ dvol_{h(t)}}\ dt\\[5mm]
\ds \le \delta \int_0^T\int_{B_{2\,R}^{h(t)}(x_0)}\ e^{-2\al_{R,x_0}(t)}\ |\Delta_{h(t)}\vec{H}|^2_{h(t)}\ dvol_{h(t)}\ dt\ds+\ \frac{C}{\delta^3\,R^4}\, e^{4\|\al\|_\infty}\ T\ E(2R,x_0,T)\ .
\end{array}
\ee
Next we estimate
\be
\label{n-ene-con-52}
\begin{array}{l}
\ds\int_0^T\lf|\int_\Sigma \lf<d\chi^{h(t)}_{R,x_0}(x)\ ,d\vec{H}\rg>_{g_{\vec{\Phi}}}\cdot \ov{\p}\lf[|\vec{H}|^2\,\p\vec{\Phi}\rg]\rg|\ dt\\[5mm]
\ds\le \frac{C}{R}\ \int_0^T \int_{B_{R}^{h(t)}(x_0)}\ e^{-\al}\ |d\vec{H}|^2_{h(t)}\ |\vec{H}|\ dvol_{h(t)}\ dt\\[5mm]
\ds+\frac{C}{R}\ \int_0^T \int_{B_{R}^{h(t)}(x_0)}\ |d\vec{H}|_{h(t)}\ |\vec{H}|^3\ dvol_{h(t)}\ dt\ .
\end{array}
\ee
First we have
\be
\label{n-ene-con-52-1}
\begin{array}{l}
\ds \frac{C}{R}\  \int_{B_{R}^{h(t)}(x_0)}\ e^{-\al}\ |d\vec{H}|^2_{h(t)}\ |\vec{H}|\ dvol_{h(t)}\\[5mm]
\ds\le  \frac{C}{R}\ e^{-\al_{R,x_0}(t)}\ \sqrt{\int_{B_{R}^{h(t)}(x_0)}\ |d\vec{H}|^4_{h(t)}\ dvol_{h(t)}}\  \sqrt{\int_{B_{R}^{h(t)}(x_0)}\ |\vec{H}|^2\ dvol_{h(t)}}\ .
\end{array}
\ee
Recall
\be
\label{n-ene-con-52-2}
\begin{array}{l}
\ds\lf[\int_{B_R^{h(t)}(x_0)}\ |d\vec{H}|^4_{h(t)}\  \ dvol_{h(t)}\rg]^{1/2}\\[5mm]
\ds\le C\ \lf[\int_{B_{2\,R}^{h(t)}(x_0)}\ |\vec{H}|^2\  \ dvol_{h(t)}\rg]^{1/4}\ \ \lf[\int_{B_{2\,R}^{h(t)}(x_0)}\ |\Delta_{h(t)} H|^2\  \ dvol_{h(t)}\rg]^{3/4}\\[5mm]
\ds+ C\ R^{-3}\,\int_{B_{2\,R}^{h(t)}(x_0)}\ |\vec{H}|^2  \ dvol_{h(t)}\ .
\end{array}
\ee
Hence
\be
\label{n-ene-con-52-3}
\begin{array}{l}
\ds \frac{C}{R}\  \int_{B_{R}^{h(t)}(x_0)}\ e^{-\al}\ |d\vec{H}|^2_{h(t)}\ |\vec{H}|\ dvol_{h(t)}\\[5mm]
\ds\le \frac{C}{R}\ e^{-\al_{R,x_0}(t)}\  \lf[\int_{B_{2\,R}^{h(t)}(x_0)}\ |\vec{H}|^2\  e^{2\al_{R,x_0}(t)}\   \ dvol_{h(t)}\rg]^{3/4}\ \ \lf[\int_{B_{2\,R}^{h(t)}(x_0)}\ |\Delta_{h(t)} H|^2\  e^{-2\al_{R,x_0}(t)}\  \ dvol_{h(t)}\rg]^{3/4}\\[5mm]
\ds+ C\ R^{-4}\ e^{-4\al_{R,x_0}(t)}\ \lf[\int_{B_{2\,R}^{h(t)}(x_0)}\ |\vec{H}|^2 \  e^{2\al_{R,x_0}(t)}\   \ dvol_{h(t)}\rg]^{3/2}\ .
\end{array}
\ee
Thus
\be
\label{n-ene-con-52-4}
\begin{array}{l}
\ds \frac{C}{R}\ \int_0^T \int_{B_{R}^{h(t)}(x_0)}\ e^{-\al}\ |d\vec{H}|^2_{h(t)}\ |\vec{H}|\ dvol_{h(t)}\ dt\\[5mm]
\ds\le \sqrt{E(2R,x_0,T)}\ \int_{B_{2\,R}^{h(t)}(x_0)}\ |\Delta_{h(t)} H|^2\  e^{-2\al_{R,x_0}(t)}\  \ dvol_{h(t)}\\[5mm]
\ds\quad+ C\ R^{-4}\ e^{4\|\al\|_\infty}\ T\ E(2R,x_0,T)^{3/2}\ .
\end{array}
\ee
Then we have
\be
\label{n-ene-con-52-5}
\begin{array}{l}
\ds R^{-1}\, \int_{B_{R}^{h(t)}(x_0)}\ |d\vec{H}|_{h(t)}\ |\vec{H}|^3\ dvol_{h(t)}\\[5mm]
\ds\le R^{-1}\, \sqrt{\int_{B_{R}^{h(t)}(x_0)}\ |d\vec{H}|^2_{h(t)}\ dvol_{h(t)}}\  \sqrt{\int_{B_{R}^{h(t)}(x_0)}\ |\vec{H}|^6\ dvol_{h(t)}}
\end{array}
\ee
We have respectively
\be
\label{n-ene-con-52-6}
\begin{array}{l}
\ds  \sqrt{\int_{B_{R}^{h(t)}(x_0)}\ |d\vec{H}|^2_{h(t)}\ dvol_{h(t)}}\ \\[5mm]
\ds\le C\,\lf[\int_{B_{2\,R}^{h(t)}(x_0)}\ |\vec{H}|^2\  \ dvol_{h(t)}\rg]^{1/4}\ \lf[\int_{B_{2\,R}^{h(t)}(x_0)}\ |\Delta_{h(t)}\vec{H}|^2\  \ dvol_{h(t)}\rg]^{1/4}\\[5mm]
\ds\ +\frac{C}{R}\,\sqrt{\int_{B_{2R}^{h(t)}(x_0)} |\vec{H}|^2\  \ dvol_{h(t)}}\ ,
\end{array}
\ee
and
\be
\label{n-ene-con-52-7}
\begin{array}{l}
\ds  \sqrt{\int_{B_{R}^{h(t)}(x_0)}\ |\vec{H}|^6\ dvol_{h(t)}}\\[5mm]
\ds\le C\ \int_{B_{2\,R}^{h(t)}(x_0)}\ |\vec{H}|^2\  \ dvol_{h(t)}\ \ \lf[\int_{B_{2\,R}^{h(t)}(x_0)}\ |\Delta_{h(t)}\vec{H}|^2\  \ dvol_{h(t)}\rg]^{1/2}\\[5mm]
\ds+ C\ \lf[R^{-4/3}\,\int_{B_{2\,R}^{h(t)}(x_0)}\ |\vec{H}|^2\  \ dvol_{h(t)}\rg]^{3/2}\ 
\end{array}
\ee
Hence
\be
\label{n-ene-con-52-8}
\begin{array}{l}
\ds R^{-1}\, \int_{B_{R}^{h(t)}(x_0)}\ |d\vec{H}|_{h(t)}\ |\vec{H}|^3\ dvol_{h(t)}\\[5mm]
\ds\le \frac{C}{R}\ \lf[\int_{B_{2\,R}^{h(t)}(x_0)}\ |\vec{H}|^2\  \ dvol_{h(t)}\rg]^{5/4}\ \lf[\int_{B_{2\,R}^{h(t)}(x_0)}\ |\Delta_{h(t)}\vec{H}|^2\  \ dvol_{h(t)}\rg]^{3/4}\\[5mm]
\ds+\frac{C}{R^3}\ \lf[\int_{B_{2\,R}^{h(t)}(x_0)}\ |\vec{H}|^2\  \ dvol_{h(t)}\rg]^{7/4}\ \lf[\int_{B_{2\,R}^{h(t)}(x_0)}\ |\Delta_{h(t)}\vec{H}|^2\  \ dvol_{h(t)}\rg]^{1/4}\\[5mm]
\ds +\frac{C}{R^2}\ \lf[\int_{B_{2\,R}^{h(t)}(x_0)}\ |\vec{H}|^2\  \ dvol_{h(t)}\rg]^{3/2}\ \lf[\int_{B_{2\,R}^{h(t)}(x_0)}\ |\Delta_{h(t)}\vec{H}|^2\  \ dvol_{h(t)}\rg]^{1/2}\\[5mm]
\ds+  \frac{C}{ R^{4}}\,\lf[\int_{B_{2\,R}^{h(t)}(x_0)}\ |\vec{H}|^2  \ dvol_{h(t)}\rg]^{2} 
\end{array}
\ee
Thus
\be
\label{n-ene-con-52-9}
\begin{array}{l}
\ds R^{-1}\, \int_0^T\int_{B_{R}^{h(t)}(x_0)}\ |d\vec{H}|_{h(t)}\ |\vec{H}|^3\ dvol_{h(t)}\ dt\\[5mm]
\ds \le C\ E(2R,x_0,T)\ \int_0^T\int_{B_{2\,R}^{h(t)}(x_0)}\ |\Delta_{h(t)}\vec{H}|^2\  e^{-2\al(R,x_0,T)}\ dvol_{h(t)}\ dt\\[5mm]
\ds\quad+ C\ R^{-4}\ e^{4\|\al\|_\infty}\ T\ E(2R,x_0,T)^{2}\ .
\end{array}
\ee
Next we estimate, using lemma~\ref{lm-n-cod-main},
\be
\label{n-ene-con-53}
\begin{array}{l}
\ds\int_0^T\lf|\int_\Sigma \lf<d\chi^{h(t)}_{R,x_0}(x)\ ,d\vec{H}\rg>_{g_{\vec{\Phi}}}\cdot \ov{\p}\lf[\vec{H}\cdot\vec{\frak{h}}^{\,0}\res \ov{\p}\vec{\Phi}  \rg]\rg|\ dt\\[5mm]
\ds\le \frac{C}{R}\ \int_0^T \int_{B_{R}^{h(t)}(x_0)}\ e^{-\al}\  \lf[|\p\vec{H}|^2_{h(t)}\ |\vec{\frak{h}}^{\,0}|_{h(t)}+|\vec{H}|\ |\p\vec{H}|^2_{h(t)}\rg]+\, |\p\vec{H}|_{h(t)}\ |\vec{H}|\ |\vec{\frak{h}}^{\,0}|_{h(t)}^2\ dvol_{h(t)}\ dt
\end{array}
\ee
These integrals are estimated as the previous ones using in particular (\ref{n-h-0-H-L4}) and (\ref{n-h-0-H-L6}) in addition to (\ref{n-ene-con-10-21}), (\ref{n-ene-con-10-22}) and (\ref{n-ene-con-52-7}).

\medskip

It remains to estimates the two last terms in the r.h.s. of (\ref{n-ene-con-7}). The penultimate term in (\ref{n-ene-con-7}) is controlled as follows
\be
\label{n-ene-con-54}
\begin{array}{l}
\ds\int_0^T\lf|\int_\Sigma\chi^{h(t)}_{R,x_0}(x)\ \pi_T\lf(\ov{\p}\vec{H}\rg)\dot{\wedge}\lf[\ov{\p}\vec{\Phi}\res(I-\ti{P}_h)\lf(\pi_{\vec{n}}\frac{\p\vec{\Phi}}{\p t}\cdot\vec{\frak{h}}_0\rg)\rg]\rg|\ dt\\[5mm]
\ds\le \int_0^T\int_{B_{R}^{h(t)}(x_0)}\ e^\al\ |\pi_T(d\vec{H})|_{h(t)}\ \lf|\pi_{\vec{n}}\frac{\p\vec{\Phi}}{\p t}\rg|\ |\vec{\frak{h}}_0|_{h(t)}\ dvol_{h(t)}\ dt\\[5mm]
\ds +\int_0^T\int_{B_{R}^{h(t)}(x_0)}\ e^\al\ |\pi_T(d\vec{H})|_{h(t)}\ \lf|\ti{P}_{h(t)}\lf(\pi_{\vec{n}}\frac{\p\vec{\Phi}}{\p t}\cdot\vec{\frak{h}}_0\rg)\rg|\ dvol_{h(t)}\ dt
\end{array}
\ee
We have in local conformal coordinates
\[
\begin{array}{l}
\ds\pi_T(\p \vec{H})=2\,\p\vec{H}\cdot\p_z\vec{\Phi}\ \p_{\ov{z}}\vec{\Phi} \, e^{-2\la}+2\,\p\vec{H}\cdot\p_{\ov{z}}\vec{\Phi}\ \p_{{z}}\vec{\Phi} \, e^{-2\la}\\[5mm]
\ds=-2 \, \vec{H}\cdot\pi_{\vec{n}}(\p_z\p_z\vec{\Phi})\ \p_{\ov{z}}\vec{\Phi} \, e^{-2\la}\ dz-2\,\vec{H}\cdot\p_z\p_{\ov{z}}\vec{\Phi}\ \p_{{z}}\vec{\Phi} \, e^{-2\la}\ dz
\end{array}
\]
Hence
\[
|\pi_T(\p \vec{H})|_{h(t)}\le e^\al\ |\vec{H}|\ |\vec{\frak h}^0|_{h(t)}+e^\al\ |\vec{H}|^2 
\]
This gives
\be
\label{n-ene-con-54-bis}
\begin{array}{l}
\ds\int_0^T\lf|\int_\Sigma\chi^{h(t)}_{R,x_0}(x)\ \pi_T\lf(\ov{\p}\vec{H}\rg)\dot{\wedge}\lf[\ov{\p}\vec{\Phi}\res(I-\ti{P}_h)\lf(\pi_{\vec{n}}\frac{\p\vec{\Phi}}{\p t}\cdot\vec{\frak{h}}_0\rg)\rg]\rg|\ dt\\[5mm]
\ds\le \int_0^T\int_{B_{R}^{h(t)}(x_0)}\ e^{2\al}\ \lf[ |\vec{H}|\ |\vec{\frak h}^0|^2_{h(t)}+ \ |\vec{H}|^2\ |\vec{\frak{h}}_0|_{h(t)}\rg]\  \lf|\pi_{\vec{n}}\frac{\p\vec{\Phi}}{\p t}\rg|\ dvol_{h(t)}\ dt\\[5mm]
\ds +\int_0^T\int_{B_{R}^{h(t)}(x_0)}\  \ e^{2\al}\ \lf[ |\vec{H}|\ |\vec{\frak h}^0|_{h(t)}+ \ |\vec{H}|^2 \rg]\  \lf|\ti{P}_{h(t)}\lf(\pi_{\vec{n}}\frac{\p\vec{\Phi}}{\p t}\cdot\vec{\frak{h}}_0\rg)\rg|\ dvol_{h(t)}\ dt
\end{array}
\ee
First we have
\be
\label{af-1}
\begin{array}{l}
\ds\int_{B_{R}^{h(t)}(x_0)}\ e^{2\al}\ \lf[ |\vec{H}|\ |\vec{\frak h}^0|^2_{h(t)}+ \ |\vec{H}|^2\ |\vec{\frak{h}}_0|_{h(t)}\rg]\  \lf|\pi_{\vec{n}}\frac{\p\vec{\Phi}}{\p t}\rg|\ dvol_{h(t)}\\[5mm]
\ds \le e^{2\al_{R,x_0}(T)}\ \sqrt{ \int_{B_{R}^{h(t)}(x_0)}\ \lf|\pi_{\vec{n}}\frac{\p\vec{\Phi}}{\p t}\rg|^2\ dvol_{h(t)}} \\[5mm] 
\ds\quad\quad\times\sqrt{ \int_{B_{R}^{h(t)}(x_0)}\ \lf[ |\vec{H}|^2\ |\vec{\frak h}^0|^4_{h(t)}+ \ |\vec{H}|^4\ |\vec{\frak{h}}_0|^2_{h(t)}\rg]\  \ dvol_{h(t)}}
\end{array}
\ee
Using 
 \be
\label{af-n-h-0-H-L6}
\begin{array}{l}
\ds\lf[\int_\Sigma\chi^{h(t)}_{R,x_0}(x)\ |\vec{\frak{h}}^{\,0}|^6_{h(t)}\  \ dvol_{h(t)}\rg]^{1/2}\le  e^{C\ (1+l_{\ast\ast}^{-2}(T))\,E(0)}\ \lf[\int_{B_{3\,R/2}^{h(t)}(x_0)}\ |\vec{H}|^6\ dvol_{h(t)}\rg]^{1/2}\\[5mm]
\ds +e^{C\ (1+l_{\ast\ast}^{-2}(T))\,E(0)}\  R^{-2}\,\lf[\int_{B_{2\,R}^{h(t)}(x_0)}\ |\vec{\frak{h}^{\,0}}|^2_{h(t)}\ dvol_{h(t)}\rg]^{3/2}
\end{array}
\ee
combined with 
\be
\label{af-n-ene-con-31-bis}
\begin{array}{l}
\ds\lf[\int_{B_R^{h(t)}(x_0)} |\vec{H}|^6\  \ dvol_{h(t)}\rg]^{1/2}\\[5mm]
\ds\le C\ \int_{B_{2\,R}^{h(t)}(x_0)}\ |\vec{H}|^2\  \ dvol_{h(t)}\ \ \lf[\int_{B_{2\,R}^{h(t)}(x_0)}\ |\Delta_{h(t)}\vec{H}|^2\  \ dvol_{h(t)}\rg]^{1/2}\\[5mm]
\ds+ C\ \lf[R^{-4/3}\,\int_{B_{2\,R}^{h(t)}(x_0)}\ |\vec{H}|^2\  \ dvol_{h(t)}\rg]^{3/2}\ 
\end{array}
\ee
gives
\be
\label{af-2}
\begin{array}{l}
\ds C(\Sigma,l_{\ast\ast}(T),E(0))^{-1}\ \sqrt{ \int_{B_{R}^{h(t)}(x_0)}\ \lf[ |\vec{H}|^2\ |\vec{\frak h}^0|^4_{h(t)}+ \ |\vec{H}|^4\ |\vec{\frak{h}}_0|^2_{h(t)}\rg]\  \ dvol_{h(t)}}\\[5mm]
\ds\le  \int_{B_{2\,R}^{h(t)}(x_0)}\ |\vec{H}|^2\  \ dvol_{h(t)}\ \ \lf[\int_{B_{2\,R}^{h(t)}(x_0)}\ |\Delta_{h(t)}\vec{H}|^2\  \ dvol_{h(t)}\rg]^{1/2}\\[5mm]
\ds+ \lf[R^{-4/3}\,\int_{B_{2\,R}^{h(t)}(x_0)}\ |\vec{\mathbb I}_{\vec{\Phi}(t)}|^2_{g_{\vec{\Phi}(t)}}\  \ dvol_{h(t)}\rg]^{3/2}\ 
\end{array}
\ee
Hence
\be
\label{af-3}
\begin{array}{l}
\ds C(\Sigma,l_{\ast\ast}(T),E(0))^{-1}\ \int_0^T\int_{B_{R}^{h(t)}(x_0)}\ e^{2\al}\ \lf[ |\vec{H}|\ |\vec{\frak h}^0|^2_{h(t)}+ \ |\vec{H}|^2\ |\vec{\frak{h}}_0|_{h(t)}\rg]\  \lf|\pi_{\vec{n}}\frac{\p\vec{\Phi}}{\p t}\rg|\ dvol_{h(t)}\ dt\\[5mm]
\ds\le E(2R,x_0,T)\ \int_0^T \sqrt{ \int_{B_{R}^{h(t)}(x_0)}\ \lf|\pi_{\vec{n}}\frac{\p\vec{\Phi}}{\p t}\rg|^2\ dvol_{h(t)}}\ \sqrt{\int_{B_{2\,R}^{h(t)}(x_0)}\ |\Delta_{h(t)}\vec{H}|^2\  \ dvol_{h(t)}}\ dt\\[5mm]
\ds+ \int_0^T \sqrt{ \int_{B_{R}^{h(t)}(x_0)}\ \lf|\pi_{\vec{n}}\frac{\p\vec{\Phi}}{\p t}\rg|^2\ dvol_{g_{\vec{\Phi}(t)}}}\ \frac{e^{-2\al_{R,x_0}(T)}}{R^2}\ E(2R,x_0,T)^{3/2}\ .
\end{array}
\ee
This implies
\be
\label{af-4}
\begin{array}{l}
\ds  C(\Sigma,l_{\ast\ast}(T),E(0))^{-1}\ \int_0^T\int_{B_{R}^{h(t)}(x_0)}\ e^{2\al}\ \lf[ |\vec{H}|\ |\vec{\frak h}^0|^2_{h(t)}+ \ |\vec{H}|^2\ |\vec{\frak{h}}_0|_{h(t)}\rg]\  \lf|\pi_{\vec{n}}\frac{\p\vec{\Phi}}{\p t}\rg|\ dvol_{h(t)}\ dt\\[5mm]
\ds\le \sqrt{E(2R,x_0,t)}\ \int_{B_{2\,R}^{h(t)}(x_0)}\ e^{-2\al_{R,x_0}(t)}\ |\Delta_{h(t)}\vec{H}|^2\  \ dvol_{h(t)}\ dt\\[5mm]
\ds+{E(2R,x_0,t)}^{3/2}\ \sqrt{\int_0^T\int_\Sigma \lf|\pi_{\vec{n}}\frac{d\vec{\Phi}}{dt}\rg|^2\  dvol_{g_{\vec{\Phi}(t)}}\ dt}\ [\sqrt{T}\ e^{2\|\al\|_\infty}\ R^{-2}+1]\ .
\end{array}
\ee 
Now we bound
\be
\label{af-5}
\begin{array}{l}
\ds\int_0^T\int_{B_{R}^{h(t)}(x_0)}\  \ e^{2\al}\ \lf[ |\vec{H}|\ |\vec{\frak h}^0|_{h(t)}+ \ |\vec{H}|^2 \rg]\  \lf|\ti{P}_{h(t)}\lf(\pi_{\vec{n}}\frac{\p\vec{\Phi}}{\p t}\cdot\vec{\frak{h}}_0\rg)\rg|\ dvol_{h(t)}\ dt\\[5mm]
\ds  \le \lf\|\ti{P}_{h(t)}\lf(\pi_{\vec{n}}\frac{\p\vec{\Phi}}{\p t}\cdot\vec{\frak{h}}_0\rg)\rg\|_{L^\infty(\Sigma)}\ E(2R,x_0,T)
\end{array}
\ee
and we use
\be
\label{n-ene-con-54-8}
\begin{array}{l}
\ds
\lf|\ti{P}_{h(t)}\lf(\frac{d\vec{\Phi}}{dt}\cdot\vec{\frak h}^{0}\rg)\rg|_{h(t)}\ \le \sum_{i=1}^Q\|q_i(t)\|^2_{L^\infty(\Sigma)} {\int_\Sigma\lf|\frac{d\vec{\Phi}}{dt}\cdot\vec{\frak h}^{0}\rg|_{h(t)}\ dvol_{h(t)}}\\[5mm]
\ds\le C_\Sigma\  l_{\ast\ast}(T)^{-4}\ e^{2\|\al\|_\infty}\ \sqrt{\int_\Sigma \lf|\pi_{\vec{n}}\frac{d\vec{\Phi}}{dt}\rg|^2\  dvol_{g_{\vec{\Phi}(t)}}}\ \sqrt{ E(0)}\ .
\end{array}
\ee
to conclude finally
\be
\label{af-6}
\begin{array}{l}
\ds\int_0^T\int_{B_{R}^{h(t)}(x_0)}\  \ e^{2\al}\ \lf[ |\vec{H}|\ |\vec{\frak h}^0|_{h(t)}+ \ |\vec{H}|^2 \rg]\  \lf|\ti{P}_{h(t)}\lf(\pi_{\vec{n}}\frac{\p\vec{\Phi}}{\p t}\cdot\vec{\frak{h}}_0\rg)\rg|\ dvol_{h(t)}\ dt\\[5mm]
\ds\le  C_\Sigma\  l_{\ast\ast}(T)^{-4}\ e^{2\|\al\|_\infty}\ \int_0^T\sqrt{\int_\Sigma \lf|\pi_{\vec{n}}\frac{d\vec{\Phi}}{dt}\rg|^2\  dvol_{g_{\vec{\Phi}(t)}}}\ dt\ \sqrt{ E(0)}\ E(2R,x_0,T)\\[5mm]
\ds\le  C_\Sigma\  l_{\ast\ast}(T)^{-4}\ e^{2\|\al\|_\infty}\ \sqrt{T}\ \sqrt{\int_0^T\int_\Sigma \lf|\pi_{\vec{n}}\frac{d\vec{\Phi}}{dt}\rg|^2\  dvol_{g_{\vec{\Phi}(t)}}}\ dt\ \sqrt{ E(0)}\ E(2R,x_0,T)
\end{array}
\ee
Combining (\ref{n-ene-con-54-bis}), (\ref{af-4}) and (\ref{af-6}) is implying
\be
\label{af-7}
\begin{array}{l}
\ds C(\Sigma,l_{\ast\ast}(T),E(0))^{-1}\ \int_0^T\lf|\int_\Sigma\chi^{h(t)}_{R,x_0}(x)\ \pi_T\lf(\ov{\p}\vec{H}\rg)\dot{\wedge}\lf[\ov{\p}\vec{\Phi}\res(I-\ti{P}_h)\lf(\pi_{\vec{n}}\frac{\p\vec{\Phi}}{\p t}\cdot\vec{\frak{h}}_0\rg)\rg]\rg|\ dt\\[5mm]
\ds\le \sqrt{E(2R,x_0,t)}\ \int_{B_{2\,R}^{h(t)}(x_0)}\ e^{-2\al_{R,x_0}(t)}\ |\Delta_{h(t)}\vec{H}|^2\  \ dvol_{h(t)}\ dt\\[5mm]
\ds+E(2R,x_0,t)^{3/2}\ \sqrt{\int_0^T\int_\Sigma \lf|\pi_{\vec{n}}\frac{d\vec{\Phi}}{dt}\rg|^2\  dvol_{g_{\vec{\Phi}(t)}}\ dt}\ [\sqrt{T}\ e^{2\|\al\|_\infty}\ R^{-2}+1]\\[5mm]
\ds +e^{2\|\al\|_\infty}\ \sqrt{T}\ \sqrt{\int_0^T\int_\Sigma \lf|\pi_{\vec{n}}\frac{d\vec{\Phi}}{dt}\rg|^2\  dvol_{g_{\vec{\Phi}(t)}}}\ dt\ \sqrt{ E(0)}\ E(2R,x_0,T)
\end{array}
\ee
Finally the last term of (\ref{n-ene-con-7}) is controlled as follows
\be
\label{n-ene-con-55}
\begin{array}{l}
\ds\int_0^T\lf|\int_\Sigma\vec{H}\cdot\ov{\p} \lf[\chi^{h(t)}_{R,x_0}(x)\ \vec{H}\rg]\wedge\lf(U^{0,1}\res g_{\vec{\Phi}}\rg)\rg|\ dt\\[5mm]
\ds\le R^{-1}\int_0^T  \int_{B_{R}^{h(t)}(x_0)}\ |\vec{H}|^2\ |U^{0,1}|_{h(t)}\ dvol_{g_{\vec{\Phi}}}\ dt\\[5mm]
\ds+\int_0^T  \int_{B_{R}^{h(t)}(x_0)}\ |\vec{H}|\ |d\vec{H}|_{h(t)}\ |U^{0,1}|_{h(t)}\ dvol_{g_{\vec{\Phi}}}\ dt
\end{array}
\ee
These two terms in the r.h.s of the inequality (\ref{n-ene-con-55}) have already been considered previously (in (\ref{n-7}) and (\ref{n-18})). Integrating between 0 and $T$ (\ref{n-ene-con-7}) and compiling all the above estimates give the desired inequality (\ref{conc-evol}) and lemma~\ref{lm-conc-evol} is proved.\hfill $\Box$

\medskip

\begin{Lm}
\label{lm-gauss-evol}
Let $T>0$ and $0<R<l_{\ast\ast}(T)/2$ and $x_0$ such that
\[
E(2R,x_0,T):=\sup_{t\in[0,T]}\int_{B^{h(t)}_{2\,R}(x_0)}|\vec{\mathbb I}_{\vec{\Phi}(t)}|^2_{g_{\vec{\Phi}(t)}}\ dvol_{g_{\vec{\Phi}(t)}}<\frac{8\pi}{3}\ .
\]
Then
\be
\label{gauss-evol}
\begin{array}{l}
\ds\lf|\int_\Sigma\chi^{h(T)}_{R,x_0}(x)\ K_{g_{\vec{\Phi}(T)}}\ dvol_{g_{\vec{\Phi}(T)}}- \int_\Sigma\chi^{h(0)}_{R,x_0}(x)\ K_{g_{\vec{\Phi}(0)}}\ dvol_{g_{\vec{\Phi}(0)}}\rg|\\[5mm]
\ds \le C(\Sigma,l_{\ast\ast}(T))\,\frac{e^{2\|\al\|_\infty}}{R^2}\, [1+E(0)]\, \sqrt{T}\,\sqrt{E(0)}\ \sqrt{W(\vec{\Phi}(0))-W(\vec{\Phi}(T))}
\end{array}
\ee
\end{Lm}
\noindent{\bf Proof of lemma~\ref{lm-gauss-evol}.} We compute

\be
\label{s-1}
\begin{array}{l}
\ds\frac{d}{dt}\int_\Sigma\chi^{h(t)}_{R,x_0}(x)\ K_{g_{\vec{\Phi}}}\ dvol_{g_{\vec{\Phi}}}=-\frac{d}{dt}\int_\Sigma\chi^{h(t)}_{R,x_0}(x)\ \lf[\Delta_{h(t)}\al-K_{h(t)}\rg]\ dvol_{h(t)}\\[5mm]
\ds =\frac{1}{R}\int_\Sigma\ \frac{d|x-x_0|_{h(t)}}{dt}\ \chi'\lf( \frac{|x-x_0|_{h(t)}}{R} \rg)\ K_{h(t)}\ dvol_{h(t)}\\[5mm]
\ds\ +\frac{d}{dt}\int_\Sigma\ \frac{1}{R}\,\chi'\lf( \frac{|x-x_0|_{h(t)}}{R} \rg)\ \lf<d|x-x_0|_{h(t)}, d\al(t)\rg>_{h(t)}\ dvol_{h(t)}
\end{array}
\ee
Recall from (\ref{der-c-1})
\[
\forall \ x\in B_{R}^{h(t)}(x_0)\quad\quad\lf|\frac{d|x-x_0|_{h(t)}}{dt}\rg|\le\ R\ \lf\|\lf|\frac{dh}{dt}\rg|_{h(t)}\rg\|_{L^\infty(B_{R}^{h(t)}(x_0))}\ ,
\]
and from  (\ref{der-conf-23-d})
\be
\label{s-2}
\begin{array}{l}
\ds\int_0^T\lf\|\lf|\frac{dh}{dt}\rg|_{h(t)}\rg\|_{L^\infty(B_{R}^{h(t)}(x_0))}\ dt\\[5mm]
\ds C_\Sigma\ \frac{e^{2\,\|\al\|_\infty}}{l^4_{\ast\ast}(T)}\ \sqrt{T}\ \lf[\int_{\Sigma}|\vec{\mathbb I}_{\vec{\Phi}(0)}|^2\ dvol_{g_{\vec{\Phi}(0)}}\rg]^{1/2}\,\lf(W(\vec{\Phi}(0))-W(\vec{\Phi}(t))\rg)^{1/2}
\end{array}
\ee
Hence
\be
\label{s-3}
\begin{array}{l}
\ds\int_0^T\lf|\frac{1}{R}\int_\Sigma\ \frac{d|x-x_0|_{h(t)}}{dt}\ \chi'\lf( \frac{|x-x_0|_{h(t)}}{R} \rg)\ K_{h(t)}\ dvol_{h(t)}\rg|\ dt\\[5mm]
\ds\le C(\Sigma,l_{\ast\ast}(T))\ e^{2\,\|\al\|_\infty}\ R^2\ \sqrt{T}\ \sqrt{E(0)}\ \sqrt{W(\vec{\Phi}(0))-W(\vec{\Phi}(t))}
\end{array}
\ee
We have now
\be
\label{s-4}
\begin{array}{l}
\ds\frac{d}{dt}\int_\Sigma\ \frac{1}{R}\,\chi'\lf( \frac{|x-x_0|_{h(t)}}{R} \rg)\ \lf<d|x-x_0|_{h(t)}, d\al(t)\rg>_{h(t)}\ dvol_{h(t)}\\[5mm]
\ds= \int_\Sigma\ \frac{1}{R^2}\,\chi''\lf( \frac{|x-x_0|_{h(t)}}{R} \rg)\ \frac{d|x-x_0|_{h(t)}}{dt}\ \lf<d|x-x_0|_{h(t)}, d\al(t)\rg>_{h(t)}\ dvol_{h(t)}\\[5mm]
\ds+\int_\Sigma\ \frac{1}{R}\,\chi'\lf( \frac{|x-x_0|_{h(t)}}{R} \rg)\ \frac{dh}{dt}\res_2\ \lf[d|x-x_0|_{h(t)}\otimes d\al(t)\rg]\ dvol_{h(t)}\\[5mm]
\ds-\int_\Sigma\ \frac{1}{R}\,\chi'\lf( \frac{|x-x_0|_{h(t)}}{R} \rg)\ \frac{d|x-x_0|_{h(t)}}{dt}\ \Delta_{h(t)}\al(t)\ dvol_{h(t)}\\[5mm]
\ds-\int_\Sigma\ \frac{1}{R^2}\,\chi''\lf( \frac{|x-x_0|_{h(t)}}{R} \rg)\ \lf<d|x-x_0|_{h(t)}, d\al(t)\rg>_{h(t)}\ dvol_{h(t)}\\[5mm]
\ds +\int_\Sigma\ \frac{1}{R}\,d^{\ast_{h(t)}}\lf[ \chi'\lf( \frac{|x-x_0|_{h(t)}}{R} \rg)\ d|x-x_0|_{h(t)}\rg]\ \frac{d\al}{dt}\ dvol_{h(t)}
\end{array}
\ee
We bound
\be
\label{s-5}
\begin{array}{l}
\ds\int_0^T\lf| \int_\Sigma\ \frac{1}{R}\,\chi'\lf( \frac{|x-x_0|_{h(t)}}{R} \rg)\ \frac{dh}{dt}\res_2\ \lf[d|x-x_0|_{h(t)}\otimes d\al(t)\rg]\ dvol_{h(t)}  \rg|\ dt\\[5mm]
\ds\le\int_0^T\lf\|\lf|\frac{dh}{dt}\rg|_{h(t)}\rg\|_{L^\infty(B_{R}^{h(t)}(x_0))}\, \lf[\int_{B_{R}^{h(t)}(x_0))}|d\al|^2_{h(t)}\ dvol_{h(t)}\rg]^{1/2}\ dt
\end{array}
\ee
Recall 
\be
\label{r-L2-2}
\|d\al\|_{L^{2}(B^h_{R}(x_0))}\le \frac{C_{\Sigma,m}}{ l_\ast^{2}(t)}\ \int_{\Sigma}|\vec{\mathbb I}|^2_{g_{\vec{\Phi}}}\ dvol_{g_{\vec{\Phi}}}
\ee
Hence
\be
\label{s-6}
\begin{array}{l}
\ds\int_0^T\lf| \int_\Sigma\ \frac{1}{R}\,\chi'\lf( \frac{|x-x_0|_{h(t)}}{R} \rg)\ \frac{dh}{dt}\res_2\ \lf[d|x-x_0|_{h(t)}\otimes d\al(t)\rg]\ dvol_{h(t)}  \rg|\ dt\\[5mm]
\ds\le\frac{C_{\Sigma,m}}{ l^2_{\ast\ast}(T)}\ {E(0)}\ \int_0^T\lf\|\lf|\frac{dh}{dt}\rg|_{h(t)}\rg\|_{L^\infty(B_{R}^{h(t)}(x_0))}\, dt\\[5mm]
\ds\le  C(\Sigma,l_{\ast\ast}(T))\ e^{2\,\|\al\|_\infty}\ \sqrt{T}\ {E(0)}^{3/2}\ \sqrt{W(\vec{\Phi}(0))-W(\vec{\Phi}(t))}
\end{array}
\ee
Next we bound
\be
\label{s-7}
\begin{array}{l}
\ds\int_0^T\lf| \int_\Sigma\ \frac{1}{R}\,\chi'\lf( \frac{|x-x_0|_{h(t)}}{R} \rg)\ \frac{d|x-x_0|_{h(t)}}{dt}\ \Delta_{h(t)}\al(t)\ dvol_{h(t)}\rg|\ dt\\[5mm]
\ds\le \int_0^T\lf\|\lf|\frac{dh}{dt}\rg|_{h(t)}\rg\|_{L^\infty(B_{R}^{h(t)}(x_0))}\, \int_{B_{R}^{h(t)}(x_0))} \lf[e^{2\al}\,K_{g_{\vec{\Phi}(t)}}-K_{h(t)}\rg]\ dvol_{h(t)}\ dt\\[5mm]
\ds\le C(\Sigma,l_{\ast\ast}(T))\ e^{2\,\|\al\|_\infty}\ [E(R,x_0,T)+R^2]\ \sqrt{T}\ \sqrt{E(0)}\ \sqrt{W(\vec{\Phi}(0))-W(\vec{\Phi}(t))}
\end{array}
\ee
We recall from lemma~\ref{lm-param-conf-var} 
\be
\label{re-param-conf-var}
\frac{d\al}{dt}=\,\frac{d\vec{\Phi}}{dt}\cdot \vec{H}+e^{-\,2\al(t)}\frac{1}{2}\,d^{\ast_{h(t)}}\lf[\frac{d\vec{\Phi}}{dt}\cdot d\vec{\Phi}(t)\rg]\ .
\ee
Inserting it in the last line of (\ref{s-4}) is giving
\be
\label{s-8}
\begin{array}{l}
\ds\int_0^T\lf| \int_\Sigma\ \frac{1}{R}\,d^{\ast_{h(t)}}\lf[ \chi'\lf( \frac{|x-x_0|_{h(t)}}{R} \rg)\ d|x-x_0|_{h(t)}\rg]\ \frac{d\al}{dt}\ dvol_{h(t)}\rg|\ dt\\[5mm]
\ds\le\frac{C}{R^2}\,\int_0^T \int_{B_{R}^{h(t)}(x_0)} \lf|\pi_{\vec{n}}\frac{d\vec{\Phi}}{dt}\rg|\ \lf|\vec{H}\rg| \ dvol_{h(t)}\ dt\\[5mm]
\ds+\frac{C}{R^2}\,\int_0^T \int_{B_{R}^{h(t)}(x_0)} \lf[|d\al|_{h(t)}+R^{-1}\rg]\ e^{-\al}\ |\vec{U}|\ dvol_{h(t)}\ dt
\end{array}
\ee
Recall $|\vec{U}|=\sqrt{2}\,e^\al\, |U^{0,1}|_{h(t)}$. Using (\ref{L2-2}) , lemma~\ref{lm-delbar-inv} and lemma~\ref{lm-L2-infty} we deduce
\be
\label{s-9}
\begin{array}{l}
\ds\int_0^T\lf| \int_\Sigma\ \frac{1}{R}\,d^{\ast_{h(t)}}\lf[ \chi'\lf( \frac{|x-x_0|_{h(t)}}{R} \rg)\ d|x-x_0|_{h(t)}\rg]\ \frac{d\al}{dt}\ dvol_{h(t)}\rg|\ dt\\[5mm]
\ds \le C\ \frac{e^{2\|\al\|_\infty}}{R^2}\ \sqrt{E(2R,x_0,T)}\,\sqrt{T}\ \sqrt{\int_0^T \int_{B_{R}^{h(t)}(x_0)} \lf|\pi_{\vec{n}}\frac{d\vec{\Phi}}{dt}\rg|^2\ \ dvol_{g_{\vec{\Phi}}}\ dt}\\[5mm]
\ds+\, \frac{C}{R^2}\,\ \int_0^T\ \|U^{0,1}\|_{L^{2,\infty}(B_{R}^{h(t)}(x_0))}\ \|d\al\|_{L^{2,1}_{h(t)}(B_{R}^{h(t)}(x_0))}\ dt+\frac{C}{R^2}\,\ \int_0^T\ \|U^{0,1}\|_{L^{2,\infty}(B_{R}^{h(t)}(x_0))}\ dt\\[5mm]
\ds \le C\ \frac{e^{2\|\al\|_\infty}}{R^2}\ \sqrt{E(2R,x_0,T)}\,\sqrt{T}\ \sqrt{W(\vec{\Phi}(0))-W(\vec{\Phi}(T))}\\[5mm]
\ds+ \frac{C}{R^2}\, [1+E(0)]\, \int_0^T \|\p^hU^{0,1}\|_{L^{1}(\Sigma)}\ dt
\end{array}
\ee
Using (\ref{n-13}) we have
\be
\label{s-10}
\begin{array}{l}
\ds\int_0^T\lf| \int_\Sigma\ \frac{1}{R}\,d^{\ast_{h(t)}}\lf[ \chi'\lf( \frac{|x-x_0|_{h(t)}}{R} \rg)\ d|x-x_0|_{h(t)}\rg]\ \frac{d\al}{dt}\ dvol_{h(t)}\rg|\ dt\\[5mm]
\le C\ \frac{e^{2\|\al\|_\infty}}{R^2}\ \sqrt{E(2R,x_0,T)}\,\sqrt{T}\ \sqrt{W(\vec{\Phi}(0))-W(\vec{\Phi}(T))}\\[5mm]
\le+C(\Sigma,l_{\ast\ast}(T))\,\frac{e^{2\|\al\|_\infty}}{R^2}\, [1+E(0)]\, \sqrt{T}\,\sqrt{E(0)}\ \sqrt{W(\vec{\Phi}(0))-W(\vec{\Phi}(T))}
\end{array}
\ee
This is ending the proof of lemma~\ref{lm-gauss-evol}.\hfill $\Box$

\medskip

Combining  lemma~\ref{lm-conc-evol} and lemma~\ref{lm-gauss-evol} and using the fact that $$|\vec{\mathbb I}_{\vec{\Phi}(t)}|^2_{g_{\vec{\Phi}}}=4|\vec{H}_{\vec{\Phi}(t)}|^2-2K_{\vec{\Phi}(t)}$$ we obtain the following lemma.

\begin{Lm}
\label{lm-energy-evol}
Let $T>0$ and $0<R<l_{\ast\ast}(T)/2$ such that
\[
E(2R,x_0,t):=\sup_{t\in[0,T]}\int_{B^{h(t)}_{2\,R}(x_0)}|\vec{\mathbb I}_{\vec{\Phi}(t)}|^2_{g_{\vec{\Phi}(t)}}\ dvol_{g_{\vec{\Phi}(t)}}<\frac{8\pi}{3}
\]
Then for any $\delta\in (0,1)$ we have
\be
\label{conc-evol-a}
\begin{array}{l}
\ds\int_{B_{R}^{h(T)}(x_0)}|\vec{\mathbb I}_{\vec{\Phi}(T)}|^2_{g_{\vec{\Phi}}}\ \chi^{h(T)}_{R,x_0}(x)\ \ dvol_{g_{\vec{\Phi}(T)}} +\int_0^T\int_{B_{R/2}^{h(t)}(x_0)}\ |\Delta_{g_{\vec{\Phi}(t)}}\vec{H}(t)|^2\  \ dvol_{g_{\vec{\Phi}(T)}}\ dt\\[5mm]
\ds\le\int_{B_{R}^{h(0)}(x_0)}|\vec{\mathbb I}_{\vec{\Phi}(0)}|^2_{g_{\vec{\Phi}}}\ \chi^{h(0)}_{R,x_0}(x)\ \ dvol_{g_{\vec{\Phi}(0)}}  + [\delta+E(2R,x_0,T)^{1/2}]\  \int_0^T\int_{B_{2\,R}^{h(t)}(x_0)}\ |\Delta_{h(t)}\vec{H}|^2\  \ dvol_{h(t)}\ dt\\[5mm]
\ds\quad+C(\Sigma,l_{\ast\ast}(T),E(0))\ \lf[\frac{e^{4\|\al\|_\infty}}{\delta^3\,R^4}\ T+\sqrt{\frac{e^{4\|\al\|_\infty}}{2\,\delta^2\,R^4}\ T}+\frac{e^{4\|\al\|_\infty}}{\delta^3}\ \lf[ \frac{T}{R^4} \rg]^{1/3}\rg]\ E(2R,x_0,T)\\[5mm]
 \ds\quad+C(\Sigma,l_{\ast\ast}(T),E(0))\ \lf[\sqrt{E(2R,x_0,T)}\ \sqrt{W(\vec{\Phi}(0)-W(\vec{\Phi}(T)}\rg]\ E(2R,x_0,T)\\[5mm]
\ds\quad+C(\Sigma,l_{\ast\ast}(T),E(0))\ \lf[ e^{2\|\al\|_\infty}\ \sqrt{T}\ \sqrt{W(\vec{\Phi}(0)-W(\vec{\Phi}(T)} \rg]\ E(2R,x_0,T)\\[5mm]
\ds\quad+ C(\Sigma,l_{\ast\ast}(T),E(0))\ \lf[\frac{e^{8\|\al\|_\infty}}{\delta}\ \lf[W(\vec{\Phi}(0)-W(\vec{\Phi}(T)\rg] \rg]\ E(2R,x_0,T)\\[5mm]
\ds\quad+C(\Sigma,l_{\ast\ast}(T), E(0))\,\frac{e^{2\|\al\|_\infty}}{R^2}\,  \sqrt{T}\,\ \sqrt{W(\vec{\Phi}(0))-W(\vec{\Phi}(T))}
\end{array}
\ee

\hfill $\Box$
\end{Lm}

 %%%%%%%%%%%%%%%%%%%%%%%%%%%%%%%%%%%%%%%%%%%%%%%%%%%%%%%%%%%%%%%%%%%%%%%%%%%%%%%%%%%%%%%%%%%%%%%%%%%%%%%%%%%%%

\section{Proof of the main Theorem}
\reset
Let $\beta<8\pi/3$ to be fixed later on. Choose $l>0$ depending only on  $l_\ast(0)$ satisfying respectively (\ref{length-control-torus}) and (\ref{length-control}) depending on the topology of $\Sigma$. Choose $R>0$  and $T$ such that
\be
\label{V.1}
E(2R,T,0):=\sup_{x_0\in\Sigma}\sup_{t\in[0,T]}\ E(2R,x_0,t):=\sup_{t\in[0,T]}\sup_{x_0\in\Sigma}\int_{B^{h(t)}_{2\,R}(x_0)}|\vec{\mathbb I}_{\vec{\Phi}(t)}|^2_{g_{\vec{\Phi}(t)}}\ dvol_{g_{\vec{\Phi}(t)}}<\beta
\ee
and
\be
\label{V.2}
\frac{e^{4\|\al(0)\|_\infty}}{R^4} \ T< \hat{C}_{\Sigma,m}^2\, \frac{l^{26}}{E(0)^8}\ .
\ee
Then, from lemma~\ref{lm-length-control-torus} and lemma~\ref{lm-length-control} we have
\be
\label{V.3}
\sup_{t\in[0,T]}\lf\|\al(t)-\al(0)\rg\|_{L^\infty(\Sigma)}\le \hat{C}_{\Sigma,m}\ .
\ee
Using lemma~\ref{lm-energy-evol} we have for any $0<\delta<1$
\be
\label{V.4}
\begin{array}{l}
\ds\sup_{t\in[0,T]}\sup_{x_0\in\Sigma}\int_{B_{R}^{h(T)}(x_0)}|\vec{\mathbb I}_{\vec{\Phi}(T)}|^2_{g_{\vec{\Phi}}}\ \chi^{h(T)}_{R,x_0}(x)\ \ dvol_{g_{\vec{\Phi}(T)}} +\sup_{x_0\in\Sigma}\int_0^T\int_{B_{R/2}^{h(t)}(x_0)}\ |\Delta_{g_{\vec{\Phi}(t)}}\vec{H}(t)|^2\  \ dvol_{g_{\vec{\Phi}(T)}}\ dt\\[5mm]
\ds\le\sup_{x_0\in\Sigma}\int_{B_{R}^{h(0)}(x_0)}|\vec{\mathbb I}_{\vec{\Phi}(0)}|^2_{g_{\vec{\Phi}}}\ \chi^{h(0)}_{R,x_0}(x)\ \ dvol_{g_{\vec{\Phi}(0)}}  + [\delta+\beta^{1/2}]\  \int_0^T\int_{B_{2\,R}^{h(t)}(x_0)}\ |\Delta_{h(t)}\vec{H}|^2\  dvol_{g_{\vec{\Phi}(t)}} \ dt\\[5mm]
\ds\quad+C(\Sigma, l,E(0))\ \lf[\frac{e^{4\|\al(0)\|_\infty}}{\delta^3\,R^4}\ T+\sqrt{\frac{e^{4\|\al(0)\|_\infty}}{2\,\delta^2\,R^4}\ T}+\frac{e^{4\|\al(0)\|_\infty}}{\delta^3}\ \lf[ \frac{T}{R^4} \rg]^{1/3}\rg]\ \beta\\[5mm]
 \ds\quad+C(\Sigma,l,E(0))\ \lf[\sqrt{\beta}\ \sqrt{W(\vec{\Phi}(0))-W(\vec{\Phi}(T))}\rg]\ \beta\\[5mm]
\ds\quad+C(\Sigma,l,E(0))\ \lf[ e^{2\|\al(0)\|_\infty}\ \sqrt{T}\ \sqrt{W(\vec{\Phi}(0))-W(\vec{\Phi}(T))} \rg]\ \beta\\[5mm]
\ds\quad+ C(\Sigma,l,E(0))\ \lf[\frac{e^{8\|\al(0)\|_\infty}}{\delta}\ \lf[W(\vec{\Phi}(0))-W(\vec{\Phi}(T))\rg] \rg]\ \beta\\[5mm]
\ds\quad+C(\Sigma,l, E(0))\,\frac{e^{2\|\al(0)\|_\infty}}{R^2}\,  \sqrt{T}\,\ \sqrt{W(\vec{\Phi}(0))-W(\vec{\Phi}(T))}
\end{array}
\ee
Since $R<l<l_{\ast\ast}(T)$ we can cover each ball $B_{2\,R}^{h(t)}(x_0)$ of radius $2R$ by a uniformly bounded number of balls of radius $R/2$ by a universal number $N_1$. We choose $N_1\ge 8$ and fixed until the end of the proof from now on. Hence summing up (\ref{V.4}) gives
\be
\label{V.5}
\begin{array}{l}
\ds E(2R,T,0) +\sup_{x_0\in\Sigma}\int_0^T\int_{B_{2R}^{h(t)}(x_0)}\ |\Delta_{g_{\vec{\Phi}(t)}}\vec{H}(t)|^2\  \ dvol_{g_{\vec{\Phi}(T)}}\ dt\\[5mm]
\ds\le\, N_1\,\sup_{x_0\in\Sigma}\int_{B_{R}^{h(0)}(x_0)}|\vec{\mathbb I}_{\vec{\Phi}(0)}|^2_{g_{\vec{\Phi}}}\ \chi^{h(0)}_{R,x_0}(x)\ \ dvol_{g_{\vec{\Phi}(0)}}\\[5mm]
\ds \quad + N_1\, [\delta+\beta^{1/2}]\ \sup_{x_0\in\Sigma}\ \int_0^T\int_{B_{2\,R}^{h(t)}(x_0)}\ |\Delta_{h(t)}\vec{H}|^2\  dvol_{g_{\vec{\Phi}(t)}} \ dt\\[5mm]
\ds\quad+C(\Sigma, l,E(0))\ \lf[\frac{e^{4\|\al(0)\|_\infty}}{\delta^3\,R^4}\ T+\sqrt{\frac{e^{4\|\al(0)\|_\infty}}{2\,\delta^2\,R^4}\ T}+\frac{e^{4\|\al(0)\|_\infty}}{\delta^3}\ \lf[ \frac{T}{R^4} \rg]^{1/3}\rg]\ \beta\\[5mm]
 \ds\quad+C(\Sigma,l,E(0))\ \lf[\sqrt{\beta}\ \sqrt{W(\vec{\Phi}(0))-W(\vec{\Phi}(T))}\rg]\ \beta\\[5mm]
\ds\quad+C(\Sigma,l,E(0))\ \lf[ e^{2\|\al(0)\|_\infty}\ \sqrt{T}\ \sqrt{W(\vec{\Phi}(0))-W(\vec{\Phi}(T))} \rg]\ \beta\\[5mm]
\ds\quad+ C(\Sigma,l,E(0))\ \lf[\frac{e^{8\|\al(0)\|_\infty}}{\delta}\ \lf[W(\vec{\Phi}(0))-W(\vec{\Phi}(T))\rg] \rg]\ \beta\\[5mm]
\ds\quad+C(\Sigma,l, E(0))\,\frac{e^{2\|\al(0)\|_\infty}}{R^2}\,  \sqrt{T}\,\ \sqrt{W(\vec{\Phi}(0))-W(\vec{\Phi}(T))}
\end{array}
\ee
We first fix $\delta>0$, universal,  such that $ N_1\,\delta=1/4$. Now we choose $\beta=\min\{8\pi/3, 16^{-1}\,N_1^{-2}\}$. With this choice we obtain
 \be
\label{V.6}
\begin{array}{l}
\ds E(2R,T,0) +2^{-1}\,\sup_{x_0\in\Sigma}\int_0^T\int_{B_{2R}^{h(t)}(x_0)}\ |\Delta_{g_{\vec{\Phi}(t)}}\vec{H}(t)|^2\  \ dvol_{g_{\vec{\Phi}(T)}}\ dt\\[5mm]
\ds\le\, N_1\,\sup_{x_0\in\Sigma}\int_{B_{R}^{h(0)}(x_0)}|\vec{\mathbb I}_{\vec{\Phi}(0)}|^2_{g_{\vec{\Phi}}}\ \chi^{h(0)}_{R,x_0}(x)\ \ dvol_{g_{\vec{\Phi}(0)}}\\[5mm]
%\ds \quad + N_1\, [\delta+\beta^{1/2}]\ \sup_{x_0\in\Sigma}\ \int_0^T\int_{B_{2\,R}^{h(t)}(x_0)}\ |\Delta_{h(t)}\vec{H}|^2\  dvol_{g_{\vec{\Phi}(t)}} \ dt\\[5mm]
\ds\quad+C(\Sigma, l,E(0))\ \lf[\frac{e^{4\|\al(0)\|_\infty}}{\delta^3\,R^4}\ T+\sqrt{\frac{e^{4\|\al(0)\|_\infty}}{2\,\delta^2\,R^4}\ T}+\frac{e^{4\|\al(0)\|_\infty}}{\delta^3}\ \lf[ \frac{T}{R^4} \rg]^{1/3}\rg]\ \beta\\[5mm]
 \ds\quad+C(\Sigma,l,E(0))\ \lf[\sqrt{\beta}\ \sqrt{W(\vec{\Phi}(0))-W(\vec{\Phi}(T))}\rg]\ \beta\\[5mm]
\ds\quad+C(\Sigma,l,E(0))\ \lf[ e^{2\|\al(0)\|_\infty}\ \sqrt{T}\ \sqrt{W(\vec{\Phi}(0))-W(\vec{\Phi}(T))} \rg]\ \beta\\[5mm]
\ds\quad+ C(\Sigma,l,E(0))\ \lf[\frac{e^{8\|\al(0)\|_\infty}}{\delta}\ \lf[W(\vec{\Phi}(0))-W(\vec{\Phi}(T))\rg] \rg]\ \beta\\[5mm]
\ds\quad+C(\Sigma,l, E(0))\,\frac{e^{2\|\al(0)\|_\infty}}{R^2}\,  \sqrt{T}\,\ \sqrt{W(\vec{\Phi}(0))-W(\vec{\Phi}(T))}
\end{array}
\ee
The values $\delta$ and $\beta$ being fixed, we now choose $R$ such that
\be
\label{V.7}
\begin{array}{l}
 \ds N_1\,\sup_{x_0\in\Sigma}\int_{B_{R}^{h(0)}(x_0)}|\vec{\mathbb I}_{\vec{\Phi}(0)}|^2_{g_{\vec{\Phi}}}\ \chi^{h(0)}_{R,x_0}(x)\ \ dvol_{g_{\vec{\Phi}(0)}}:=8^{-1}\,\beta\ .
\end{array}
\ee
Finally, we consider the largest time $T_{\beta,R,\delta}>0$ such that respectively
\be
\label{V.7-b}
\frac{e^{4\|\al(0)\|_\infty}}{R^4} \ T_{\beta,R,\delta}\le \hat{C}_{\Sigma,m}^2\, \frac{l^{26}}{E(0)^8}\ .
\ee
as well as
\be
\label{V.8}
C(\Sigma, l,E(0))\ \lf[\frac{e^{4\|\al(0)\|_\infty}}{\delta^3\,R^4}\ T_{\beta,R,\delta}+\sqrt{\frac{e^{4\|\al(0)\|_\infty}}{2\,\delta^2\,R^4}\ T_{\beta,R,\delta}}+\frac{e^{4\|\al(0)\|_\infty}}{\delta^3}\ \lf[ \frac{T_{\beta,R,\delta}}{R^4} \rg]^{1/3}\rg]\le 1/8\ ,
\ee
we have also
\be
\label{V.9}
C(\Sigma,l,E(0))\ \lf[\sqrt{\beta}\ \sqrt{W(\vec{\Phi}(0))-W(\vec{\Phi}(T_{\beta,R,\delta})}\rg]\le1/8\ ,
\ee
as well as
\be
\label{V.10}
C(\Sigma,l,E(0))\ \lf[ e^{2\|\al(0)\|_\infty}\ \sqrt{T_{\beta,R,\delta}}\ \sqrt{W(\vec{\Phi}(0))-W(\vec{\Phi}(T_{\beta,R,\delta})} \rg]\le1/8\ ,
\ee
moreover
\be
\label{V.11}
C(\Sigma,l,E(0))\ \lf[\frac{e^{8\|\al(0)\|_\infty}}{\delta}\ \lf[W(\vec{\Phi}(0))-W(\vec{\Phi}(T_{\beta,R,\delta})\rg] \rg]\le1/8\ ,
\ee
and finally
\be
\label{V.12}
C(\Sigma,l, E(0))\,\frac{e^{2\|\al(0)\|_\infty}}{R^2}\,  \sqrt{T_{\beta,R,\delta}}\,\ \sqrt{W(\vec{\Phi}(0))-W(\vec{\Phi}(T_{\beta,R,\delta}))}\le\beta/8\ .
\ee
Under the assumptions (\ref{V.7})...(\ref{V.12}) , for any $T<T_{\beta,R,\delta}$ we obtain the a-priori estimate
 \be
\label{V.13}
\begin{array}{l}
\ds E(2R,T,0) +2^{-1}\,\sup_{x_0\in\Sigma}\int_0^T\int_{B_{2R}^{h(t)}(x_0)}\ |\Delta_{g_{\vec{\Phi}(t)}}\vec{H}(t)|^2\  \ dvol_{g_{\vec{\Phi}(T)}}\ dt\le\frac{3\beta}{4}\ .
\end{array}
\ee
From \cite{PalRiv} (Appendix A.2) we have the existence of $T>0$ such that the flow admits a smooth solution at least until $T$ for the given smooth initial data $\vec{\Phi}(0)$. Assume the flow has a maximal existence time
$T_{max}<T_{\beta,R,\delta}$. Thanks to the a-priori estimate (\ref{V.13}) we have 
\[
\sup_{T<T_{max}}E(2R,T,0)<\frac{3\beta}{4}\ .
\]
Thanks to (\ref{V.3}), $\al$ is uniformly bounded in $L^\infty([0,T_{max})\times \Sigma)$. Since $l_{\ast\ast}(T_{max})\ge l>0$ the constant Gauss curvature metric $h(t)$ is controlled uniformly in $[0,T_{max})$ in any norm. Hence the parametric Willmore flow equation is uniformly parabolic and extends to a smooth immersion $\vec{\Phi}(T_{max})$ at the time $T_{max}$ (see section 5 of \cite{PalRiv}). One can ``restart'' the flow from 
$\vec{\Phi}(T_{max})$ and we get a contradiction. Hence the flow exists at least up to $T:=T_{\beta,R,\delta}$. 

\medskip

Assume equality holds in either (\ref{V.7-b}), (\ref{V.8}), (\ref{V.10}) or (\ref{V.12}) then the theorem is proved. Now we assume that the inequality is strict in (\ref{V.7-b}), (\ref{V.8}), (\ref{V.10}) or (\ref{V.12}).
We then have equality in either (\ref{V.9}) or (\ref{V.11}). If equality holds in (\ref{V.9}) we take initially $\beta$ a bit smaller satisfying
\be
\label{V.99}
\beta<\frac{1}{8^2\, C^2(\Sigma,l,E(0))\ W(\vec{\Phi}(0))}\ .
\ee
and hence the possibility that equality holds in (\ref{V.9}) is excluded. Now we have to deal with the case where equality holds in (\ref{V.11}) and modify the argument as follows.

We keep the value $\beta$ satisfying
\be
\label{beta}
\beta:=\min\lf\{ \frac{8\pi}{3}, \frac{1}{16\,N_1^{2}},\frac{1}{2\, 8^2\, C(\Sigma,l,E(0))\ W(\vec{\Phi}(0))}\rg\}\ .
\ee
We keep the value of $l$ and the value of $\delta$. Let $0<\rho<l$ such that
\be
\label{rho}
\sup_{x_0\in\Sigma}\int_{B_{\rho}^{h(0)}(x_0)}|\vec{\mathbb I}_{\vec{\Phi}(0)}|^2_{g_{\vec{\Phi}}}\ \chi^{h(0)}_{R,x_0}(x)\ \ dvol_{g_{\vec{\Phi}(0)}}=N_1^{-k}\,\frac{\beta}{8}
\ee
where
\be
\label{k}
\lf[8\, C(\Sigma,l,E(0))\ \frac{e^{8\|\al(0)\|_\infty}}{\delta}\ W(\vec{\Phi}(0))\rg]+1= k
\ee
where $[8\, C(\Sigma,l,E(0))\ \frac{e^{8\|\al(0)\|_\infty}}{\delta}\ W(\vec{\Phi}(0))]$ is the integer part of $8\, C(\Sigma,l,E(0))\ \frac{e^{8\|\al(0)\|_\infty}}{\delta}\ W(\vec{\Phi}(0))$. We shall restrict exclusively to $T>0$ such that
\be
\label{V.14}
\frac{e^{4\|\al(0)\|_\infty}}{\rho^4} \ T< \hat{C}_{\Sigma,m}^2\, \frac{l^{26}}{E(0)^8}\ .
\ee
We have the a-priori bound
\be
\label{V.15}
\begin{array}{l}
\ds E(2\rho,T) +2^{-1}\,\sup_{x_0\in\Sigma}\int_0^T\int_{B_{2\rho}^{h(t)}(x_0)}\ |\Delta_{g_{\vec{\Phi}(t)}}\vec{H}(t)|^2\  \ dvol_{g_{\vec{\Phi}(T)}}\ dt\\[5mm]
\ds\le\, N_1\,\sup_{x_0\in\Sigma}\int_{B_{\rho}^{h(0)}(x_0)}|\vec{\mathbb I}_{\vec{\Phi}(0)}|^2_{g_{\vec{\Phi}}}\ \chi^{h(0)}_{\rho,x_0}(x)\ \ dvol_{g_{\vec{\Phi}(0)}}\\[5mm]
%\ds \quad + N_1\, [\delta+\beta^{1/2}]\ \sup_{x_0\in\Sigma}\ \int_0^T\int_{B_{2\,R}^{h(t)}(x_0)}\ |\Delta_{h(t)}\vec{H}|^2\  dvol_{g_{\vec{\Phi}(t)}} \ dt\\[5mm]
\ds\quad+C(\Sigma, l,E(0))\ \lf[\frac{e^{4\|\al(0)\|_\infty}}{\delta^3\,\rho^4}\ T+\sqrt{\frac{e^{4\|\al(0)\|_\infty}}{2\,\delta^2\,\rho^4}\ T}+\frac{e^{4\|\al(0)\|_\infty}}{\delta^3}\ \lf[ \frac{T}{\rho^4} \rg]^{1/3}\rg]\  N_1^{-k}\  \beta\\[5mm]
 \ds\quad+C(\Sigma,l,E(0))\ \lf[\sqrt{\beta}\ \sqrt{W(\vec{\Phi}(0))-W(\vec{\Phi}(T))}\rg]\  N_1^{-k}\ \beta\\[5mm]
\ds\quad+C(\Sigma,l,E(0))\ \lf[ e^{2\|\al(0)\|_\infty}\ \sqrt{T}\ \sqrt{W(\vec{\Phi}(0))-W(\vec{\Phi}(T))} \rg]\  N_1^{-k}\  \beta\\[5mm]
\ds\quad+ C(\Sigma,l,E(0))\ \lf[\frac{e^{8\|\al(0)\|_\infty}}{\delta}\ \lf[W(\vec{\Phi}(0))-W(\vec{\Phi}(T))\rg] \rg]\  N_1^{-k}\ \beta\\[5mm]
\ds\quad+C(\Sigma,l, E(0))\,\frac{e^{2\|\al(0)\|_\infty}}{\rho^2}\,  \sqrt{T}\,\ \sqrt{W(\vec{\Phi}(0))-W(\vec{\Phi}(T))}
\end{array}
\ee
As before we fix the maximal value we denote $T_1$ such that simultaneously 
\be
\label{V.14-a}
\frac{e^{4\|\al(0)\|_\infty}}{\rho^4} \ T_1< \hat{C}_{\Sigma,m}^2\, \frac{l^{26}}{E(0)^8}\ .
\ee
as well as the following five conditions.
\be
\label{V.8-a}
C(\Sigma, l,E(0))\ \lf[\frac{e^{4\|\al(0)\|_\infty}}{\delta^3\,\rho^4}\ T_1+\sqrt{\frac{e^{4\|\al(0)\|_\infty}}{2\,\delta^2\,\rho^4}\ T_1}+\frac{e^{4\|\al(0)\|_\infty}}{\delta^3}\ \lf[ \frac{T_1}{\rho^4} \rg]^{1/3}\rg]\le 1/8\ ,
\ee
we have also
\be
\label{V.9-a}
C(\Sigma,l,E(0))\ \lf[\sqrt{\beta}\ \sqrt{W(\vec{\Phi}(0))-W(\vec{\Phi}(T_1))}\rg]\le1/8\ ,
\ee
as well as
\be
\label{V.10-a}
C(\Sigma,l,E(0))\ \lf[ e^{2\|\al(0)\|_\infty}\ \sqrt{T_1}\ \sqrt{W(\vec{\Phi}(0))-W(\vec{\Phi}(T_{1}))} \rg]\le1/8\ ,
\ee
moreover
\be
\label{V.11-a}
C(\Sigma,l,E(0))\ \lf[\frac{e^{8\|\al(0)\|_\infty}}{\delta}\ \lf[W(\vec{\Phi}(0))-W(\vec{\Phi}(T_{1})\rg] \rg]\le1/8\ ,
\ee
and finally
\be
\label{V.12-a}
C(\Sigma,l, E(0))\,\frac{e^{2\|\al(0)\|_\infty}}{R^2}\,  \sqrt{T_{1}}\,\ \sqrt{W(\vec{\Phi}(0))-W(\vec{\Phi}(T_{1}))}\le N_1^{-k}\ \beta/8\ .
\ee
Arguing as above we deduce the existence of the flow until $T_1$ and the a-priori estimate holds true
\be
\label{V.13-a}
E(2\rho,T_1)=\le \frac{3}{4}\, N_1^{-k}\ \beta\ .
\ee
 Assume equality holds in either (\ref{V.14-a}), (\ref{V.8-a}), (\ref{V.10-a}) or (\ref{V.12-a}) then the theorem is proved. Now we assume that the inequality is strict in (\ref{V.14-a}), (\ref{V.8-a}), (\ref{V.10-a}) or (\ref{V.12-a}) .
We then have equality in either (\ref{V.9}) or (\ref{V.11}). Because of the choice of $\beta$, (\ref{V.9}) cannot be an equality. Hence we are left with the case where
\be
\label{V.11-a-a}
C(\Sigma,l,E(0))\ \lf[\frac{e^{8\|\al(0)\|_\infty}}{\delta}\ \lf[W(\vec{\Phi}(0))-W(\vec{\Phi}(T_{1})\rg] \rg]=1/8\ ,
\ee
If this happens we first observe that, because of (\ref{V.13-a}) we have
\be
\label{rho-bis}
\sup_{x_0\in\Sigma}\int_{B_{\rho}^{h(0)}(x_0)}|\vec{\mathbb I}_{\vec{\Phi}(0)}|^2_{g_{\vec{\Phi}}}\ \chi^{h(0)}_{R,x_0}(x)\ \ dvol_{g_{\vec{\Phi}(0)}}\le N_1^{-k}\, {\beta}
\ee
We now consider the a-priori estimate given by lemma~\ref{lm-energy-evol} but between $T_1$ and $T>T_1$. Assuming 
\be
\label{V.17}
E(2\rho,T,T_1):=\sup_{x_0\in\Sigma}\sup_{t\in[T_1,T]}\ E(2R,x_0,t):=\sup_{t\in[0,T]}\sup_{x_0\in\Sigma}\int_{B^{h(t)}_{2\,R}(x_0)}|\vec{\mathbb I}_{\vec{\Phi}(t)}|^2_{g_{\vec{\Phi}(t)}}\ dvol_{g_{\vec{\Phi}(t)}}< N_1^{-k+1}\ \beta\ ,
\ee
we have in particular
\[
E(2\rho,T,0)=\max\lf\{ E(2\rho,T_1,0) ,E(2\rho,T,T_1)\rg\}<\beta\ .
\]
Hence for $T_1<T$ satisfying
\be
\label{V.17-a}
\frac{e^{4\|\al(0)\|_\infty}}{\rho^4} \ T\le \hat{C}_{\Sigma,m}^2\, \frac{l^{26}}{E(0)^8}\ ,
\ee
we have
\be
\label{V.18}
\sup_{t\in[0,T]}\lf\|\al(t)-\al(0)\rg\|_{L^\infty(\Sigma)}\le \hat{C}_{\Sigma,m}\ ,
\ee
where $\hat{C}_{\Sigma,m}$ is exactly the \underbar{same constant} as in (\ref{V.3}). Bounding $\|\al\|_\infty$ by $\|\al(0)\|_\infty+\hat{C}_{\Sigma,m}$ in  the a-priori estimate given by lemma~\ref{lm-energy-evol} but between $T_1$ and $T>T_1$ satisfying (\ref{V.17}) and (\ref{V.17-a}) exactly as we did (with the same constants) to obtain (\ref{V.6}), Since $E(0)>E(\vec{\Phi}(T_1))$ we obtain the a-priori estimate assuming (\ref{V.17})
\be
\label{V.19}
\begin{array}{l}
\ds E(2\rho,T,T_1) +2^{-1}\,\sup_{x_0\in\Sigma}\int_{T_1}^T\int_{B_{2\rho}^{h(t)}(x_0)}\ |\Delta_{g_{\vec{\Phi}(t)}}\vec{H}(t)|^2\  \ dvol_{g_{\vec{\Phi}(t)}}\ dt\\[5mm]
\ds\le\, N_1\,\sup_{x_0\in\Sigma}\int_{B_{\rho}^{h(T_1)}(x_0)}|\vec{\mathbb I}_{\vec{\Phi}(T_1)}|^2_{g_{\vec{\Phi}}}\ \chi^{h(T_1)}_{\rho,x_0}(x)\ \ dvol_{g_{\vec{\Phi}(0)}}\\[5mm]
%\ds \quad + N_1\, [\delta+\beta^{1/2}]\ \sup_{x_0\in\Sigma}\ \int_0^T\int_{B_{2\,R}^{h(t)}(x_0)}\ |\Delta_{h(t)}\vec{H}|^2\  dvol_{g_{\vec{\Phi}(t)}} \ dt\\[5mm]
\ds\quad+C(\Sigma, l,E(0))\ \lf[\frac{e^{4\|\al(0)\|_\infty}}{\delta^3\,\rho^4}\ (T-T_1)+\sqrt{\frac{e^{4\|\al(0)\|_\infty}}{2\,\delta^2\,\rho^4}\ (T-T_1)}+\frac{e^{4\|\al(0)\|_\infty}}{\delta^3}\ \lf[ \frac{(T-T_1)}{\rho^4} \rg]^{1/3}\rg]\  N_1^{-k+1}\  \beta\\[5mm]
 \ds\quad+C(\Sigma,l,E(0))\ \lf[\sqrt{\beta}\ \sqrt{W(\vec{\Phi}(T_1))-W(\vec{\Phi}(T))}\rg]\  N_1^{-k+1}\ \beta\\[5mm]
\ds\quad+C(\Sigma,l,E(0))\ \lf[ e^{2\|\al(0)\|_\infty}\ \sqrt{T}\ \sqrt{W(\vec{\Phi}(T_1))-W(\vec{\Phi}(T))} \rg]\  N_1^{-k+1}\  \beta\\[5mm]
\ds\quad+ C(\Sigma,l,E(0))\ \lf[\frac{e^{8\|\al(0)\|_\infty}}{\delta}\ \lf[W(\vec{\Phi}(T_1))-W(\vec{\Phi}(T))\rg] \rg]\  N_1^{-k+1}\ \beta\\[5mm]
\ds\quad+C(\Sigma,l, E(0))\,\frac{e^{2\|\al(0)\|_\infty}}{\rho^2}\,  \sqrt{T-T_1}\,\ \sqrt{W(\vec{\Phi}(T_1))-W(\vec{\Phi}(T))}
\end{array}
\ee
It is fundamental in the argument that the constant $C(\Sigma,l,E(0))$ in (\ref{V.19}) are exactly the same as in (\ref{V.15}). Choosing the maximal value of $T=T_2>T_1$ such that, in addition to (\ref{V.17-a}), we have the following five conditions being satisfied
\be
\label{V.8-a-b}
C(\Sigma, l,E(0))\ \lf[\frac{e^{4\|\al(0)\|_\infty}}{\delta^3\,\rho^4}\ (T_2-T_1)+\sqrt{\frac{e^{4\|\al(0)\|_\infty}}{2\,\delta^2\,\rho^4}\ (T_2-T_1)}+\frac{e^{4\|\al(0)\|_\infty}}{\delta^3}\ \lf[ \frac{T_2-T_1}{\rho^4} \rg]^{1/3}\rg]\le 1/8\ ,
\ee
we have also
\be
\label{V.9-a-b}
C(\Sigma,l,E(0))\ \lf[\sqrt{\beta}\ \sqrt{W(\vec{\Phi}(T_1))-W(\vec{\Phi}(T_2))}\rg]\le1/8\ ,
\ee
as well as
\be
\label{V.10-a-b}
C(\Sigma,l,E(0))\ \lf[ e^{2\|\al(0)\|_\infty}\ \sqrt{T_2-T_1}\ \sqrt{W(\vec{\Phi}(T_1))-W(\vec{\Phi}(T_{2}))} \rg]\le1/8\ ,
\ee
moreover
\be
\label{V.11-a-b}
C(\Sigma,l,E(0))\ \lf[\frac{e^{8\|\al(0)\|_\infty}}{\delta}\ \lf[W(\vec{\Phi}(T_1))-W(\vec{\Phi}(T_{2})\rg] \rg]\le1/8\ ,
\ee
and finally
\be
\label{V.12-a-b}
C(\Sigma,l, E(0))\,\frac{e^{2\|\al(0)\|_\infty}}{R^2}\,  \sqrt{T_{1}}\,\ \sqrt{W(\vec{\Phi}(T_1))-W(\vec{\Phi}(T_{2}))}\le N_1^{-k+1}\ \beta/8\ .
\ee
Similarly as before, since $N_1$ is taken from the beginning of the argument to be larger than $8$, we obtain the a-priori estimate
\be
\label{V.19-aaa}
 E(2\rho,T,T_1)\le \frac{3}{4}\ N_1^{-k+1}\ \beta\ .
\ee
Then, similarly as before, this implies that the solution exists until $T_2$.

\medskip

Assuming now there is an equality in either (\ref{V.17-a}), (\ref{V.8-a-b}), (\ref{V.10-a-b}), (\ref{V.11-a-b}) or (\ref{V.12-a-b}), the theorem is proved. Because of the assumption (\ref{V.99}) we cannot have an equality in (\ref{V.9-a-b}) and hence, the only remaining possibility which has to be considered is
\be
\label{V.20}
8\ C(\Sigma,l,E(0))\ \lf[\sqrt{\beta}\ \sqrt{W(\vec{\Phi}(T_1))-W(\vec{\Phi}(T_2))}\rg]=1
\ee
Combined with (\ref{V.11-a-a}), (\ref{V.20})  is implying
 \be
\label{V.21}
8\ C(\Sigma,l,E(0))\ \lf[\sqrt{\beta}\ \sqrt{W(\vec{\Phi}(0))-W(\vec{\Phi}(T_2))}\rg]=2
\ee
By iterating the argument we would either prove the theorem or obtain an increasing  sequence of times $0<T_1<T_2\cdots <T_k$ such that
 \be
\label{V.22}
8\ C(\Sigma,l,E(0))\ \lf[\sqrt{\beta}\ \sqrt{W(\vec{\Phi}(0))-W(\vec{\Phi}(T_i))}\rg]=i
\ee
The identity (\ref{V.22}) for $i=k$ contradicts (\ref{k}) hence the theorem must have been proved in one of the preceding steps. This concludes the proof of the main theorem~\ref{th-I.1}.\hfill $\Box$


\begin{thebibliography}{99}
 \bibitem{Aub}  Aubin, Thierry Nonlinear analysis on manifolds. Monge-Amp\`ere equations. Grundlehren der mathematischen Wissenschaften. Vol. 252. New York: Springer-Verlag. 1982.
\bibitem{Bav} Bavard, Christophe
Disques extr\'emaux et surfaces modulaires.[Extremal disks and modular surfaces]
Ann. Fac. Sci. Toulouse Math. (6)5(1996), no.2, 191-202.
\bibitem{BeRi} Bernard, Yann; Rivi\`ere, Tristan
Energy quantization for Willmore surfaces and applications
Ann. of Math. (2) 180 (2014), no. 1, 87-136.
\bibitem{Bers} Bers, Lipman Spaces of degenerating Riemann surfaces.
Ann. of Math. Studies, No. 79,
Princeton Univ. Press, Princeton, N.J., 1974, pp. 43-55
 \bibitem{ChLi}  Cheng, Shiu Yuen; Li, Peter Heat kernel estimates and lower bound of eigenvalues Comment. Math. Helv. 56 (1981), no. 3, 327-338.
 \bibitem{Cro} Croke, Christopher B. Some isoperimetric inequalities and eigenvalue estimates
Ann. Sci. \'Ecole Norm. Sup. (4) 13 (1980), no. 4, 419-435.
\bibitem{Gra1} Grafakos, Loukas { Classical Fourier analysis.} Third edition. Graduate Texts in Mathematics, 249. Springer, New York, 2014.
\bibitem{HaJo} Habermann, Lutz; Jost, J\"urgen Riemannian metrics on Teichm\"uller space
Manuscripta Math. 89 (1996), no. 3, 281-306.
 \bibitem{Hel} H\'elein Fr\'ed\'eric Harmonic maps, conservation laws and moving frames
Cambridge Tracts in Math., 150
Cambridge University Press, Cambridge, 2002.
\bibitem{Hub} Huber, A. On the Isoperimetric Inequality on Surfaces of Variable Gaussian Curvature  Ann. of Math.,
Vol. 60, 1954, p. 237-247
 \bibitem{Jost} Jost, J\"urgen
Compact Riemann surfaces. An introduction to contemporary mathematics. Third edition.  Universitext
Springer-Verlag, Berlin, 2006. 
 \bibitem{KS1} Kuwert, Ernst; Sch\"atzle, Reiner The Willmore flow with small initial energy
J. Differential Geom. 57 (2001), no. 3, 409-441.
 \bibitem{KS2} Kuwert, Ernst; Sch\"atzle, Reiner  Gradient flow for the Willmore functional
Comm. Anal. Geom. 10 (2002), no. 2, 307-339.
\bibitem{KS3}  Kuwert, Ernst; Sch\"atzle, Reiner The Willmore functional
CRM Series, 13
Edizioni della Normale, Pisa, 2012, 1-115.
\bibitem{LaRi} Laurain, Paul; Rivi\`ere, Tristan Optimal estimate for the gradient of Green's function on degenerating surfaces and applications
Comm. Anal. Geom. 26 (2018), no. 4, 887-913.
 \bibitem{LaRi-w} Laurain, Paul; Rivi\`ere, Tristan Energy quantization of Willmore surfaces at the boundary of the moduli space
Duke Math. J. 167 (2018), no. 11, 2073-2124.
\bibitem{MN} Marques Fernando Cod\'a and Neves Andre, Min-max theory and the Willmore conjecture, Ann. of Math. 179 (2014), 683-782.
\bibitem{Mas} Masur, Howard Extension of the Weil-Petersson metric to the boundary of Teichmuller space
Duke Math. J. 43 (1976), no. 3, 623-635.
\bibitem{MiRi} Michelat, Alexis; Rivi\`ere, Tristan Pointwise expansion of degenerating immersions of finite total curvature
J. Geom. Anal. 33 (2023), no. 1,Paper No. 24, 91 pp.
\bibitem{PalRiv} Palmurella, Francesco; Rivi\`ere, Tristan The parametric approach to the Willmore flow
Adv. Math. 400 (2022), Paper No. 108257, 48 pp.
\bibitem{Riv-w} Rivi\`ere Tristan Analysis aspects of Willmore surfaces
Invent. Math. 174 (2008), no. 1, 1-45
\bibitem{Riv-notes}   Rivi\`ere Tristan   Conformally Invariant Variational Problems {\it unpublished notes} (2012)   https://people.math.ethz.ch/~riviere/lecture-notes
 \bibitem{Pa-Ci} Rivi\`ere, Tristan Weak immersions of surfaces with $L^2-$bounded second fundamental form
IAS/Park City Math. Ser., 22
American Mathematical Society, Providence, RI, 2016, 303-384.
\bibitem{Riv-14} Rivi\`ere, Tristan
Variational principles for immersed surfaces with $L^2$-bounded second fundamental form.(English summary)
J. Reine Angew. Math.695 (2014), 41-98.
\bibitem{Riv-min} Rivi\`ere, Tristan Willmore minmax surfaces and the cost of the sphere eversion
J. Eur. Math. Soc. (JEMS) 23 (2021), no. 2, 349-423.
  \bibitem{Ri1} Rivi\`ere, Tristan Lower semi-continuity of the index in the viscosity method for minimal surfaces. Int. Math. Res. Not. IMRN 2021, no. 8, 5651-5675. 
  \bibitem{RuTo}  Rupflin, Melanie; Topping, Peter M. Teichm\"uller harmonic map flow into nonpositively curved targets. J. Differential Geom. 108 (2018), no. 1, 135-184. 
  \bibitem{Sim}  Simon, Leon Existence of surfaces minimizing the Willmore functional, Commun. Analysis Geom. 1(2) (1993) 281-326.  
  \bibitem{Tro}  Tromba, Anthony J. Teichm\"uller theory in Riemannian geometry. Lecture notes prepared by Jochen Denzler. Lectures in Mathematics ETH Z\"urich. Birkh\"auser Verlag, Basel, 1992.
 \end{thebibliography}
\end{document}